\documentclass[11pt,a4paper]{amsart}
\usepackage[T1]{fontenc}
\usepackage[utf8]{inputenc}
\usepackage[margin=25mm]{geometry}
\usepackage{amsmath,amssymb,amsthm,mathtools}
\usepackage{mathrsfs}
\usepackage{booktabs,longtable,array}
\usepackage[expansion=false]{microtype}
\usepackage[colorlinks=true,linkcolor=blue,citecolor=blue,urlcolor=blue]{hyperref}
\allowdisplaybreaks[2]
\makeatletter
\renewcommand{\paragraph}{\@startsection{paragraph}{4}%
  {\normalparindent}{\z@}{-\fontdimen2\font}{\normalfont}}
\makeatother

\numberwithin{equation}{section}
\newtheorem{theorem}{Theorem}[section]
\newtheorem{proposition}[theorem]{Proposition}
\newtheorem{lemma}[theorem]{Lemma}
\newtheorem{corollary}[theorem]{Corollary}
\newtheorem*{mainresult}{Main theorem}
\theoremstyle{remark}
\newtheorem{remark}[theorem]{Remark}
\theoremstyle{plain}
\newcommand{\eb}{\epsilon_{\rm b}}
\newcommand{\OO}{\mathcal O}
\newcommand{\sfq}{\mathrm{sf}}
\newcommand{\eps}{\varepsilon}
\newcommand{\Rea}{\operatorname{Re}}

\newcommand{\rad}{\operatorname{rad}}
\newcommand{\Res}{\operatorname*{Res}}
\newcommand{\Src}[2]{\textup{\cite[#1~#2]{QRH}}}
\newcommand{\Prev}[2]{\textup{\cite[#1~#2]{LiuArXiv}}}
\newcommand{\LiuLoc}[1]{\textup{[#1]}}

\newcommand{\pos}[1]{\left(#1\right)_+}
\hypersetup{pdftitle={Zero-free half-planes for Hecke L-functions over number fields and Quasi-Riemann Hypothesis}}
\title[Zero-free half-planes over number fields]{Zero-free half-planes for Hecke $L$-functions over number fields and Quasi-Riemann Hypothesis}
\author[Baiying Liu]{Baiying Liu}
\address{Department of Mathematics\\
Purdue University\\
West Lafayette, IN, 47907, USA}
\email{liu2053@purdue.edu}
\date{20 October 2026}
\begin{document}
\begin{abstract}
In \cite{LiuArXiv}, varying the length and moment parameters in the
cubic-theta argument of \cite{QRH} gave the zero-free half-plane
$\Rea s>B_{\mathrm{new}}=0.874957019420098946\ldots$ for finite-order
Hecke $L$-functions over $\mathbb Q(\sqrt{-3})$ and for Dirichlet
$L$-functions, and its optimality in the specified exponent
comparison system.
In this paper, we extend the argument
to finite-order Hecke
$L$-functions over cyclotomic fields containing $\mu_{12}$ and, by
class-field descent, to every pair $(K,\eta)$ whose cyclic class field
has metabelian normal closure over $\mathbb Q$, including every
finite-order character over every abelian number field. Then,
adapting the averaged numerator estimate developed in \cite{XS26} for Hecke $L$-functions over $\mathbb Q(\sqrt{-3})$
and combining it with a limiting contradiction argument, we obtain the boundary $B_\flat=0.874886647039964\ldots$, the root near
$0.8749$ of $133488B^4-376596B^3+384534B^2-167541B+26231$; it is the
exact infimum of the refined comparison system, unattained by finite
chains starting above it. Following the general-field strategy of
\cite{CT26}, we extend the analytic argument to all number fields
$F\supset\mu_{12}$. Descent from $K(\zeta_{12})$ gives the same bound
for finite-order Hecke $L$-functions over arbitrary number fields $K$,
giving a uniform Quasi-Riemann Hypothesis over number fields.
The extended analytic argument uses the
exceptional-theta constructions and local calculations of the
preprint~\cite{Artin}.
\end{abstract}
\maketitle

\tableofcontents

\section*{Introduction}
\label{sec:introduction}

The Quasi-Riemann Hypothesis asks whether there is a fixed $B<1$ such
that $\zeta(s)\ne0$ for $\Rea s>B$. Its uniform version for a family of
$L$-functions asks for the same $B$ independently of conductor and
height, with principal poles at $s=1$ allowed. For a number field $K$,
the corresponding question concerns the Dedekind zeta function
$\zeta_K$ and the finite-order Hecke $L$-functions $L_K(s,\chi)$;
see \cite[Chapter~5]{IK} for background.

Recently, \cite{QRH} proved, by a cubic-theta
argument over $\mathbb Q(\sqrt{-3})$, that every finite-order Hecke
$L$-function over $\mathbb Q(\sqrt{-3})$ and every Dirichlet
$L$-function is zero-free for $\Rea s>7/8$. In \cite{LiuArXiv},
varying the length parameters $b,\ell$ and the moment parameter
$\kappa$ of that argument, we proved the zero-free half-plane
$\Rea s>B_{\mathrm{new}}$, where
\[
 B_{\mathrm{new}}=\frac{18\xi_{\mathrm{new}}+7}{18\xi_{\mathrm{new}}+9}
 =0.874957019420098946\ldots,
\]
and $\xi_{\mathrm{new}}$ is the unique root in $(0.38858,0.38859)$ of
$1188\xi^3-2160\xi^2+171\xi+190$, and that
$B_{\mathrm{new}}$ is optimal within the specified exponent
comparison system. The results of \cite{LiuArXiv} used below are recalled in Part~\ref{part:prelim}.

The present paper extends that argument beyond $\mathbb Q(\sqrt{-3})$.
Cureton and Thorner~\cite[Theorem 1.1, Corollary 1.2 and Sections 5.1 and 14]{CT26}
obtained, for every number field $K$, the region $\Rea s>31/32$ for
$\zeta_K$ and finite-order Hecke $L$-functions, by an argument over
$K(\zeta_{12})$ followed by abelian descent. We follow this
route with the varying-parameter moments of \cite{LiuArXiv} and extend
the averaged numerator estimate of Xu and Sang~\cite{XS26} from
$\mathbb Q(\sqrt{-3})$ to the fields considered here, obtaining
\[
 \begin{gathered}
 B_\flat=0.874886647039964\ldots<B_{\mathrm{new}},\\
 133488B_\flat^4-376596B_\flat^3+384534B_\flat^2-167541B_\flat+26231=0.
 \end{gathered}
\]

\begin{mainresult}
For every number field $K$ and every finite-order Hecke character $\chi$
of $K$,
\[
 L_K(s,\chi)\ne0\qquad(\Rea s>B_\flat),
\]
with the pole of the principal character at $s=1$ permitted. In
particular $\zeta_K(s)\ne0$ for $\Rea s>B_\flat$, $s\ne1$.
\end{mainresult}
This is Theorem~\ref{vi:thm:main}; for cyclotomic
base fields containing $\mu_{12}$ and their class-field descendants it
is Theorem~\ref{thm:v-main}.

A subsequent refinement of the zero-free boundary over arbitrary
number fields is developed in~\cite{LiuHexagon}, using different method, based on the work in this paper. 

\paragraph{\it The cyclotomic argument.}
Part~\ref{part:V} replaces $\mathbb Q(\sqrt{-3})$ by a cyclotomic field
$F=\mathbb Q(\zeta_m)\supset\mu_{12}$, and the cubic theta function by
the exceptional theta representation of
Kazhdan and Patterson~\cite{KP,KPerr}, with the allowed generators, local
Whittaker factors and radial functional equation of \cite{Artin}.
Sections~\ref{v:sec:unmarked-assembly}--\ref{v:sec:M2-assembly} give
derivations of the cyclotomic inverse and unweighted product moments,
the latter for $13/18\le\kappa\le1$ as in \Prev{Part}{IV}.
Section~\ref{v:sec:compensated-assembly} combines them into
the compensated comparison, whose length and count exponents are
those of that part; with the scalar inequalities and the
iteration of Section~\ref{sec:prev-iv}, it gives
$B_\infty=B_{\mathrm{new}}$.

\paragraph{\it Pointwise and averaged numerator comparisons.}
The distinction between $B_{\mathrm{new}}$ and the bound of
Xu and Sang~\cite{XS26} is visible in the exponent of a nonprincipal
row class. To compare the formulas, use balanced lengths
\[
 l_x=\frac{1-b-\ell}{2},\qquad l_y=\frac{1+b-\ell}{2},\qquad
 h=\frac{1+b+3\ell}{2},\qquad B(\ell)=\frac{11}{12}-\frac\ell4,
\]
and put $C_0=-1/4+b/6+5\ell/4$. A class at norm scale $U=Z^d$
has zero-class parameter $a=(1+\delta)/2$ and weighted mean
prime-polynomial amplitude $g=\delta x$, $0\le x\le1/2$
(the mean is denoted $\bar g$ in Section~\ref{sec:prev-limit}).
Write $R=R_\kappa(\delta,x)$ for the counting exponent of the class,
so its cardinality is at most $U^{R+\eps}$ up to polynomial height
factors; its explicit formula is \eqref{eq:v-refined-count}.
Here $\kappa$ is the moment parameter, and
$c_\kappa=1/(3\kappa)$.

The argument giving $B_{\mathrm{new}}$ reflects the numerator on
$\Rea w=1-a$ and bounds the reflected value pointwise, at conductor
cost $U^{a-1/2}=U^{\delta/2}$, with small contour losses suppressed.
At the endpoint $d=h$, its relative exponent
\eqref{eq:balanced-extension-exponent} is
\[
 E_{\rm pt}=C_0+(1/2+\ell)\delta+\ell\delta x-h(1-R).
\]
Varying $b,\ell,\kappa$ and the justifying zero-free bound within
this pointwise-numerator system gives $B_{\mathrm{new}}$;
Theorem~\ref{thm:iv-system-optimality} concerns this system.

Xu and Sang~\cite[Sections~5--6]{XS26} treat finite-order Hecke
$L$-functions only over $\mathbb Q(\sqrt{-3})$, with Dirichlet
$L$-functions obtained as consequences. They use an amplified first
moment of the reflected numerator over the same class $\mathscr B$.
In the notation of Proposition~\ref{v:prop:refined-numerator}, the estimate is
\[
 \sum_{u\in\mathscr B}
 |L(1/2+\theta+it,\bar\psi_u^*)|
 \ll U^{(3R+1)/4-\theta/2+\eps}(1+H)^{A'},
\]
where $\psi_u^*$ is the primitive character inducing the numerator,
$H$ bounds the numerator and prime-polynomial heights, and the
shift $\theta$ obeys the available prime-length and moment-capacity
conditions. This comes from a fourth moment of the family, an approximate
functional equation and H\"older's inequality (over $F$, the sextic
fourth moment of Remark~\ref{v:rem:sextic-fourth-moment}, a case of
Theorem~\ref{v:thm:M2-assembled}), with prime amplification
selected from the factors that define the class. Reflecting from
$\Rea w=1/2-\theta$ contributes $U^\theta$; the resulting replacements
in the full exponent are
\[
 U^{R+\delta/2}\longmapsto U^{(3R+1)/4+\theta/2},\qquad
 Y^{-a}\longmapsto Y^{-1/2-\theta},\qquad Y=Z^{l_y}.
\]
For $1/8\le\delta\le5/6$ the admissible full-capacity choice is
$\theta=g/(6\kappa)$, and the new endpoint
\eqref{eq:v-refined-endpoint} becomes
\[
 \begin{aligned}
 E_{\rm av}={}&C_0+\frac\delta2
 +\delta x\left\{\ell+c_\kappa\left(\frac h4-\frac{l_y}{2}\right)\right\}
 -\frac{3h}{4}(1-R).
 \end{aligned}
\]
Thus the improvement uses a different numerator estimate and a different
exponent, while retaining the counting exponent $R_\kappa$.
The optimum for $E_{\rm pt}$ does not bound this averaged-numerator
system. These displays describe the main endpoint; the shorter rows,
limited prime supply, local errors and contour losses also require
control. Over $\mathbb Q(\sqrt{-3})$, Xu and Sang use two parameter
choices justified in order:
\[
 \beta_*\le\frac78
 \ \Longrightarrow\ \kappa=\frac34
 \ \Longrightarrow\ \beta_*\le\frac{69991}{80000}
 \ \Longrightarrow\ \kappa=\frac{29991}{40000}
 \ \Longrightarrow\ \beta_*\le\frac{8748867}{10^7},
\]
where $\beta_*$ is the supremum of real parts of nontrivial zeros in
the fixed-field finite-order Hecke family. Each application uses the
zero-free condition $\beta_*\le(1+\kappa)/2$, supplied by the preceding
bound in the displayed chain.

\paragraph{\it Extending the averaged-numerator method.}
Our extension of this argument has an analytic component and a
scalar component. Over a fixed cyclotomic field $F\supset\mu_{12}$,
Proposition~\ref{v:prop:refined-numerator} derives the amplified numerator
estimate from the ordinary product moment, using the
number-field approximate functional equation and the actual prime-sum
factors of the detected class. Sections~\ref{v:sec:refined-local}--\ref{v:sec:refined-high}
assemble the amplitude-dependent contours and local errors into the
refined remainder. The inverse moment, the row count and the direct
bound used in the preceding comparison are retained. In particular,
we adapt the numerator argument of \cite{XS26}; we do not use its
zero-free conclusion as an input over $F$.

For the scalar component, we vary the parameters throughout
$\kappa_\flat\le\kappa\le3/4$. The critical class is
$(\delta_c(\kappa),1/2)$, where $E_{\rm av}$ is independent of $b$.
Its vanishing determines $\ell^*(\kappa)$, and the equation
$\partial_\delta E_{\rm av}=0$ determines $b^*(\kappa)$.
Lemma~\ref{v:lem:refined-family} proves that this tangent choice controls
the full endpoint rectangle and all remaining admissibility conditions;
choosing $\ell<\ell^*(\kappa)$ sufficiently close gives a strict margin.
Put
\[
 F_\flat(\kappa)=\frac{11}{12}-\frac{\ell^*(\kappa)}4,
 \qquad
 F_\flat(\kappa_\flat)=\frac{1+\kappa_\flat}{2}=B_\flat.
\]
This fixed-point equation gives the quartic for $B_\flat$ displayed
above. To justify the moment parameter, suppose
$\beta=\beta_*(F)>B_\flat$. The preceding pointwise comparison gives
$\beta\le B_\infty$, so $\kappa=2\beta-1$ lies in
$(\kappa_\flat,3/4)$ and satisfies the zero-free condition required by
M2 with equality.
Lemma~\ref{v:lem:refined-critical} gives $F_\flat(\kappa)<\beta$.
We can therefore choose $\ell<\ell^*(\kappa)$ close enough that
$B(\ell)<\beta$, and Lemma~\ref{v:lem:refined-closure} gives
$\beta\le B(\ell)<\beta$, a contradiction.
Theorem~\ref{v:thm:refined-limit} thus obtains the limiting
bound in one contradiction step from $B_\infty$.
\paragraph{\it Optimality within the refined comparison system.}
Theorem~\ref{v:thm:refined-optimality} proves that every permitted finite
chain starting from a bound $B_{\rm prev}\ge B_\flat$ yields only
boundaries $B\ge B_\flat$. If the starting bound is strictly larger
than $B_\flat$, no finite chain attains $B_\flat$. Conversely, for every
$\epsilon>0$, a finite chain of admissible refined steps starting from
$7/8$ yields a boundary satisfying
\[
 B_\flat<B<B_\flat+\epsilon.
\]
Thus $B_\flat$ is the exact infimum of this comparison system.
These are scalar optimality statements for the explicitly specified
system of Section~\ref{v:sec:refined-optimality}. The limiting
contradiction argument above gives the endpoint conclusion, even though a finite chain starting above
$B_\flat$ cannot attain it. The analytic dependencies of these comparisons are specified in
Remark~\ref{v:rem:analytic-inputs}.

As an independent check in exact rational arithmetic,
Appendix~\ref{app:refined-chain} gives a chain of four verified stages:
the first two recover over $F$ the boundary $B_1=8748867/10^7$ of Xu and
Sang, and the last reaches $B_\sharp=874886647041/10^{12}$, with
$B_\sharp-B_\flat<1.04\times10^{-12}$. By Section~\ref{sec:v-optimality}, the
pointwise-numerator system has the same optimum $B_\infty$ as the
system of \cite[Theorem~34.4]{LiuArXiv}; the refined
boundary lies outside it. Section~\ref{sec:v-descent} proves
the descent to every pair $(K,\eta)$ whose cyclic class
field becomes abelian over some cyclotomic field containing $\mu_{12}$,
and characterizes this class arithmetically.

\paragraph{\it Other parameter refinement.}
Following Xu and Sang, Newton~\cite[Theorem~1]{Newton26} gives an exact rational parameter certificate for
\[
 B_{\mathrm N}=\frac{43744332481}{50000000000}=0.87488664962
\]
for the same finite-order Hecke family over $\mathbb Q(\sqrt{-3})$
and its Dirichlet consequences. The certificate controls the continuous endpoint rectangle
by a rational polynomial identity and a negative-discriminant argument;
it does not introduce a stronger mean-value estimate. In
\cite[Theorem~2]{Newton26}, Newton also
obtains the barrier
\[
 B_{\mathrm{lim}}=0.874886649616631419\ldots
\]
within the stated endpoint admissibility system with
$29991/40000\le\kappa\le3/4$ and a strictly positive endpoint margin.

The smaller boundary $B_\flat<B_{\mathrm{lim}}$ here uses a
varying justifying boundary and moment parameter, with the condition
$\beta_*(F)\le(1+\kappa)/2$ retained at every application. The tangent
family of Lemma~\ref{v:lem:refined-family} reaches
$\kappa_\flat=2B_\flat-1=0.749773294079928\ldots<29991/40000$;
in the contradiction argument of Theorem~\ref{v:thm:refined-limit}, the choice is
$\kappa=2\beta_*(F)-1$. Newton's fixed-range barrier and the optimality
statement of Theorem~\ref{v:thm:refined-optimality} therefore concern
different admissible moment ranges. This distinction concerns the scalar parameter systems; the analytic
inputs for our refinement are specified in
Remark~\ref{v:rem:analytic-inputs}.

\paragraph{\it Arbitrary number fields.}
Following the general-field strategy of \cite{CT26}, we check in
Section~\ref{vi:sec:generality} of Part~\ref{part:VI} that the analytic argument
of Part~\ref{part:V} uses only $\mu_{12}\subset F$ and
standard facts valid for every number field: total imaginarity,
$S$-unit and class-group finiteness, Hilbert reciprocity, Kummer
theory, the Kazhdan--Patterson theory over fields containing $\mu_3$,
and the large sieves of Goldmakher--Louvel and
Blomer--Goldmakher--Louvel, which are stated for general number fields
containing the relevant roots of unity. Running the argument over
$F=K(\zeta_{12})$ and descending by class-field factorization gives the
main theorem.

\paragraph{\it Structure of the paper.}
Part~\ref{part:prelim} fixes the notation of \Prev{Section}{2} and
recalls the results of \Prev{Parts}{I--IV} used later.
Part~\ref{part:V} gives the cyclotomic extension and the class-field
descent; Part~\ref{part:VI} treats arbitrary number fields.
Appendix~\ref{app:inputs-index} indexes the cited results and their
roles.

\subsection*{Acknowledgments}
Exact symbolic checks used SymPy~\cite{SymPy}.
OpenAI Codex (GPT 6.1 Sol Medium-Ultra) assisted with drafting, checking
references, and the search for improved bounds. Anthropic Claude
(Opus 5.5 Medium) assisted with auditing. The Lean formalization is in progress, based on the Lean files in \cite{LeanExtension} and in \cite{CTLean}.

\part{Notation and results from the moment-parameter extension}
\label{part:prelim}

In this part, we fix the notation and recall, with precise locators,
the results of \cite{LiuArXiv} used later; no new results are proved.
Section~\ref{sec:arithmetic} gives the arithmetic conventions over
$F=\mathbb Q(\sqrt{-3})$, which Part~\ref{part:V} adapts to cyclotomic
fields in Section~\ref{sec:v-scope}. Section~\ref{sec:prev-iv} recalls two
kinds of results, as explained at the start of that section: analytic
estimates over $\mathbb Q(\sqrt{-3})$, which serve as models for the
constructions over $F$ in Sections~\ref{sec:v-packet}--\ref{v:sec:compensated-assembly},
and scalar relations between explicit real functions, which are used as
stated in Sections~\ref{v:sec:compensated-assembly}--\ref{sec:v-optimality}.

\section{Notation and preliminaries}
\label{sec:arithmetic}

In this section, we recall the conventions of \Prev{Section}{2} for the
ground field $F=\mathbb Q(\sqrt{-3})$ of \cite{QRH,LiuArXiv}, in which
Section~\ref{sec:prev-iv} is formulated: ideal and asymptotic notation,
the length exponent $N=Z^r$ of \eqref{eq:length-exponent-convention},
character presentations and primitive induction
\eqref{eq:primitive-induction}, the family zero supremum $\beta_*$ of
\eqref{eq:beta-star}, and the sextic residue symbols with the zero
convention \eqref{eq:zeros}. Part~\ref{part:V} retains them, with the
change of ground field of Section~\ref{sec:v-scope}.

Set
\[
 \omega=e^{2\pi i/3},\quad
 \OO=\mathbb Z[\omega],\quad \lambda=\sqrt{-3},\quad
 q_a=|a|^2,\quad q_{\mathfrak a}=N\mathfrak a,\quad
 \varpi(a)=a/|a|.
\]
$\OO$ is norm-Euclidean. Ideals coprime to $3$ have unique
primary generators $a\equiv1\pmod3$. Ideal sums use these
generators; sums over elements may also include units.
Unless stated otherwise, ideal sums run over nonzero integral ideals.
A sum over $x\bmod c$ runs over residue classes in $\OO/(c)$.
For an integral ideal $\mathfrak a$, $N\mathfrak a=|\OO/\mathfrak a|$.
Write $v_p(a)$ for its prime valuation, $\rad(a)=\prod_{p\mid a}p$,
and $\mu$ for the ideal M\"obius function. For a nonzero integral ideal,
\[
 \mu(\mathfrak a)=
 \begin{cases}
  1,&\mathfrak a=\OO,\\
  (-1)^k,&\mathfrak a=\mathfrak p_1\cdots\mathfrak p_k
          \text{ with distinct prime ideals }\mathfrak p_j,\\
  0,&\mathfrak p^2\mid\mathfrak a\text{ for some prime ideal }\mathfrak p.
 \end{cases}
\]
When an ideal is represented by a primary generator $n$, we write
$\mu(n)=\mu((n))$. The notation $\sfq$ in a sum
means squarefree, $1_E$ is the indicator of a condition $E$, and
$(a,b)=1$ means that the corresponding ideals are coprime.
A prime outside $S$ is called good.

Unless specified otherwise, $Z\ge2$ is a real variable tending to infinity;
$\eps,\eps_1>0$ are arbitrarily small losses. For $a\in\mathbb C$
and $b\ge0$, $a\ll b$ means $|a|\le Cb$ for some $C>0$, and
\[
 a\ll_\eps b\quad\Longleftrightarrow\quad
 |a|\le C_\eps b,
\]
where $C_\eps>0$ may depend on $\eps$ but on no other varying quantity;
$\ll_S$ is analogous for the fixed prime set $S$.
For $a,b\ge0$, $a\gg b$ means $b\ll a$, and $a\asymp b$ means
\[
 cb\le a\le Cb\qquad\text{for some }c,C>0
\]
(not that $a/b$ is constant).
These constants are independent of the varying quantities, subject to
any indicated dependence. We write $x_+=\max(x,0)$.

A summation length $N>0$ is often expressed through its
\emph{length exponent} $r$:
\begin{equation}\label{eq:length-exponent-convention}
 N=Z^r,\qquad r=\frac{\log N}{\log Z},
\end{equation}
so $q_n\asymp N$ becomes $q_n\asymp Z^r$ and
$Z^{r_1}Z^{r_2}=Z^{r_1+r_2}$; ``length $r$'' means norm parameter $Z^r$.

Auxiliary letters are used independently in separately stated
lemmas. For $c\in\mathbb R$,
$\int_{(c)}F(s)\,ds=\lim_{T\to\infty}\int_{c-iT}^{c+iT}F(s)\,ds$.
The test functions restricting norms are smooth, supported
on fixed compact subintervals of $(0,\infty)$, with bounded logarithmic
derivatives. Constants $A,J$ in height estimates denote finite
nonnegative exponents or derivative orders fixed before the height cutoff;
their values may differ from one estimate to another.

Fix a finite set $S$ containing the prime ideals above $6$, all prime
ideals dividing the defining moduli of the fixed characters, and all
prime ideals of norm at most a fixed cutoff $P_0>0$. Unless stated otherwise, ideal sums are restricted to ideals coprime
to $S$. These data are independent
of $Z$ and the summation variables, but may depend on a fixed
character $\eta$; the fixed moduli, characters, ray group,
exclusions and normalizations are denoted collectively by $\mathcal A$.

\begin{center}
\small
\begin{tabular}{@{}p{.24\textwidth}p{.69\textwidth}@{}}
\toprule Symbol & Meaning and scope\\\midrule
$\beta_*,\Delta,\kappa$ & Family zero supremum \eqref{eq:beta-star};
$\Delta=\beta_*-7/8$; $\kappa$ is the moment parameter, equal to
$2\beta_*-1$ in a supremum contradiction and an independent parameter
$13/18\le\kappa\le1$ from Section~\ref{sec:prev-moment} on.\\
$\alpha=5/6,\ \varpi(a)$ & Fixed amplification slope;
unit phase $a/|a|$ (written $\alpha(a)$ in \cite{QRH}).\\
$\eb,\ \eps$ & Fixed zero-class width (written $e$ in \cite{QRH});
arbitrarily small exponent loss.\\
$\mathfrak a(r)$ & Amplified moment exponent
$\max\{1,(1+5r)/6\}$ (written $e(r)$ in \cite{QRH}).\\
$Z,U=Z^d,T_1$ & Real variable tending to infinity; dyadic norm parameter;
zero-detection height.\\
$l_x,l_y,b,h,\ell$ & Averaging exponents; $b=l_y-l_x$;
combined frequency exponent; total auxiliary prime length.\\
$z_{\rm cap},z_0$ & Upper bound for total prime length in
\cite[(9.2)]{LiuArXiv}; real contour value $z_0=17/50$ in
\Prev{Section}{11}.\\
$\mathfrak r$ & Coprimality modulus in \cite[(9.2)]{LiuArXiv}
(written $R$ in \cite{QRH}); distinct from a counting exponent $R$.\\
$\mathcal R,\mathcal M$ & Reciprocity phase; Mellin transform
of the Fourier-transformed weight.\\
$H_\eta,J_\eta,f_\eta$ & Nonvanishing correction; normalized average;
principal-residue Mellin integral.\\
$\Gamma(s)$ & Euler gamma function.\\
\bottomrule
\end{tabular}
\end{center}
\noindent Other indices and length parameters are local to the
statement in which they are introduced. A \emph{zero class} is a
set of elements satisfying the same pair of zero-location inequalities
\cite[(5.1)]{LiuArXiv}.

\paragraph{\it Hecke characters, moduli and conductors.}
For a nonzero integral ideal (modulus) $\mathfrak f$, the ray class group is
$\mathrm{Cl}_{\mathfrak f}(F)=I(\mathfrak f)/P_1(\mathfrak f)$, where $I(\mathfrak f)$ consists of
fractional ideals coprime to $\mathfrak f$ and $P_1(\mathfrak f)$
of principal ideals $(a)$ with $a\equiv1\bmod\mathfrak f$.
A finite-order Hecke character is a
homomorphism from this finite group to roots of unity;
its conductor $\mathfrak f_*$ is the smallest modulus through which
it factors. A \emph{character presentation} specifies a defining modulus
$\mathfrak f$ with $\mathfrak f_*\mid\mathfrak f$, together with extension
by zero on integral ideals not coprime to $\mathfrak f$; it
is primitive if $\mathfrak f=\mathfrak f_*$ and imprimitive otherwise.
A character $\psi$ with defining modulus $\mathfrak f$
is induced by a unique primitive character $\psi^*$ of conductor
$\mathfrak f_*\mid\mathfrak f$: the natural surjection gives
\begin{equation}\label{eq:primitive-induction}
 \pi:\mathrm{Cl}_{\mathfrak f}(F)\longrightarrow
       \mathrm{Cl}_{\mathfrak f_*}(F),\qquad
 \psi=\psi^*\circ\pi.
\end{equation}
Thus $\psi(\mathfrak a)=\psi^*(\mathfrak a)$ if
$(\mathfrak a,\mathfrak f)=1$, where $\psi^*$
does not factor through the ray class group of any proper divisor
of $\mathfrak f_*$; but $\psi(\mathfrak a)=0$ if $(\mathfrak a,\mathfrak f)\ne1$,
even when $(\mathfrak a,\mathfrak f_*)=1$.
So the defining modulus records a coprimality restriction,
including primes absent from the conductor; see \cite[Section 4]{QRH}.

Initially,
\[
 L_F(s,\eta)=\sum_{\mathfrak a\ne0}
       \eta(\mathfrak a)(N\mathfrak a)^{-s}
   =\prod_p(1-\eta(p)q_p^{-s})^{-1},
 \qquad \Rea s>1;
\]
outside this region the notation refers to its meromorphic
continuation. A norm twist multiplies a character by
$\mathfrak a\mapsto(N\mathfrak a)^{it}$, $t\in\mathbb R$.

Let $\mathcal Z_F$ be the set of nontrivial zeros of all primitive
finite-order Hecke $L$-functions over $F=\mathbb Q(\sqrt{-3})$, of
arbitrary conductor; poles are excluded. As in \cite{LiuArXiv}, the
family zero supremum is
\begin{equation}\label{eq:beta-star}
 \beta_*:=\sup\bigl(\{1/2\}\cup\{\Rea\rho:\rho\in\mathcal Z_F\}\bigr).
\end{equation}
Its analogue $\beta_*(F)$ for the ground fields of Parts~\ref{part:V}
and~\ref{part:VI} is defined in Section~\ref{sec:v-scope}.

For a prime $p\notin S$, the sextic residue character is
\[
 \chi_p(a)=(a/p)_6,\qquad
 \chi_p(a)\equiv a^{(q_p-1)/6}\pmod p\quad(p\nmid a),
 \qquad \chi_p(a)=0\quad(p\mid a).
\]
For a primary element $c$, meaning $c\equiv1\pmod3$, write
$c=\prod_p p^{v_p(c)}$ with primary prime generators $p$, and put
$\chi_c(a)=\prod_{p\mid c}\chi_p(a)^{v_p(c)}$.
Character powers, for every integer exponent (including zero, negative
exponents and multiples of six), are taken on invertible residue
classes and extended by zero:
\begin{equation}\label{eq:zeros}
 \chi_c(a)^6=1_{(a,c)=1}.
\end{equation}
Passing to a primitive inducing character therefore requires
retaining the coprimality indicator for the original defining modulus. Coefficients and
such conditions are independent of the averaged variable,
except for the displayed residue symbols.

The additive normalization is
\[
 e(z)=\exp\!\left(2\pi i\operatorname{Tr}_{F/\mathbb Q}(z/\lambda)\right),
 \qquad d\mu(z)=\frac2{\sqrt3}\,d\xi_1\,d\xi_2.
\]
Here $z=\xi_1+i\xi_2$ with $\xi_1,\xi_2\in\mathbb R$;
for $z\in\mathbb C$ the trace in $e(z)$ means
$z/\lambda+\overline{z/\lambda}$. The pairing $(z_1,z_2)\mapsto e(z_1z_2)$
makes $\OO$ self-dual. For squarefree $c$ define
\[
 \gamma_j(c)=q_c^{-1/2}\sum_{v\bmod c}\chi_c(v)^j e(v/c).
\]
Fixed-numerator symbols $n\mapsto\chi_n(u)$ are finite-order
ray characters \Src{Lemma}{4.1}. Prime Gauss identities and
sextic reciprocity provide root-of-unity-valued functions $\mathcal R,G$
of fixed ray classes:
$\mathcal R$ is the bicharacter on residue square
classes modulo $4$ of \Src{Lemma}{4.4}, and on coprime
primary pairs
\[
 \chi_b(a)=\mathcal R(a,b)\chi_a(b).
\]
The same lemma defines $G$ on all primary indices, with $|G|=1$,
depending only on the ray class of its argument in a fixed
finite ray quotient; $G(c)=\overline{\chi_c(4)}\gamma_3(c)$ for
squarefree $c$, while on
nonsquarefree indices $G$ is that finite-ray function, not a
Gauss sum, and
\[
 G(AB)=G(A)G(B)\mathcal R(A,B)
\]
on the admissible primary elements \Src{Lemmas}{4.2--4.4}.
They cancel the reciprocity factors linking distinct primes,
so the coefficients of the theta average of \Prev{Section}{4}
form an Euler product.

\section{Results from the moment-parameter extension}
\label{sec:prev-iv}

In this section, we quote, with locators and without proofs, results from
Parts~I--IV (chiefly Part~IV) of \cite{LiuArXiv} used in
Parts~\ref{part:V} and~\ref{part:VI}, notably the limiting bound
$\beta_*\le B_{\mathrm{new}}$ of Theorem~\ref{thm:iv-limit}, the extended
moment Proposition~\ref{prop:iv-plain}, the count
Proposition~\ref{prop:iv-count}, and the optimality
Theorem~\ref{thm:iv-system-optimality}. Here $F=\mathbb Q(\sqrt{-3})$ and
$\beta_*$ is \eqref{eq:beta-star}. The analytic statements (the modified
average, the detector and the principal term, the initial $7/8$ count and
the row bounds, the balanced and free-length direct estimates, the
extended product moment and the variable-parameter count,
Sections~\ref{sec:prev-average}--\ref{sec:iv-count}) serve as models for
the cyclotomic estimates of Part~\ref{part:V}. The scalar
statements (the exponent algebra such as \eqref{eq:general-extension-exponent},
\eqref{eq:prev-intermediate} and the endpoint saving of
Lemma~\ref{lem:prev-initial-endpoint}, the balanced geometry, the
endpoint and margin inequalities, the fixed-point identity, the iteration
and the optimality of the comparison system,
Sections~\ref{sec:iv-geometry}--\ref{sec:iv-further-parameters}) are
relations between explicit real functions, independent of the ground
field, and are used as stated.

Throughout this section, bracketed source locators such as
\LiuLoc{Theorem~27.1} and \LiuLoc{(18.4)} refer to that paper.
References to the present paper retain their usual local numbering.

Some letters have a second, local meaning in this section: $p(\xi)$ is the cubic \eqref{eq:iv-cubic} and $F(\kappa)$ the
boundary function of Section~\ref{sec:iv-count-obstruction}, while $p$ and
$F$ otherwise denote a prime and the ground field; $m$ is a detector length
except in Lemma~\ref{lem:prev-free-reflection}, where it is an outer range
exponent, and in the average over $m,v$, where it is the summation element;
$d$ is the norm exponent of $U=Z^d$, $d_{\mathcal I}$ the length of a subset
of rescaled factors, and $d\bmod v$ a residue class; $A$ is a height
exponent, the theta index in $\mathcal F_{v,A}$, or part of the normalizer
$A_{\mathcal T}$, while $A(y)$, $A_j$ and $A_0$ are polynomial coefficients or
local exponents; the bare letter $\zeta$ is a range width, whereas
$\zeta_F^S$ and $\zeta$ with an argument denote zeta functions; $M'$ is the
outer range exponent $1-\ell-2d_{\mathcal I}$ of the rescaled theta sums of
\LiuLoc{Section~16} in Section~\ref{sec:prev-free-direct}, and an enlarged
moment range, like $m_0+q_0$ in Proposition~\ref{prop:iv-plain}, in
Section~\ref{sec:iv-count-obstruction};
and $c=2/(3\kappa)$ is unrelated to the theta index $c$ in $cn^3$.

\subsection{The limiting bound and the comparison exponent}
\label{sec:prev-limit}
In this subsection, we recall the limiting bound
$\beta_*\le B_{\mathrm{new}}$ of Theorem~\ref{thm:iv-limit} and the
comparison exponent $E_B(d)$ of \eqref{eq:general-extension-exponent}, with its
balanced form \eqref{eq:balanced-extension-exponent}. The continuation
criterion \LiuLoc{Proposition~3.1} converts a strict comparison between the
direct estimate of the normalized theta average and its nonprincipal
contributions into a zero-free half-plane. Theorem~\ref{thm:iv-limit} is
the bound extended in Section~\ref{sec:v-scope}, and
\eqref{eq:general-extension-exponent} is the exponent of the contour
assembly of Section~\ref{sec:v-detail-30}.

Throughout,
\[
 \alpha=5/6,\quad\delta=2a-1,\quad x=\bar g/\delta\in[0,1/2],\quad
 b=l_y-l_x,\quad M=l_x+l_y,\quad h=1-l_x+\ell,
\]
where $\alpha$ is the amplification slope, $a$ the zero-class parameter of \LiuLoc{(5.1)},
$\bar g$ the weighted prime-sum mean of
\LiuLoc{(10.1)}, $l_x,l_y$ the averaging exponents and $\ell$
the total selected prime length. Except in Proposition~\ref{prop:iv-plain}, where it
is written $m_0+q_0$, the unprimed letter $M$ always denotes $l_x+l_y$.

\begin{theorem}[{\LiuLoc{Theorem~27.1}}]\label{thm:iv-limit}
Let
\begin{equation}\label{eq:iv-explicit-root}
 \xi_{\mathrm{new}}=\frac{20}{33}+\frac{\sqrt{1391}}{33}
 \cos\!\left(\frac13\arccos\frac{28470}{1391^{3/2}}
                         -\frac{2\pi}{3}\right),
\end{equation}
where $\arccos\in[0,\pi]$. Equivalently, $\xi_{\mathrm{new}}$ is the
unique root in $(0.38858,0.38859)$ of
$1188\xi^3-2160\xi^2+171\xi+190=0$. Then
\[
 \beta_*\le B_{\mathrm{new}}:=\frac{18\xi_{\mathrm{new}}+7}{18\xi_{\mathrm{new}}+9}
 =0.874957019420098946\ldots.
\]
Every finite-order Hecke $L$-function over $F$, and every Dirichlet
$L$-function, including $\zeta$, is zero-free for $\Rea s>B_{\mathrm{new}}$;
principal poles are allowed.
\end{theorem}

The comparison exponent is that of \LiuLoc{Proposition~18.2}. Let
$B\le\beta_*$ with $437/500\le B\le7/8$ \LiuLoc{(17.1)}, assume
positive lengths and the joint-error and moment hypotheses of
\LiuLoc{Sections~8--9}, and let $C(s)=s-5/6+l_x/3+\ell/6$ be the principal
exponent. For a class with element-norm exponent $d$, count exponent
$R$, zero-class parameter $a$ and prime-sum mean $\bar g$, and for each
prescribed loss $\eps>0$ and sufficiently small fixed contour and
zero-class widths, the normalized nonprincipal contribution is bounded by
$Z^{C(B)+E_B(d)+\eps}(1+T_1)^A$, where
\begin{equation}\label{eq:general-extension-exponent}
 E_B(d)=a-B+h(z_0-1/6)-a l_y-\ell/2+\bar g\ell
             +d(R+\delta/2-z_0),\quad z_0=17/50
\end{equation}
\LiuLoc{(18.4)}. For balanced lengths $l_x+l_y+\ell=1$ the
direct comparison gives $B=B(\ell)=11/12-\ell/4$
\LiuLoc{(15.2)}, and then
\begin{equation}\label{eq:balanced-extension-exponent}
 E_B(h)=-\frac14+\frac b6+\frac{5\ell}{4}
 +(1/2+\ell)\delta+\ell\bar g-h(1-R)
\end{equation}
\LiuLoc{(18.5)}.

\subsection{The modified average, the detector and the principal term}
\label{sec:prev-average}
In this subsection, we recall the objects of \LiuLoc{Sections~4--8}: the
excluded group $\Theta$, the modified average $I_{\eta,\rm mod}$ with its
exact Mellin representation, the principal term $f_\eta$ compared with
the normalized average $J_\eta$, and the zero classes of the detector. Their cyclotomic counterparts are constructed in
Sections~\ref{sec:v-detail-26}--\ref{sec:v-detail-29}.

\paragraph{\it Ray classes and the excluded group.}
Fix a finite ray class group $\mathcal T$ through which the reciprocity
function $\mathcal R$, the function $G$ and the fixed root-of-unity
factors of the primary-generator convention factor
\LiuLoc{Section~4.3}. Write $\widehat{\mathcal T}$ for its character
group and $1_{\mathcal T}$ for its identity class; $p\in1_{\mathcal T}$
means that the primary generator of $p$ has trivial image in
$\mathcal T$. The finite group of excluded Hecke characters is
$\Theta=\langle\eta,\widehat{\mathcal T}\rangle$. Let $b_*$ generate the
squarefree product of the primes in $S$. In Part~\ref{part:V}, $\Theta$
denotes the corresponding fixed finite group over $F$, specified where
it is used.

\paragraph{\it The modified average.}
The smoothed average $I_\eta(X,Y,Z)$ of \LiuLoc{(4.5)} is a
weighted average over $m,v$ of cubic theta sums of index $cn^3$, with
averaging lengths $X=Z^{l_x}$ and $Y=Z^{l_y}$ \LiuLoc{(7.1)}. For
$K$ disjoint prime windows with smooth weights $W_i$ and norms
$q_p\asymp P_i=Z^{\ell_i}$, $\sum_i\ell_i=\ell$, the allowed set
$\mathcal P_i$ consists of the primes $p\notin S$ with $p\in1_{\mathcal T}$
in the $i$th window.

The modified average \LiuLoc{(7.2)} is
\[
 I_{\eta,\rm mod}(Z)=\sum_{(p_i)}\prod_iW_i(q_{p_i}/P_i)
 \sum_{J}(-1)^{|J|}q_{p_J}^{-3/2}\overline{\eta(p_{J^c})}
 I_{\eta;p_{J^c}}(X/q_{p_J},Y/q_{p_J},Zq_{p_{J^c}}),
\]
where $I_{\eta;D}$ inserts the condition $D\mid cn^3$, the tuple sums
have $p_i\in\mathcal P_i$, and $p_J=\prod_{i\in J}p_i$.

Poisson summation
\LiuLoc{(4.9)} expresses the unmodified average $I_\eta$ as a
triple Mellin integral over sixth-power-free frequencies $u$.

Its integrand contains the scalar quotient
$\zeta_F^S(6z)L_F^S(w,\chi_\bullet(u))/L_F^S(s,\eta\overline{\chi_\bullet(u)})$
and a holomorphic correction $\mathcal H_{\eta,u}=\prod_{p\notin S}H_p$.
Here $L_F^S(s,\eta)=L_F(s,\eta)\prod_{p\in S}(1-\eta(p)q_p^{-s})$ is the
$L$-function with the Euler factors at $S$ removed \LiuLoc{Section~3},
$\zeta_F^S=L_F^S(\cdot,1)$, $\chi_\bullet(u)$ is the character
$n\mapsto\chi_n(u)$, and $H_p=P_p(1-V_p)(1-W_p)/(1-D_p)$ is the local
factor of \LiuLoc{(4.10)}, obtained from the complete local
valuation series $P_p$ by removing the scalar Euler factors
$V_p=q_p^{-6z}$, $W_p=\chi_p(u)q_p^{-w}$ and
$D_p=\eta(p)\overline{\chi_p(u)}q_p^{-s}$. For the modified average the
correction is replaced by
\[
 \mathfrak H_{\eta,u,Z}=\sum_{(p_i)}\prod_i\bigl[W_i(q_{p_i}/P_i)
 q_{p_i}^{z-1}G_{p_i}\bigr]\prod_{p\notin S\cup\{p_i\}}H_p
\]
\LiuLoc{(8.2)}, where $G_p$ is the selected local factor of
\LiuLoc{(8.1)}; replacing $\mathcal H$ by $\mathfrak H$ in
\LiuLoc{(4.9)} gives the exact modified average.

The finite Fourier coefficients $\mathcal F_{v,A}(H)$ of its arithmetic
coefficients \LiuLoc{Section~4.4} are sums over $d\bmod v$ of Gauss sums
modulo $b_*A$, so the direct bound from these sums is
$|\mathcal F_{v,A}(H)|\ll q_vq_A$.

\paragraph{\it The principal term.}
Only $u=1$ gives a principal numerator, with poles at $w=1$, $z=1/6$.
The principal correction is
\[
 H_\eta(s)=\mathcal H_{\eta,1}(s,1,1/6)
\]
\LiuLoc{(4.13)}; by \LiuLoc{(4.12) and (4.14)} it is
holomorphic and within $1/2$ of $1$ for
$\Rea s>7/8$ after a fixed choice of $P_0$. With the principal exponent
$C(s)$ of Section~\ref{sec:prev-limit}, the normalized average
$J_\eta=I_{\eta,\rm mod}/(c_SA_{\mathcal T})$ \LiuLoc{(8.5)} is
compared with
\[
 f_\eta(Z)=\frac1{2\pi i}\int_{(2)}Z^{C(s)}e^{(s-5/6)^2}
 \frac{H_\eta(s)}{L_F^S(s,\eta)}\,ds .
\]

In the proof of \LiuLoc{Proposition~17.2}, the principal term is
extracted by moving $w$ to $\Rea w=19/20$ and then, in its residue, $z$ to
$\Rea z=33/200$. Since the raw exponent is
$l_x/2+\Rea s-1+h\Rea z+l_y(\Rea w-1)$, the two unresidued integrals save
$l_y/20$ and $h(1/6-33/200)=h/600$ \LiuLoc{(17.7)}.

Here
$c_S=\widehat W_1(1)\mathcal M(1/6)(\Res_{v=1}\zeta_F^S(v))^2/6>0$ is the
scalar residue at $w=1$, $z=1/6$ \LiuLoc{(4.13)}; in it, $W_1$ is
the fixed test function restricting the averaging variable $v$ in
\LiuLoc{(4.5)} (not the prime-window weight $W_i$ with $i=1$),
$\widehat W_1(w)=\int_0^\infty W_1(r)r^{w-1}\,dr$, and $\mathcal M$ is the
Mellin transform of the Fourier transform of the $m$-weight
\LiuLoc{Section~4.4}. Further,
$A_{\mathcal T}(Z)=(-1)^KZ^{-\ell/6}\prod_i\sum_{p\in\mathcal P_i}
W_i(q_p/P_i)q_p^{-5/6}$ is the prime normalizer, with
$|A_{\mathcal T}|^{-1}\ll_\eps Z^\eps$ \LiuLoc{(8.4)--(8.5)}.

\paragraph{\it Zero classes and the detector.}
For sixth-power-free $u$ with $q_u\asymp U=Z^d$ the detector of
\LiuLoc{Section~5} uses the fixed lowest threshold $a_0=51/100$ and a fixed
width $\eb$. By pigeonhole, each such $u$ lies in a
zero class, determined by the real parts of the zeros of its associated
characters, with parameter $a\in a_0+\eb\mathbb Z_{\ge0}$
\LiuLoc{(5.1)}; the lowest class has $\delta=2a_0-1=1/50$.

A selected zero forces simultaneous large values
\LiuLoc{(5.5)} of two detecting character polynomials $M_r$ and
$S_m$ of length exponents $r,m$, with $r+m\ge t$ for the detecting
cutoff $t\in[1,3/2]$; the height cutoff is $T_1=Z^\tau\le U^{1/100}$. For
a fixed zero class and prime-sum size class, $\mathcal B$ denotes the set
of detected elements $u$ with $q_u\asymp U$. The prime-sum size classes have
exponents $g_i$ with weighted mean $\bar g$ \LiuLoc{(10.1)}.

\subsection{The initial \texorpdfstring{$7/8$}{7/8} comparison}
\label{sec:prev-initial}
In this subsection, we recall the $7/8$ comparison of
\LiuLoc{Sections~7--12}: the count Proposition~\ref{prop:prev-initial-count},
the endpoint saving $-E_*\ge49/440640$ of
Lemma~\ref{lem:prev-initial-endpoint}, and the remaining norm ranges of
\LiuLoc{Section~11.3} and \LiuLoc{Lemma~19.2}, including the intermediate
bound \eqref{eq:prev-intermediate}. The first application of the
continuation criterion (Section~\ref{sec:v-detail-31}) reproduces this
comparison, and Proposition~\ref{prop:prev-initial-count} also enters the
counting estimates of Section~\ref{sec:v-detail-29}.

The initial comparison has geometry
\begin{equation}\label{eq:prev-initial-geometry}
 (b,\ell,l_x,l_y,h)=(1/8,1/6,17/48,23/48,13/16),\qquad
 C(s)=s-11/16
\end{equation}
\LiuLoc{(7.1), (8.7)}, so $C(7/8)=3/16$ and $\ell/h=8/39>1/5$.
Suppose $\Delta=\beta_*-7/8>0$ and take $\kappa=3/4+2\Delta$.

\begin{proposition}[{\LiuLoc{Proposition~10.1}}]\label{prop:prev-initial-count}
Assume $\Delta\in(0,1/24]$ and $\kappa=3/4+2\Delta$. Let $\mathcal B$
consist of the elements in one zero class and one prime-sum size class,
with $q_u\asymp U=Z^d$, $a>a_0$, $0<\delta\le5/6$,
$0\le\bar g\le\delta/2$ and $x=\bar g/\delta$. Assume the detector height
hypotheses and that the available prime length $\ell/d$ exceeds $7/37$
by a fixed positive amount. Put
$D_x=3-17x/9$ and $P_x=(2-8x/9)(1-x)$. For $t\in[1,3/2]$ let
\[
 R_{\rm short}(t)=1-\delta+\frac{\delta P_x}{D_x}(3/2-t),\qquad
 R_{\rm long}(t)=1-\delta+(5/6-\delta)(t-1).
\]
Then, for every $\eps>0$, the length decrements and the maximum individual
prime length can be chosen so that
$\#\mathcal B\ll U^{\max(R_{\rm short}(t),R_{\rm long}(t))+\Delta/4+\eps}(1+T_1)^A$.
With no selected prime sums, $t=1$ gives
$\#\mathcal B\ll U^{1-2\delta/3+\eps}(1+T_1)^A$ without any hypothesis on
the available prime length.
\end{proposition}

In its proof the two short counts cross at
$r_*(t)=[(2-8x/9)t-5x/9]/D_x$, the case $\kappa=3/4$ of
\eqref{eq:iv-crossing}, and $1/3\le t-r_*(t)\le1/2$, so $m\ge1/3-O(\eps)$ on that branch. Restoring
the actual unweighted prime length at $\kappa=3/4+2\Delta$ increases the
counting exponent by at most
\[
 2\bar g(1-2m)\Bigl(\frac29-\frac1{9/2+12\Delta}\Bigr)
 =\frac{24\bar g(1-2m)\Delta}{(9/2)(9/2+12\Delta)}
 \le\frac{40}{243}\Delta+O(\eps)<\frac\Delta4+O(\eps) .
\]
The count without selected prime sums comes from balancing $1-\delta r$
and $1-2\delta(1-r)$ at $r=2/3$.

\begin{lemma}[{\LiuLoc{Lemma~11.1}}]\label{lem:prev-initial-endpoint}
Let $1/50<\delta<5/6$, $0\le x\le1/2$, choose $t$ by balancing
$R_{\rm short}=R_{\rm long}=R_*$, and let $E_*=E(h)$ be the exponent
\eqref{eq:general-extension-exponent} relative to $C(7/8)$ at the
geometry \eqref{eq:prev-initial-geometry}, with $R=R_*$. Then
$-E_*\ge49/440640$. Let $E_{\rm actual}(h)$ be the resulting exponent
relative to $C(7/8)$ including the count loss $\Delta/4$, and let
$\eps_{\rm total}$ be the sum of the remaining contour, length,
size-partition and smoothing losses. Since
$C(\beta_*)-C(7/8)=\Delta$,
\[
 E_{\rm actual}(h)-\Delta\le-\frac{49}{440640}-\frac{51}{64}\Delta
 +\eps_{\rm total},\qquad \frac{51}{64}=1-\frac h4
\]
\LiuLoc{(11.3)}.
\end{lemma}

In the discussion following \LiuLoc{Lemma~11.1}, the lowest zero class
$a=a_0$, with the direct count $R=1$, gives $E(h)\le-7/1200$.

The remaining norm ranges at the initial geometry are those of
\LiuLoc{Section~11.3}, with $d_{\min}=1/100$ and $0<\zeta<1/48$:
the small range $d<d_{\min}$ saves $63/800$; the intermediate range
$d_{\min}\le d\le1/2$ uses $R=76/75-2\delta/3$ and gives
$E(d)\le E(1/2)\le-49/14400$; the moderate range $1/2\le d\le h$ has
slope $R+\delta/2-z_0\in[4/25,2]$ in $d$, so the endpoint
$d=h$ controls it; the range $h<d\le h+\zeta$ costs at most $2\zeta$;
and $d>h+\zeta$ is handled by absolute convergence on
$\Rea z=z_{\mathrm{new}}$ large. These savings are positive, independent
of the target, and assembled in \LiuLoc{(12.1)--(12.3)}.

For general balanced lengths these ranges are treated by
\LiuLoc{Lemma~19.2}. Assume the balanced hypotheses \LiuLoc{(15.3)},
$B\le\beta_*\le7/8$ and $0\le\delta\le5/6$, and let $d_{\min}=1/100$; the
following hold with arbitrarily small fixed losses:
\begin{enumerate}
\item the small range $q_u\ll Z^{d_{\min}}$ has normalized exponent,
relative to $C(\beta_*)$, at most $13h/75-l_y/2+2d_{\min}$;
\item for $d_{\min}\le d\le1/2$ and $\delta\ge1/50$, the intermediate range uses $R=76/75-2\delta/3$, the count
without selected prime factors with an extra $1/75$ to cover the lowest
class, and its slope $101/150-\delta/6$ in $d$ is positive;
\item the moderate range has slope at least $4/25$ and at most $2$;
\item the lowest zero class has $R=1$, $\bar g\le1/100$ and endpoint
exponent at most $-1/4+b/6+5\ell/4+(1/2+3\ell/2)/50$;
\item extending the moderate range to $h<d\le h+\zeta$ increases the
exponent by at most $2\zeta$, provided $\ell/(h+\zeta)>7/37$;
\item the far range $q_u>Z^{h+\zeta}$, on $(\Rea s,\Rea w,\Rea z)=(2,2,z_{\mathrm{new}})$,
has exponent
$l_x/2+1+l_y+(h+\zeta)(1+\eps)-\zeta z_{\mathrm{new}}$.
\end{enumerate}
Substituting $R=76/75-2\delta/3$, $a=(1+\delta)/2$, $\bar g\le\delta/2$,
$B=11/12-\ell/4$ and $d=1/2$ in \eqref{eq:general-extension-exponent}
gives the intermediate bound
\begin{equation}\label{eq:prev-intermediate}
 E_B(1/2)\le\frac{(225\ell+50-75b)\delta+78\ell-49b-73}{300}.
\end{equation}

\subsection{The direct estimate with free total length}
\label{sec:prev-free-direct}
In this subsection, we recall the direct estimate
\eqref{eq:prev-balanced-low} of \LiuLoc{Proposition~16.2}, with exponent
$l_x/2+b/12$ for balanced lengths, and its extension to free
$M=l_x+l_y$: the reflection bound Lemma~\ref{lem:prev-free-reflection}
and the direct majorant $B_{\rm low}$ of Proposition~\ref{prop:prev-free-low},
\eqref{eq:prev-B-low}. The balanced estimate is the model for the direct
bound of Section~\ref{sec:v-detail-27}; Lemma~\ref{lem:prev-free-reflection}
and Proposition~\ref{prop:prev-free-low} enter the estimate system of
Section~\ref{sec:v-detail-33}.

For balanced lengths the direct estimate is \LiuLoc{Proposition~16.2}:
if $l_x+l_y+\ell=1$, $b=l_y-l_x>0$, $l_x>\ell>0$, $l_y>\ell$,
$l_y-\ell-11b/6>0$, $1-3\ell>0$ and $\ell<1/5$, then
\begin{equation}\label{eq:prev-balanced-low}
 |I_{\eta,\rm mod}(Z)|\ll_{\eta,\eps}Z^{l_x/2+b/12+\eps}
\end{equation}
\LiuLoc{(16.2)}.

Its proof applies the Cauchy--Schwarz separation
of the average into an additive and a theta factor, the additive
estimate with parameter $P_a=q_{b_*}^{-1}Z^b$ (where $l_y-\ell-11b/6>0$
gives $P_a^2/Y'\ll P_a^{1/6}$), and the functional-equation estimate
\LiuLoc{Lemma~16.1}: for a subset of rescaled factors of total length
$d_{\mathcal I}$, with $M'=1-\ell-2d_{\mathcal I}$, the rescaled theta sums
$B^{\mathcal I}_{a,\sigma}$ have squared norm
$\ll Z^{M'+\pos{d_{\mathcal I}+5\ell-1}/4+\eps}$ \LiuLoc{(16.1)}.
The net exponent of a subset is at most $-7d_{\mathcal I}/8$.

The optimality argument allows $M=l_x+l_y$ to vary independently of
$\ell$. Assume $l_x,l_y>\ell$, $b>0$ and $h>0$, and use
$C(s)=s-5/6+l_x/3+\ell/6$ \LiuLoc{(22.1)}.

\begin{lemma}[{\LiuLoc{Lemma~22.1}}]\label{lem:prev-free-reflection}
Under the common-weight and coefficient hypotheses of \Src{Lemma}{14.3},
the rescaled theta sums $B^{\mathcal I}_{a,\sigma}$ of \LiuLoc{Section~16},
with outer element-range exponent $m$ and
remaining prime-length exponent $l$ have squared norm
$\ll Z^{m+k(m,l)+\eps}$, where
\[
 k(m,l)=\max\Bigl\{0,\frac{1+3l-2m}{4},m+l-1,\frac{3l-1}{2}\Bigr\}.
\]
\end{lemma}

The proof treats four cases in the local length variables of the
functional-equation estimate: the exponents $O$ and $H$ of the norms of the
powerful and squarefree parts of the outer element, the exponents
$S_0,B_0,A_0,N_0$ of the local divisor decompositions, the retained prime length $z_a\le l$, the
retained dual length exponents $v,\ell_b,e_\lambda$ and the available dual
length exponent $T_d$, with $v+3\ell_b+e_\lambda\le T_d$ and the term
$\max(H,v+\ell_b)$ retained. In the balanced setting of
\LiuLoc{Lemma~16.1}, where $m+l-1=-3d_{\mathcal I}\le0$ for a subset of
rescaled factors of total length $d_{\mathcal I}$, one has
$T_d\le H-3d_{\mathcal I}\le H$, so the maximum is $H$.

\begin{proposition}[{\LiuLoc{Proposition~22.2}}]\label{prop:prev-free-low}
The modified average satisfies
$|I_{\eta,\rm mod}(Z)|\ll Z^{L_{\rm low}+\eps}$ with
\[
 L_{\rm low}=\frac{l_x}{2}+\max\{0,b/12,b-l_y/2\}+\frac12k(M,\ell),
\]
and the direct comparison $L_{\rm low}\le C(B)$ is equivalent to
$B\ge B_{\rm low}$, where
\begin{equation}\label{eq:prev-B-low}
 B_{\rm low}=L_{\rm low}+5/6-l_x/3-\ell/6
\end{equation}
\LiuLoc{(22.4)}. For $M+\ell=1$,
$B_{\rm low}=11/12-\ell/4+\pos{8b+3\ell-3}/12+\pos{5\ell-1}/8$.
\end{proposition}

In the proof, a subset of rescaled factors of total length $d_{\mathcal I}$
changes the exponent by at most $-3d_{\mathcal I}/8$, so the empty subset controls the
direct exponent.

\subsection{The extended unweighted product moment}
\label{sec:prev-moment}
In this subsection, we recall the analytic input of \LiuLoc{Part~IV}:
Proposition~\ref{prop:iv-plain} extends the
unweighted fourth moment \Src{Lemma}{18.1} from $3/4\le\kappa\le1$ to
$13/18\le\kappa\le1$. It is the model for the unweighted product
moment of Section~\ref{v:sec:M2-assembly}. Its proof extends the argument of \Src{Section}{18},
rather than its printed parameter range.

\begin{proposition}[{\LiuLoc{Proposition~28.1}}]\label{prop:iv-plain}
Retain the notation, character families, fixed twists, common
coprimality conditions and bounded length ranges of the unweighted
moment \LiuLoc{(9.4)}. Let $n_1,n_2$ be the ideal-length
exponents, $z_{\rm p}$ the total prime-length exponent, and $M=m_0+q_0$
the combined range exponent of the moment estimate, where $m_0$ bounds the
element range and $q_0$ bounds the radical of its variable residue-symbol
factors. In this proposition only, $M$ is not $l_x+l_y$. Let
$\Theta=\langle\eta,\widehat{\mathcal T}\rangle$ be the finite group of
excluded Hecke characters of \LiuLoc{Section~4.3}, and assume that the
prime polynomials have disjoint underlying prime supports. Then
\LiuLoc{(9.4)} holds uniformly for $13/18\le\kappa\le1$, in the
cases
\[
 \begin{array}{c|c|c}
 &\text{length condition}&\text{excluded inducing characters}\\ \hline
 z_{\rm p}=0&n_1,n_2\text{ arbitrary bounded}&\text{principal}\\
 z_{\rm p}>0&n_1+n_2+6\kappa z_{\rm p}\le M&\Theta.
 \end{array}
\]
If $z_{\rm p}>0$ and $\kappa<1$, assume $\beta_*\le(1+\kappa)/2$.
The upper bound $\ell_{\rm mesh}$ on the individual prime lengths depends
only on the bounded real ranges and $\eps$, not on the number of
prime-polynomial factors. Constants and finite smooth-seminorm and polynomial-height orders may
depend on the fixed number of prime-polynomial factors and the
arithmetic data; they are uniform over the character families, moduli,
common coprimality restrictions and fixed outer summation parameters.
\end{proposition}

The cyclotomic counterpart of this estimate is Theorem~\ref{v:thm:M2-assembled}.

\subsection{The variable-parameter count}
\label{sec:iv-count}
In this subsection, we recall the variable-parameter count
\LiuLoc{Proposition~29.1}, Proposition~\ref{prop:iv-count}, which gives
$\#\mathcal B\ll U^{R_\kappa+\eps}(1+T_1)^A$ with $R_\kappa$ of
\eqref{eq:iv-balanced-count}, together with the crossing
\eqref{eq:iv-crossing} from its proof. It is the model for the counting
estimates of Section~\ref{sec:v-detail-29}.

The count combines the extended moment of
Section~\ref{sec:prev-moment} with the inverse moment
\LiuLoc{(9.1)}, its amplified form \LiuLoc{(9.3)} and
the simultaneous large values \LiuLoc{(5.5)} of the zero
detector of \LiuLoc{Section~5}. The detecting cutoff is $t$, and $r,m$
are the norm-length exponents of the two detecting polynomials. Put
\[
 c=\frac{2}{3\kappa},\qquad
 D_{\kappa,x}=3-(1+c)x,\qquad
 P_{\kappa,x}=(2-cx)(1-x)
\]
\LiuLoc{(29.1)}.

\begin{proposition}[{\LiuLoc{Proposition~29.1}}]\label{prop:iv-count}
Retain the detector, character-family, coefficient, smooth-weight and
height hypotheses of \LiuLoc{Proposition~18.1}, with the simultaneous
large values for one common character. Assume
\[
 13/18\le\kappa\le1,\qquad
 \beta_*\le(1+\kappa)/2,\qquad
 \ell/d\ge1/5+\rho_{\rm sup},\quad\rho_{\rm sup}>0.
\]
For $0<\delta\le\alpha=5/6$ and $0\le x\le1/2$, put
\begin{equation}\label{eq:iv-balanced-count}
 \mathcal J_\kappa=(\alpha-\delta)D_{\kappa,x}
                         +\delta P_{\kappa,x},\quad
 t_\kappa=1+\frac{\delta P_{\kappa,x}}{2\mathcal J_\kappa},\quad
 R_\kappa=1-\delta+
 \frac{(\alpha-\delta)\delta P_{\kappa,x}}{2\mathcal J_\kappa}.
\end{equation}
For each $\eps>0$, sufficiently fine fixed prime lengths give
\[
 \#\mathcal B\ll U^{R_\kappa+\eps}(1+T_1)^A,
\]
where $\mathcal B$ is the set of detected elements of the class,
$U=Z^d$, $T_1$ is the zero-detection height cutoff and $A\ge0$ is a
fixed finite height exponent.
\end{proposition}

The proof of \LiuLoc{Proposition~29.1} uses the inverse moment for
$r\ge r_\kappa$ and the unweighted moment for $r\le r_\kappa$, where,
with $m=t-r$, the maximal selected prime-length exponents are
$z_M=(1-r)/2$ and $z_S=(1-2m)/(6\kappa)$, the savings divided by
$\delta$ are $F_I(r)=x+(1-x)r$ and
$F_S(r)=2(t-r)+(c/2)x[1-2(t-r)]$, and their crossing is
\begin{equation}\label{eq:iv-crossing}
 r_\kappa(t,x)=
 \frac{(2-cx)t+(c/2-1)x}{3-(1+c)x}
\end{equation}
\LiuLoc{(29.3)}. For $1\le t\le3/2$ and $0\le x\le1/2$ it
satisfies
\[
\begin{gathered}
 \partial_t r_\kappa=\frac{2-cx}{D_{\kappa,x}}>0,\qquad
 \partial_x r_\kappa=\frac{(2-c)(2t-3)}{2D_{\kappa,x}^2}\le0,\\
 r_\kappa(3/2,x)=1,\qquad
 \frac{6-c}{10-2c}\le r_\kappa\le1,\qquad
 \frac13\le t-r_\kappa\le\frac12,
\end{gathered}
\]
and the admissible prime lengths satisfy
$z_M\le(4-c)/(20-4c)<1/5$ and $z_S\le c/12<1/5$. The short and long
counts
\[
 R_{\kappa,\rm short}(t)=1-\delta+
 \frac{\delta P_{\kappa,x}}{D_{\kappa,x}}(3/2-t),\qquad
 R_{\rm long}(t)=1-\alpha+(\alpha-\delta)t
\]
balance at $t=t_\kappa\in[1,3/2]$ with common value $R_\kappa$.

\subsection{Balanced lengths and the critical class}
\label{sec:iv-geometry}
In this subsection, we recall the one-parameter family
\eqref{eq:iv-geometry} of \LiuLoc{(30.2)}, whose averaging lengths
balance the direct estimate \LiuLoc{Proposition~16.2} with the critical
class of Proposition~\ref{prop:iv-count}, and Lemma~\ref{lem:iv-geometry}:
at $(\delta,x)=(\xi,1/2)$ one has $R_\kappa=2/3$ and $E_{B(\xi)}(h)=0$. The
family \eqref{eq:iv-geometry} is the parameter choice of
Sections~\ref{sec:v-detail-27} and~\ref{sec:v-detail-31}.

For $\xi$ in the interval
$I=[38858/100000,38859/100000]$ of \LiuLoc{(30.1)}, define
\begin{equation}\label{eq:iv-geometry}
\begin{gathered}
 \kappa(\xi)=\frac{18\xi^2-57\xi+20}
                       {6(36\xi^2-78\xi+25)},\qquad
 \ell(\xi)=\frac{5-6\xi}{9+18\xi},\qquad
 B(\xi)=\frac{18\xi+7}{18\xi+9},\\
 b(\xi)=\frac{8(3\xi-1)(12\xi^2-40\xi+15)}
                    {9(2\xi+1)(24\xi^2-20\xi+5)},\\
 l_x(\xi)=\frac{1-b(\xi)-\ell(\xi)}2,\qquad
 l_y(\xi)=\frac{1+b(\xi)-\ell(\xi)}2,\qquad
 h(\xi)=\frac{1+b(\xi)+3\ell(\xi)}2,\\
 C_\xi(s)=s-\frac{4+b(\xi)}6,
 \qquad B(\xi)=\frac{11}{12}-\frac{\ell(\xi)}4
\end{gathered}
\end{equation}
\LiuLoc{(30.2)}. At fixed $\xi$, write $b,\ell,\kappa,h$ for
$b(\xi),\ell(\xi),\kappa(\xi),h(\xi)$. For $\xi_0=(49-\sqrt{921})/48$,
these are the parameters of \LiuLoc{Theorem~21.1} and
$\kappa(\xi_0)=3/4$; numerically $\xi_0=0.3885837\ldots\in I$.

\begin{lemma}[{\LiuLoc{Lemma~30.1}}]\label{lem:iv-geometry}
For the length parameters \eqref{eq:iv-geometry}, at
$(\delta,x)=(\xi,1/2)$,
\[
 R_\kappa=2/3,\qquad E_{B(\xi)}(h)=0,\qquad
 \partial_\delta E_{B(\xi)}(h)=0.
\]
The equation $R_\kappa(\delta,1/2)=2/3$ is equivalent to
\begin{equation}\label{eq:iv-critical-polynomial}
 (216\kappa-18)\delta^2-(468\kappa-57)\delta
                           +150\kappa-20=0,
\end{equation}
and the solution $\delta=\xi\in I$ is its smaller root.
\end{lemma}

\subsection{The endpoint and margin inequalities}
\label{sec:iv-certificate}
In this subsection, we recall the endpoint inequality
$E_{B(\xi)}(h)\le0$ of Lemma~\ref{lem:iv-endpoint}, proved uniformly along
the family \eqref{eq:iv-geometry} by the exact quadratic identity
\eqref{eq:iv-Q}, and the margin inequalities of Lemma~\ref{lem:iv-reserves},
proved by Bernstein-basis bounds on $I$. Both lemmas are used in the
variable scalar comparison of Section~\ref{sec:v-detail-31}, and
Lemma~\ref{lem:iv-reserves} also in Section~\ref{sec:v-detail-27}.

\begin{lemma}[{\LiuLoc{Lemma~31.1}}]\label{lem:iv-endpoint}
For every $\xi\in I$, use \eqref{eq:iv-geometry} and the count exponent
$R=R_{\kappa(\xi)}$ of \eqref{eq:iv-balanced-count}. Then
\begin{equation}\label{eq:iv-endpoint}
 E_{B(\xi)}(h)\le0
 \qquad(0\le\delta\le5/6,\quad0\le x\le1/2).
\end{equation}
Equality holds only at $(\delta,x)=(\xi,1/2)$.
\end{lemma}

The proof of \LiuLoc{Lemma~31.1} is the following verification. Put
$y=1/2-x$, $c=2/(3\kappa(\xi))$,
\[
 D=(5-c)/2+(1+c)y,\qquad P=1-c/4+2y+cy^2,
 \qquad \mathcal J=(\alpha-\delta)D+\delta P,
\]
and $C_0=-1/4+b/6+5\ell/4$. By \eqref{eq:balanced-extension-exponent},
$E=E_{B(\xi)}(h)$ satisfies
\begin{equation}\label{eq:iv-Q}
 Q=-2\mathcal J E
 =-2C_0\mathcal J+(b+2\ell y)\delta\mathcal J
                    -h(\alpha-\delta)\delta P
\end{equation}
\LiuLoc{(31.2)}. Write $Q=A(y)\delta^2+B_q(y)\delta+C_q(y)$ and
$\mathscr D(y)=4A(y)C_q(y)-B_q(y)^2=\sum_{j=1}^4d_jy^j$, the constant
coefficient being identically zero. The coefficients $A_j$ of $A(y)$
and $d_j$ of $\mathscr D(y)$ are explicit rational functions of $\xi$,
and the exact Bernstein-basis bounds on $I$ give
\begin{center}
\begin{tabular}{c|c@{\qquad}c|c}\toprule
quantity&lower bound&quantity&lower bound\\\midrule
$A_0$&$47/100$&$d_1$&$41/1000$\\
$A_1$&$121/100$&$d_2$&$198/1000$\\
$A_2$&$86/100$&$d_3$&$267/1000$\\
$A_3$&$29/100$&$d_4$&$76/1000$\\\bottomrule
\end{tabular}
\end{center}
Therefore
\[
\begin{gathered}
 A(y)>0\quad(y\ge0),\qquad
 \mathscr D(y)>0\quad(y>0),\\
 4A(y)Q=(2A(y)\delta+B_q(y))^2+\mathscr D(y)\ge0,\qquad
 \mathcal J>0\quad\Longrightarrow\quad E\le0,\\
 Q(\delta,0)=A_0(\delta-\xi)^2.
\end{gathered}
\]

\begin{lemma}[{\LiuLoc{Lemma~31.2}}]\label{lem:iv-reserves}
All parameter choices in \eqref{eq:iv-geometry} satisfy
\[
\begin{gathered}
 b>0,\quad0<\ell<1/5,\quad l_x-b-\ell>0,
 \quad l_y-\ell>11b/6,\\
 h<1,\qquad5\ell>h,\qquad437/500<B<7/8.
\end{gathered}
\]
The estimates for the lowest zero class and the intermediate and small
norm ranges of \LiuLoc{Lemma~19.2} have uniform positive exponent
savings; the detecting and rounded count formulas have uniformly
positive denominators.
\end{lemma}

The proof of \LiuLoc{Lemma~31.2} gives the explicit lower bound
$1/25000$ on $I$ for each expression in the list of margins
\LiuLoc{(31.3)}, which includes $\kappa-13/18$, $4/5-\kappa$,
$\kappa'-1000B'$ and $-p'$, where primes denote derivatives in $\xi$ of
the functions \eqref{eq:iv-geometry} and $p$ is the cubic \eqref{eq:iv-cubic} below,
together with the side margins for the remaining ranges of
Section~\ref{sec:prev-initial}:
\[
\begin{gathered}
 -C_0-(1/2+3\ell/2)/50,\qquad l_y/2-13h/75-1/50-(7/8-B),\\
 -\bigl[(225\ell+50-75b)\delta+78\ell-49b-73\bigr]/300
 \quad(\delta=1/50,\ 19/25),
\end{gathered}
\]
with $C_0$ as above: the savings of the lowest
zero class and of the small range, and (being
affine in $\delta$) a bound for the intermediate exponent
\eqref{eq:prev-intermediate} on $1/50\le\delta\le19/25$.

\subsection{The one-step bound and its iteration}
\label{sec:iv-iteration}
In this subsection, we recall the one-step bound
Proposition~\ref{prop:iv-bootstrap},
$\beta_*\le(1+\kappa(\xi))/2\Rightarrow\beta_*\le B(\xi)$, and its iteration
\LiuLoc{Section~33}, which converges to $B_{\mathrm{new}}$ through the
fixed-point identity \eqref{eq:iv-fixed-point-identity} and so proves
Theorem~\ref{thm:iv-limit}. The zero-free hypotheses pass from step to
step in the same way in Section~\ref{v:sec:M2-assembly}, and the
fixed-point identity enters the supremum argument of
Section~\ref{v:sec:compensated-assembly}.

\begin{proposition}[{\LiuLoc{Proposition~32.1}}]\label{prop:iv-bootstrap}
For $\xi\in I$ and $\kappa(\xi)\le3/4$,
\[
 \beta_*\le\frac{1+\kappa(\xi)}2
 \quad\Longrightarrow\quad \beta_*\le B(\xi).
\]
\end{proposition}

The iteration of \LiuLoc{Section~33} repeats this implication. Put
\begin{equation}\label{eq:iv-cubic}
 p(\xi)=1188\xi^3-2160\xi^2+171\xi+190
\end{equation}
\LiuLoc{(33.1)}. It has a unique zero $\xi_{\mathrm{new}}\in I$,
given by \eqref{eq:iv-explicit-root} \LiuLoc{Lemma~33.1}, and $p'<0$ on
$I$. Define
\begin{equation}\label{eq:iv-final-parameters}
 \kappa_{\mathrm{new}}=\kappa(\xi_{\mathrm{new}}),\quad
 B_{\mathrm{new}}=B(\xi_{\mathrm{new}}),\quad b_{\mathrm{new}}=b(\xi_{\mathrm{new}}),\quad
 \ell_{\mathrm{new}}=\ell(\xi_{\mathrm{new}})
\end{equation}
\LiuLoc{(33.2)}. The identity
\begin{equation}\label{eq:iv-fixed-point-identity}
 \kappa(\xi)-(2B(\xi)-1)
 =-\frac{p(\xi)}{18(2\xi+1)(36\xi^2-78\xi+25)}
\end{equation}
\LiuLoc{(33.3)} gives $\kappa_{\mathrm{new}}=2B_{\mathrm{new}}-1$.

On $I$ one has $\kappa'>0$, $B'>0$ and $\kappa'-1000B'>0$. Let
$\xi(\kappa)$ be the inverse of $\xi\mapsto\kappa(\xi)$ on $I$ and
$g(\kappa)=2B(\xi(\kappa))-1$. Then $0<g'<1/500$ and, by
\eqref{eq:iv-fixed-point-identity},
\[
 \kappa_{\mathrm{new}}<\kappa\le3/4
 \quad\Longrightarrow\quad
 \kappa_{\mathrm{new}}<g(\kappa)<\kappa.
\]

Starting from the $7/8$ bound of \cite{QRH}, the proof of
\LiuLoc{Theorem~27.1} sets $\kappa_0=3/4$, $\xi_n=\xi(\kappa_n)$,
$B_n=B(\xi_n)$ and $\kappa_{n+1}=g(\kappa_n)=2B_n-1$.
Proposition~\ref{prop:iv-bootstrap} gives
$\beta_*\le(1+\kappa_n)/2\Rightarrow\beta_*\le B_n=(1+\kappa_{n+1})/2$, so
\[
 \kappa_n\downarrow\kappa_{\mathrm{new}},\qquad
 B_n\downarrow B_{\mathrm{new}},\qquad
 \beta_*\le\inf_{n\ge0}B_n=B_{\mathrm{new}}.
\]

At each step the preceding bound supplies the zero-free
hypothesis of the moment estimate; only the inequalities
$\beta_*\le B_n$ pass to the limit, so the analytic constants and
prime intervals need not be uniform in $n$. The limiting parameters are
\[
\begin{array}{c|c@{\qquad}c|c}
 B_{\mathrm{new}}&0.874957019420098946\ldots&
 \kappa_{\mathrm{new}}&0.749914038840197892\ldots\\
 b_{\mathrm{new}}&0.123405685419541002\ldots&
 \ell_{\mathrm{new}}&0.166838588986270882\ldots,
\end{array}
\]
and $7884B_{\mathrm{new}}^3-18819B_{\mathrm{new}}^2+14643B_{\mathrm{new}}-3686=0$
\LiuLoc{(33.4)}.

\subsection{The critical counting exponent}
\label{sec:iv-count-obstruction}
In this subsection, we recall the first step of the optimality argument of
\LiuLoc{Section~34}, which minimizes the supplied counting exponent at the
critical class: by Lemma~\ref{lem:iv-count-obstruction} the minimum is
$R_\kappa(\xi(\kappa),1/2)=2/3$, and \LiuLoc{Lemma~34.3} turns it into the
constraints $B\ge L_0$, $B\ge L_K$, $B\ge H_{\xi(\kappa)}$, the last through
the identity \eqref{eq:prev-critical-necessity}. These are used for the
cyclotomic critical count of Section~\ref{sec:v-detail-34} and in
Sections~\ref{sec:v-detail-33}, \ref{sec:v-detail-36}
and~\ref{sec:v-detail-37}.

For
$13/18\le\kappa\le1$ let $\xi(\kappa)$ be the smaller root of
\eqref{eq:iv-critical-polynomial}; it lies in $I_0=[19/50,39/100]$, and
$\xi(\kappa)$ inverts $\kappa(\xi)$ of \eqref{eq:iv-geometry} on $I_0$.
Put $F(\kappa)=(18\xi(\kappa)+7)/(18\xi(\kappa)+9)$. The length
exponents $r,m$ of the two detecting polynomials and the detecting
cutoff $t$ range over
\[
 t-1/2\le r\le t,\qquad0\le m\le1/2,
 \qquad r+m\ge t,\qquad1\le t\le3/2
\]
\LiuLoc{(34.5)}. At $\delta=\xi(\kappa)$, $x=1/2$,
$c=2/(3\kappa)$ and combined range parameter one, the optimized
exponents of the inverse, unweighted and amplified long estimates are
\[
\begin{aligned}
 e_I(r)&=1-\delta(1+r)/2 &&(r\le1),\\
 e_S(m)&=1-\delta\bigl[2m+(c/4)(1-2m)\bigr] &&(0\le m\le1/2),\\
 e_L(r)&=1-\alpha+(\alpha-\delta)r &&(r\ge1)
\end{aligned}
\]
\LiuLoc{(34.6)--(34.7)}.

\begin{lemma}[{\LiuLoc{Lemma~34.2}}]\label{lem:iv-count-obstruction}
For these moment estimates and admissible length pairs, the smallest
uniform upper-bound exponent obtainable at the critical class is
\[
 \inf_{1\le t\le3/2}\ \sup_{(r,m)}
 \min\{e_{\mathrm{inv}}(r),e_S(m)\}
 =R_\kappa(\xi(\kappa),1/2)=2/3,
\]
where $e_{\mathrm{inv}}=e_I$ for $r\le1$ and $e_L$ for $r\ge1$. Upper
bounds on the selected prime lengths, estimates without selected prime
factors and H\"older interpolation between these estimates cannot
lower it.
\end{lemma}

The proof of \LiuLoc{Lemma~34.2} uses two admissible splits. The crossing
split $m=t-r_\kappa(t,1/2)$ imposes the short count
$R_{\kappa,\rm short}(t)$ of Section~\ref{sec:iv-count}; the split
$(r,m)=(t,0)$ imposes $R_{\rm long}(t)=e_L(t)$, because
\[
 e_S(0)-e_L(t)\ge-5/12+\delta(3/2-c/4)\ge-5/12+(19/50)(5/4)=7/120>0;
\]
if a discrete presentation omits $m=0$, one passes to the scalar limit.
The short/long derivation of Proposition~\ref{prop:iv-count} gives the
reverse bound, so the supplied exponent is
$\max\{R_{\kappa,\rm short}(t),R_{\rm long}(t)\}$, minimized at
$t=t_\kappa$.

Enlarging the combined range exponent to $M'\ge1$ does not
help: beyond the short length bounds the costs have slopes $2-\delta>0$
and $6\kappa-\delta>0$ in the selected length, and for $r>1$ the
unamplified inverse cost $r(1-\delta)$ is at least $e_L(r)$.

For $M,\ell\in\mathbb R$ put
\[
 L_0=5/6+M/12-\ell/6,\qquad
 L_K=1/3+7M/12+\ell/3,\qquad
 H_\xi=1+\xi-(1+\xi)M/2+\xi\ell
\]
\LiuLoc{(34.1)}. By \LiuLoc{Lemma~34.3}, a construction in the
comparison system with $B<B_{\mathrm{new}}$, $13/18\le\kappa\le1$ and
critical count exponent at least $2/3$ must satisfy
$B\ge L_0$, $B\ge L_K$ and $B\ge H_{\xi(\kappa)}$; the first two come
from the full direct majorant \eqref{eq:prev-B-low}, and the third from
the identity
\begin{equation}\label{eq:prev-critical-necessity}
 E_B(h)=H_\xi-B+h(R-2/3)\qquad\text{at }(d,\delta,x)=(h,\xi,1/2)
\end{equation}
\LiuLoc{(34.8)}, in which all $z_0$ terms cancel. The proof of
\LiuLoc{Lemma~34.3} also shows that the auxiliary choices cannot evade
these constraints:
\begin{enumerate}
\item raising the lower zero-class threshold to $\delta_f\ge\xi$ replaces
the critical class by its trivial count, with
$E_{B,\rm floor}(h)=H_{\delta_f}-B+h/3$ and
$H_{\delta_f}-H_\xi=(\delta_f-\xi)(1-M/2+\ell)\ge0$, so it requires the
stronger $B\ge H_\xi+h/3$, while a threshold below $\xi$ retains the
critical class;
\item global reciprocal contours and absolute contours give positive
endpoint exponents relative to $C(B)$, so they do not improve the
comparison;
\item principal contour changes affect only the savings of the error
terms: the residues $w=1$, $z=1/6$ and the exponent of the principal
term are fixed, and new residue extractions or estimates change the
system.
\end{enumerate}

\subsection{The optimality theorem for the specified comparison system}
\label{sec:iv-system-optimality}
In this subsection, we recall that the three constraints of
Section~\ref{sec:iv-count-obstruction} admit the positive-weight dual bound
$B\ge F(\kappa)$ of Proposition~\ref{prop:iv-dual-optimality}, which gives
the optimality Theorem~\ref{thm:iv-system-optimality},
$\inf B=B_{\mathrm{new}}$ over the admissible parameter choices, and
Corollary~\ref{cor:iv-all-constructions} for iterated constructions
\LiuLoc{Section~34}. Their cyclotomic analogues are in
Sections~\ref{sec:v-detail-35} and~\ref{sec:v-detail-36} of
Section~\ref{v:sec:strong-parameter-barrier}, and
Theorem~\ref{thm:iv-system-optimality} is also used in
Section~\ref{sec:v-metabelian}.

\begin{proposition}[{\LiuLoc{Proposition~34.1}}]\label{prop:iv-dual-optimality}
Let $13/18\le\kappa\le1$ and $M,\ell,B\in\mathbb R$. If
\[
 B\ge L_0,\qquad B\ge L_K,\qquad B\ge H_{\xi(\kappa)},
\]
then $B\ge F(\kappa)$, with equality at
$\ell=(5-6\xi(\kappa))/(9+18\xi(\kappa))$ and $M=1-\ell$. If moreover
$B\le(1+\kappa)/2$, then $\kappa\ge\kappa_{\mathrm{new}}$ and
$B\ge B_{\mathrm{new}}$. Equality is attained by
\eqref{eq:iv-final-parameters}.
\end{proposition}

The proof uses the weights
$(w_0,w_K,w_H)=(4+18\xi,2,3)/(9+18\xi)$, $\xi=\xi(\kappa)$, which are
positive, sum to one and eliminate $M$ and $\ell$:
$B\ge w_0L_0+w_KL_K+w_HH_\xi=F(\kappa)$. The fixed-point identity
\eqref{eq:iv-fixed-point-identity} then converts $F(\kappa)\le B\le(1+\kappa)/2$
into $p(\xi)\le0$, that is, $\kappa\ge\kappa_{\mathrm{new}}$.

The condition $B\le(1+\kappa)/2$ alone does not justify applying the
moment estimate; this requires a previously proved bound
$\beta_*\le\widehat B$ for the family. Then
\[
 \kappa\ge2\widehat B-1,\qquad B\ge F(\kappa)\ge F(2\widehat B-1),
\]
and $F(2B_{\mathrm{new}}-1)=B_{\mathrm{new}}$ \LiuLoc{Section~34}.

Let $\mathscr A$ be the set of parameter choices of \LiuLoc{Section~34.2}:
with all exponent losses set to zero,
\[
 \begin{gathered}
 l_x,l_y>\ell>0,\quad b=l_y-l_x>0,\quad
 M=l_x+l_y,\quad h=1-l_x+\ell>0,\\
 13/18\le\kappa\le1,\qquad B\le(1+\kappa)/2,\\
 B\ge B_{\mathrm{low}},\qquad
 E_B(d)\le0\quad\text{for every permitted zero and amplitude class},
 \end{gathered}
\]
where $B_{\mathrm{low}}$ is the full direct majorant
\eqref{eq:prev-B-low} and $E_B(d)$ is
\eqref{eq:general-extension-exponent}. The supplied counts are those of
Section~\ref{sec:iv-count-obstruction}; the contour choices, lower
zero-class thresholds and norm ranges are those examined in the proof of
\LiuLoc{Lemma~34.3}, recalled after \eqref{eq:prev-critical-necessity}; and
there is no additional bound on the total selected prime length.

\begin{theorem}[{\LiuLoc{Theorem~34.4}}]\label{thm:iv-system-optimality}
For the admissible parameter choices $\mathfrak C\in\mathscr A$,
\[
 \inf_{\mathfrak C\in\mathscr A}B(\mathfrak C)=B_{\mathrm{new}}.
\]
Every such choice satisfies $B(\mathfrak C)\ge B_{\mathrm{new}}$; equality is
attained when all exponent losses are set to zero, using
\eqref{eq:iv-final-parameters}.
This includes free $M$, all detecting cutoffs in $[1,3/2]$, the
specified estimates with and without selected prime factors, and the
auxiliary choices examined in the proof of \LiuLoc{Lemma~34.3}.
\end{theorem}

\begin{corollary}[{\LiuLoc{Corollary~34.5}}]\label{cor:iv-all-constructions}
At each improving step of this system, allow arbitrary admissible
choices of the length parameters, the zero-detector truncation exponent
$t$, the moment parameter $\kappa$ and the other auxiliary parameters.
If each moment hypothesis comes from the preceding proved boundary,
every resulting bound and every real limit is at least $B_{\mathrm{new}}$.
The sequence of Section~\ref{sec:iv-iteration} attains this infimum.
\end{corollary}

\subsection{Scope of the optimality statement}
\label{sec:iv-further-parameters}
In this subsection, we recall from \LiuLoc{Section~34.3} the scope of
Theorem~\ref{thm:iv-system-optimality}: it concerns
the uniform counting and comparison exponents of the
specified estimates, and stronger analytic inputs can change the
obstruction. The scope statement of Section~\ref{sec:v-detail-37} is
modeled on it.

Rational parameters give the same
infimum, then not attained \LiuLoc{Section~34.3}; a one-step
bound with rational parameters is \LiuLoc{Proposition~34.6}. The total
length is free, balance occurs at equality and the detecting
cutoff is minimized; $\delta,x,d$ are class
coordinates requiring uniform control.

A stronger analytic estimate can
change the obstruction: with the present $D,P$,
$\partial_\alpha R=\delta^2P^2/(2\mathcal J^2)>0$,
so a proved amplification exponent below $5/6$ would improve the
upper-bound exponents. Larger admissible prime lengths in the moment
estimates, a joint moment estimate for the detecting polynomials,
better functional-equation or Gram estimates, or restrictions on actual
critical classes could likewise alter the comparison. The averaged
numerator estimate of Section~\ref{v:sec:refined} is a change of this
kind.

\part{Cyclotomic fields and class-field descent}
\label{part:V}

In this part, we give an analytic extension of the argument of
\cite{LiuArXiv} to finite-order Hecke $L$-functions over a cyclotomic
field $F\supset\mu_{12}$, using the analytic inputs specified in
Remark~\ref{v:rem:analytic-inputs}. Section~\ref{sec:v-scope} states the
main result, Theorem~\ref{thm:v-main}, and the conventions over
$F$. Sections~\ref{sec:v-packet}--\ref{sec:v-detail-6} construct the
exceptional theta function, its radial functional equation and the sieve
and Poisson inputs; Sections~\ref{v:sec:unmarked-assembly}--\ref{v:sec:M2-assembly}
derive the inverse moments and the product moment, the counterparts over
$F$ of \cite[(9.1), (9.3)]{LiuArXiv} and of Proposition~\ref{prop:iv-plain}.
Section~\ref{v:sec:compensated-assembly} combines them with the scalar
results of Section~\ref{sec:prev-iv} to reach
$B_\infty=B_{\mathrm{new}}$ (Theorem~\ref{v:thm:original-target}),
and Section~\ref{v:sec:refined} lowers the boundary to $B_\flat$
(Theorem~\ref{v:thm:refined-limit}) by an averaged numerator
estimate. Section~\ref{sec:v-optimality} proves the lower bound
$B\ge B_\infty$ in the pointwise-numerator comparison system
(Theorem~\ref{thm:v-optimality}), attained under Theorem~\ref{v:thm:original-target}, and
Section~\ref{sec:v-descent} descends by class field theory to the pairs
$(K,\eta)$ of \eqref{eq:v-class}. Part~\ref{part:VI} shows that the
argument uses only $\mu_{12}\subset F$.

\section{The extended Hecke family and its boundary}
\label{sec:v-scope}
In this section, we state an extension of Theorem~\ref{thm:iv-limit} to the
cyclotomic fields of \cite{Artin} and the class fields
\eqref{eq:v-class}, stated as Theorem~\ref{thm:v-main}, and fix
the conventions of Part~\ref{part:V}: the boundary \eqref{eq:v-boundary},
the Hilbert-symbol orientation \eqref{eq:v-hilbert-orientation}, and the
allowed generators, reciprocity \eqref{eq:v-generator-unit-reciprocity}
and finite unit characters $\Upsilon_S$ of
Lemmas~\ref{v:lem:allowed-reciprocity} and~\ref{v:lem:unit-products},
which are used throughout Part~\ref{part:V}. Section~\ref{v:sec:compensated-assembly}
reproduces $B_\infty=B_{\mathrm{new}}$ and the limiting parameters
\eqref{eq:iv-final-parameters}; Section~\ref{v:sec:refined} refines this to
$B_\flat=0.874886647039964\ldots$ by an averaged numerator estimate adapted
from \cite{XS26} and a fixed-point limit argument, and Theorem~\ref{thm:v-main} is stated with $B_\flat$.
Sections~\ref{sec:v-packet}--\ref{v:sec:refined} give
the analytic derivation;
Sections~\ref{sec:v-optimality} and \ref{sec:v-descent} give
the scalar obstruction and the arithmetic descent implication.
Part~\ref{part:VI} removes the cyclotomic restriction and
extends $B_\flat$ to every number field.
Remark~\ref{v:rem:analytic-inputs} specifies the analytic dependencies
and the scope of the scalar and arithmetic arguments.

Cureton and Thorner~\cite[Theorem 1.1, Corollary 1.2 and Sections 5.1 and 14]{CT26}
obtain $\Rea s>31/32$ for $\zeta_K$ and finite-order Hecke functions
over every number field $K$, working over $K(\zeta_{12})$ and descending by abelian
factorization; $K/\mathbb Q$ need not be abelian.
The regions
$\Rea s>B_\infty$, $\Rea s>B_\flat$ require the varying-parameter moment estimates
of Sections~\ref{sec:v-packet}--\ref{v:sec:refined} and do not follow
by substituting $B_\infty$ or $B_\flat$ into the
fixed-parameter bounds of that report. We cite its version~2;
its region is $\Rea s>31/32$.

Throughout Part~II put
\begin{equation}\label{eq:v-boundary}
 \xi_\infty=\xi_{\mathrm{new}},\qquad
 B_\infty=B_{\mathrm{new}}=
 \frac{18\xi_\infty+7}{18\xi_\infty+9}
 =0.874957019420098946\ldots,
\end{equation}
where $\xi_{\mathrm{new}}$ is defined in Theorem~\ref{thm:iv-limit}.
The parameters $b_\infty,\ell_\infty,\kappa_\infty$ are those of
\eqref{eq:iv-final-parameters}; in particular
$\kappa_\infty=2B_\infty-1$. Section~\ref{v:sec:refined} gives the
limiting boundary
\[
 \begin{gathered}
 B_\flat=0.874886647039964\ldots,\\
 133488B_\flat^4-376596B_\flat^3+384534B_\flat^2-167541B_\flat+26231=0,
 \end{gathered}
\]
in Theorem~\ref{v:thm:refined-limit}, the exact infimum of its
refined comparison system; Appendix~\ref{app:refined-chain} adds the
explicit verified rational stage
\[
 B_\sharp=\frac{874886647041}{10^{12}}=0.874886647041,
 \qquad B_\flat<B_\sharp<B_\infty .
\]
For a number field $K$ and a finite-order Hecke character $\eta$ of its
idele class group $C_K=\mathbb A_K^\times/K^\times$, let $M_\eta/K$
be the cyclic extension corresponding to $\ker\eta$ under class field
theory. Thus $[M_\eta:K]=\operatorname{ord}\eta$, and $M_1=K$.
Fix an algebraic closure $\overline{\mathbb Q}$ and take all number
fields, class fields, normal closures and composita inside it. Define
\begin{equation}\label{eq:v-class}
\begin{split}
 \mathscr C=\{(K,\eta):\;&K\subset\overline{\mathbb Q},\quad
     [K:\mathbb Q]<\infty,\quad\eta\text{ a finite-order Hecke character},\\
 &\exists m\ge1:\ F=\mathbb Q(\zeta_m)\supset\mu_{12},\quad
                       M_\eta F/F\text{ abelian}\}.
\end{split}
\end{equation}
Here and below an abelian extension is finite and Galois.

\begin{remark}[Analytic inputs and proof structure]\label{v:rem:analytic-inputs}
For each fixed cyclotomic $F\supset\mu_{12}$, the analytic argument
uses the following results, with their displayed coefficients and natural zeros:
\begin{enumerate}
\item Proposition~\ref{v:prop:checked-inputs} constructs the normalized
theta datum from the cited representation-theoretic results.
Lemma~\ref{v:lem:theta-characters} gives its coefficient characters,
and Lemma~\ref{v:lem:raw-theta} proves the norm-sheet identity with
\eqref{eq:v-exact-S-vector}. Proposition~\ref{v:prop:radial} derives
the compact radial functional equation and its uniform bounds.
\item Theorems~\ref{v:thm:unmarked-inverse} and
\ref{v:thm:marked-inverse-new} state the inverse moments;
Theorem~\ref{v:thm:M2-assembled} states the product moment.
Their finite coefficient groups contain the unit characters of
Lemma~\ref{v:lem:allowed-reciprocity}. The zero-free hypothesis in the
product moment for $\kappa<1$ is retained in each application.
\item Section~\ref{v:sec:compensated-assembly} applies those moments
to the compensated theta average \eqref{v:cmp:1}, the cyclotomic counterpart of
the modified average of Section~\ref{sec:prev-average}. Its direct bound,
nonzero principal term and remainder are respectively
\eqref{v:cmp:4}, \eqref{v:cmp:8} and \eqref{v:cmp:13}.
Lemma~\ref{v:lem:arithmetic-independent-mesh} and the separate
principal and nonprincipal margins fix the numerical saving before
the target-dependent arithmetic and height data.
\item Section~\ref{v:sec:refined} replaces the pointwise numerator
estimate \eqref{v:cmp:11} by the amplified first moment of
Proposition~\ref{v:prop:refined-numerator}, which uses M2 together
with the standard approximate functional equation of
Lemma~\ref{v:lem:refined-afe}, and proves the refined remainder
\eqref{eq:v-refined-remainder}; Lemma~\ref{v:lem:refined-family} and
Theorem~\ref{v:thm:refined-limit} give the limiting boundary
$B_\flat$. This step adapts \cite{XS26}.
\end{enumerate}
The derivations below use the stated theta,
reciprocity, Poisson, Hecke and quadratic/cubic sieve inputs; the
zero-free conclusion of \cite{Artin} is not an input. The local and
global theta calculations of that preprint remain source inputs; the
moment estimates are stated and derived below. The arguments specify
the cited local results, the uniformity conditions
and the zero-free conditions used in each application.
The inverse and product moment estimates in item~(2) are the analytic
assertions whose derivations are given in the cited sections of this
paper. The condition
$\beta_*(F)\le(1+\kappa)/2$ is an intermediate condition of M2:
Section~\ref{sec:v-detail-31} verifies it in the initial and pointwise
comparisons, and Theorem~\ref{v:thm:refined-limit} verifies it with
equality by setting $\kappa=2\beta_*(F)-1$ in the limiting contradiction.
Thus it is discharged within the proof of the main result.
The arithmetic descent and the scalar optimality arguments do not
use the moment estimates themselves.

\end{remark}

\begin{remark}[Order of choices and uniformity]
\label{v:rem:parameter-ledger}
Fix $F$ and a contradiction $\beta=\beta_*(F)>B$. In both the
pointwise and refined comparisons, choose, in order:
\begin{enumerate}
\item the admissible geometry and moment parameter, numerical margins,
row cutoffs and meshes;
\item $K$ and the window lengths $\ell_i$, $\sum_i\ell_i=\ell$, then
the principal saving and buffer $e$, normalization loss, bins and contour rules;
\item the target $\eta$ and final exceptional set $S$, then the
$S$-dependent generators, coefficient groups, tests and normalizers;
\item the internal derivative orders and height exponents, then
$\tau_\eta>0$, the external tail order $N$ satisfying
\[
 N\tau_\eta>D_\eta-C(\beta)+m_{\rm save},
\]
and finally the onset $Z_0$.
\end{enumerate}
All fixed data precede $Z\to\infty$; only
$P_i=Z^{\ell_i}$ and $T_1=Z^{\tau_\eta}$ vary by their prescribed rules.
The principal $s$-contour remains in $\Rea s>\beta$
(Remark~\ref{v:rem:no-circularity}).

For fixed $F$ and scalar comparison, choices (1)--(2) and the resulting
$\sigma_{\rm cont}<\beta_*(F)$ are common to all targets.
Arithmetic data, constants, derivative orders and onsets may depend
on $F$ and $\eta$. The comparison parameters may vary with $F$;
the final boundary $B_\flat$ is field independent because its scalar
inequalities are field independent.
\end{remark}

\begin{theorem}[Hecke bound and its descent]\label{thm:v-main}
For every $(K,\eta)\in\mathscr C$,
\[
 L_K(s,\eta)\ne0\qquad(\Rea s>B_\flat),\qquad
 B_\flat=0.874886647039964\ldots,
\]
with the pole of the principal character at $s=1$ permitted.
In particular this holds for every finite-order Hecke character over
 every finite abelian extension $K/\mathbb Q$.
The numerical boundary and the scalar admissibility inequalities are
independent of $K$ and $\eta$. The comparison parameters may depend
on $\beta_*(F)$ and the contradiction threshold, but are common to
all targets over the chosen auxiliary field $F$; implied constants
and thresholds may depend on the fixed fields and characters.
No assertion is made on $\Rea s=B_\flat$.
Within the pointwise-numerator comparison system specified in
Section~\ref{sec:v-optimality}, $B_\infty$ is the infimum over the
admissible parameter choices and successive constructions; the lower bound for that
infimum is the separate algebraic statement of
Theorem~\ref{thm:v-optimality}. The boundary $B_\flat<B_\infty$ is
obtained with the additional averaged numerator estimate of
Section~\ref{v:sec:refined}, which lies outside that system. For the
refined comparison system of Section~\ref{v:sec:refined-optimality}
$B_\flat$ is the exact infimum, Theorem~\ref{v:thm:refined-optimality}.
\end{theorem}
\begin{proof}
The cyclotomic assertion is Theorem~\ref{v:thm:refined-limit}.
Proposition~\ref{prop:v-descent} proves descent, and
Corollary~\ref{cor:v-abelian} gives the abelian-field case.
Proposition~\ref{prop:v-metabelian} identifies the arithmetic scope.
Theorem~\ref{thm:v-optimality} supplies the scalar lower bound for the
pointwise-numerator system;
Theorem~\ref{v:thm:original-target} and the limiting construction
supply the reverse inequality for that infimum.
\end{proof}

The conventions of Section~\ref{sec:arithmetic} remain in force
except for the ground field: there and in \cite{LiuArXiv},
$F=\mathbb Q(\sqrt{-3})$, $\OO=\mathbb Z[\omega]$; here
$\OO_F$ is the ring of integers of the fixed cyclotomic field $F$.
Absolute ideal norms, valuations, radicals, M\"obius coefficients,
$1_E$, $\sfq$, and the length convention
\eqref{eq:length-exponent-convention} have the same meanings.
Character presentations retain their defining moduli and zero extensions
as in \eqref{eq:primitive-induction}--\eqref{eq:zeros}; these moduli
are not replaced by primitive conductors in intermediate sums.
The local additive characters and Haar measures below replace the
quadratic-field normalization of Section~\ref{sec:arithmetic}.

For fixed $F=\mathbb Q(\zeta_m)\supset\mu_{12}$ set
\[
 \beta_*(F)=\sup\bigl(\{1/2\}\cup
 \{\Rea\rho:L_F(\rho,\nu)=0,\ \nu\text{ primitive of finite order}\}\bigr),
\]
where only nontrivial zeros are included and principal poles are excluded.
The argument summarized in Remark~\ref{v:rem:analytic-inputs}
gives $\beta_*(F)\le B_\infty$, and with Section~\ref{v:sec:refined}
$\beta_*(F)\le B_\flat$, without assuming the zero-free conclusion
of \cite{Artin}.

Choose a finite set $S$ containing the archimedean places, the places
above $6$, the support of the different, the fixed character conductors and the exceptional places of the theta
representation, and enough finite places that $R=\OO_{F,S}$ is principal.
Include also the nondistinguished places of a nonzero pure tensor in
Proposition~\ref{v:prop:checked-inputs}, all fixed finite coefficient
conductors, all primes below the fixed cutoff required for
\eqref{v:cmp:7}, and those below the contour cutoff of
Lemma~\ref{v:lem:H-strip}. These cutoffs are independent of the unitary
target phases; the latter may depend on the already chosen buffer $e$.
In a comparison application the geometry, prime lengths and buffer
are fixed first, as in Remark~\ref{v:rem:parameter-ledger}; for each
target $\eta$ this prescription then fixes the final set $S$.
These finite choices precede the definitions of the allowed generators,
$B_S$, the Gauss sums, $\Lambda,G_0$ and local tests, and $S$ is fixed
throughout a fixed-target argument.
Write $S_f$ for the finite places in $S$. Normalize local absolute
values by $|\varpi_v|_v=(N\mathfrak p_v)^{-1}$ at finite places and
$|z|_v=|z|^2$ at complex places. Put
\[
 F_S=\prod_{v\in S}F_v,\quad
 |y|_S=\prod_{v\in S}|y_v|_v,\quad
 G=F_S^1=\ker|\cdot|_S,\quad
 U_S=R^\times,\quad\Gamma=U_S^6.
\]

An exterior ideal is a nonzero integral ideal prime to $S_f$;
$q_{\mathfrak a}=N\mathfrak a$ and
$\mu_F(\mathfrak a)$ denote its absolute norm and ideal M\"obius function. We also write $\mu(a)=\mu_F(a)$ when the field is fixed.
Indices such as $a,n,f$ denote allowed generators of exterior ideals
modulo $\Gamma$, with the fixed unit classes retained. The generators
are multiplicative modulo $\Gamma$ as in \cite[Section 3]{Artin}.
The notation $q_a$, $(a,b)$, $d\mid a$ and $a\ \sfq$ refers to the
corresponding exterior ideals; $\operatorname{rad}(a)$ is their squarefree
radical. The sextic residue symbol is $\chi_a(b)$; every power,
including the zeroth, is zero when $(a,b)\ne1$.
Use the trace additive character and the dual lattice normalization
of Sections~3 and~6 of the same source; Haar measures are fixed throughout.
A multiplicative norm section is written $\ell_T$, with
$|\ell_T|_S=T$ and finite components one. Fixed local classes and
finite character groups are chosen before any norm decomposition.

We use \textup{(M1)} for the inverse estimate with prime-polynomial factors of
Theorem~\ref{v:thm:marked-inverse-new} and \textup{(M2)} for the
unweighted product estimate of Theorem~\ref{v:thm:M2-assembled}.
Their common-character and prime-support conditions remain part of the
notation. Write $d(\mathfrak a)=\sum_{\mathfrak d\mid\mathfrak a}1$ for the
ideal divisor function. In the moment estimates, $f$ is an auxiliary
ideal index and a prime-polynomial factor is a normalized sum over one
prescribed prime annulus. Auxiliary norm exponents are local to their
stated estimates; in particular, $F=N+V$ in
Section~\ref{sec:v-detail-9} is a real exponent, not the ground field.

Several letters of Part~\ref{part:prelim} are reused locally in
Part~\ref{part:V}. In Sections~\ref{v:sec:unmarked-assembly}
and~\ref{v:sec:marked-assembly}, $z_0$ is the cap on the total
selected-prime length, the counterpart of $z_{\rm cap}$ in
\cite[(9.2)]{LiuArXiv}, not the contour value $17/50$ of
\eqref{eq:general-extension-exponent}; the margin conditions called
canonical there are the counterparts of the inequalities in that display; and the
fine-scale mesh $\eta$ is unrelated to the character $\eta$. In
Section~\ref{v:sec:M2-assembly}, $M=m+q$ is the combined range exponent
$m_0+q_0$ of Proposition~\ref{prop:iv-plain}, not $l_x+l_y$. The letter
$\kappa$ also denotes local indices, such as an element of $\mathcal V_S$
in Section~\ref{sec:v-detail-2} and the constants $\kappa_i,\kappa_c$ of
Sections~\ref{v:sec:unmarked-assembly}--\ref{v:sec:marked-assembly}; the
moment parameter is always written $\kappa$ without subscript. The amplified exponent
$e(r)=\max\{1,(1+5r)/6\}$ is $\mathfrak a(r)$ of the table in
Section~\ref{sec:arithmetic}. Elsewhere $G=F_S^1$, $\Gamma=U_S^6$ and
$R=\OO_{F,S}$ replace the uses of these letters in
Part~\ref{part:prelim}, and $u=g/\delta$ in
Section~\ref{v:sec:compensated-assembly} is the variable $x$ of
Section~\ref{sec:prev-iv}.
The parameter functions $b(\xi),\ell(\xi),\kappa(\xi),B(\xi)$ agree
with those of \eqref{eq:iv-geometry}. A norm twist is multiplication by $q_{\mathfrak a}^{iT}$,
$T\in\mathbb R$; its cost is measured by fixed polynomial powers of
$1+|T|$.

Fix an embedding $\jmath:F\hookrightarrow\mathbb C$, and let
$\jmath_6$ denote its restriction to $\mu_6(F)$. At an exterior prime,
$\chi_p(u)=\jmath_6(\widetilde{\bar u^{(q_p-1)/6}})$, where the tilde
is the unique lift to $\mu_6(F)$. Our Hilbert symbols have orientation
\begin{equation}\label{eq:v-hilbert-orientation}
 (a,b)_{n,v}=\frac{\operatorname{Art}_v(b)(a^{1/n})}{a^{1/n}},
 \qquad n\in\{2,3,6\},
\end{equation}
where a uniformizer acts by arithmetic Frobenius on unramified
extensions. At a tame place this gives
\[
 (\varpi^a u,\varpi^b v)_n
 =\widetilde{\left((-1)^{ab}\bar u^{\,b}\bar v^{-a}
                                  \right)^{(q-1)/n}}.
\]
We apply $\jmath_6$ when writing Hilbert symbols as complex scalars.
The genuine central embedding for the cubic cover is instead
\begin{equation}\label{eq:v-genuine-embedding}
 \iota(\zeta)=\overline{\jmath_6(\zeta)}\quad(\zeta\in\mu_3(F)).
\end{equation}
These conventions are retained in the local Whittaker calculation
and the global theta identity below.

\begin{lemma}[Allowed generators, reciprocity and unit characters]
\label{v:lem:allowed-reciprocity}
Put $\mathcal V_S=F_S^\times/F_S^{\times6}$ and
$H_S(x,y)=\prod_{v\in S}(x_v,y_v)_{6,v}$, with the Hilbert-symbol
orientation of \cite[Section 3, equations (3.3)--(3.5)]{Artin}.
The image $E_S$ of $U_S/U_S^6$ is injected, $E_S^\perp=E_S$,
and one may choose an isotropic complement $B_S$ satisfying
\[
 \mathcal V_S=E_S\oplus B_S,\qquad E_S\cap B_S=\{1\}.
\]
An exterior fractional ideal has an allowed generator $a$ with
$[a]\in B_S$, unique modulo $\Gamma=U_S^6$. Products of allowed
generators remain allowed. For allowed $a,b$ and $\epsilon\in U_S$,
\begin{equation}\label{eq:v-generator-unit-reciprocity}
 \chi_a(b)=\chi_b(a),\qquad
 \chi_a(\epsilon b)=\chi_a(\epsilon)\chi_b(a),\qquad
 \chi_a(\epsilon)=H_S(\epsilon,a)^{-1}.
\end{equation}
The first two identities include the shared-prime zeros.
Choose a set $\mathcal U_6$ of representatives for $U_S/U_S^6$.
Every nonzero $h\in R$, modulo $\Gamma$, has a unique form
$h=\epsilon a_h$, $\epsilon\in\mathcal U_6$, with $a_h$ allowed.
The characters
\[
 \xi_\epsilon(\mathfrak a)=\chi_a(\epsilon),\qquad
 \Upsilon_S=\langle\xi_\epsilon:\epsilon\in\mathcal U_6\rangle
\]
form a fixed finite group of finite-order Hecke characters ramified
only in $S$; all have order dividing six.
\end{lemma}
\begin{proof}
These are the constructions of \cite[Section 3, Lemmas 3.2--3.3]{Artin}.
For completeness, at exponent $\ell=2,3$ the Hilbert pairing is
nondegenerate and alternating, and
\[
 \dim_{\mathbb F_\ell}F_S^\times/F_S^{\times\ell}
       =2(|S_f|+[F:\mathbb Q]/2),\qquad
 \dim_{\mathbb F_\ell}U_S/U_S^\ell=|S_f|+[F:\mathbb Q]/2.
\]
The unit map is injective: its kernel would give an everywhere
unramified Kummer extension split at $S$, hence a quotient of the
trivial class group of $R$. Global reciprocity makes its image
isotropic. Symplectic linear algebra gives a complement for each
$\ell$; combining them gives $B_S$ and the displayed direct sum.
Multiplication by an $S$-unit removes the $E_S$ component of any
generator; two choices differ by $E_S\cap B_S$, hence by $U_S^6$.
For coprime allowed generators, the exterior Hilbert symbols and
global product formula give
$\chi_a(b)^{-1}\chi_b(a)H_S(a,b)=1$. Isotropy of $B_S$
proves reciprocity. Shared-prime zeros and multiplicativity give
\eqref{eq:v-generator-unit-reciprocity}. Finally $\xi_\epsilon$
is the Kummer character of $F(\epsilon^{1/6})/F$ in this orientation;
exterior ramification is absent. Finiteness follows from the finite
group $U_S/U_S^6$.
\end{proof}

\begin{lemma}[Finite unit characters in quadratic products]
\label{v:lem:unit-products}
With the notation of Lemma~\ref{v:lem:allowed-reciprocity}, put
$\mathcal E=U_S/U_S^6$ and $N_S=|\mathcal E|=|B_S|$.
For $b\in B_S$ write $\xi_\epsilon(b)=H_S(\epsilon,b)^{-1}$.
For $w:\mathcal E\to\mathbb C$ define
\[
 \widehat w(b)=N_S^{-1}\sum_{\epsilon\in\mathcal E}
                         w(\epsilon)\xi_\epsilon(b).
\]
Then
\begin{equation}\label{eq:v-unit-projection}
 \begin{aligned}
 N_S^{-1}\sum_\epsilon\xi_\epsilon(b)
                   \overline{\xi_\epsilon(c)}&=1_{b=c},\\
 N_S^{-1}\sum_\epsilon w(\epsilon)
       \overline{\xi_\epsilon(b_U)}\xi_\epsilon(b_V)
       &=\widehat w(b_Vb_U^{-1}).
 \end{aligned}
\end{equation}
Here $b_U=[U],b_V=[V]\in B_S$ for allowed generators.
For the full Gauss sums and complete correlations of
Section~\ref{sec:v-detail-21}, including imprimitive moduli,
\begin{equation}\label{eq:v-unit-Gauss-correlation}
 \begin{aligned}
 \mathcal G(U,\epsilon h)
   &=\overline{\xi_\epsilon(U)}\mathcal G(U,h),\\
 \mathcal G(U,\epsilon h)\overline{\mathcal G(V,\epsilon h)}
   &=\overline{\xi_\epsilon(U)}\xi_\epsilon(V)
                 \mathcal G(U,h)\overline{\mathcal G(V,h)},\\
 \mathcal F(U,V;\epsilon j)
   &=\xi_\epsilon(U)\overline{\xi_\epsilon(V)}
                 \mathcal F(U,V;j).
 \end{aligned}
\end{equation}
All identities preserve the natural zeros.
\end{lemma}
\begin{proof}
Nondegeneracy of $H_S$, and $E_S^\perp=E_S$,
$\mathcal V_S=E_S\oplus B_S$, make the pairing
$\mathcal E\times B_S\to\mu_6$ perfect. Character orthogonality
proves \eqref{eq:v-unit-projection}. In the Gauss sum substitute
$x=\epsilon^{-1}y$; this is a bijection modulo $U$ even when
$h$ or $y$ is a nonunit. For the correlation, whose defining
congruence is $Vx-Uy=j\pmod{UV}$, substitute
$(x,y)=(\epsilon x',\epsilon y')$. Multiplicativity on the unit
$\epsilon$ proves \eqref{eq:v-unit-Gauss-correlation}.
\end{proof}

The matrix $(N_S^{-1/2}\xi_\epsilon(b))_{\epsilon,b}$ is unitary,
so the weighted projector in
\eqref{eq:v-unit-projection} has operator norm $\le
\max_\epsilon|w(\epsilon)|$. Alternatively, for the
sector sums $T_\epsilon$,
\begin{equation}\label{eq:v-unit-finite-Cauchy}
 \left|\sum_\epsilon c_\epsilon T_\epsilon\right|^2
 \le\left(\sum_\epsilon|c_\epsilon|^2\right)
                       \sum_\epsilon|T_\epsilon|^2.
\end{equation}
This bounds a signed sum; no unit projector is deleted. The origin
in orbit Poisson is kept separately, with its compact-quotient measure
and without a factor $N_S$. If a fixed finite character group
$\Theta\supset\Upsilon_S$, then
\begin{equation}\label{eq:v-unit-exception-invariance}
 \psi\in\Theta\quad\Longleftrightarrow\quad
                  \psi\xi_\epsilon\in\Theta.
\end{equation}
Principal-only exclusion is not invariant under such twists;
any additional principal sector must therefore be estimated explicitly.

\section{The exceptional representation and its local factors}
\label{sec:v-packet}
In this section, Proposition~\ref{v:prop:checked-inputs} gives the covering datum,
nonzero Fourier functional and local normalizations of the theta
representation over $F$, from \cite[Sections 3--6 and Appendix A]{Artin},
the applicable Kazhdan--Patterson results and \cite{KPerr}; its exterior
Whittaker table obeys the recurrence \eqref{eq:v-whittaker-recurrence}.
Lemma~\ref{v:lem:theta-characters} then identifies the coefficient
characters $\Lambda,G_0$ through $a_0(n)=\overline{\Lambda(n)}\gamma_2(n)$
and $a_0(n)\gamma_1(n)=\mu_F(n)G_0(n)$. These are the inputs of
Proposition~\ref{v:prop:radial} in Section~\ref{sec:v-detail-2} and of
the inverse transforms of Section~\ref{sec:v-detail-10}.
Subsection~\ref{v:ss:classical-theta-comparison} explains why no exact
identification with the classical spherical theta function is required,
and records the treatment of bad places and normalization.

\begin{proposition}[The exceptional theta representation]
\label{v:prop:checked-inputs}
For the rank-two, degree-three Kazhdan--Patterson cover with $c=2$,
quotient by the split scalar subgroup and take the genuine datum
\[
 \omega_v(\zeta s_v(\operatorname{diag}(x^3,y^3))s_v(zI))
       =\iota(\zeta)|x^3/y^3|_v^{1/6}.
\]
There is an exceptional automorphic representation $\Theta_F$ with
a nonzero global Whittaker functional and one-dimensional local
Whittaker spaces. Normalize its exterior spherical vectors by $W_p(1)=1$.
Here and below $\gamma_j(\mathfrak a)$ denotes the normalized Gauss sum of
the $j$th sextic-symbol power with the trace additive character of $F$,
not the quadratic-field sum of Section~\ref{sec:arithmetic}.
With the trace character and ordinary torus section,
\[
 W_p(a(p^eu))=\chi_p(u)^{2e}
 \begin{cases}
 q_p^{-k},&e=3k,\ k\ge0,\\
 q_p^{-k-1/2}\gamma_2(p),&e=3k+1,\ k\ge0,\\
 0,&e<0\text{ or }e\equiv2\pmod3.
 \end{cases}
\]
At a complex place the normalized inducing exponents are $\pm1/6$.
The two first-root derivative vectors have Whittaker functions
$|y|_v\omega_{\pm1}(y)K_{1/3}(c_v\sqrt{|y|_v})$, up to fixed
nonzero constants, with $\omega_{\pm1}(y)=(y/|y|_v^{1/2})^{\pm1}$.
At $S_f$ choose the direct cutoff on one fixed sufficiently small unit
neighborhood. Its Weyl rotations form a finite family with compact
upper support and $O(|y|_v^\delta)$ small-norm bounds for every
$\delta<1/3$. The global Whittaker expansion factors with these
local normalizations. These normalizations follow from
\cite[Section 5 and Appendix A]{Artin} and the
Kazhdan--Patterson results cited there.
\end{proposition}
\begin{proof}
For the Kazhdan--Patterson cover take rank $2$, covering degree $3$
and $c=2$. Its torus cocycle is
\[
 \sigma_{c,v}(\operatorname{diag}(a,b),
                    \operatorname{diag}(a',b'))
       =(a,b')_{3,v}(ab,a'b')_{3,v}^{c}.
\]
Its scalar action is obtained through
\eqref{eq:v-genuine-embedding}. On the quotient torus
$a(y)=\operatorname{diag}(y,1)$ the effective cocycle is
\begin{equation}\label{eq:v-effective-torus-cocycle}
 C_v(x,y)=\iota((x,y)_{3,v})^2
          =\jmath_6((x,y)_{6,v})^2.
\end{equation}
At an exterior place, therefore,
\[
 C_p(\varpi,u)=\chi_p(u)^{-2},\quad
 C_p(u,\varpi)=\chi_p(u)^2,\quad
 C_p(u,u')=C_p(\varpi,\varpi)=1.
\]
The root groups have their canonical splittings; the compact splitting
agrees on unit diagonals, integral roots and the Weyl element, as in
\cite[Proposition~0.1.3]{KP} and \cite[Appendix~A]{Artin}.

Proposition~0.1.1 of \cite{KP} gives scalar centrality exponent
$2-1+2\cdot2c=9$. On two scalar sections,
\[
 \sigma_{2,v}(zI,z'I)=(z,z')_{3,v}^{9}=1.
\]
Thus the ordinary scalar section is split and central. The rational
splitting constructed by the product formula agrees with this section
on rational scalars, so its quotient is a rationally split cover of
$\operatorname{PGL}_2$. With $s_v$ the ordinary torus section, prescribe
\[
 \omega_v\bigl(\zeta s_v(\operatorname{diag}(x^3,y^3))s_v(zI)\bigr)
       =\iota(\zeta)|x^3/y^3|_v^{1/6}.
\]
Cube arguments make all remaining Hilbert symbols one. Changing a
presentation multiplies both cube coordinates by the same scalar
cube, hence preserves the ratio. The datum is therefore well-defined,
genuine, and trivial on the split scalar subgroup. On a coroot cube,
\[
 \omega_v(s_v(\operatorname{diag}(t^3,t^{-3})))=|t|_v,
\]
the exceptional condition in \cite[Section I.2]{KP}.

For global compatibility, a rational ratio that is a cube at every
completion is a cube in $F$: otherwise its nontrivial Kummer extension
would have an unramified nonsplit prime by Chebotarev. The idele
product formula then makes $\omega$ trivial on the intersection with
the rational torus lift. This verifies the compatibility condition of
\cite[Section II.1]{KP}, and the scalar central action is trivial
and thus unitary.

Theorem~I.2.9 (p.~72), Corollary~I.3.6 (p.~79) and
Theorems~II.2.1 (p.~118) and II.2.5 (p.~122) give the
local exceptional factors, automorphic residue, nonzero global Fourier
functional and all-place Whittaker uniqueness. The tame uniqueness
criterion is $2(c+1)=6\equiv0\pmod3$; wild places are covered by the
conclusion of II.2.5.
The 1985 correction replaces the coefficient limit in
Theorems~II.2.2--II.2.3 by
\[
 c(\alpha)=\lim_{S'}C_{S'}(\alpha)\mathcal T(S'),\qquad
 \mathcal T(S')=\prod_{v\in S'_f}\prod_{j=1}^{r-1}
                             (1+q_v^{-1}+\cdots+q_v^{-j}),
\]
and inserts $\mathcal T(S')^{-1}$ in the finite-$S'$ expansion;
here $r$ is the rank in \cite{KP}, hence $r=2$ here.
The nonvanishing proof of II.2.5(a), pp.~123--126, uses the
nonzero residue, irreducibility and the finite-$S'$ Fourier expansion;
the uniqueness argument in part (b), pp.~126--127, uses that
nonvanishing and the local mirabolic filtration. The correction is
on p.~203 of~\cite{KPerr}. Neither requires the withdrawn
spherical-vanishing assertion on pp.~119--120. With the corrected
coefficient normalization, II.2.5(a),(b) give respectively the
nonzero global functional and all-place local uniqueness used here.
The set $S'$ is auxiliary; the arithmetic set $S$ is unchanged.

To normalize the global functional $\lambda$, use a pure tensor
$\varphi^\circ$ with $\lambda(\varphi^\circ)\ne0$ whose finitely
many nondistinguished finite places were included in the preliminary
choice of $S$ in Section~\ref{sec:v-scope}. This tensor is chosen
from the automorphic representation before the $S$-dependent
generators and coefficient characters are defined. Varying
one exterior component gives a nonzero local Whittaker functional,
normalized at $\varphi_v^\circ$. Local uniqueness and induction on
the number of varied components give
\[
 \lambda(\pi(g)\varphi)
       =\mathcal W_{S,\varphi_S}(g_S)
                   \prod_{v\notin S}W_{v,\varphi_v}(g_v),
       \qquad W_{v,\varphi_v^\circ}(1)=1.
\]
The finite $S$-functional is a finite sum of local products, since the
local twisted Jacquet quotients are finite-dimensional.
This uses the nonzero functional and finitely many component changes;
it does not use the spherical vanishing assertion withdrawn by
\cite{KPerr}, or the uncorrected infinite coefficient product.

The local table also follows from three initial
Jacquet values and a spherical recurrence. Fix an exterior place,
write $q=q_p$, $\varpi=p$, and $g_j=\gamma_j(p)$, and give
$\mathcal O_p$ additive volume one. In the vector-valued realization of
the Weyl-conjugate normalized induction, let $e_0$ be its unit-fixed
torus vector and put
\[
 e_j=q^{-j/6}\sigma^-(a(\varpi^j))e_0,\qquad e_{j+3}=e_j.
\]
Indeed $a(\varpi^3)$ acts in $\sigma^-$ by $q^{1/2}$.
The ordinary torus cocycle and the chosen orientation give
\[
 \sigma^-(a(\varpi^ju))e_0
       =q^{j/6}\chi_p(u)^{2j}e_j\quad(u\in\mathcal O_p^\times).
\]

The spherical section $f^\circ$, with $f^\circ(k)=e_0$, belongs to
the exceptional subrepresentation of this induction by
\cite[Theorem I.2.9(a),(f)]{KP}. Set
\[
 J_e=\int_{F_p}^{\mathrm{st}}
       f^\circ(wn(x)a(\varpi^e))\overline{\psi_p(x)}\,dx.
\]
The superscript denotes the stabilized integrals over expanding balls.
Stabilization follows by splitting large valuation shells into unit
cosets on which the multiplicative factor is constant and the additive
character has zero integral. The exact Bruhat decomposition, followed
by $x=\varpi^eH$, gives, for $e\ge0$,
\[
 J_e=q^{-2e/3}\left(e_{-e}+
   \sum_{l\ge1}q^{l/3}e_{-e+2l}
    \int_{\mathcal O_p^\times}\chi_p(u)^{-2l}
             \overline{\psi_p(\varpi^{e-l}u)}\,du\right).
\]
Here the three factors in the common power are the change of measure
$q^{-e}$, the modular factor $q^{e/2}$, and the torus factor
$q^{-e/6}$. For $0\le e\le2$, the terms $1\le l\le e$ vanish by
unit-character orthogonality, and those $l>e+1$ vanish by additive
orthogonality. The remaining integral at $l=e+1$ is respectively
$q^{-1/2}g_4$, $q^{-1/2}g_2$, and $-q^{-1}$. Consequently
\[
 \begin{split}
 J_0&=e_0+q^{-1/6}g_4e_2,\\
 J_1&=q^{-1/2}g_2e_0+q^{-2/3}e_2
       =q^{-1/2}g_2J_0,\qquad J_2=0,
 \end{split}
\]
using $g_2g_4=1$. Thus $J_0\ne0$; taking a scalar coordinate
nonzero on $J_0$ and using local uniqueness yields
$W_0=1$, $W_1=q^{-1/2}g_2$, and $W_2=0$.

For the recurrence evaluate every shell: put
\[
 I_{e,l}=\int_{\mathcal O_p^\times}\chi_p(u)^{-2l}
               \overline{\psi_p(\varpi^{e-l}u)}\,du.
\]
Unit orthogonality and the final additive residue digit give
\[
 I_{e,l}=\begin{cases}
 1-q^{-1},&1\le l\le e,\quad3\mid l,\\
 0,&1\le l\le e,\quad3\nmid l,\\
 q^{-1/2}g_{-2l},&l=e+1,\quad3\nmid l,\\
 -q^{-1},&l=e+1,\quad3\mid l,\\
 0,&l>e+1.
 \end{cases}
\]
The zeroth character retains its zero at the origin, so its
unnormalized Gauss sum is $-1$. Let $k=\lfloor e/3\rfloor$ and
$A_e=I_{e,e+1}$. For $l=3a\le e$, the vector $e_{-e+2l}$ equals
$e_{-e}$ and
\[
 1+(1-q^{-1})\sum_{a=1}^{k}q^a=q^k.
\]
Thus the full stabilized integral is
\[
 J_e=q^{-2e/3}\bigl(q^k e_{-e}+q^{(e+1)/3}A_e e_{e+2}\bigr).
\]
Since $A_{e+3}=A_e$ and $e_{j+3}=e_j$, replacing $e$ by $e+3$
multiplies the parenthesis by $q$ and the exterior factor by $q^{-2}$.
Consequently
\begin{equation}\label{eq:v-whittaker-recurrence}
 J_{e+3}=q^{-1}J_e,\qquad W_{e+3}=q^{-1}W_e\quad(e\ge0).
\end{equation}
The second identity follows by applying the scalar coordinate
nonzero on $J_0$ and normalizing at one, in agreement with
\cite[Appendix~A, (A.9)--(A.11)]{Artin}.
Spherical upper-root invariance forces $W_e=0$ for $e<0$, and
the torus law and unit-diagonal invariance give
$W(a(\varpi^eu))=\chi_p(u)^{2e}W_e$. This proves the table,
including its phase and unit factor. Since $\chi_p(-1)=1$, the sign of
the Weyl representative or of the additive Gauss argument is irrelevant.

At a finite place $v\in S_f$ choose a smooth vector $f_v$ with
$\lambda_v(f_v)=1$. Local constancy gives a small compact-open unit
subgroup $U_v$ on which $W_{f_v}(a(y))=1$; shrink it into all the
required local kernels and sixth powers. Additive Fourier inversion
gives a Schwartz--Bruhat function $\phi_v$ with
\[
 \int\phi_v(x)\psi_v(yx)\,dx=1_{U_v}(y).
\]
For $f_{v,U}=\int\phi_v(x)\pi_v(n(x))f_v\,dx$, root equivariance gives
\[
 W_{f_{v,U}}(a(y))=W_{f_v}(a(y))1_{U_v}(y)=1_{U_v}(y).
\]
This is the direct cutoff. For a Weyl-rotated vector,
upper-root invariance gives compact upper support. Its small-norm
bound is the local estimate \cite[Section~5.2, (5.13)]{Artin};
the Jacquet-module recurrence yielding it is recorded below.

For the finite small-norm bound, set $z=a(\varpi_v^3)$. The
covering torus modulo its central cube subgroup is finite because
$F_v^\times/F_v^{\times3}$ is finite, including at wild places.
Its fixed-central-character torus constituents are finite-dimensional.
The two-cell Jacquet calculation
\cite[Proposition~I.2.1, p.~61; Theorem~I.2.9(e), p.~72]{KP}
gives the possible unnormalized exponents $1/2\pm1/6$.
Thus on the ordinary Jacquet module the eigenvalues of $z$ have
absolute values $q_v^{-1}$ or $q_v^{-2}$. Finite dimensionality
and exactness give an annihilating polynomial $P$ with those root
moduli. For any fixed vector $\phi$, ordinary Jacquet coinvariants give
\[
 P(\pi(z))\phi=\sum_{i=1}^h(\pi(n(x_i))-1)\phi_i.
\]
Apply $\lambda_v\pi(a(y))$. At sufficiently small $|y|_v$, every
$\psi_v(yx_i)$ is one and the right side vanishes. The cube cocycles
are trivial, so
\[
 P(E)W_\phi(y)=0,\qquad EW_\phi(y)=W_\phi(y\varpi_v^3).
\]
A small open subgroup of cube units fixing $\phi$ makes the unit
variable range over finitely many classes. On each valuation
progression the recurrence has solutions equal to root powers times
polynomials in the valuation. Consequently
\[
 |W_\phi(y)|\ll |y|_v^{1/3}(1+|\log|y|_v|)^A
                  \ll_\delta |y|_v^\delta\quad(0<\delta<1/3).
\]
This derives \cite[pp.~21--22, (5.13)]{Artin} for the fixed finite
family of direct and Weyl-rotated vectors; the small-norm threshold
is uniform over the finitely many unit classes.

At complex places the ordinary
principal series has parameters $\pm1/6$, giving the Bessel order $1/3$;
first-root derivatives supply the two angular types. The weighted
constant terms vanish by the pure-norm Jacquet-exponent calculation
in Lemma~\ref{v:lem:raw-theta}. These establish
the raw reflection; Proposition~\ref{v:prop:radial} proves the arbitrary
radial extension used in the moment recursions.
Finally Lemma~\ref{v:lem:primary-kummer-transport} verifies the primary
quadratic and cubic sieve hypotheses and transports their primitive
characters to the actual exterior symbols. No zero-free conclusion
of \cite{Artin} or strengthened inverse-count estimate is used here.
\end{proof}

By Proposition~\ref{v:prop:checked-inputs}, at an unramified finite place
\[
 W_{\mathfrak p}(\varpi_{\mathfrak p}^{3j})=q_{\mathfrak p}^{-j},\quad
 W_{\mathfrak p}(\varpi_{\mathfrak p}^{3j+1})
 =q_{\mathfrak p}^{-j-1/2}\gamma_2(\mathfrak p),\quad
 W_{\mathfrak p}(\varpi_{\mathfrak p}^{3j+2})=0\qquad(j\ge0).
\]
These Gauss factors determine
the coefficient characters of the two inverse transforms of
Section~\ref{sec:v-detail-10}; the following construction also fixes
the complex angular types needed in Section~\ref{sec:v-detail-2}.

\begin{lemma}[The theta coefficient characters]\label{v:lem:theta-characters}
There are a unitary Hecke character $\Lambda$ of complex angular degree
$1$ or $-1$ at each infinite place and a finite-order Hecke character
$G_0$, both unramified outside $S$, such that, for squarefree exterior $n$,
\[
 a_0(n)=\overline{\Lambda(n)}\gamma_2(n),\qquad
 a_0(n)\gamma_1(n)=\mu_F(n)G_0(n),\qquad
 \mu_F(n)\gamma_{-1}(n)=a_0(n)\overline{G_0(n)}.
\]
The identities retain every common natural zero. The character $\Lambda$
is unitary and need not have finite order.
\end{lemma}
\begin{proof}
We use \cite[Lemma 3.4, equations (3.9)--(3.14) and its
proof]{Artin}, with the same allowed generators and trace character.
Put $E_0=\mathbb Q(\omega)$, $\omega^2+\omega+1=0$. If
$\pi_+(\mathfrak a)$ is the unique generator of an $E_0$-ideal prime
to $3$ satisfying $\pi_+(\mathfrak a)\equiv1\pmod3$, set
\[
 \Lambda_{E_0}(\mathfrak a)
      =\frac{\jmath(\pi_+(\mathfrak a))}{|\jmath(\pi_+(\mathfrak a))|},\qquad
 \Lambda=\Lambda_{E_0}\circ N_{F/E_0}.
\]
Here $\jmath|_{E_0}$ is the residue embedding fixed above, also
used for the cubic character and Jacobi sum. The conjugate genuine
embedding $\iota$ is used only on the covering center.
The primary generators multiply, so $\Lambda_{E_0}$ is the unitary
Hecke character of conductor dividing $(3)$ with this ideal-character
infinity type.

For an exterior prime $p$ above $\mathfrak q$ of
$E_0$, of residue degree $f$, the cubic residue character is the
residue-norm lift of that at $\mathfrak q$. The Gauss--Jacobi identity
and Davenport--Hasse lifting formula give
\[
 \begin{split}
 \gamma_2(p)^3
 &=q_p^{-1/2}J_{\mathbb F_{q_p}}(\rho,\rho),\\
 J_{\mathbb F_{q_p}}(\rho,\rho)
 &=(-1)^{f-1}(-\jmath(\pi_+(\mathfrak q)))^f
   =-\jmath(\pi_+(\mathfrak q))^f,
 \end{split}
\]
and therefore $\gamma_2(p)^3=-\Lambda(p)$. The base Jacobi identity
$J(\rho,\rho)=-\jmath(\pi_+(\mathfrak q))$ and the lifting hypotheses
($\rho,\rho,\rho^2$ all nontrivial) are verified in the cited proof.
At every complex place the local norm to $E_0$ is a coordinate or
its conjugate, so the angular exponent of $\Lambda$ there is
$\pm1$, with no radial component.

For the quadratic phase write
$\mathfrak w_v(a)=\operatorname{WeilIndex}_{\psi_v}(ax^2)/
\operatorname{WeilIndex}_{\psi_v}(x^2)$ for $a\in F_v^\times$, where
$\psi=\prod_v\psi_v$ is the trace character of $\mathbb A_F/F$. It
depends only on the square class of $a$, and Weil's product formula and
the Hilbert-symbol identity \cite[Chapter~II, \S\S28 and~30]{Weil64} are
\begin{equation}\label{eq:v-weil-index}
 \prod_v\mathfrak w_v(a)=1\quad(a\in F^\times),\qquad
 \mathfrak w_v(a)\mathfrak w_v(b)
       =(a,b)_{2,v}\mathfrak w_v(ab).
\end{equation}
Put $\mathfrak w_S=\prod_{v\in S}\mathfrak w_v$.

\emph{The local evaluation.} Let $p$ be the allowed generator of an
exterior prime $\mathfrak p$, $q=q_p$. At an exterior place $w$ the
residue characteristic is odd and $\OO_w$ is self-dual for $\psi_w$,
since $S$ contains the places above $2$ and the support of the
different. The local Gauss integral for the Weil index
\cite[Chapter~II, \S27]{Weil64} gives $\mathfrak w_w(u)=1$ for every
unit $u$ of $F_w$, and at $w=\mathfrak p$
\[
 \mathfrak w_{\mathfrak p}(p^{-1})
 =q^{-1/2}\sum_{x\bmod p}\psi_{\mathfrak p}(x^2/p)
 =q^{-1/2}\sum_{x\bmod p}\chi_p(x)^3\psi_{\mathfrak p}(x/p).
\]
In the second equality, $\chi_p^3$ is the quadratic character of
$\OO/\mathfrak p$, so $x$ has $1+\chi_p(x)^3$ square roots, and the
complete sum of the nontrivial additive character $x\mapsto
\psi_{\mathfrak p}(x/p)$ vanishes. For $x\in R$ the global character
is trivial on $x/p\in F$, and $x/p$ is integral at every exterior place
other than $\mathfrak p$; hence $e(x/p)=\prod_{v\in S}\psi_v(x/p)
=\psi_{\mathfrak p}(-x/p)$. Since $\chi_p(-1)=1$, the last sum is
$\gamma_3(p)$. Applying the product formula in \eqref{eq:v-weil-index} to
$p^{-1}\in F^\times$, every exterior place other than $\mathfrak p$
contributes $1$, so
\begin{equation}\label{eq:v-gamma3-weil}
 \gamma_3(p)=\mathfrak w_S(p^{-1})^{-1}.
\end{equation}

\emph{A character of $B_S$.} Each $\mathfrak w_v$ factors through
$F_v^\times/F_v^{\times2}$, hence through $\mathcal V_S$. For $a,b\in B_S$,
\eqref{eq:v-weil-index} gives
$\mathfrak w_S(a)\mathfrak w_S(b)=\mathfrak w_S(ab)\prod_{v\in S}(a,b)_{2,v}$,
and $\prod_{v\in S}(a,b)_{2,v}=H_S(a,b)^3=1$ by isotropy of $B_S$.
Thus $\lambda(b)=\mathfrak w_S(b^{-1})^{-1}$ is a character of $B_S$.
It is trivial on squares, so $\lambda^2=1$. Extend $\lambda$ to
$\mathcal V_S=E_S\oplus B_S$ by $\lambda|_{E_S}=1$.

\emph{Extension to a Hecke character.} For every $v\in S_f$ the
subgroup $F_v^{\times6}$ is open of finite index in $F_v^\times$, so
there is $e_v\ge1$ with $1+\mathfrak p_v^{e_v}\subset F_v^{\times6}$.
Put $\mathfrak m=\prod_{v\in S_f}\mathfrak p_v^{e_v}$. Every infinite
place of $F$ is complex, so no sign conditions occur there.
For an exterior ideal $\mathfrak a$ with allowed generator $a$ put
$G'(\mathfrak a)=\lambda([a])$. This is well defined, since $a$ is
unique modulo $\Gamma=U_S^6\subset F_S^{\times6}$, and multiplicative,
since products of allowed generators are allowed
(Lemma~\ref{v:lem:allowed-reciprocity}). Let $\alpha\in F^\times$ satisfy
$\alpha\equiv1\pmod{^\times\mathfrak m}$, prime to $S_f$. Then
$\alpha_v\in F_v^{\times6}$ for every $v\in S$, so $[\alpha]=1$ in
$\mathcal V_S$. If $a=\epsilon\alpha$, $\epsilon\in U_S$, is the
allowed generator of $\alpha R$, then $[a]=[\epsilon]\in E_S\cap
B_S=\{1\}$, so $G'(\alpha R)=1$. Hence $G'$ is a character of the
finite ray class group of modulus $\mathfrak m$: a Hecke character of
order dividing $2$, unramified outside $S$, with
$G'(\mathfrak p)=\gamma_3(p)$ at every exterior prime by
\eqref{eq:v-gamma3-weil}.

Since $4\in U_S$, \eqref{eq:v-generator-unit-reciprocity} gives
$\chi_a(4)=H_S(4,a)^{-1}$; this is $\xi_\epsilon(\mathfrak a)$ for the
representative $\epsilon\in\mathcal U_6$ of $4U_S^6$, a character in
$\Upsilon_S$, which we denote $\xi_4$; it is again trivial on the principal ideals
$\alpha R$ above, because $H_S(4,\cdot)$ is trivial on $F_S^{\times6}$.
Put $G_0=G'\,\overline{\xi_4}$. Then $G_0$ is a character of the ray class
group of modulus $\mathfrak m$, of order dividing six, unramified outside
$S$, and
\[
 G_0(p)=\overline{\chi_p(4)}\gamma_3(p)
\]
at every exterior prime. Finite-field sextic duplication gives
\[
 \gamma_1(p)\gamma_2(p)
 =\overline{\chi_p(4)}\gamma_3(p)\gamma_2(p)^3
 =-\Lambda(p)G_0(p).
\]

For coprime squarefree indices the CRT multiplier for
$\gamma_1\gamma_2$ is
$\chi_{c_1}(c_2)^3\chi_{c_2}(c_1)^3=1$, by allowed-generator
reciprocity. Consequently, on squarefree exterior ideals,
\[
 a_0(n)\gamma_1(n)=\mu_F(n)G_0(n),\qquad
 \mu_F(n)\gamma_{-1}(n)=a_0(n)\overline{G_0(n)},
 \qquad a_0(n)=\overline{\Lambda(n)}\gamma_2(n).
\]
For the second identity use
$\gamma_{-1}(n)=\overline{\gamma_1(n)}$, which follows from
$\chi_n(-1)=1$, and $|a_0(n)|=1$. The identities persist
after multiplication by any common zero-extended coprimality or
character factor; no character is divided by at a nonunit.

The finite-order coefficient characters, $G_0$, and the finite unit
characters $\Upsilon_S$ lie in one fixed finite group; $G_0$ and
$\Upsilon_S$ are characters of the ray class group of modulus
$\mathfrak m$ constructed above.
The two Gauss transformations send $\nu\mapsto\nu G_0\mapsto
\nu G_0\overline{G_0}=\nu$, with the separately retained unit
twists. For a fixed finite-order target $\eta$, its conductor lies in
$S$ and its finite character in this group by the
preliminary choice of Section~\ref{sec:v-scope}. Its complex local components are trivial,
so the angular types of $\overline\Lambda\overline G_0\eta$
remain $\pm1$.
\end{proof}

\section{Theta functions and the radial functional equation}
\label{sec:v-detail-2}
In this section, Lemma~\ref{v:lem:raw-theta} proves the unsmoothed
norm-sheet theta identity \eqref{eq:v-raw-theta-identity}, the local
Weyl formula \eqref{eq:v-local-weyl-formula} records the row-prime rules,
and Proposition~\ref{v:prop:radial} extends the norm-sheet reflection of
Proposition~\ref{v:prop:checked-inputs} to compactly supported smooth
radial tests, with finite seminorm bounds and polynomial dependence on
norm twists; these justify the profile separation
and discarded-tail estimates in
Sections~\ref{v:sec:unmarked-assembly} and \ref{v:sec:M2-assembly}.

Write
$a(y)=\operatorname{diag}(y,1)$ and $n(x)=\bigl(\begin{smallmatrix}1&x\\0&1\end{smallmatrix}\bigr)$,
using the sections of Proposition~\ref{v:prop:checked-inputs}.
For the ordinary torus section, \eqref{eq:v-effective-torus-cocycle}
gives the $S$-scalar
\[
 \mathcal C_S(\nu,y)=H_S([\nu],[y])^2.
\]
It is one on the direct support $[\nu y]=1$ by alternation. In the
reflected term it is $\mathcal C_S(\nu,y^{-1})$, which becomes
$H_S([\nu],[v])^2$ under $v=\nu/y$.

Fix a finite-order target $\eta$ and put
$\vartheta=\bar\Lambda\bar G_0\eta$. For $m\ne0$, choose an allowed
$m_0$ generating $(m)_S$ and write $m=\epsilon_m m_0$;
$\upsilon=[\epsilon_m]\in U_S/\Gamma$. Let $L_\upsilon$ be the
unitary $S$-character, trivial on $U_S$, characterized on allowed $A$ by
\[
 L_\upsilon(A)=\vartheta((A)_S)\chi_A(\epsilon_m),
 \qquad L_\upsilon=L_S H_S(\epsilon_m,\cdot)^{-1},
\]
where $L_S(A)=\vartheta((A)_S)$, as in
\cite[Section 3, Lemma 3.4 and Section 6]{Artin}.
Choose its theta vector with unrotated finite cutoff and matching
complex angular type; denote its automorphic function by
$\theta_\upsilon$ and its Weyl-rotated $S$-Whittaker function by
$f_{w,\upsilon}$. Normalize the direct function by
\begin{equation}\label{eq:v-exact-S-vector}
 f_\upsilon(y)=1_{\mathcal U}(y_f)
       \prod_{v\mid\infty}|y_v|_v\omega_{\upsilon,v}(y_v)
                 K_{1/3}(c_v\sqrt{|y_v|_v}).
\end{equation}
Here $\mathcal U\subset\prod_{v\in S_f}F_v^{\times6}$ is a fixed
compact open subgroup contained in the kernels of all finite
components of $L_\upsilon$, and $\omega_{\upsilon,v}$ has angular
degree $1$ or $-1$, matching $L_\upsilon$. The scalar in
\eqref{eq:v-exact-S-vector} is one; all tensor coefficients after
Weyl rotation are included in $f_{w,\upsilon}$. Counting measure
on $\Gamma$ induces the quotient sheet measure.

Put
$X_x=\{v\in F_S^\times:|v|_S=x\}$,
with multiplicative sheet measure $d_x^\times v$. Define the direct profile by
\begin{equation}\label{eq:v-direct-profile-definition}
 V(x)=\int_{X_x}1_{\{[v]=1\}}L_\upsilon(v)^{-1}
                         f_\upsilon(v)\,d_x^\times v.
\end{equation}
The complex components of $L_\upsilon$ are pure angular characters,
with no norm factor: $\Lambda$ is the norm pullback constructed in
Lemma~\ref{v:lem:theta-characters}, and finite-order complex characters
are trivial. Hence the angular factors in this integral cancel.
With the fixed Haar measures, let $c_{V,S}>0$ be the product of the finite
unit-neighborhood and angular volumes. Then
\[
 \widehat V(t)=c_{V,S}\prod_{v\mid\infty}
 2^{2t+1}c_v^{-2t-2}\Gamma(5/6+t)\Gamma(7/6+t),
 \qquad\Re t>-5/6.
\]
Hence $V$ is independent of the row-unit class $\upsilon$.

For $p\mid(m)_S$ set $j_p\equiv v_p(m)\pmod6$, $0\le j_p<6$, and define
\begin{equation}\label{eq:v-local-digit-weights}
 \alpha_{p,j}(z)=
 \begin{cases}
 \tau_{p,j}^{-1}\chi_p(z)^{-j},&j\ne0,\ z\ne0,\\
 0,&j\ne0,\ z=0,\\
 1-q_p^{-1},&j=0,\ z=0,\\
 -q_p^{-1},&j=0,\ z\ne0,
 \end{cases}
 \qquad
 \tau_{p,j}=\sum_{z\ne0}\chi_p(z)^{-j}e_p(z/p).
\end{equation}
Then $\sum_z\alpha_{p,j}(z)e_p(tz/p)=\chi_p(t)^j$,
including $t=0$ with its natural zero. Define, for $\nu\in F^\times$,
\[
 \mathcal C_{m,p}(\nu)=
 \begin{cases}
 \displaystyle\sum_{z\bmod p}\alpha_{p,j_p}(z)
                  W_p(a(\nu)w n(z/p)),&p\mid(m)_S,\\
 W_p(a(\nu)),&p\nmid(m)_S.
 \end{cases}
\]
For $\kappa\in\mathcal V_S$, set
\begin{equation}\label{eq:v-raw-profile-definition}
 H_{\upsilon,\kappa}(x)=\int_{X_x}
 1_{B_S}(\kappa[v]^{-1})H_S(\kappa,[v])^2
                 L_\upsilon(v)f_{w,\upsilon}(v)\,d_x^\times v.
\end{equation}
Finite tensor decompositions are included in $f_{w,\upsilon}$, so this
is a fixed finite family indexed by $\upsilon,\kappa$.
For $\nu$ modulo $\Gamma$, define all the reflected data by
\begin{equation}\label{eq:v-raw-dual-data}
 a_m^\vee(\nu)=L_\upsilon(\nu)^{-1}
                         \prod_{p\notin S}\mathcal C_{m,p}(\nu),
 \qquad \rho_{m,\nu}=|\nu|_S,\qquad \alpha_\nu=[\nu]\in\mathcal V_S.
\end{equation}
For fixed $m,\nu$ the product has finitely many nonunit factors, and
it vanishes unless $v_p(\nu)\ge0$ for $p\nmid(m)_S$ and
$v_p(\nu)\ge-2$ for $p\mid(m)_S$. Thus every nonzero frequency has
an integral exterior numerator over a divisor of
$\operatorname{rad}((m)_S)^2$.

\begin{lemma}[The unsmoothed norm-sheet theta identity]\label{v:lem:raw-theta}
With the normalized theta datum above, let
\[
 \mathcal T_m(l)=\sum_{\substack{c\ \sfq,n}}
 \frac{\gamma_2(c)\vartheta(cn^3)\chi_{cn^3}(m)}{\sqrt{q_c}q_n}
                          V(q_{cn^3}l).
\]
Here $V$ is \eqref{eq:v-direct-profile-definition}, with
the Haar normalization and finite cutoff stated there. Then
\begin{equation}\label{eq:v-raw-theta-identity}
 \mathcal T_m(l)=\sum_{\nu\in F^\times/\Gamma}
       a_m^\vee(\nu)
       H_{\upsilon,\alpha_\nu}(\rho_{m,\nu}/l).
\end{equation}
Both sides converge absolutely after the sheet integration. The
frequency-unit class is summed independently of the fixed row-unit
class $\upsilon$.
\end{lemma}
\begin{proof}
Fix $m\ne0$, write $L=L_\upsilon$, and let $\Omega_l$ be a
fundamental set for $X_l/\Gamma$, with counting measure on $\Gamma$.
Let $g(y,\mathbf z)$ have components $a(y)$ at $S$,
$n(z_p/p)$ at row primes and identity elsewhere, and put
\[
 \mathscr D_m^+(y)=\sum_{\mathbf z}\Bigl(\prod_{p\mid m}
       \alpha_{p,j_p}(z_p)\Bigr)\theta_\upsilon(g(y,\mathbf z)),
 \qquad
 \mathscr D_m^-(y)=\sum_{\mathbf z}\Bigl(\prod_{p\mid m}
       \alpha_{p,j_p}(z_p)\Bigr)\theta_\upsilon(wg(y,\mathbf z)).
\]
Rational automorphy gives $\mathscr D_m^+=\mathscr D_m^-$.
For these digit-averaged theta functions, the
change $z'_p\equiv\gamma z_p\pmod p$, $\gamma\in\Gamma$, gives
\[
 g(\gamma y,\mathbf z')=a(\gamma)g(y,\mathbf z)k^{-1},\qquad
 wg(\gamma y,\mathbf z')=a(\gamma^{-1})wg(y,\mathbf z)k^{-1}.
\]
Here $k_S=1$ and the exterior components of $k$ are compact and fix
the original spherical vector. Since
$\alpha_{p,j}(\gamma z)=\alpha_{p,j}(z)$, rational automorphy proves
$\Gamma$-invariance of both digit averages. Their constant terms
descend as well: rational diagonal conjugation rescales the unipotent
variable by a rational scalar, whose adelic modulus is one.
The weights $1_{B_S}([y])L(y)^{-1}$ are $\Gamma$-invariant.

At a complex place the constant-term Jacquet characters are pure
norm characters: in this normalization their possible powers are
$|y_v|_v^{1/3}$ and $|y_v|_v^{2/3}$, with an absent constituent given
coefficient zero. This holds for the chosen vector and its Weyl
rotation. The class weight is rotation invariant, since
$\mathbb C^\times/\mathbb C^{\times6}=1$, whereas $L(y)^{-1}$ has
angular degree $-m_v$, $m_v\in\{1,-1\}$. Thus each weighted constant
term contains $\int_0^{2\pi}e^{-im_v\theta}\,d\theta=0$.
The Whittaker function instead contributes $e^{im_v\theta}$, which
cancels $L(y)^{-1}$ and survives angular integration.

We justify all exchanges absolutely. If $F_\upsilon^{\rm abs}$ and
$F_{w,\upsilon}^{\rm abs}$ are the sums of absolute values of the
fixed tensor presentations, their sheet integrals satisfy
$H^{\rm abs}(x)\ll_\sigma x^{-\sigma}$ for $\sigma>1/2$.
At an exterior prime,
\[
 \sum_{e\ge0}|W_p(a(p^e))|q_p^{-\sigma e}
 =\frac{1+q_p^{-\sigma-1/2}}{1-q_p^{-3\sigma-1}}.
\]
The exterior product away from $m$ is bounded by
$\zeta_F^S(\sigma+1/2)\zeta_F^S(3\sigma+1)$.

At a row prime the total digit mass is at most $q_p^{1/2}$ when
$j_p\ne0$, and less than $2$ when $j_p=0$. A nonzero digit shifts
the spherical numerator valuation by two; its weighted local sum
therefore costs at most an additional $q_p^{2\sigma}$.
The finitely many row-prime factors give a finite constant depending
on $m$ and $\sigma$. Tonelli now bounds the direct and reflected
integrated absolute Fourier series by
$O_{m,\sigma}(l^{-\sigma})$ and $O_{m,\sigma}(l^\sigma)$,
respectively. Hence digit evaluation and signed orbit
unfolding are legitimate for this fixed row.

In the direct expansion the cutoff is supported on $[\nu y]=1$.
Together with $[y]\in B_S$, this forces $[\nu]\in B_S$.
If $A$ is the allowed generator of the exterior ideal of $\nu$,
then $\nu/A\in U_S$; since $E_S\cap B_S=1$, it lies in $\Gamma$.
The direct cocycle is one on this support. Substitution
$v=\gamma Ay$ unfolds $\Omega_l$ into $X_{q_A l}$, giving
$L(A)V(q_A l)$ under the specified finite $S$ normalization.
The spherical coefficient is supported at $A=cn^3$ and equals
$\gamma_2(c)/(\sqrt{q_c}q_n)$. The digit identities give
\[
 L(A)\prod_{p\mid m}\chi_p(A)^{j_p}
   =\vartheta(A)\chi_A(m),
\]
including shared-prime zeros. This is $\mathcal T_m(l)$.

In the reflected expansion the exterior coefficient is
$\prod_{p\notin S}\mathcal C_{m,p}(\nu)$, invariant under
$\nu\mapsto\gamma\nu$, $\gamma\in\Gamma$:
the spherical and digit-averaged local monomials have only sextic
unit characters. For each representative $\nu_0$, absolute unfolding
is justified by
\[
 \sum_{\gamma\in\Gamma}\int_{\Omega_l}
 F_{w,\upsilon}^{\rm abs}(\gamma\nu_0/y)\,d_l^\times y
 =\int_{X_{|\nu_0|_S/l}}
 F_{w,\upsilon}^{\rm abs}(v)\,d_{|\nu_0|_S/l}^\times v.
\]
Indeed, the sets $\gamma\Omega_l^{-1}$ partition $X_{1/l}$.
Under $v=\gamma\nu_0/y$, the class weight, character and cocycle become
\[
 [y]=[\nu_0][v]^{-1},\qquad
 L(y)^{-1}=L(\nu_0)^{-1}L(v),\qquad
 \mathcal C_S(\gamma\nu_0,y^{-1})=H_S([\nu_0],[v])^2.
\]
The last equality uses $[\gamma]=1$ and alternation.
Multiplicative sheet measure has Jacobian one, so this orbit
contributes
\[
 L(\nu_0)^{-1}\prod_{p\notin S}\mathcal C_{m,p}(\nu_0)
 H_{\upsilon,[\nu_0]}(|\nu_0|_S/l).
\]
Summing the absolutely convergent orbit expansion proves the identity.
The finite $S$ tensor coefficients were included in $f_{w,\upsilon}$,
so there is no additional scalar or unit-class multiplicity.
\end{proof}

For a unit $t$ put
\[
 \mathcal R_{p,j,e}(t)=\sum_{z\ne0}\alpha_{p,j}(z)
                 W_p(a(p^{e-2}t)w n(z/p)).
\]
For $H\ne0$, direct matrix multiplication gives
\begin{equation}\label{eq:v-bruhat-identity}
 wn(H)=HI\,n(-H^{-1})a(H^{-2})\bar n(H^{-1}),\qquad
 \bar n(x)=\begin{pmatrix}1&0\\x&1\end{pmatrix}.
\end{equation}
The scalar $HI$ is split and trivial in the quotient. Substitution
in the chosen Kubota cocycle gives the same lift identity: the
remaining symbols have cube arguments, powers of $-1$, or two
powers of $H$; the last case is one by alternation.
For $H=z/p$ the final factor is compact. The unipotent argument
is initially $-p^{e-1}t/z$; replacing $z$ by $-z$ gives the
positive argument below, since $\alpha_{p,j}(-z)=\alpha_{p,j}(z)$.
Thus the local formula is
\begin{equation}\label{eq:v-local-weyl-formula}
 \mathcal R_{p,j,e}(t)=W_p(a(p^e))\chi_p(t)^{4+2e}
 \sum_{z\ne0}\alpha_{p,j}(z)\chi_p(z)^4e_p(p^{e-1}t/z).
\end{equation}
Indeed the cocycle of $a(p^{e-2}t)a(p^2z^{-2})$ contributes
$\chi_p(t)^4\chi_p(z)^{4(e-2)}$; the remaining spherical unit
factor is $\chi_p(t)^{2e}\chi_p(z)^{-4e}$.

Thus for $e=0$ the unit exponent is $8-j$ modulo six; for $e>0$
it is $4+2e$, and the digit sum vanishes unless $j=4$.
The spherical factor further forces $e\equiv0,1\pmod3$.
For $j=1$ only $e=0$ survives, with exponent one and magnitude one.
For an inserted full-index factor $1_{p\mid A}=1-\chi_p(A)^0$,
the inactive coefficient is $q_p^{-1}$ and the active coefficient
has magnitude $q_p^{-1/2}$, exponent two and numerator valuation zero.
This proves the local rules used in
Lemma~\ref{v:lem:marked-terminal-typed} and the width tables below.

At a complex place let \(v^\circ\) be the compact-invariant spherical vector.
For \(z\in\mathbb C\), write
\[
 E(z)=\begin{pmatrix}0&z\\0&0\end{pmatrix},\qquad
 F(z)=\begin{pmatrix}0&0\\z&0\end{pmatrix},\qquad
 w=\begin{pmatrix}0&-1\\1&0\end{pmatrix}.
\]
Since \(E(z)-F(\bar z)\in\mathfrak{su}(2)\),
\[
 d\pi(E(z))v^\circ=d\pi(F(\bar z))v^\circ.
\]
Also \(\pi(w)v^\circ=v^\circ\) and
\(\operatorname{Ad}(w)E(z)=-F(z)\). Hence \(w\) exchanges the two
first-root derivative vectors of angular types \(\pm1\), up
to fixed nonzero constants. The source Whittaker functions are
\[
 W_{v,m}(a(y))=|y|_v\omega_m(y)K_{1/3}(c_v|y|_v^{1/2}),
 \qquad m=\pm1.
\]
Pairing the rotated vector with the reversed angular weight therefore gives
the same radial Bessel factor as before rotation.

Let \(r_F=[F:\mathbb Q]/2\), and denote the direct and reflected norm-sheet
profiles by \(V\) and \(H_{\upsilon,\alpha}\). The finite index
\(\alpha\) retains the actual local power classes. The finite-place support
and small-norm bounds in \cite[Section 5]{Artin} imply that their fixed
finite-place Mellin factors are bounded and holomorphic on each closed strip
\(\Rea t>-\delta\), for any fixed \(\delta<1/3\). Thus
\[
 m_{\upsilon,\alpha}(t)
 =\frac{\widehat H_{\upsilon,\alpha}(t)}{\widehat V(-t)}
 =E_{\upsilon,\alpha}(t)A_1^t
 \left[\frac{\Gamma(5/6+t)\Gamma(7/6+t)}
 {\Gamma(5/6-t)\Gamma(7/6-t)}\right]^{r_F}.
\]

For the finite factors, the direct finite test is $1_U$ on a product of
unit neighborhoods, with Mellin transform the nonzero constant
$\operatorname{vol}(U)$. For each reflected finite test $h_v$, the
source bound gives
\[
 h_v(y)=0\quad(|y|_v>C_v),\qquad
 |h_v(y)|\ll_{\delta}|y|_v^\delta\quad(|y|_v\le1)
\]
for every fixed $\delta<1/3$. Hence
\[
 e_v(t)=\int_{F_v^\times}h_v(y)|y|_v^t\,d^\times y
\]
converges absolutely and locally uniformly for $\Re t>-\delta$, and
is bounded on closed substrips independently of
$\Im t$, since $||y|_v^{i\Im t}|=1$; likewise for
finite sums of products of local tests. Division by the direct
finite factor creates no pole or zero denominator.

At a complex place, put $x=|y|_v$. After angular integration the
local Mellin factor is, for $\Re t>-5/6$,
\[
 \begin{split}
 \int_0^\infty xK_{1/3}(c_v\sqrt{x})x^t\,\frac{dx}{x}
 &=2^{2t+1}c_v^{-2t-2}
                   \Gamma(5/6+t)\Gamma(7/6+t).
 \end{split}
\]
This follows by $u=c_v\sqrt{x}$ and the Mellin integral of $K_{1/3}$.
The reversed angular weight and the Weyl relation above give the same
radial factor on the reflected side, up to a fixed constant. Taking
the quotient at $t$ and $-t$ gives the displayed Gamma ratio and a
positive real exponential factor $A_1^t$. Reciprocal Gamma functions
are entire and the numerator poles start at $t=-5/6$, so the quotient
is holomorphic on every strip considered. Here $E$ is
the bounded finite-place factor above, including the
fixed normalization constants. Stirling gives
\[
 |m_{\upsilon,\alpha}(\sigma+i\tau)|
 \ll_\sigma(1+|\tau|)^{4r_F\sigma}.
\]
The sign \(-t\) in the denominator is essential.

Use Lemma~\ref{v:lem:raw-theta} before Gaussian convolution.
For an integral exterior ideal $\mathfrak a$, define $a_m(\mathfrak a)$
to be its displayed direct coefficient when $\mathfrak a=c n^3$ and
zero otherwise. In the following Dirichlet series the dual index $b$
is the frequency $\nu\in F^\times/\Gamma$ and all its data are
those of \eqref{eq:v-raw-dual-data}. Put
\[
 D_m(s)=\sum_{\mathfrak a}a_m(\mathfrak a)q_{\mathfrak a}^{-s},\qquad
 D^\vee_{m,\alpha}(s)=
       \sum_{\nu:\alpha_\nu=\alpha}a_m^\vee(\nu)\rho_{m,\nu}^{-s}.
\]
Mellin integration of \eqref{eq:v-raw-theta-identity} on the respective
absolute-convergence half-planes gives
\[
 \widehat V(s)D_m(s)
    =\widehat{\mathcal T_m}(s)
    =\sum_\alpha\widehat H_{\upsilon,\alpha}(-s)
                         D^\vee_{m,\alpha}(-s),
 \qquad D_m(s)=\sum_\alpha m_{\upsilon,\alpha}(-s)
                         D^\vee_{m,\alpha}(-s).
\]

The theta integral is entire: its direct and reflected expressions
give rapid decay at the two endpoints for a fixed nonzero row.
The direct and reflected Dirichlet series give polynomial bounds
on the two boundary lines. Inside the strip, division by the
nonvanishing direct Gamma product gives at most
$\exp(C|\Im s|)$ growth. To interpolate, divide by $(s+c)^N$
with $-c$ outside the strip and $N$ larger than both boundary
orders, apply the maximum principle to
$e^{\varepsilon s^2}D_m(s)/(s+c)^N$ on rectangles, and then let
$\varepsilon\downarrow0$. The horizontal sides vanish before this
limit. Hence $D_m$ has polynomial growth throughout the strip,
which justifies shifting compact smooth Mellin inversion. For a
scale $D>0$ this gives
\[
 \sum_n a_m(n)W(q_n/D)
 =\sum_{\nu\in F^\times/\Gamma}a_m^\vee(\nu)
                  \Phi_{W,\upsilon,\alpha_\nu}(D\rho_{m,\nu}).
\]
All row-dependent powers remain in the displayed reflected scale.

\begin{proposition}[Functional equation for radial test functions]\label{v:prop:radial}
For \(W\in C_c^\infty(I)\), \(I\Subset(0,\infty)\), the actual reflected
coefficients and coprimality indicators of the theta functions in the source are unchanged; their
Gaussian profiles may be replaced by
\[
 \Phi_{W,\upsilon,\alpha}(x)=\frac1{2\pi i}\int_{(A)}
 \widehat W(-t)m_{\upsilon,\alpha}(t)x^{-t}\,dt.
\]
For every \(A>0\), integer \(j\ge0\), and \(0<a_0<1/3\), there is
an integer \(J\) such that
\[
 |(x\partial_x)^j\Phi_{W,\upsilon,\alpha}(x)|
 \ll\|W\|_{C^J(I)}\min(x^{a_0},x^{-A}).
\]
The constants may depend on the fixed field and finite family of theta functions, but are uniform
in the row and in the fixed profile indices. For each fixed
$A,j,a_0$ and fixed field, replacing $W(x)$ by $W(x)x^{iT}$ costs at
most $(1+|T|)^J$, after increasing $J$ if necessary. The degree $J$
may depend on the decay order $A$, derivative order $j$, $a_0$ and $F$.
\end{proposition}
\begin{proof}
Use the theta Mellin identity above (not inverse-Gaussian convolution
of an arbitrary compact smooth function). Integration by parts in
\(\widehat W\) gives arbitrary polynomial vertical decay with finite
seminorm control. Shift the profile contour to \(-a_0\) for small \(x\)
and to \(A\) for large \(x\), choosing \(\delta>a_0\). The multiplier
bounds make the horizontal edges vanish and permit Euler derivatives under
the integral. The coefficient identity first holds row by row; uniform
bounds follow from the fixed multiplier family, rather than from a uniform
row-dependent entire-continuation estimate.
\end{proof}

\begin{lemma}[Identification of the theta data used in the moments]
\label{v:lem:theta-moment-identification}
The direct and reflected coefficients used in the inverse moments below
are those of Proposition~\ref{v:prop:checked-inputs} and
\eqref{eq:v-raw-theta-identity}, with the same trace character,
allowed generators and finite local data. Passing to radial tests
changes their profiles, not their arithmetic coefficients or their
zero sets. The terminal quadratic and cubic estimates are
\eqref{v:eq:transport-quadratic} and \eqref{v:eq:transport-cubic}
for the literal exterior symbols, after the separations specified in
Lemma~\ref{v:lem:hybrid-sieve}.
\end{lemma}
\begin{proof}
Fix the exterior spherical normalization $W_p(1)=1$, the two
archimedean derivative vectors and the finite family of cutoff
vectors at $S_f$ as in Proposition~\ref{v:prop:checked-inputs}.
The direct coefficient is $a_m(\mathfrak a)$ defined above, including
its zero value outside the displayed cube classes. In the reflected
coefficient retain the frequency $\nu$, its finite index
$\alpha_\nu$, its scale $\rho_{m,\nu}$ and all the data of
\eqref{eq:v-raw-dual-data}. The Mellin identity preceding
Proposition~\ref{v:prop:radial} is an identity for these coefficients.
Its multiplication by $\widehat W$ and inversion therefore does not
replace the finite Gauss correction or the angular character by a
norm function. Unit representatives and primes of $S$ remain fixed
throughout this operation.

For the subsequent sieve application, the quadratic family is the
primitive family $\Xi_{\mathfrak a}$ of
Lemma~\ref{v:lem:primary-kummer-transport}; its conductor outside $S$
is exactly $\mathfrak a$. The cubic family is the primitive family
$\omega_{\varepsilon_{\mathfrak a}}\zeta_{\mathfrak a}$ of
\eqref{eq:v-cubic-character}. On each fixed unit-class slice its extra
factor is a column twist independent of the varying row. The respective
source estimates are \cite[Theorem~1.1]{GL} and
\cite[Theorem~1.3, $n=3$]{BGL}, applied to the full families before
row deletion. Shared-prime zeros agree with the primitive conductors;
auxiliary decomposition ideals introduce no zeros. Finally the
coprimality expansion and quadratic estimate in
Lemma~\ref{v:lem:hybrid-sieve} precede the cubic estimate. Thus its
remaining coefficient is independent of the varying cubic row.
This order, rather than a bound for arbitrary coefficients depending
on both indices, supplies the terminal sieve application.
\end{proof}

\subsection{Relation to the classical cubic theta function}
\label{v:ss:classical-theta-comparison}

\indent{\it Representation and choice of vector}.
The classical cubic theta function over
$E_0=\mathbb Q(\omega)=\mathbb Q(\sqrt{-3})$ is a classical
realization of the rank-two cubic exceptional theta construction;
see \cite[Introduction, pp.~36--37]{KP}. Its explicit Fourier
coefficients and cusp expansions are recorded in
\cite[Section~5.1 and Appendix~A]{DR24}. An automorphic representation
contains many vectors. The classical spherical function and the
angular derivative vectors with finite cutoffs used here are not
asserted to be the same normalized automorphic function.

Our argument requires the direct and reflected expansions of one
fixed choice of theta vectors, with matching arithmetic coefficients
and measures. These are Proposition~\ref{v:prop:checked-inputs},
Lemma~\ref{v:lem:raw-theta}, and
Lemma~\ref{v:lem:theta-moment-identification}. Equality with the
classical spherical function, including all of its cusp constants,
is not an additional input to these statements. The standing
assumption $\mu_{12}\subset F$ also excludes $F=E_0$ itself; the
comparison with $E_0$ explains the theta construction, rather than
specializing every assertion of this part verbatim to $E_0$.

\indent{\it Exterior coefficients and additive characters}.
At exterior primes the normalization is $W_p(1)=1$ and
\eqref{eq:v-whittaker-recurrence}. Multiplying the local values,
including their unit factors, gives the direct coefficient
$\gamma_2(c)/(\sqrt{q_c}q_n)$ at $A=cn^3$, with $c$ squarefree;
$c$ and $n$ may share primes. The trace additive character is fixed
throughout this construction.

For comparison, in the conjugated classical expansion write
$\check e(z)=\exp(2\pi i\operatorname{Tr}_{E_0/\mathbb Q}z)$ and
$\lambda=\sqrt{-3}$. With primary generators, the infinity-cusp
coefficient in the normalization of \cite[(5.7), (5.16)]{DR24} is
\[
 d_0(\lambda^{-3}cn^3)=3^{5/2}|n|\widetilde g_3(c),\qquad
 \widetilde g_3(c)=q_c^{-1/2}\sum_{x\bmod c}(x/c)_3\check e(x/c).
\]
Consequently, since $|cn^3|=\sqrt{q_c}\,q_n^{3/2}$,
\[
 \frac{d_0(\lambda^{-3}cn^3)}{3^{5/2}|cn^3|}
       =\frac{\widetilde g_3(c)}{\sqrt{q_c}q_n}.
\]
This explains the different powers in the classical Fourier and
adelic Whittaker conventions. If instead the Gauss sum is defined
using $e(z)=\check e(z/\lambda)$, the change of additive character
introduces the factor $\chi_c(\lambda)^{-2}$ at primes away from $6$.
One must also retain conjugation and the frequency change
$\mu\mapsto-\mu$ in comparing classical expansions. This coefficient
comparison does not identify the bad-place vectors.

\indent{\it Bad finite places}.
No tame spherical formula is imposed at $v\in S_f$. The proof of
Proposition~\ref{v:prop:checked-inputs} constructs an actual smooth
vector with $W_{v}(a(y))=1_{U_v}(y)$ by additive Fourier inversion
and upper-unipotent translation. The compact open subgroup $U_v$
is chosen inside the sixth powers and the required character kernels.
Its Weyl rotation is retained as a local function, not replaced by
one. It has compact upper support and is $O(|y|_v^\delta)$ at zero
for every fixed $\delta<1/3$, also at places above $3$.

The finite tensor coefficients of these rotated vectors occur in
$f_{w,\upsilon}$ and hence in the profiles
\eqref{eq:v-raw-profile-definition}. The indices $\upsilon$ and
$\kappa$, the Hilbert factor $H_S(\kappa,[v])^2$, and the finite
character values remain in those profiles. Thus the proof needs
bounds for a fixed finite family of actual bad-place transforms,
not equality with the classical coefficients at the prime above $3$.

\indent{\it Cusp expansions and constant terms}.
In the classical argument one uses explicit expansions at several
cusps. Here rational automorphy under the Weyl element gives the
identity of digit-averaged theta functions used in
Lemma~\ref{v:lem:raw-theta}. Its full reflected expansion retains every
frequency class and the actual transformed $S$-vectors. This serves
in place of an identification with the classical three-cusp table.
The constant terms are not omitted on grounds of cuspidality:
they are radial at each complex place, whereas the integrating
character has angular degree $1$ or $-1$. Their weighted angular
integrals vanish, while the matching Whittaker terms survive.

\indent{\it Normalization and radial profiles}.
The direct $S$-Whittaker function is fixed with scalar one in
\eqref{eq:v-exact-S-vector}; all scalars after Weyl rotation remain
in the reflected function. The direct norm-sheet profile satisfies
\[
 \widehat V(t)=c_{V,S}\prod_{v\mid\infty}
 2^{2t+1}c_v^{-2t-2}\Gamma(5/6+t)\Gamma(7/6+t),
 \qquad c_{V,S}>0.
\]
The finite direct Mellin factor is a nonzero volume. The reflected
finite Mellin factors are holomorphic on $\Re t>-\delta$ and bounded
on closed substrips, by their support and small-norm estimates.
The radial multiplier is $\widehat H_{\upsilon,\alpha}(t)/\widehat V(-t)$;
no division by a possibly vanishing reflected factor is made.

For the Gaussian average, exact smoothing in
\cite[Section~5.4, Lemma~5.5]{Artin} uses
\[
 \widehat K(-s)=\frac{e^{s^2}}{\widehat V(s)},\qquad
 \int_0^\infty K(T)V(uT)\,\frac{dT}{T}
      =\frac{e^{-(\log u)^2/4}}{2\sqrt\pi}.
\]
The field-dependent volumes and Bessel scales are included in
$\widehat V$, so this is an exact normalization, not a comparison
of the sizes of two profiles. Compact radial tests are treated by
Proposition~\ref{v:prop:radial}. The corrected global Whittaker
normalization of \cite{KPerr} is used in
Proposition~\ref{v:prop:checked-inputs}; the classical constant
$3^{5/2}/2$ need not be imported.

Finally the principal Mellin term is normalized by the explicitly
computed $c_SA_{\mathcal T}$ in Section~\ref{sec:v-detail-28}.
The required nonzero factors are proved there, using the chosen
additive tests and principal Euler correction. Changing the theta
vector does not permit an arbitrary change in these factors: the
direct bound, reflected bounds and principal calculation must all
refer to the same coefficients and tests.

\indent{\it Scope of this comparison}.
A full classical comparison would additionally match the covering
splittings, genuine embeddings, vector at $\lambda$, cusp
representatives, additive characters and Haar measures. Such a
comparison is not supplied here and is unnecessary for the stated
route. What remains necessary is the analytic validity of the
constructed theta identity and the moment estimates for its actual
coefficients.

\section{Quadratic and cubic sieve estimates}
\label{v:sec:primary-kummer-transport}
In this section, Lemma~\ref{v:lem:primary-kummer-transport} places the exterior quadratic
and cubic symbols in the primitive Hecke families of \cite{GL,BGL},
giving \eqref{v:eq:transport-quadratic} and \eqref{v:eq:transport-cubic};
these bound the terminal sums of
Sections~\ref{v:sec:unmarked-assembly} and \ref{v:sec:marked-assembly}
after the shared-prime coprimality indicators are separated.

Here $F\supset\mu_{12}$, $S$ (containing the places above $6$),
$R=\mathcal O_{F,S}$ (principal) and the allowed generators, multiplicative
modulo $U_S^6$ with exact sextic reciprocity, are as in
Section~\ref{sec:v-scope} and \cite[Section~3]{Artin}; $a$ denotes the
allowed generator of $\mathfrak a$. Character values at primes of their
conductor are zero, and $\chi_{\mathfrak a}(b)^j=0$ at a shared prime
even if $6\mid j$.

\begin{lemma}[The two primary sieve norms for the actual exterior symbols]
\label{v:lem:primary-kummer-transport}
For $A,B\geq1$, squarefree exterior indices, arbitrary complex coefficients
$\lambda_{\mathfrak b}$, and $\epsilon>0$, one has
\begin{align}
 \sum_{\substack{\mathfrak a\ \sfq\\q_{\mathfrak a}\leq A}}
 \left|\sum_{\substack{\mathfrak b\ \sfq\\q_{\mathfrak b}\leq B}}
       \lambda_{\mathfrak b}\chi_{\mathfrak a}(b)^3\right|^2
 &\ll_{F,S,\epsilon}(AB)^\epsilon(A+B)
              \sum_{\mathfrak b}|\lambda_{\mathfrak b}|^2,
 \label{v:eq:transport-quadratic}\\
 \sum_{\substack{\mathfrak a\ \sfq\\q_{\mathfrak a}\leq A}}
 \left|\sum_{\substack{\mathfrak b\ \sfq\\q_{\mathfrak b}\leq B}}
       \lambda_{\mathfrak b}\chi_{\mathfrak a}(b)^{-2}\right|^2
 &\ll_{F,S,\epsilon}(AB)^\epsilon
       \bigl(A+B+(AB)^{2/3}\bigr)
              \sum_{\mathfrak b}|\lambda_{\mathfrak b}|^2.
 \label{v:eq:transport-cubic}
\end{align}
The constants may also depend on the fixed choices of generators and
finite local data. Separate row restrictions are allowed by deletion;
separate column restrictions and bounded column phases are allowed by
changing the coefficients. Neither assertion permits a coefficient
which depends on a varying opposite sieve index.
\end{lemma}
\begin{proof}
We first verify the conductor condition in the quadratic theorem~\cite{GL}. For squarefree exterior
$\mathfrak a$, let $\Xi_{\mathfrak a}$ be the primitive Hecke character
of $F(\sqrt a)/F$. It is independent of the allowed representative,
since a change by $U_S^6$ is a square. All infinite places of $F$ are
complex, so its infinite type is trivial. At an exterior prime dividing
$\mathfrak a$, the valuation of $a$ is odd and the residue characteristic
is odd; the character is tamely ramified with conductor exponent one.
It is unramified at every other exterior prime. Fix an ordinary ideal
$\mathfrak q_S$, supported at $S$, with exponents $m_v$ so large that
\[
 1+\mathfrak p_v^{m_v}\mathcal O_v\subset F_v^{\times2}.
\]
Every local quadratic character is trivial on this subgroup. Consequently
\[
 \operatorname{cond}(\Xi_{\mathfrak a})
    =\mathfrak a\mathfrak s_{\mathfrak a},
       \qquad \mathfrak s_{\mathfrak a}\mid\mathfrak q_S,
 \qquad
 \Xi_{\mathfrak a}(\mathfrak b)=\chi_{\mathfrak b}(a)^3.
\]
The last identity holds also at a shared prime: that prime is in the
actual primitive conductor. Exact allowed-generator reciprocity gives
$\Xi_{\mathfrak a}(\mathfrak b)=\Xi_{\mathfrak b}(\mathfrak a)$.

Let $B_2$ be the image of the source allowed local subgroup $B$ in
$F_S^\times/F_S^{\times2}$, and assign to an exterior ideal the class
$[a]_{S,2}\in B_2$. Multiplicativity modulo $U_S^6$ makes this a
homomorphism on fractional exterior ideals with finite image.
Choose the ray modulus deep enough that $x\equiv1$ there implies
$[x]=1$ in $\mathcal V_S$. If $a=\epsilon x$ is the allowed generator
of $(x)$, then $[a]=[\epsilon]\in B_S\cap E_S=\{1\}$ by
Lemma~\ref{v:lem:allowed-reciprocity}. Thus its quadratic class is
trivial, which proves factorization through fixed ray data. If coprime squarefree
$\mathfrak a_1,\mathfrak a_2$ have the same class, $a_1/a_2$ is a square
at every place of $S$. Thus the primitive inducing character of
$\Xi_{\mathfrak a_1}\overline{\Xi_{\mathfrak a_2}}$ is trivial at $S$,
is tamely ramified with exponent one at every prime of
$\mathfrak a_1\mathfrak a_2$, and is unramified elsewhere. Its global
conductor is exactly $\mathfrak a_1\mathfrak a_2$.

These are the three defining conditions of the quadratic Hecke family in
\cite{GL}, Definition~1: order dividing two, finite-class
reciprocity, and the exact conductor product for coprime equal-class
indices. Apply their Theorem~1.1 to the full family. This proves
\eqref{v:eq:transport-quadratic}; unwanted rows or classes may then be
deleted. A single class slice need not itself satisfy the family conditions.
The product in the conductor condition is the primitive inducing character;
no artificial shared-prime zero is removed from the displayed matrix.

For the cubic bound, fix once and for all the $n=3$ Kummer construction
of Blomer--Goldmakher--Louvel, Appendix~II. Namely, choose its ray modulus
$\mathfrak c$ supported at $S$, large enough for local cube powers,
its representatives $E$ of $H_{\mathfrak c}\otimes\mathbb Z/3\mathbb Z$,
and its generators $m_E$ as $R$-ideals. For each exterior ideal
$\mathfrak a$, fix one admissible decomposition
\[
 \mathfrak a=(x)E\mathfrak g^3,
            \qquad x\equiv1\pmod{\mathfrak c},
\]
and a generator $z$ of the fractional $R$-ideal $\mathfrak g$.
There is an $S$-unit $\varepsilon_{\mathfrak a}$ such that
\begin{equation}\label{eq:v-kummer-decomposition}
 a=\varepsilon_{\mathfrak a}\,x m_E z^3.
\end{equation}

Let $\zeta_{\mathfrak a}$ denote the character of that appendix
\emph{after} its final primitive-induction convention. This convention
is stated also at the beginning of \cite[Section~2]{BGL}; it replaces
the initial residue-symbol function by its underlying primitive
Hecke character. Kummer theory and \eqref{eq:v-kummer-decomposition}, with one consistent Artin
orientation, give the identity of primitive global cubic characters
\begin{equation}\label{eq:v-cubic-character}
 \Psi_{\mathfrak a}
     =\omega_{\varepsilon_{\mathfrak a}}\zeta_{\mathfrak a},
 \qquad
 \Psi_{\mathfrak a}(\mathfrak b)=\chi_{\mathfrak b}(a)^2.
\end{equation}
Here $\omega_{\varepsilon}$ is the cubic Kummer character of the $S$-unit
$\varepsilon$, ramified only at $S$. If the primary source uses the
opposite Artin convention, conjugate both characters in \eqref{eq:v-cubic-character}; the same
sieve norm results. The identity is between the primitive global
characters, not just residue symbols away from $\mathfrak g$;
a prime dividing the auxiliary $\mathfrak g$ acquires
no additional zero: choose the decomposition with
$(\mathfrak g,\mathfrak b)=1$ to evaluate away from the conductor,
as Appendix~II prescribes, and then use primitive induction. At every exterior prime dividing squarefree
$\mathfrak a$, the valuation of $a$ is one and the residue characteristic
is not three. Both sides are tamely ramified there with conductor
exponent one, so their natural shared-prime zeros agree on all exterior
$\mathfrak b$. Outside $S$, $\omega_{\varepsilon}$ is unramified and
has modulus one.

Partition the outer indices $\mathfrak a$ by the finitely many classes
$\varepsilon_{\mathfrak a}\in U_S/U_S^3$. On one part replace
$\lambda_{\mathfrak b}$ by
$\lambda_{\mathfrak b}\omega_{\varepsilon}(\mathfrak b)$, which
preserves the squared coefficient mass, and apply \cite[Theorem~1.3]{BGL}
for $n=3$ to its full fixed family, and then delete the unwanted outer
indices. Sum the finitely many bounds. This proves the cubic bound for
$\Psi_{\mathfrak a}(\mathfrak b)$. Exact allowed-generator reciprocity
identifies it with $\chi_{\mathfrak a}(b)^2$; conjugating the character
and the coefficients gives \eqref{v:eq:transport-cubic}. This also shows
that a fixed inverse orientation, fixed finite character twist, or finite
class partition changes only fixed constants. No norm-dependent choice
of family is involved. Both starred index families in the primary BGL
theorem are squarefree and prime to $S$; unrestricted cube factors must
be stratified before applying this theorem.
\end{proof}

In the terminal hybrid estimate the two norms are invoked only after
the shared-prime coprimality indicators are expanded and the positive
quadratic bound has made the remaining coefficient independent of the
opposite cubic index. Angular factors $\Lambda$, which are neither finite Hecke
characters nor necessarily functions of the norm, are retained as
row-only, tuple-only, or dual-only bounded coefficients.
An inner angular phase independent of the opposite index is permissible
in either primary sieve; the remaining bounded dual
coefficient is discarded only after the quadratic positive estimate.
This does not authorize angular prime coefficients in M2.
All constants concern one fixed cyclotomic field and fixed $S$;
no uniformity in the field or its discriminant is claimed.

\section{Coprimality conditions and cube completion}
\label{sec:v-detail-4}
In this section, Lemma~\ref{v:lem:hybrid-sieve} applies Lemma~\ref{v:lem:primary-kummer-transport}
to terminal sums with an unrestricted cube factor and moving
coprimality conditions, expanded before positivity so that the two
sieve indices stay independent; the resulting bound
\eqref{eq:v-hybrid-sieve-detailed} has the factor
$(K+NB)B\{N+L+(NL)^{2/3}\}$, with one factor $B$ from the cube index.
We then record the cube-completion normalization and the table
relocating excluded prime supports as natural zeros. The hybrid bound and the
cube-completion normalization are used in
Sections~\ref{sec:v-detail-9} and \ref{sec:v-detail-14}.

\begin{lemma}[The quadratic--cubic mean square]\label{v:lem:hybrid-sieve}
Let $K,N,B,L\ge1$. Let $R,n,P$ be squarefree exterior ideals and
let $b$ be an arbitrary exterior ideal, with respective norm bounds
$K,N,L,B$. Suppose $|a(P)|\le1$ and $|\beta(n,b)|\le1$, where
$a$ is independent of $R,n,b$ and $\beta$ is independent of $R,P$.
Then
\begin{equation}\label{eq:v-hybrid-sieve-detailed}
 \begin{split}
 &\sum_{\substack{R\ \sfq\\q_R\le K}}
 \left|\sum_{\substack{P\ \sfq\\q_P\le L}}
       \frac{a(P)}{\sqrt{q_P}}1_{(P,R)=1}
       \sum_{\substack{n\ \sfq,b\\q_n\le N,\ q_b\le B}}
       \beta(n,b)\chi_R(nb)^3\chi_P(n)^{-2}1_{(P,b)=1}
                                                        \right|^2\\
 &\hspace{15mm}\ll_{F,S,\epsilon}(KNBL)^\epsilon
             (K+NB)B\{N+L+(NL)^{2/3}\}.
 \end{split}
\end{equation}
The assertion permits arbitrary separate restrictions on $P$ and
on $(n,b)$, and fixed bounded phases with the same dependence.
\end{lemma}
\begin{proof}
Put $X=2KNBL$ and
$\mathcal C(U,V)=U+V+(UV)^{2/3}$. Throughout the proof a factor
$X^\delta$ may be inserted for any prescribed $\delta>0$; the
finitely many occurrences and logarithmic decompositions are absorbed
in the final $X^\epsilon$ by choosing $\delta$ first sufficiently small.
We use $d(a)^A\ll_{A,\delta}q_a^\delta$ and
$\#\{a:q_a\le T\}\ll_F T$.

Expand $1_{(P,R)=1}=\sum_{D\mid P,R}\mu_F(D)$, and put
$P=DP'$, $R=DR'$. At each fixed $R$, Cauchy--Schwarz over its
squarefree divisors costs only $d(R)\ll X^\delta$. Therefore the
left side of \eqref{eq:v-hybrid-sieve-detailed} is at most
\begin{equation}\label{eq:v-hybrid-first-divisor}
 X^\delta\sum_{\substack{D\ \sfq\\q_D\le\min(K,L)}}q_D^{-1}
                   \mathcal H_D,
\end{equation}
where $\mathcal H_D$ has row bound $K/q_D$, prime-product bound
$L/q_D$, and no remaining $R',P'$ coprimality indicator. Its coefficients
have the form
\[
 a_D(P')=\mu_F(D)a(DP')1_{(P',D)=1},\quad
 \beta_D(n,b)=\beta(n,b)\chi_D(nb)^3\chi_D(n)^{-2}1_{(D,b)=1}.
\]
Then $|a_D|,|\beta_D|\le1$ with the same independence, the restriction
$(R',D)=1$ can be dropped after squaring, and all zero-extended symbols
retain their zeros.

Write uniquely
\[
 b=gt^2,\quad g\ \sfq,\quad c=(n,g),\quad n=cm,\quad g=ch.
\]
Then $c,m,h$ are pairwise coprime and squarefree, and
\[
 \chi_{R'}(nb)^3=\chi_{R'}(mh)^3 1_{(R',ct)=1}.
\]
Decompose $g,t,c$ into norm dyads $G,T,C$, so $GT^2\ll B$,
$C\ll\min(G,N)$. First use Cauchy--Schwarz over $t$; its cost
is $O(T)$. For each fixed $t$, retain the row restriction $(R',t)=1$
and expand $1_{(R',c)=1}=\sum_{r\mid R',c}\mu_F(r)$.
Cauchy--Schwarz at each $R'$ again costs only $X^\delta$.
Write $R'=rR''$, $c=rc'$. The coefficient $\chi_r(mh)^3$ is
independent of $R''$ and can be included in the column coefficient.
The remaining row restrictions may now be dropped in the positive norm.

For fixed $D,r,t$, the quadratic column is $u=mh$, which is
squarefree because $(m,h)=1$, and
$q_u\ll NG/C^2$. Its coefficient is a sum over factorizations
$u=mh$ and over $c'$ with $q_{c'}\ll C/q_r$. The quadratic sieve
\eqref{v:eq:transport-quadratic}, followed by Cauchy--Schwarz over
these factorizations and over $c'$, bounds this row norm by
\begin{equation}\label{eq:v-hybrid-quadratic-stage}
 X^\delta\left(\frac K{q_Dq_r}+\frac{NG}{C^2}\right)
 \frac C{q_r}\sum_{h,c',m}|F_{D,r,h,t}(c'm)|^2.
\end{equation}
Here the factors of modulus at most one originating from $\beta_D$
are discarded only after this positive estimate, and
\[
 F_{D,r,h,t}(j)=
 \sum_{P'\ \sfq}\frac{a_D(P')\chi_{P'}(r)^{-2}
                         1_{(P',ht)=1}}{\sqrt{q_{P'}}}
                              \chi_{P'}(j)^{-2}.
\]

To obtain this expression use the exact identity
\[
 \chi_{P'}(rc'm)^{-2}1_{(P',rc'ht^2)=1}
 =\chi_{P'}(r)^{-2}\chi_{P'}(c'm)^{-2}1_{(P',ht)=1}.
\]
Thus the remaining coefficient of $P'$ is independent of the cubic
index $j=c'm$, and no $j$-dependent restriction has been discarded
from a linear polynomial.

For fixed $h$, the ideal $j=c'm$ is squarefree, has norm $\ll N/q_r$,
and has at most $d(j)$ representations of the required kind. Applying
\eqref{v:eq:transport-cubic} with outer index $j$ gives
\[
 \sum_{c',m}|F_{D,r,h,t}(c'm)|^2
 \ll X^\delta\mathcal C(N/q_r,L/q_D).
\]
The orientation is justified by the zero-preserving reciprocity
identity $\chi_{P'}(j)^{-2}=\chi_j(P')^{-2}$ of
Lemma~\ref{v:lem:allowed-reciprocity}. The squared coefficient mass is
at most $\sum_{q_{P'}\le L/q_D}q_{P'}^{-1}\ll1+\log L$.

The number of $h$ is $O(G/C)$. Summing over the $O(T)$ values of
$t$, including its earlier Cauchy factor, shows that a fixed $r$
contributes at most
\[
 X^\delta\frac{GT^2}{q_r}
 \left(\frac K{q_Dq_r}+\frac{NG}{C^2}\right)
                    \mathcal C(N/q_r,L/q_D)
\]
to $\mathcal H_D$ in these dyads. Bound $GT^2\ll B$,
$G/C^2\ll B$, and use monotonicity of $\mathcal C$.
Summation over $r$ costs at most
$\sum_{q_r\ll C}q_r^{-1}\ll1+\log C$. Consequently
\[
 \mathcal H_D\ll X^\delta B(K+NB)\mathcal C(N,L),
\]
with the logarithmic decompositions absorbed in $X^\delta$.
Finally insert this in \eqref{eq:v-hybrid-first-divisor} and use
$\sum_{q_D\le\min(K,L)}q_D^{-1}\ll1+\log(KL)$.
This proves \eqref{eq:v-hybrid-sieve-detailed} after distributing the
initial loss shares; the unrestricted cube index contributes one
factor $B$, not $B^2$.
\end{proof}

For the cube-completion identity, whose cube-M\"obius notation is
specified in Section~\ref{sec:v-detail-9}, fix $Z\ge2$, an initial column
length $Z^N$, and a smooth annular test $W$. The completed theta test is \(y^{1/2}W(y)\). If its free cube index is
\(b_0\) and its cube-M\"obius divisor is \(h_0\), putting \(B=h_0b_0\)
turns the combined coefficient into
\[
 \frac{\sqrt{q_B}}{\sqrt{Z^N}}\,
       \sum_{h_0\mid B\text{ in the chosen branch}}\mu_F(h_0),
\]
up to its retained fixed characters, natural coprimality indicators and bounded profiles.
The divisor sum is at most
$d(B)\ll_\epsilon q_B^\epsilon$ in absolute value, uniformly in the row and
residual squarefree column. On $q_B\asymp Z^\ell$, with $r=N-3\ell$,
\[
 q_B/Z^N=(q_B/Z^\ell)Z^{-r-2\ell}\asymp Z^{-r-2\ell},
\]
and the bounded ratio $q_B/Z^\ell$ remains in the common profile. This
gives the normalization of the transformed sum without a free
arbitrary column coefficient.

For surviving inverse-moment prime factors, the following table
anticipates the indices defined in Section~\ref{sec:v-detail-10}.
There $C,J,d_2,v'$ are the squarefree factors of the new label,
$q_0,t',g'/d_2$ are its fixed extra coprimality factors, $d_k$ is
the first Poisson divisor, and $k''$ is the final dual frequency.
The old support exclusions become natural zeros:
\[
 \begin{array}{c|l}
 \text{excluded prime support}&\text{location in the transformed sum}\\\hline
 C,J,d_2,v'&f_{\rm new}=JCd_2v'\\
 q_0,t',g'/d_2&\rho'(n)=\rho(n)1_{(n,\operatorname{rad}(q_0t'g'/d_2))=1}\\
 d_k,d_2,k''&k_{\rm new}=d_kd_2k''
 \end{array}
\]

A weight selecting \(p\mid n\) then vanishes automatically when \(p\)
is in one of these supports; this removes the earlier
prime-selection restriction exactly, without deleting terms from a
linear polynomial. The identity
\[
 \chi_n(d_k)\chi_n(k'')\bar\chi_n(d_2)\chi_n(CJ)^4\chi_n(v')^4
 =\chi_n(k_{\rm new})\chi_n(f_{\rm new})^4
\]
uses \(1+4\equiv-1\pmod6\) at \(d_2\), and both sides have the same
nonunit zeros. The complete
uniform proofs are in Sections~\ref{v:sec:marked-assembly}
and \ref{v:sec:M2-assembly}.

\section{Poisson summation over unit orbits with coprimality conditions}
\label{sec:v-detail-5}
In this section, Proposition~\ref{v:prop:orbit-poisson} proves Poisson
summation with coprimality conditions on the $\Gamma$-orbits of
Section~\ref{sec:v-scope}: the nonzero dual orbits unfold to one common
norm kernel $\mathcal K_\Phi=\Psi_{\widehat\Phi}$ and the dual origin
gives a separate compact-quotient integral; no row-sector projector is
used (Sections~\ref{sec:v-detail-10}, \ref{sec:v-detail-22}). The unit
normalization \eqref{eq:v-poisson-unit-normalization}, from
Lemma~\ref{v:lem:allowed-reciprocity}, and the imprimitive formula
\eqref{eq:v-imprimitive-orbit-poisson} are recorded after it.

Let \(\lambda\) be a primitive exterior residue character modulo
\(\mathfrak m=mR\), of order dividing six, extended by zero on nonunits;
it is a periodic function on \(R/mR\), not a global \(S\)-character.
Let \(D\) be a complementary exterior
radical, \((D,\mathfrak m)=1\), and let \(\Phi\) be an actual row test.

\begin{proposition}[The kernel in Poisson summation over unit orbits]\label{v:prop:orbit-poisson}
Define $\Psi_\Phi(x)=\int_G\Phi(g\ell_x)\,dg$ and
$\mathcal K_\Phi=\Psi_{\widehat\Phi}$, with finite components
of the norm section equal to one, and assume \(\Phi(0)=0\). The original
row expression is \(\sum_{[k]\ne0}\lambda(k)1_{(k,D)=1}
\Psi_\Phi(q_k/K)\). On the nonzero \(\Gamma\)-orbits, the
transform with coprimality conditions has the exact nonzero part
\[
 \sum_{d\mid D}\mu_F(d)\lambda(d)
 \frac{K\tau(\lambda)}{q_dq_{\mathfrak m}}
 \sum_{[h]\ne0}\bar\lambda(h)
 \Psi_{\widehat\Phi}\left(\frac{Kq_h}{q_dq_{\mathfrak m}}\right).
\]
The dual-origin term vanishes unless the primitive modulus is one; in that
case it is
\[
 K\operatorname{vol}(G/\Gamma)\widehat\Phi(0)
            \sum_{d\mid D}\frac{\mu_F(d)}{q_d}.
\]
The norm kernel is independent of the chosen generator. The finite
unit factors in \eqref{eq:v-poisson-unit-normalization} below are
retained separately.
\end{proposition}
\begin{proof}
For \(t\in G/\Gamma\), expand the complementary coprimality indicator and apply ordinary
residue-class Poisson to
\(\sum_{k'\in R}\lambda(k')\Phi(dk'/(t\ell_K))\).
Its covolume and scaling give \(K/(q_dq_{\mathfrak m})\); the finite
Fourier sum is \(\tau(\lambda)\bar\lambda(h)\). Its test argument is
\(t\ell_Kh/(dm)\). The order-six condition gives
\(\lambda(\gamma)=1\) for \(\gamma\in\Gamma\), so the nonzero dual
orbits unfold the compact quotient into \(G\). Translation by the
norm-one part of \(h/(dm)\) removes its shape and yields the
displayed norm kernel; absolute Schwartz bounds justify the
unfolding. At \(h=0\) one integrates the constant Fourier value over
the compact quotient (not periodized over its infinite stabilizer);
character orthogonality kills it for a nonprincipal
primitive modulus. Since \(\Phi(0)=0\), the nonzero-row convention is unchanged.
\end{proof}

The dual sum runs over $(R\setminus\{0\})/\Gamma$,
not only over allowed ideal generators; by Lemma~\ref{v:lem:allowed-reciprocity},
\begin{equation}\label{eq:v-poisson-unit-normalization}
 \sum_{[h]\ne0}\bar\lambda(h)\mathcal K_\Phi(Aq_h)
 =\sum_{\epsilon\in\mathcal U_6}\bar\lambda(\epsilon)
       \sum_{\mathfrak h}\bar\lambda(a_h)
                              \mathcal K_\Phi(Aq_{\mathfrak h}).
\end{equation}
For every full, possibly imprimitive Gauss sum,
\[
 \mathcal G(u,\epsilon h)
       =\overline{\chi_u(\epsilon)}\mathcal G(u,h),
\]
so normalization can introduce $\xi_\epsilon$ into an entire column.
These factors are retained, and every fixed coefficient group in M1
and M2 is closed under $\Upsilon_S$ before any bin is selected.
They act on the whole product, so within one transformed
side the two ordinary factors and both centered rectangles carry the
same twist; compatibility with the inverse and product moment recursions
is \eqref{eq:v-inverse-unit-return} and
Lemma~\ref{v:m2:lem:common-coefficient}, respectively.

For an arbitrary periodic function $f$ modulo an exterior ideal
$\mathfrak m$, the formula used in the imprimitive case is instead
\begin{equation}\label{eq:v-imprimitive-orbit-poisson}
 \sum_{[k]\ne0}f(k)\Psi_\Phi(q_k/K)
 =\frac K{q_{\mathfrak m}}
     \sum_{[h]\ne0}\widehat f_{\mathfrak m}(h)
             \mathcal K_\Phi(Kq_h/q_{\mathfrak m})
   +\frac K{q_{\mathfrak m}}\widehat f_{\mathfrak m}(0)
            \operatorname{vol}(G/\Gamma)\widehat\Phi(0),
\end{equation}
where $f(\gamma k)=f(k)$ for $\gamma\in\Gamma$, $\Phi(0)=0$, and
$\widehat f_{\mathfrak m}(h)=\sum_{x\bmod\mathfrak m}
 f(x)e(hx/\mathfrak m)$ with the sign paired with the infinite transform.
This follows directly from residue-class Poisson (no primitive
Gauss identity is used); the last term is the compact-quotient origin,
and the nonzero dual sum is reindexed by the exact identity
\begin{equation}\label{eq:v-general-unit-normalization}
 \sum_{[h]\ne0}\widehat f_{\mathfrak m}(h)\mathcal K_\Phi(Aq_h)
 =\sum_{\epsilon\in\mathcal U_6}\sum_{\mathfrak h}
       \widehat f_{\mathfrak m}(\epsilon a_h)
                              \mathcal K_\Phi(Aq_{\mathfrak h}).
\end{equation}
Indeed, $\Gamma$-invariance of $f$ implies that of its finite Fourier
transform, and Lemma~\ref{v:lem:allowed-reciprocity} partitions the
nonzero $\Gamma$-orbits into the displayed unit classes and allowed
ideal generators. For general $f$ the coefficient
$\widehat f_{\mathfrak m}(\epsilon a_h)$ is retained in full;
no scalar unit factor is extracted without a separate covariance law.
For a primitive character that law gives
\eqref{eq:v-poisson-unit-normalization}; for the imprimitive character
sums used in M2 it is the full Gauss-sum scaling law above, with all
unit twists retained in the arithmetic coefficient.
On retained polynomial norm ranges the kernel
and its Euler derivatives are $\ll(\log Z)^{O(1)}$
near zero (Section~\ref{sec:v-detail-8}); this is included
in $Z^\epsilon$.

In the first inverse-moment transform the primitive characters come from the
residual squarefree columns and the nonzero common-prime exponents, and
the complementary coprimality indicator supplies \(D\). The first
product-moment transform keeps the full imprimitive Gauss sums, with the
same unrestricted orbit unfolding; its finite class, unit and Gauss
information stays in the explicit arithmetic coefficient. A projector
\(b(k)\) could introduce new coefficient characters and is therefore
excluded. The kernel is used in Sections~\ref{v:sec:unmarked-assembly} and \ref{v:sec:M2-assembly}.

\section{Compact support of the dual norm kernel}
\label{sec:v-detail-6}
In this section, Lemma~\ref{v:lem:compact-dual} shows that for a
positive row lift of compact complex additive Fourier support the dual
norm kernel vanishes beyond a fixed $C_\Phi$, so the nonzero dual
frequencies satisfy $q_h\le C_\Phi q_dq_{\mathfrak m}/K$ (used in
Sections~\ref{v:sec:unmarked-assembly} and \ref{sec:v-detail-22});
theta-reflection tails are controlled by Proposition~\ref{v:prop:radial}.
We then construct the nonnegative band-limited lift to which the lemma
is applied.

\begin{lemma}[Compact dual norm support]\label{v:lem:compact-dual}
Let \(\Phi\) be a product Schwartz--Bruhat row lift whose complex additive
Fourier transforms are compactly supported. Suppose its finite-place factors
run through a fixed finite family of locally constant compactly supported
tests, including any fixed local class indicators. There is a fixed \(C_\Phi\)
such that the nonzero dual norm-sheet kernel satisfies
\[
 \mathcal K_\Phi(x)=0\qquad(x>C_\Phi).
\]
In primitive Poisson summation with coprimality conditions with row scale \(K\), primitive modulus
\(\mathfrak m\), and divisor \(d\) from the complementary coprimality indicator, every
nonzero retained dual frequency satisfies
\[
 q_h\le C_\Phi\,q_dq_{\mathfrak m}/K.
\]
The dual origin remains a separate compact-quotient integral.
\end{lemma}
\begin{proof}
Choose compact sets \(C_v\subset F_v\) with
\[
 \operatorname{supp}\widehat\Phi\subseteq\prod_{v\in S}C_v,
 \qquad C_\Phi=\prod_{v\in S}\sup_{y\in C_v}|y|_v.
\]
For \(g\in F_S^1\), one has \(|g\ell_x|_S=x\). Hence
\[
 x>C_\Phi\quad\Longrightarrow\quad
 \widehat\Phi(g\ell_x)=0\quad(g\in F_S^1)
 \quad\Longrightarrow\quad\mathcal K_\Phi(x)=0.
\]
A nonzero frequency \(h\) has norm argument
\(x=Kq_h/(q_dq_{\mathfrak m})\), whence
\(q_h\le C_\Phi q_dq_{\mathfrak m}/K\). The term \(h=0\) is the
compact-quotient integral.
\end{proof}

At a complex place take the nonnegative lift \(|\check b|^2\), with \(b\)
a compact smooth bump near zero: its Fourier transform, the convolution
of \(b\) with its reflected conjugate, is compactly supported, and the
lift is bounded below on the required fixed ball.
Finite-place unit-neighborhood indicators avoid the additive origin and
remain Schwartz--Bruhat after fixed class indicators. This lift is used in Section~\ref{v:sec:unmarked-assembly}.

The lemma applies to this fresh lift only, not to a row test multiplied
by a norm phase or another kernel unless its Fourier support is
rechecked. Pure row phases cancelling in a positive square may be
removed first; column norm modes stay in their profiles.

\section{The inverse moment}
\label{v:sec:unmarked-assembly}
In this section, Proposition~\ref{v:prop:radial},
Lemma~\ref{v:lem:primary-kummer-transport} and Section~\ref{sec:v-detail-5}
give Theorem~\ref{v:thm:unmarked-inverse}, the inverse moment
$\sum_{q_k\le C_0H}|D^{-1/2}\sum_n\mu_F(n)\nu(n)\chi_n(k)W(q_n/D)|^2
\ll H(HD)^\eps$ for $H\ge\max(2,D^{1+c})$. Sixth-power amplification then
gives $\beta_*(F)\le11/12$ (Theorem~\ref{v:thm:unmarked-seed}),
the initial bound for Section~\ref{v:sec:compensated-assembly}, and the
bound \eqref{eq:v-general-amplification} for general rows
(Lemma~\ref{v:lem:general-amplification}), used in
Section~\ref{sec:v-detail-29}. We first record how the fixed coefficient
character changes under the two transforms, \eqref{eq:v-inverse-unit-return}.

Write \(a_0(n)=\bar\Lambda(n)\gamma_2(n)\) on squarefree exterior ideals,
with the unitary character \(\Lambda\) and the finite-order character
\(G_0\) of Lemma~\ref{v:lem:theta-characters}. By that lemma,
\[
 a_0(n)\gamma_1(n)=\mu_F(n)G_0(n),\qquad
 \mu_F(n)\gamma_{-1}(n)=a_0(n)\bar G_0(n),
\]
where \(G_0\) is a character of the ray class group of the modulus
\(\mathfrak m\) supported on \(S_f\) constructed in its proof, of order
dividing six (the CRT phases of a single normalized Gauss sum are
retained). At an exterior prime,
\[
 \gamma_1(p)\gamma_2(p)
   =\bar\chi_p(4)\gamma_3(p)\gamma_2(p)^3,
 \quad\gamma_2(p)^3=-\Lambda(p),
 \quad G_0(p)=\bar\chi_p(4)\gamma_3(p).
\]
Hence $a_0(p)\gamma_1(p)=-G_0(p)$. For coprime squarefree
$c_1,c_2$, the CRT factor for $\gamma_j(c_1c_2)$ is
$\chi_{c_1}(c_2)^j\chi_{c_2}(c_1)^j$; in
$\gamma_1\gamma_2$ the total exponent is three in each orientation, so
by allowed-generator reciprocity the factor is
$\chi_{c_1}(c_2)^3\chi_{c_2}(c_1)^3=\chi_{c_2}(c_1)^6=1$. This
proves the first identity for every squarefree $n$.
Since $\mu_{12}\subset F$, $\chi_p(-1)=1$, so
$\gamma_{-1}(n)=\overline{\gamma_1(n)}$; conjugating the first
identity and using $|a_0(n)|=1$ gives the second. In the two
successive transforms the fixed character therefore changes by
\[
 \nu\longmapsto\nu G_0\longmapsto\nu G_0\bar G_0=\nu.
\]

For the unit factors: in the first transform the
primitive row character is
\[
 \psi_1=\chi_{u_1}\bar\chi_{u_2}
                    \prod_{p\mid R_1}\chi_p^{t_p}.
\]
Writing its frequency as $\epsilon_1 h_0$ gives
\[
 \bar\psi_1(\epsilon_1)=
 \overline{\xi_{\epsilon_1}(u_1)}\xi_{\epsilon_1}(u_2)
              \prod_{p\mid R_1}\bar\chi_p(\epsilon_1)^{t_p}.
\]
The last product is a frozen outer scalar, so both first inverse
polynomials have the twist $\overline{\xi_{\epsilon_1}}$.
Similarly, for
$\psi_2=\bar\chi_{z_1}\chi_{z_2}$,
$\bar\psi_2(\epsilon_2)=
 \xi_{\epsilon_2}(z_1)\overline{\xi_{\epsilon_2}(z_2)}$.
The complete two-transform change is therefore
\begin{equation}\label{eq:v-inverse-unit-return}
 \nu\longmapsto\nu G_0\overline{\xi_{\epsilon_1}}
 \longmapsto\nu\overline{\xi_{\epsilon_1}}\xi_{\epsilon_2}.
\end{equation}

Common-divisor factors have unit modulus and cancel between the
conjugate columns. Lemma~\ref{v:lem:unit-products} preserves all
imprimitive zeros. Induction is simultaneous over one fixed finite
group containing $\nu,G_0$ and $\Upsilon_S$. The angular character
$\Lambda$ is unitary (finite order is not assumed) and occurs only
in the cube completion and terminal reflection, not as a new fixed
coefficient after the two transforms. All Kummer symbol powers,
including zero powers, retain their nonunit zeros.

\subsection{Uniform kernels and the full independent profile boxes}
\label{sec:v-detail-8}
In this subsection, we prove the Euler-derivative bounds
$|\mathcal D^jK(T)|\ll_j(1+|\log T|)^{d_F}$ for the norm-sheet kernel and
the common Fourier density \eqref{v:eq:H2-common-density}, which
separates all independent norm coordinates before any arithmetic index
is evaluated, at the cost $(1+\log Z)^{d_F}$ with $d_F$ independent of the
derivative order. Lemma~\ref{v:lem:row-column-test-domains} then keeps
the column tests compactly supported and the row lifts band-limited.
Both are used in Sections~\ref{sec:v-detail-10} and
\ref{sec:v-detail-11}.

For a fresh unrestricted positive row lift \(\Phi\), put
\[
 K(T)=\int_{F_S^1}\widehat\Phi(t\ell_T)\,dt,\qquad
 \mathcal D=T\partial_T,\qquad d_F=[F:\mathbb Q]/2+|S_f|-1.
\]
For every \(j\),
\[
 |\mathcal D^jK(T)|\ll_j(1+|\log T|)^{d_F}\quad(0<T\le1),
 \qquad \mathcal D^jK(T)\ll_{A,j}T^{-A}\quad(T\ge1).
\]
These follow by writing the complex radii logarithmically and the finite
radii as valuation shells. On each fixed norm sheet the nonnegative
logarithmic variables form a simplex of dimension \([F:\mathbb Q]/2-1\);
the finite shell vectors contribute degree \(|S_f|\). Negative complex
logs are summed with Schwartz decay. A norm derivative acts on one
complex factor and changes its seminorm, not this dimension, so the
exponent \(d_F\) is independent of \(j\). For the band-limited lift
of Lemma~\ref{v:lem:compact-dual}, the upper norm support is compact;
the additive origin keeps its separate compact-quotient term.

Choose $\varphi\in C_c^\infty((1/2,2))$ with
$\sum_{j\in\mathbb Z}\varphi(2^{-j}x)=1$ for $x>0$.
For each nonzero dual frequency insert $\varphi(q_h/H_j)$,
$H_j=2^j$, in the signed source sum before taking an absolute Fourier
measure. Compact dual support and the polynomial bounds for all
row, column and divisor norms leave $O(1+\log Z)$ nonempty dyads
(also for a frequency with an explicit, polynomially bounded denominator).
On each dyad $\log(q_h/H_j)$ ranges over a fixed interval, with
uniform Euler derivatives. Reflection tails without compact support
are first truncated by the arbitrary decay in
Proposition~\ref{v:prop:radial}, and the retained polynomial range is
partitioned likewise. The total fixed power of $\log Z$ is
absorbed in a reserved subsidiary loss; the band-limited row lift is unchanged.

Let \(u\) be a full independent logarithmic box of normalized ideal
norms. Each coordinate is \(u_i=\log(q_i/Y_i)\) for its actual fixed-ratio
dyad \(q_i/Y_i\in[c_i,C_i]\), with constants independent of \(Z\).
The coarse \(\eta\)-centres are used only in the exponent calculation;
replacing \(Y_i\) there by its centre pays the displayed affine edge
allowance, not a derivative-dependent support loss. Choose a smooth
cutoff \(\chi\) equal to one on this actual fixed box. Every inverse-moment profile has the form
\[
 P_Z(u)=\chi(u)K(Z^\zeta e^{L(u)})e^{R(u)}
             \prod_j W_j(c_j e^{L_j(u)})e^{i\vartheta\cdot A u},
\]
where \(L,R,L_j\) are fixed linear forms, \(\zeta\) is bounded, and
real inverse-root ratios are in \(e^{R(u)}\). For every required weight
\(A_0\), integration by parts gives a common Fourier density satisfying
\begin{equation}\label{v:eq:H2-common-density}
 P_Z(u)=\int e^{i\xi\cdot u}\widehat P_Z(\xi)\,d\xi,\qquad
 \int(1+|\xi|)^{A_0}|\widehat P_Z(\xi)|\,d\xi
 \ll (1+\log Z)^{d_F}(1+|\vartheta|)^J\mathfrak S_J.
\end{equation}

Indeed, let $n=\dim u$, $|\zeta|\le C_0$, $L_Z=1+\log Z$.
Since $\operatorname{supp}\chi$ is fixed, the chain rule and the
Euler-derivative bounds for $K$ give, for every multi-index $\alpha$,
\[
 \|\partial^\alpha P_Z\|_{L^1(\mathbb R^n)}
 \ll_\alpha L_Z^{d_F}(1+|\vartheta|)^{|\alpha|}
                 \mathfrak S_{|\alpha|}.
\]
On this box,
$|\log(Z^\zeta e^{L(u)})|\le C_0\log Z+O(1)$ when the
kernel argument is $\le1$; above $1$ the derivatives are bounded
by rapid decay. For each annular factor,
\[
 \sup_{c>0,u}\left|(y\partial_y)^jW(y)
                       \big|_{y=ce^{L_j(u)}}\right|
 \le\sup_{y>0}|(y\partial_y)^jW(y)|,
\]
so the constants are uniform in $c_j>0$.

Choose an integer $m$ with $2m>A_0+n$. Applying $(1-\Delta)^m$
before Fourier transformation yields
\[
 |\widehat P_Z(\xi)|
 \ll (1+|\xi|^2)^{-m}
       \sum_{|\alpha|\le2m}\|\partial^\alpha P_Z\|_1.
\]
Since
$\int_{\mathbb R^n}(1+|\xi|)^{A_0}(1+|\xi|^2)^{-m}\,d\xi<\infty$,
this proves \eqref{v:eq:H2-common-density} with $J=2m$, increasing
$J$ if required by the fixed tests. Increasing a Fourier order changes
$J$ and the constant, but not the exponent $d_F$ of $L_Z$, and
introduces no factor $Z^{J\eta}$.

Fixed coordinates give outer unitary
phases; residual columns give pure norm modes in their annular tests.
Fixed finite-place tests form one finite family. Hard restrictions
remain in the outer index set and are not differentiated. Tails are
removed on the signed source before any absolute Fourier measure is taken.

\begin{lemma}[Row and column test functions]
\label{v:lem:row-column-test-domains}
Proposition~\ref{v:prop:radial} is applied to compactly supported
smooth column tests on $(0,\infty)$. The Schwartz--Bruhat row lifts
of Lemma~\ref{v:lem:compact-dual} enter through Poisson summation
and their norm-sheet kernels.
\end{lemma}
\begin{proof}
The original column test is $W\in C_c^\infty(I)$, and cube completion
uses $y^{1/2}W(y)$ on the same annulus. At each later step include
the entire Poisson kernel in $P_Z$ of
\eqref{v:eq:H2-common-density}. Choose
$\chi_i\in C_c^\infty(\mathbb R)$ equal to one on each coordinate
projection of $\operatorname{supp}P_Z$. Fourier inversion gives exactly
\[
 P_Z(u)=\int\widehat P_Z(\xi)
              \prod_i\bigl(\chi_i(u_i)e^{i\xi_i u_i}\bigr)\,d\xi.
\]
Every residual column factor is therefore
\[
 W_{i,\xi}(y)=\chi_i(\log y)y^{i\xi_i},\qquad
 \operatorname{supp}W_{i,\xi}\subset
                \exp(\operatorname{supp}\chi_i)\Subset(0,\infty),
\]
possibly times a fixed compact annular profile or a bounded real
power. Its $C^J$ seminorm is bounded by a fixed test seminorm times
$(1+|\xi_i|)^J$, integrable against the common weighted Fourier
measure. Hence each subsequent theta reflection has a compact column
test, with no Poisson kernel inside the column square.
The row lift is unchanged, so its compact
additive Fourier support is preserved; a fresh positive row majorant
is inserted only when enlarging a positive row norm
(Lemma~\ref{v:lem:compact-dual}).
\end{proof}

\subsection{The normalized inverse family and its direct estimates}
\label{sec:v-detail-9}
In this subsection, we define the normalized inverse family
$\mathcal E(N,V,M;\rho)$, whose bound $\mathcal E\ll Z^{F+\eps}$ under the
margin \eqref{v:eq:unmarked-canonical-margin} is the target of the
recursion, and prove its direct bounds. Through
Lemma~\ref{v:lem:primary-kummer-transport} and the terminal reflection of
Lemmas~\ref{v:lem:terminal-separated-coefficients} and
\ref{v:lem:summed-terminal}, the energy bound
\eqref{eq:v-terminal-total-energy} closes the short terminal range; the
remaining range is handled in Section~\ref{sec:v-detail-10}.

First partition the inverse row sum by its finite unit class:
\begin{equation}\label{eq:v-inverse-row-unit-sector}
 k=\epsilon k_0,\quad\epsilon\in\mathcal U_6,\quad k_0\text{ allowed},
 \qquad \nu_\epsilon(n)=\nu(n)\xi_\epsilon(n),\quad
 \nu(n)\chi_n(k)=\nu_\epsilon(n)\chi_n(k_0).
\end{equation}
Each positive sector norm is estimated and the finitely many
bounds summed. On a sector, $L_\upsilon$ and the finite
profile indices are fixed before the arithmetic sums. Reciprocity
steps below use the allowed row $k_0$; we suppress $\epsilon$ in
$\nu_\epsilon$ and write $k$ for $k_0$. The coefficient group is fixed
once to contain $\nu,G_0,\Upsilon_S$ and their conjugates; the
allowed-row restriction is used for the terminal coefficient identities.

Before inverse row Poisson, enlarge the positive sector square to the
unrestricted $\Gamma$-orbit row norm with the same
$\nu_\epsilon$; since the bounds are simultaneous over the
finite coefficient group, the enlarged norm has the same required
bound. Orbit Poisson is applied to the unrestricted norm, and a later
terminal reflection is again partitioned into unit sectors, so no
sector Poisson formula without its projector is used, and no sector
projector enters an unrestricted M2 row Poisson formula. Dual
frequency-unit classes are summed
independently (Proposition~\ref{v:prop:orbit-poisson}).

All ideal-factor dyads start at norm one, so label
lengths are $\ge0$. For \(N,V,M\ge0\), a fixed coprimality indicator
\(\rho(n)=1_{(n,R_\rho)=1}\), \(Q=\log_Z q_{R_\rho}\), and
\(0\le w(f)\ll d(f)^L\) depending only on \(f\), define
\[
 \begin{aligned}
 \mathcal E(N,V,M;\rho)
 &=Z^{-N-V}\sum_{\substack{f\ \sfq\\q_f\asymp Z^V}}w(f)\\
 &\quad\times\sum_{0<q_k\le C Z^M}
 \left|\sum_{n\ \sfq}a_0(n)\nu(n)\rho(n)
     \chi_n(k)\chi_n(f)^4 W(q_n/Z^N)\right|^2.
 \end{aligned}
\]
Put \(F=N+V\). Here \(\nu\) is fixed, of finite order, with conductor in \(S\);
no arbitrary residual column coefficient is permitted, while complex annular
tests and unitary norm modes are permitted with finite seminorm and
polynomial height bounds. We prove \(\mathcal E\ll Z^{F+\eps}\) in
bounded ranges whenever
\begin{equation}\label{v:eq:unmarked-canonical-margin}
 F-M-Q\ge c>0.
\end{equation}
Then \(4F-3M\ge M+4c\) automatically. This recursion contains no selected-prime factors.

\begin{proposition}[Canonical inverse family bound]
\label{v:prop:canonical-family}
Fix $F$, $S$, the finite coefficient group, the annular test family and
bounded ranges for $N,V,M,Q$. For every $c>0$ and $\eps>0$,
\[
 \mathcal E(N,V,M;\rho)\ll_{F,S,c,\eps,W}Z^{F+\eps},
 \qquad F=N+V,
\]
uniformly over all families with $F-M-Q\ge c$, all fixed coprimality
indicators $\rho$ and all label weights $0\le w(f)\ll d(f)^L$, with finite
test seminorm and polynomial norm-twist height orders.
\end{proposition}
The proof occupies the rest of this section and
Sections~\ref{sec:v-detail-10}--\ref{sec:v-detail-11}:
\begin{itemize}
\item the bounded-row and short terminal ranges ($V<2d$, $\ell_h<2d$) are
treated below by \eqref{eq:v-terminal-total-energy};
\item each nonterminal block ($V+\ell\ge d$) is bounded by two child
families of the same shape, whose row cap is lower by more than $d$ and
whose two margins lose at most $128\eta$
(Section~\ref{sec:v-detail-10});
\item backward induction over the fixed depth
$D=1+\lceil(M_{\max}+2)/d\rceil$ closes the recursion
(Section~\ref{sec:v-detail-11}).
\end{itemize}
Theorem~\ref{v:thm:unmarked-inverse} is the case reached from the
original M\"obius square after one Poisson transform (its proof).
The inputs are
Proposition~\ref{v:prop:checked-inputs}, Proposition~\ref{v:prop:radial}
and Lemma~\ref{v:lem:primary-kummer-transport}.

For \(q_k\le C\), absolute column counting and the
label count give
\[
 \mathcal E\ll Z^{N+\eps}\le Z^{F+\eps};
\]
this covers every negative-but-nonempty formal
row scale, whose fixed window lies in \(q_k\le C\).
Bounded column scales below one are replaced by one with their exact
normalization conversion (a fixed compact dilation
family); only this replacement changes \(F\).

For a zero-extended multiplicative presentation $\Psi$ define
\begin{equation}\label{eq:v-completed-theta-definition}
 \mathcal T_W(X;\Psi)=
 \sum_{\substack{c\ \sfq,b}}
 \frac{\gamma_2(c)\bar\Lambda(cb^3)\Psi(cb^3)}{\sqrt{q_c}q_b}
 \left(\frac{q_cq_b^3}{X}\right)^{1/2}
                       W\!\left(\frac{q_cq_b^3}{X}\right).
\end{equation}
This is finite for fixed $X$; its radial test is
$W^\#(y)=y^{1/2}W(y)\in C_c^\infty(I)$. Below
$\mathcal T(X;k)=\mathcal T_W(X;\Psi_k)$, with the fixed label
and unit sector retained in $\Psi_k$.

For the other direct bound use the exact cube-M\"obius completion
\[
 Z^{-N/2}\sum_{n\ \sfq}a_0(n)\Psi_k(n)W(q_n/Z^N)
 =\sum_h\frac{\mu(h)\bar\Lambda(h)^3\Psi_k(h)^3}{q_h}
                   \mathcal T(Z^N/q_h^3;k),
\]
with \(\Psi_k(n)=\nu(n)\rho(n)\chi_n(k)\chi_n(f)^4\) and the
completed test \(y^{1/2}W(y)\). Indeed, expanding
\eqref{eq:v-completed-theta-definition} with $B=hb$, the
coefficient is
\[
 \frac{\sqrt{q_B}}{Z^{N/2}}\gamma_2(c)
 \bar\Lambda(cB^3)\Psi_k(cB^3)
 W(q_cq_B^3/Z^N)\sum_{h\mid B}\mu_F(h),
\]
and the divisor sum is $1_{B=1}$, retaining all shared-prime zeros.
In the short-range direct bound, \(\ell_h<2d\) and
\(N_*=N-3\ell_h\); the triangle inequality on this divisor dyad costs
\(\sum_{q_h\asymp Z^{\ell_h}}q_h^{-1}\ll Z^\eps\).

On long completion dyads put \(B_{\mathrm{cube}}=hb\) (the theta cube
factor times the completion divisor); its coefficient
is \(\sqrt{q_{B_{\rm cube}}}/\sqrt{Z^N}\) times a divisor-bounded
partial sum over \(h\mid B_{\rm cube}\). Choose a fine cube scale $Y_b\ge1$ with
\[
 \widehat\ell=\log_ZY_b,\quad|\widehat\ell-\ell|\le\eta,
 \qquad Y_r=Z^N/Y_b^3,\quad\widehat r=N-3\widehat\ell.
\]
Rename the combined factor $b$ below. The exact squared block
normalization and the coarse residual exponent are
\[
 Z^{-N-V}Y_b=Y_r^{-1}Y_b^{-2}Z^{-V},\qquad r=N-3\ell.
\]

Before assigning the cube coefficient, put
$\nu_f(n)=\nu(n)\rho(n)\chi_n(f)^4$, so that
$\Psi_k(n)=\nu_f(n)\chi_n(k)$. For $b=hb_0$, multiplicativity
and allowed-generator reciprocity give, including natural zeros,
\[
 \Psi_k(h)^3\Psi_k(cb_0^3)
 =\nu_f(c)\nu_f(b)^3\chi_c(k)\chi_b(k)^3
 =\nu_f(c)\nu_f(b)^3\chi_c(k)\chi_k(b)^3.
\]
Here $k=k_0$ as in
\eqref{eq:v-inverse-row-unit-sector}; the
cube unit phase $\xi_\epsilon(b)^3=\chi_b(\epsilon)^3$ stays inside
$\nu_f(b)^3$. Likewise $\bar\Lambda(h)^3\bar\Lambda(b_0)^3=\bar\Lambda(b)^3$.
Thus the row dependence of the cube
factor is $\chi_k(b)^3$, which belongs to the quadratic matrix before
any coefficient is frozen. The bounded ratio $(q_b/Y_b)^{1/2}$,
$\bar\Lambda(b)^3\nu_f(b)^3$, and the partial M\"obius divisor
sum form the divisor-bounded coefficient $\beta(b,f)$,
independent of the residual column and current row;
every coprimality condition is retained in the displayed symbols or
separate indicators.

Separate the
two-coordinate completion weight on a full fixed box:
\[
 y_b=q_b/Y_b,\qquad y_n=q_n/Y_r,\qquad W(y_b^3y_n).
\]
Its common Fourier density is chosen before \(b,n\); cube norm modes
enter \(\beta(b,f)\), and column modes enter a fresh annular \(W_i\).
Its weighted moments have finite seminorm and polynomial incoming-height
bounds by \eqref{v:eq:H2-common-density}, and this measure is
retained with the later outer Fourier parameters; so no cube-dependent
column coefficient is hidden in the nine-coordinate first profile.
No coprimality between the cube and squarefree factors is assumed. If \(V\ge2d\), the
uncompleted polynomial is the same block with
\(\ell=0,b=1,Y_b=1,Y_r=Z^N\). Hence every remaining block has \(V+\ell\ge d\):
the long completion cutoff is \(2d\), so the combined cube
length is \(\ge2d-O(1/\log Z)>d\) after the fixed threshold.

To remove fixed-support primes from the varying squarefree row,
for fixed $f$ and $R_\rho$ write uniquely
\begin{equation}\label{eq:v-terminal-row-freeze}
 k=k_{\rm pow}dR,\qquad d\mid\operatorname{rad}(fR_\rho),
 \qquad d,R\ \sfq,\qquad(dR,k_{\rm pow})=(R,fR_\rho)=1.
\end{equation}
All three factors are chosen allowed; their fixed unit sector is
retained in \eqref{eq:v-inverse-row-unit-sector}. Freeze
$k_{\rm pow}$ and $d$ before reflection. If their lengths
are $O$ and $L$, respectively, then
\[
 H\le M-O-L+O(\eta),\qquad
 \#\{k_{\rm pow}\}\ll Z^{O/2+\eps},\qquad
 \#\{d\}\le d(\operatorname{rad}(fR_\rho))\ll Z^\eps,
\]
so the outer counting exponent is $O/2$, and
every varying row prime has local digit one. This is
necessary: if $v_p(k)=v_p(f)=1$, its digit would be five, and its
reflected interaction with a squarefree column would have exponent
$(8-5)-4\equiv5\pmod6$, rather than three.

\begin{lemma}[Separated coefficients in a terminal reflected block]
\label{v:lem:terminal-separated-coefficients}
Fix the row-unit sector of \eqref{eq:v-inverse-row-unit-sector},
the supported data in \eqref{eq:v-terminal-row-freeze}, and one
local reflection branch. In particular $L_\upsilon$ is fixed before
the varying allowed row $R$. Let $D$ be its active frozen radical,
$a_D$ its supported numerator, $\epsilon$ its frequency-unit class,
and write the frequency as
\[
 \nu=\frac{\epsilon a_D n b^3}{(RD)^2},\qquad
 n,R\ \sfq,\quad(R,D)=(RD,nb)=1,\quad
 \operatorname{rad}(a_D)\mid D.
\]
At $p\mid D$ put $e_p=v_p(a_D)$ and write the local monomial as
$s_p\chi_p(t)^{c_p}$. With $L=L_\upsilon$, the exterior coefficient
and $L(\nu)^{-1}$ have the exact factorization
\begin{equation}\label{eq:v-terminal-coefficient-factorization}
 \frac{\Omega\rho(R)\beta(n,b)}{\sqrt{q_n}q_b}
                              \chi_R(nb)^3,
\end{equation}
where $|\rho|,|\beta|\le1$, $|\Omega|\le\prod_{p\mid D}|s_p|$,
and $\beta$ is independent of $R$.
\end{lemma}
\begin{proof}
Let $I$ be the product of the inactive frozen scalars, $|I|\le1$.
The explicit factors are
\begin{align*}
 \Omega={}&I L(\epsilon a_D)^{-1}L(D)^2
   \prod_{p\mid D}s_p
      \chi_p\!\left(\epsilon\frac{a_D}{p^{e_p}}
                               (D/p)^{-2}\right)^{c_p},\\
 \rho(R)={}&L(R)^2
   \prod_{p\mid R}\omega_p
      \chi_p\!\left(\epsilon a_D D^{-2}(R/p)^{-2}\right)
   \prod_{p\mid D}\chi_p(R)^{-2c_p},\\
 \beta(n,b)={}&L(n)^{-1}L(b)^{-3}
   \prod_{p\mid D}\chi_p(n)^{c_p}\chi_p(b)^{3c_p}
   \prod_{q\mid n}\gamma_2(q)
       \chi_q\!\left(\epsilon a_D(n/q)D^{-2}\right)^2.
\end{align*}
Here $\omega_p=\mathcal R_{p,1,0}(1)$ has modulus one.
Every displayed character argument is a unit on the stated support.
The digit-one formula at $p\mid R$ contributes
$\chi_R(n)\chi_R(b)^3$ in addition to $\rho$.
At $q\mid n$, put $a=v_q(b)$. The spherical valuation is $1+3a$,
so the unit exponent is $2\bmod 6$ and the cube-unit part has
exponent six; hence the formula holds also when $n,b$ share
primes. The remaining $R$-dependence from these places is
$\chi_n(R)^{-4}$, and reciprocity gives
\[
 \chi_R(n)\chi_n(R)^{-4}\chi_R(b)^3=\chi_R(nb)^3.
\]
The inverse $L$-factor splits multiplicatively as displayed;
all other factors have modulus one, and vanishing supports are retained
by zero extension.
\end{proof}

Fixing separately the finite classes of $R,n,b$ also fixes
$[\nu]=[\epsilon a_D][n][b]^3[R]^{-2}[D]^{-2}$.
so the remaining profile depends only on the three norms and is
separated on one full independent logarithmic box before
the sieve. At a digit-four active prime,
\[
 \mathcal R_{p,4,e}(t)=\tau_{p,4}^{-1}W_p(a(p^e))
       \chi_p(t)^{4+2e}\bigl(q_p1_{e\ge1}-1\bigr)
\]
separates the small term from the divisibility term. Partition
$1_{p\mid nb^3}=1_{p\mid n}+1_{p\nmid n}1_{p\mid b}$
before extracting its forced squarefree or cube factor.
The small term contributes $q_p^{-1/2}$ to the
coefficient; the other two force $p$ into the squarefree
or the cube numerator. These three
local contributions are summed in Lemma~\ref{v:lem:summed-terminal}.

\begin{lemma}[Terminal reflection with all numerator valuations summed]
\label{v:lem:summed-terminal}
Use the normalized theta and radial data of
Propositions~\ref{v:prop:checked-inputs} and \ref{v:prop:radial}.
For a fixed tuple in \eqref{eq:v-terminal-row-freeze}, partition
the remaining row into fixed-ratio annuli $q_R\asymp K=Z^H$.
Put $X=Z^{N_*}$. Its reflected squared norm is
\begin{equation}\label{eq:v-summed-terminal-energy}
 \ll Z^\eps\left(K+\frac{K^2q_D^2}{Xq_{D_{0,4}}}\right),
\end{equation}
where $D_{0,4}$ is the product of active frozen primes of digit
zero or four. Summing the supported tuples, the energy is
\begin{equation}\label{eq:v-terminal-total-energy}
 \ll Z^\eps\left(Z^M+Z^{2M-O/2+Q+V-N_*}\right).
\end{equation}
\end{lemma}
\begin{proof}
The coprimality modulus is inserted by the digit-zero identity
\eqref{eq:v-local-digit-weights} at every prime of $R_\rho$ not
already present in the row or label, retaining its unit indicator
without changing the moving digit-one row.

For a fixed active numerator $a_D$, put
$Y=K^2q_D^2/(Xq_{a_D})$. Lemma~\ref{v:lem:terminal-separated-coefficients}
gives a sum with coefficient $\Omega\rho(R)\beta(n,b)$ and
profile $\Phi(q_nq_b^3/Y)$, with its exact $q_R/K$ coordinate.
Taking $P=1$ in Lemma~\ref{v:lem:hybrid-sieve} gives, for
$q_n\asymp U$, $q_b\asymp B$,
\[
 \sum_{R\ \sfq,\ q_R\le CK}
       \left|\sum_{n,b}\beta(n,b)\chi_R(nb)^3\right|^2
       \ll(KUB)^\eps(K+UB)UB.
\]
Here $n$ is squarefree, $b$ is unrestricted, $\beta$ is
independent of $R$; no coprimality between $n,b$ is imposed.
Separating the inverse roots and the smooth norm profile on the three
full independent boxes, the weighted Fourier moments are
$\ll m(UB^3/Y)$ times one fixed seminorm, where $m(t)=\min(t^{1/4},t^{-J})$
(factor out $m$, then apply the Euler derivative
bounds of Proposition~\ref{v:prop:radial}). Hence
\[
 \|S_{a_D;U,B}\|_2
 \ll(KUB)^\eps\left(K^{1/2}B^{-1/2}+U^{1/2}\right)m(UB^3/Y).
\]

The dyadic sums, for $J>1/2$, satisfy uniformly for $Y>0$
\[
 \sum_{U,B}B^{-1/2}m(UB^3/Y)\ll1,\qquad
 \sum_{U,B}U^{1/2}m(UB^3/Y)\ll Y^{1/2}
\]
(sum over $U$ by the geometric
series for $m(t)$ and $t^{1/2}m(t)$, then use
$\sum_B B^{-1/2}\ll1$ and $\sum_B B^{-3/2}\ll1$ over dyadic $B\ge1$).
The small sieve powers are absorbed by a smaller input
$\eps$ and a larger $J$. Since $Y\le K^2q_D^2/X\le Z^{O(1)}$,
\begin{equation}\label{eq:v-fixed-numerator-norm}
 \|S_{a_D}\|_2\ll Z^\eps
       \left(K^{1/2}+\frac{Kq_D}{\sqrt{Xq_{a_D}}}\right).
\end{equation}

Multiplication by the row phase $\rho_{a_D}(R)$ is a contraction, so
Minkowski sums all $a_D$ using the local series
\[
 U_{p,j}=\sum_{e\ge0}|s_{p,j,e}|,\qquad
 V_{p,j}=\sum_{e\ge0}|s_{p,j,e}|q_p^{-e/2}.
\]
For $P=q_p$, the exact local formulas give
\[
 \begin{array}{c|c|c}
 j&U_{p,j}&V_{p,j}\\\hline
 0\text{ active}&P^{-1/2}&P^{-1/2}\\
 1,2,3,5&1&1\\
 4&1+2P^{-1/2}&\le3P^{-1/2}.
 \end{array}
\]
For the last row, the nonzero magnitudes are
\[
 |s_{p,4,0}|=P^{-1/2},\quad
 |s_{p,4,1+3a}|=(1-P^{-1})P^{-a}\ (a\ge0),\quad
 |s_{p,4,3a}|=(1-P^{-1})P^{1/2-a}\ (a\ge1).
\]
Summing the geometric series,
\[
 V_{p,4}=P^{-1/2}+(1-P^{-1})
              \frac{P^{-1/2}+P^{-2}}{1-P^{-5/2}}\le3P^{-1/2}.
\]
Thus
$\prod U_{p,j}\ll q_D^\eps$ and
$\prod V_{p,j}\ll q_D^\eps q_{D_{0,4}}^{-1/2}$,
which prove \eqref{eq:v-summed-terminal-energy}.

The primewise norm allowance is
\[
 \frac{q_R^2q_D^2}{q_{D_{0,4}}}
       \le\frac{q_k^2}{q_{k_{\rm pow}}}
                      q_fq_{\operatorname{rad}R_\rho}.
\]
Indeed, a varying row prime has charge two; a prime supported only on
$fR_\rho$ has digit zero or four and charge $\le1$;
a powerful row prime has charge $\le2$, bounded by its full
valuation; a supported valuation-one row prime has charge
$\le2\le2+v_p(f)+1_{p\mid R_\rho}$.
On the fixed-ratio annulus $K^2\asymp q_R^2$.
Summing $O(Z^{O/2+\eps})$ powerful parts, divisor-many supported
parts, and the $Z^{o(1)}$ active subsets gives
\eqref{eq:v-terminal-total-energy}.
\end{proof}

In the short terminal range $V<2d$, $\ell_h<2d$ one has
$N_*=N-3\ell_h$, $F=N+V$ and $F-M-Q\ge c$.
The two exponents in \eqref{eq:v-terminal-total-energy} are at most
\[
 M\le F-c,\qquad
 2M-O/2+Q+V-N_*
       \le F-2c+2V+3\ell_h\le F-2c+10d.
\]
Choose $d<c/100$. The fixed-label estimate is uniform and the label
count cancels $Z^{-V}$, proving the short terminal bound.

\begin{remark}[The terminal reflection in one display]
\label{v:rem:terminal-reflection-display}
For a formal statement, the terminal step combines the following data.
\begin{itemize}
\item \emph{Fixed data:} a row-unit sector \eqref{eq:v-inverse-row-unit-sector}; a
frozen tuple $(k_{\rm pow},d)$ of \eqref{eq:v-terminal-row-freeze};
one local reflection branch with active frozen radical $D$, numerator
$a_D$ and frequency-unit class $\epsilon$; and a varying squarefree
allowed row $R$ with $q_R\asymp K$, $(R,fR_\rho)=1$.
\item \emph{Reflection:} the reflected frequency is
$\nu=\epsilon a_Dnb^3/(RD)^2$, with $n$ squarefree, $b$ unrestricted
and $(R,D)=(RD,nb)=1$. The reflected coefficient factors as
\eqref{eq:v-terminal-coefficient-factorization}, with the row entering
only through $\rho(R)\chi_R(nb)^3$. The reflected profile is
$\Phi(q_nq_b^3/Y)$, where
\[
 Y=\frac{K^2q_D^2}{Xq_{a_D}},\qquad X=Z^{N_*},
 \qquad |(x\partial_x)^j\Phi(x)|\ll\min(x^{1/4},x^{-A})
\]
(Proposition~\ref{v:prop:radial}).
\item \emph{Sieve:} with $P=1$, Lemma~\ref{v:lem:hybrid-sieve} reduces to
the quadratic estimate \eqref{v:eq:transport-quadratic} after the
coprimality expansion. For $q_n\asymp U$ and $q_b\asymp B$ it gives
$\ll(KUB)^\eps(K+UB)UB$.
\end{itemize}
The output is
\[
 \text{(reflected energy)}\ll Z^\eps\Bigl(K+\frac{K^2q_D^2}{Xq_{D_{0,4}}}\Bigr),
\]
that is, \eqref{eq:v-summed-terminal-energy}. Summed over the frozen
data this is \eqref{eq:v-terminal-total-energy}. The only analytic
input is the theta functional equation with vanishing weighted constant
terms (Proposition~\ref{v:prop:checked-inputs},
Lemma~\ref{v:lem:theta-moment-identification}). The only sieve input
is the optimal quadratic family sieve of \cite{GL}; the cubic sieve,
with its $(AB)^{2/3}$ term, is not used in the terminal range.
\end{remark}

\subsection{Two Poisson transforms and their positive-norm estimates}
\label{sec:v-detail-10}
In this subsection, we treat the nonterminal range $V+\ell\ge d$:
Proposition~\ref{v:prop:radial}, the orbit formula of
Section~\ref{sec:v-detail-5} (summand
\eqref{eq:v-inverse-first-full-summand}) and a second Gauss conversion
(summand \eqref{eq:v-inverse-second-full-summand}) return the block to
the same inverse coefficient family, with new row $k'=d_kd_2k''$, label
$f'=JCd_2v'$ and a common Fourier measure. The block is bounded by child
families at exponent $\lambda_c+C_c+2F'$; Section~\ref{sec:v-detail-11}
closes the recursion.

A nonterminal block has \(V+\ell\ge d\). Expand its square, write
\(n_i=A_i C u_i\), where \(A_i=(n_i,\operatorname{rad}(b_1b_2))\),
and \(C\) is the remaining gcd. The residual \(u_i\) are squarefree,
mutually coprime and prime to the fixed factors. At each cube prime put
\(a_i=1_{p\mid A_i}\), \(\pi=v_p(b_1b_2)\bmod2\),
\(t_p=a_1-a_2+3\pi\pmod6\), and
\(R_1=\prod_{t_p\ne0}p\). Define also the integral frequency factor
\[
 E=\prod_{\substack{p\mid b_1b_2\\v_p(b_1b_2)\ \mathrm{odd}}}
                         p^{a_1+a_2}.
\]
The primitive modulus is
\(m_1=u_1u_2R_1\); the complementary coprimality indicator is
\(C\operatorname{rad}(b_1b_2)/R_1\).
Apply the unrestricted Poisson identity with coprimality conditions, with
divisor \(d_k\) of that indicator and kernel argument
\(Z^M q_h/(q_{d_k}q_{m_1})\).

The first principal modulus forces $u_1=u_2=1$, $\pi=0$,
$a_1=a_2$, so the two squarefree columns agree and
$b_1b_2$ is a square. With $b_i=gt_i^2$, $g$ squarefree,
the cube-pair count at the cube scale is
\[
 \sum_{q_g\ll Y_b}O\bigl((Y_b/q_g)^{1/2}\bigr)^2
                  \ll Y_b^{1+\epsilon}.
\]
The column, label and row counts are $\ll Y_rZ^\epsilon$,
$Z^{V+\epsilon}$ and $Z^{M+\epsilon}$; with the squared
normalization $Z^{-N-V}Y_b$ the product is
\[
 \ll Z^{M+\epsilon}\frac{Y_b^2Y_r}{Z^N}
       = Z^{M+\epsilon}/Y_b
       \le Z^{M+\epsilon}\le Z^{F-c+\epsilon}.
\]
Divisor weights are absorbed in the same loss.

The complete M\"obius
extension also restores the case $u_1=u_2=1$, $R_1=1$. Let
$\mathcal D_1$ be its original principal block, $\mathcal O_1$ its
dual origin, and $\mathcal T_1^{\rm pr}$ its signed nonzero dual sum;
modulus-one Poisson gives
$\mathcal D_1=\mathcal O_1+\mathcal T_1^{\rm pr}$.
Hence, with $\mathcal T_1^{\rm all}$ the full extended nonzero
expression below, the original first square is
$\mathcal O_1+\mathcal T_1^{\rm all}$ (equivalently,
$\mathcal D_1-\mathcal T_1^{\rm pr}+\mathcal T_1^{\rm all}$).
Artificial diagonal pairs with $u_1=u_2\ne1$ cancel in the complete
M\"obius sum, and the primitive Gauss and column factors at $1$ agree
with the modulus-one identity.

The principal restoration has
$q_h\le Cq_{d_k}/Z^M$ and $y=hf^2$, contained in the same compact support
and primewise row enclosure used below; the column and cube-pair counts
apply to both principal terms. At each $d_k$,
$(Z^M/q_{d_k})\sum_{[h]\ne0}|K_1(Z^Mq_h/q_{d_k})|\ll1$;
the divisor-many $d_k$ cost $Z^\epsilon$, and $1\le Z^M$.
The origin is $O(Z^{M+\epsilon})$. Hence both
principal terms are $\ll Z^{M+\epsilon}/Y_b$.
Their signed grouping precedes the first absolute Fourier measure.

For the nonprincipal part the exact first Gauss identity yields
\(\mu(u)(\nu G_0)(u)\), with row \(y=hf^2E\).
Write \(b_1b_2=q_0^2J_2^2s\), where \(s\) is the odd-parity
radical and \(J_2\) is the product of the even-parity cube primes
with \(a_1=a_2=1\); the remaining half-exponents define \(q_0\), and
\(J=sJ_2\) is squarefree.

At a cube prime let
$e_1,e_2$ be the two residual CRT exponents. If $t_p\ne0$, the
Gauss CRT rule gives
\[
 e_1\equiv4a_1+1+t_p,\qquad e_2\equiv4a_2+1-t_p\pmod6;
\]
if $t_p=0$, $e_i\equiv4a_i$ (the $4a_i$ from
the cubic Gauss coefficients, the $1$ from the residual first
Gauss conversion, $\pm t_p$ from the common primitive character).
With $e_y\equiv-\pi(a_1+a_2)\pmod6$,
\[
\begin{array}{ccc|c|cc|c|c}
 \pi&a_1&a_2&t_p&e_1&e_2&e_y&e_i-e_y\pmod6\\\hline
 0&0&0&0&0&0&0&0\\
 0&1&0&1&0&0&0&0\\
 0&0&1&5&0&0&0&0\\
 0&1&1&0&4&4&0&4\\
 1&0&0&3&4&4&0&4\\
 1&1&0&4&3&3&5&4\\
 1&0&1&2&3&3&5&4\\
 1&1&1&3&2&2&4&4
\end{array}
\]
The factor $\bar\chi_u(E)$ in $y=hf^2E$ supplies
$e_y$, leaving exponent $4v_p(J)$; the
common factor supplies $\chi_u(C)^4$ and the Poisson divisor
$\chi_u(d_k)$. Hence the remaining coefficient is
\[
 \chi_u(CJ)^4\chi_u(d_k)1_{(u,q_0)=1}.
\]
This includes nonunits: primes in $J$ keep their
zero-extended fourth powers, all other cube primes occur in
$q_0$, and no exponent $\equiv0$ is replaced by $1$ on a nonunit.

The primewise inequality
\[
 \frac{q_{R_1}q_Eq_J}{q_{A_1}q_{A_2}}
                  \le q_{b_1b_2}^{2}
\]
gives a column-independent outer ball for \(y\), since the local
left exponent is
\[
 1_{t_p\ne0}+\pi(a_1+a_2)+v_p(J)-a_1-a_2
 =\begin{cases}2,&\pi=1,\\-a_1a_2,&\pi=0,\end{cases}
\]
which is at most $2v_p(b_1b_2)$. At each cube prime
of total valuation \(t\ge1\), the exponent of \(q_0J\) is
\((t+\pi)/2\), and every odd-parity prime belongs to \(R_1\).
Thus \(q_{q_0J}\le(q_{b_1b_2}q_{R_1})^{1/2}\) (counting) and
\(q_{q_0}\le q_{b_1b_2}^{1/2}\) (coprimality modulus), and every cube
support lies in \(\operatorname{rad}(q_0J)\).

For logarithmic
centers $B,\delta,j,R,t$ of $C,d_k,J,R_1,t'$, use the complete
M\"obius identity
\[
 1_{(u_1,u_2)=1}=\sum_{t'\mid u_1,u_2}\mu_F(t'),
 \qquad u_i=t'x_i.
\]
Before any off-coprime extension the full nonprincipal
summand, at fixed $A_i,C,b_i,f,d_k$ and one unit class, is
\begin{equation}\label{eq:v-inverse-first-full-summand}
 \begin{aligned}
 &\frac{Z^{M-N-V}Y_b\Theta_0}
             {q_{d_k}\sqrt{q_{R_1}q_{u_1}q_{u_2}}}
       c(u_1;y)\overline{c(u_2;y)}\\
 &\quad\times W_1\!\left(\frac{q_{A_1}q_Cq_{u_1}}{Y_r}\right)
       \overline{W_2\!\left(\frac{q_{A_2}q_Cq_{u_2}}{Y_r}\right)}
       K_1\!\left(\frac{Z^Mq_h}
                       {q_{d_k}q_{R_1}q_{u_1}q_{u_2}}\right),\\
 &c(u;y)=\mu_F(u)\nu_1(u)\rho(u)\bar\chi_u(y)
                   \chi_u(CJ)^4\chi_u(d_k)1_{(u,q_0)=1},\\
 &\nu_1=\nu G_0\overline{\xi_{\epsilon_1}},\qquad y=hf^2E.
 \end{aligned}
\end{equation}
Here $\Theta_0$ is the product of the frozen $A_i,C,b_i,f,d_k$
coefficients and the normalized primitive Gauss factor on $R_1$.
The eight local CRT states remove its $u_i$-dependence;
its remaining frequency and unit phases are outer scalars of
modulus at most a fixed divisor power, and fixed support indicators
(with their zeros) are kept in this product or in $c$.
$K_1$ is the norm-sheet kernel of the first row lift.
The common profile ($K_1$, the two $W$ factors and the smooth
inverse-root ratios) is evaluated on
the nine independent coordinates below, before positivity.
A common enclosing frequency ball is imposed on the genuine source,
its raw tail removed by the radial bounds.

The first transformed formula is extended to artificial noncoprime
pairs with the same finite coefficients, kernel and inverse-root
normalization; the full M\"obius sum annihilates them (primitive
Poisson was applied only on the coprime locus). After this
signed identity, and before positivity, set
\[
 H_c=2r-2B-M+4\ell+2V+\delta-j,
 \quad H_{\rm use}=H_c+32\eta+\tau,\quad
 s_i=r-A_i-B-t,
\]
\[
 \kappa_i=M-2r-2\ell-V-\delta+A_i+B-R/2.
\]

The common first profile uses the nine norm coordinates
\(A_1,A_2,C,d_k,R_1,t',h,x_1,x_2\), before remaining ideals are selected.
The prefactor
$Z^{M-N-V}Y_b/(q_{d_k}\sqrt{q_{R_1}q_{u_1}q_{u_2}})$
is extracted coarsely as $Z^{(\kappa_1+\kappa_2)/2}$, with the
actual-scale error allocated in Section~\ref{sec:v-detail-11};
there is no second root for $t'$.
Let $Y_{A_i},Y_C,Y_t$ denote the actual fine factor scales and put
\[
 Y_{s_i}=Y_r/(Y_{A_i}Y_CY_t),\qquad
 \widehat s_i=\widehat r-\widehat A_i-\widehat B-\widehat t.
\]
These scales are fixed before the ideals
are evaluated; the $s_i$ are coarse exponents only.
The pure side polynomials are
\[
 P_i(\lambda,y)=\sum_{x\ \sfq}\mu(x)\nu_1(x)
  \bar\chi_x(y)\chi_x(CJ)^4\chi_x(d_k)
  \rho(x)1_{(x,q_0t')=1}W_i(q_x/Y_{s_i}).
\]

All Fourier norm modes are in \(W_i\) or in outer phases; e.g.\
\(|q_h^{iv}|=|q_y^{iv}q_f^{-2iv}q_E^{-iv}|=1\).
At fixed \(\lambda,y\), \(f^2E\mid(y)\), so the weighted fibre of the
old label is divisor-bounded. Bound it on the genuine ball,
then take the common absolute Fourier measure and Cauchy over
\((\lambda,y,\xi)\). Apart from the root-extraction factor
$Z^{15\eta/2}$ allocated in Section~\ref{sec:v-detail-11}, this bounds
the block by the geometric mean of
\[
 Z^{\kappa_i}\int\sum_{\lambda,y}
       \Psi_\Phi(q_y/Z^{H_{\rm use}})|P_i(\lambda,y)|^2\,d|\mu_1|(\xi).
\]
The old divisor weight is removed before the fresh lift, and no first kernel
remains in either column square.

In this fresh square put \(x_j=g'z_j\), take its Poisson divisor
$d_2\mid g'$, and insert the complete residual indicator
\[
 1_{(z_1,z_2)=1}=\sum_{v'\mid z_1,z_2}\mu_F(v'),
 \qquad z_j=v'n_j.
\]
Before the $v'$ extraction, put $K_i=Z^{H_{\rm use}}$,
$H_o=CJ$, $r_o=\operatorname{rad}(q_0t')$, and define
\[
 \begin{aligned}
 B(g')&=\mu_F(g')\nu_1(g')\rho(g')1_{(g',r_o)=1}
                            \chi_{g'}(H_o)^4\chi_{g'}(d_k),\\
 L_2(z)&=a_0(z)\nu_2(z)\rho(z)1_{(z,r_o)=1}
        \chi_z(H_o)^4\chi_z(d_k)\bar\chi_z(d_2)\chi_z(k''),\\
 \nu_2&=\nu_1\bar G_0\xi_{\epsilon_2}.
 \end{aligned}
\]
Up to the already fixed positive outer measure, that summand is
\begin{equation}\label{eq:v-inverse-second-full-summand}
 \begin{aligned}
 &\mu_F(d_2)|B(g')|^2 L_2(z_1)\overline{L_2(z_2)}
                 1_{(z_1z_2,g')=1}
                     \frac{K_i}{q_{d_2}\sqrt{q_{z_1}q_{z_2}}}\\
 &\quad\times K_2\!\left(\frac{K_iq_{k''}}
                              {q_{d_2}q_{z_1}q_{z_2}}\right)
        W_i(q_{g'}q_{z_1}/Y_{s_i})
                  \overline{W_i(q_{g'}q_{z_2}/Y_{s_i})}.
 \end{aligned}
\end{equation}

$K_2$ comes from the fresh positive row lift (no old divisor weight).
Reciprocity cancels the two primitive
cross factors; all remaining column phases are in $L_2$, and
$1_{(z_1z_2,g')=1}$ is kept before the $v'$ sum.
The enclosing ball for $k'=d_kd_2k''$ and the smooth
kernel agree on the genuine locus and its extension, which, as in
the first step, is fixed before the signed sum is estimated. The
second Gauss identity returns $\nu G_0$ to
$\nu$, up to the retained unit characters, so the transformed
coefficient stays in the fixed finite group containing
$\nu,G_0,\Upsilon_S$. The new row, label and coprimality modulus are
\[
 k'=d_kd_2k'',\qquad f'=JCd_2v',\qquad
 \gamma=(q_0,t',g'/d_2),\quad
 \rho_\gamma(n)=\rho(n)1_{(n,\operatorname{rad}(q_0t'g'/d_2))=1}.
\]

The common factor at $v'$ is retained before projection.
With $\nu_2$ as in
\eqref{eq:v-inverse-unit-return}, the returned coefficient
before extracting $z=v'n$ is
\[
 A(z)=a_0(z)\nu_2(z)\rho_\gamma(z)
         \chi_z(d_k k'')\bar\chi_z(d_2)\chi_z(CJ)^4.
\]
For squarefree coprime $v',n$ its extraction gives
\[
 A(v'n)=A(v')a_0(n)\nu_2(n)\rho_\gamma(n)
    \chi_n(d_k k'')\bar\chi_n(d_2)\chi_n(CJv')^4.
\]
Hence $|A(v')|^2$ stays outside the second column norm, and
$A(v')=0$ unless $(v',g'CJk'')=1$ (zero extensions and
$\rho_\gamma$). With the earlier squarefree supports,
$J,C,d_2,v'$ are pairwise coprime and squarefree.
Only on this genuine support do we use $|A(v')|\le1$
and form the positive projected majorant.
The zero-preserving coefficient identity is
\[
 \chi_n(d_k)\chi_n(k'')\bar\chi_n(d_2)\chi_n(CJ)^4\chi_n(v')^4
       =\chi_n(k')\chi_n(f')^4,
\]
so there is no new moving column character. The common second profile uses
\(g',d_2,v',k'',n_1,n_2\) before \(\gamma,f',k'\) are fixed.

For centers \(g,\theta,v\) of \(g',d_2,v'\), put
\[
 N_c=r-A_i-B-t-g-v,\quad V_c=B+\theta+v+j,\quad F_c=N_c+V_c,
\]
\[
 M_c=M-4\ell-2A_i-2t-2g+2\theta-2V+j,
\]
\[
 \lambda_c=H_c-\theta-(s_i-g),\qquad
 C_c=\ell+R/2-j+t+g-\theta.
\]
The genuine structural set projects once to distinct $\gamma$, before
new row or label selection, with containing balls
\[
 q_{q_0}\ll Y_bY_{R_1}^{1/2}/Y_J,\qquad q_{t'}\ll Y_t,
 \qquad q_{g'/d_2}\ll Y_g/Y_\theta,
\]
($Y$ the actual fine scales), of divisor-weighted total
$\ll Z^{C_c+11\eta/2+\epsilon_1}$. The coprimality
modulus is bounded separately by the stronger
$q_{q_0}\le q_{b_1b_2}^{1/2}\ll Y_b$; $C_c$ is a counting exponent only.

At fixed \(\gamma,f',k'\), split \(f'\) into \(J,C,d_2,v'\), then
\(J=sJ_2\), \(b_1b_2=q_0^2J_2^2s\); all choices are divisor-bounded,
\(d_k\mid f'\operatorname{rad}(q_0)\), and \(g',k''\) are determined
when integral (the erased \(f,h\) are not counted again).
Hence the fibre is \(\le A_\gamma w_\gamma(f')\) with
\(w_\gamma(f')\ll d(f')^L\) depending only on \(f'\).

Normalize the single projection measure as
\(d\nu(\gamma)=Z^{-C_c}\sum A_\gamma\delta_\gamma\), of mass
$\ll Z^{11\eta/2+\epsilon_1}$ (allocated in
Section~\ref{sec:v-detail-11}). With $Y_g,Y_v,Y_\theta,Y_J$ the
remaining factor scales, the child column and label scales are
\[
 Y_N=Y_{s_i}/(Y_gY_v),\qquad Y_V=Y_CY_\theta Y_vY_J.
\]
After at most a bounded subunit conversion, let $N',V',F'=N'+V'$
be the actual nonnegative child parameters and
$\mathcal E^{\rm act}_{a,\gamma,\omega}$ the normalized family
at these scales. Cauchy on $(\gamma,f',k',\omega)$
and the common second absolute Fourier measure give, apart from the
root and projection errors listed in Section~\ref{sec:v-detail-11},
\[
 |\mathcal E_i^{\rm off}|
 \ll Z^{\lambda_c+C_c+F'+o(1)}
 \prod_{a=1}^2\left(\int\mathcal E^{\rm act}_{a,\gamma,\omega}\,
                  d\nu(\gamma)d|\mu_2|(\omega)\right)^{1/2}.
\]
Hence $\mathcal E^{\rm act}_{a,\gamma,\omega}\ll Z^{F'+\eps'}$
gives $Z^{\lambda_c+C_c+2F'+o(1)+\eps'}$, with one projection count:
one $F'$ restores the normalization, the other bounds the energy.
Relative to the coarse identity the change is $2(F'-F_c)\le26\eta$
(proved below). The second principal block and the diagonal
reintroduced by the signed extension are treated in
\eqref{eq:v-inverse-second-principal-bound}, with nominal exponent
$F-B-j$.

\subsection{Error terms, bounded ranges and finite induction}
\label{sec:v-detail-11}
In this subsection, the losses of Section~\ref{sec:v-detail-10} are
allocated within \eqref{v:eq:unmarked-canonical-margin}, with internal
orders chosen backward through the finite depth via
\eqref{v:eq:H2-common-density}: each nonterminal step lowers the row cap
by more than $d$ and loses at most $128\eta$ in either canonical margin,
and the second principal block obeys
\eqref{eq:v-inverse-second-principal-bound}. This gives Theorems~\ref{v:thm:unmarked-inverse} and \ref{v:thm:unmarked-seed}
(zero-freeness for $\Rea s>11/12$) and
Lemma~\ref{v:lem:general-amplification}.

With \(D_c=\ell+A_i+t+g-\theta+V\), the exact identities are
\[
 F_c=F-D_c-2\ell+j,\quad M_c=M-2D_c-2\ell+j,
 \qquad \kappa_i+\lambda_c+C_c+2F_c=F.
\]
For a canonical family with common pre-loss margin $c>0$, fix
$c_*=c/2$ before all mesh choices. Choose $d<c_*/400$ and put
$D=1+\lceil(M_{\max}+2)/d\rceil$.
\paragraph{Actual fine scales and coarse exponents.}
Every arithmetic annulus retains its actual fine scale $Y_a\ge1$,
with fixed support ratios independent of $Z$. Write
$\widehat a=\log_ZY_a$ and let $a$ be its coarse center, so that
$|\widehat a-a|\le\eta$; coarse centers are used only in affine
exponent estimates, never for $Y_a$ inside a test.
The parent parameters $N,V,M$ are actual. In a cube block put
\[
 \widehat r=N-3\widehat\ell,
 \qquad |\widehat r-r|\le3\eta.
\]
Products and quotients of fine scales are actual,
with fixed support constants left in the profiles, so
the child column and label have logarithmic scales
\[
 \widehat N=\widehat r-\widehat A_i-\widehat B
              -\widehat t-\widehat g-\widehat v,
 \qquad
 \widehat V=\widehat B+\widehat\theta+\widehat v+\widehat j.
\]
Then $|\widehat N-N_c|\le8\eta$,
$|\widehat V-V_c|\le4\eta$, and $\widehat V\ge0$ (fine
factor scales start at one).

If $Z^{\widehat N}<1$ and
its annulus is nonempty, then $Z^{\widehat N}$ is bounded below by a
fixed positive support endpoint; only then do we dilate the scale to
one. The dilation lies in a fixed compact interval, retains a fixed
annular test family, and costs $O(1/\log Z)\le\eta$ in logarithmic
scale after the threshold. Let $N'\ge0$ be the
resulting column length, $V'=\widehat V$,
$F'=N'+V'$. Then
\[
 N'\ge\widehat N,\qquad
 |F'-F_c|\le13\eta.
\]
If $N_c<0<\widehat N$, the child theorem is applied at
$\widehat N$ without dilation; a negative coarse
center alone never selects a bounded endpoint.

\paragraph{Enclosing the actual rows.}
Use the slightly larger first row lift
\[
 H_{\rm use}=H_c+32\eta+\tau.
\]
In the first support inequality the center errors are $\le6\eta$
(two squarefree-column norms), $4\eta$ (square of the cube pair) and
$\eta,2\eta,\eta$ ($d_k,C^2,J$), total $14\eta$ (the parent label
scale is actual). Fixed support constants and the radial cutoff
are absorbed in the rest of $32\eta+\tau$ after the threshold,
so the first lift contains every genuine first row.

For the second transform each $z_i$ has logarithmic scale
\[
 \widehat r-\widehat A_i-\widehat B-\widehat t-\widehat g,
\]
within $7\eta$ of $s_i-g$. In its support inequality the two $z_i$,
$d_k$ and $d_2^2$ have total error $\le14\eta+\eta+2\eta=17\eta<32\eta+\tau$,
the first row allowance. Hence, after
absorbing fixed support constants, the new row ball has cap
\[
 M'\le M_c+\eta.
\]
A negative actual row cap is stopped by the bounded-row estimate only
when its actual ball is bounded, not on a negative formal center alone.

\paragraph{The second principal block at its actual scales.}
Isolate the positive diagonal $x_1=x_2$ in each fresh first-side
square. For $K_i=Z^{H_{\rm use}}$, a fixed outer tuple and first
Fourier parameter, its raw contribution $\mathcal D_i^{\rm raw}$ is at most
\[
 \sum_{[y]\ne0}\Psi_\Phi(q_y/K_i)
       \sum_{x\ \sfq}|W_i(q_x/Y_{s_i})|^2
 \ll Z^{\epsilon_1}K_iY_{s_i},
\]
since all arithmetic coefficients (also with selected-prime
divisibility weights) are bounded by a fixed divisor power; the row
bound follows from ideal counting and the lift's sheet
bounds, the column bound from the fixed-ratio annulus, also at a
nonempty bounded subunit column scale. The
seminorm factor is polynomial in the incoming
Fourier height and is integrated by
\eqref{v:eq:H2-common-density}.

For the principal terms reintroduced by the complete
M\"obius extension, let
$\mathcal O_i$ be the diagonal's dual origin and
$\mathcal T_i^{\rm pr}$ its signed nonzero dual sum.
Since $\Phi(0)=0$, the additive origin contributes zero, and
modulus-one Poisson with coprimality divisor
$d_2\mid g'=x$ gives
\[
 \mathcal D_i^{\rm raw}=\mathcal O_i+\mathcal T_i^{\rm pr}.
\]
At the origin the divisor coefficient is
\[
 K_i\sum_{d_2\mid g'}\frac{\mu_F(d_2)}{q_{d_2}}
   =K_i\prod_{p\mid g'}(1-q_p^{-1}),
\]
times the fixed compact-quotient integral and
$|B(g')W_i(q_{g'}/Y_{s_i})|^2$. Hence
\[
 |\mathcal O_i|\ll Z^{\epsilon_1}K_iY_{s_i},\qquad
 |\mathcal T_i^{\rm pr}|\ll Z^{\epsilon_1}K_iY_{s_i},
\]
the signed principal-frequency sum being grouped before
the second absolute Fourier measure (no termwise absolute bound is
claimed).

The complete M\"obius sum annihilates every artificial pair
$z_1=z_2\ne1$, and the diagonal $z_1=z_2=1$ of the full
extended nonzero-frequency expression $\mathcal T_i^{\rm all}$ is
$\mathcal T_i^{\rm pr}$ (the second summand's Gauss and column
factors at $1$ are $1$). So the off-diagonal source is
$\mathcal T_i^{\rm all}-\mathcal T_i^{\rm pr}$, and the
fresh square is
\[
 \mathcal D_i^{\rm raw}
   +\mathcal T_i^{\rm all}-\mathcal T_i^{\rm pr}
       =\mathcal O_i+\mathcal T_i^{\rm all}.
\]
Thus the principal allowance below covers the raw origin and
the diagonal-restoration correction, and only $\mathcal T_i^{\rm all}$
is passed to the positive child bound. For the band-limited lift of
Lemma~\ref{v:lem:compact-dual}, on $z_1=z_2=1$,
\[
 q_{k'}\le C_\Phi q_{d_k}q_{d_2}^{2}/K_i,
\]
the second support inequality with $q_{z_1}q_{z_2}=1$; so the
common enclosing ball retains the whole signed
principal sum.

For the outer structural count in this principal block, choose $C,J,t'$ at their actual scales; then
\[
 q_{q_0}\ll Y_bY_{R_1}^{1/2}/Y_J.
\]
The $O(Y_J)$ choices of $J$ cancel the denominator
$Y_J$. Given $q_0,J$, split $J=sJ_2$, form
$b_1b_2=q_0^2J_2^2s$; splitting into $b_1,b_2$ and choosing
$A_1,A_2$ and $d_k\mid C\operatorname{rad}(b_1b_2)$ cost
fixed divisor powers, other restrictions only
reduce the count, and $f,h$ (erased by the
first-stage divisor-bounded fibre) are not counted again. Hence
\[
 \#\{\text{outer structural tuples}\}_{\rm weighted}
 \ll Z^{\epsilon_1}Y_CY_tY_bY_{R_1}^{1/2}.
\]

The first-side prefactor (bounded inverse-root ratios
left in the common profile) is
\[
 \mathscr A_i=
 \frac{Z^{M-N-V}Y_bY_{A_i}Y_C}
      {Y_{d_k}Y_{R_1}^{1/2}Y_r},
 \qquad Y_{s_i}=\frac{Y_r}{Y_{A_i}Y_CY_t}.
\]
Put
\[
 K_{\rm fine}=
 \frac{Z^{2V-M}Y_r^2Y_b^4Y_{d_k}}{Y_C^2Y_J}.
\]
Since $|\log_ZK_{\rm fine}-H_c|\le14\eta$,
$K_i\le Z^{46\eta+\tau}K_{\rm fine}$.
Multiplying the diagonal bound, the structural count and the
prefactor, and integrating the common Fourier measure,
\begin{equation}\label{eq:v-inverse-second-principal-bound}
 \begin{split}
 \mathcal D_i
 &\ll Z^{\epsilon_1+\pi}
      \mathscr A_i K_iY_{s_i}Y_CY_tY_bY_{R_1}^{1/2}\\
 &=Z^{\epsilon_1+\pi}
      \frac{Z^{M-N-V}K_iY_b^2Y_C}{Y_{d_k}}\\
 &\ll
 \frac{Z^{F+46\eta+\tau+\pi+\epsilon_1}}{Y_CY_J}
 \le Z^{F+46\eta+\tau+\pi+\epsilon_1},
 \end{split}
\end{equation}
since $Y_rY_b^3=Z^N$ and $Y_C,Y_J\ge1$.
Its coarse exponent is $F-B-j$; the origin and the signed principal
restoration satisfy the same bound up to a constant. This
allowance lies in the $256\eta$ edge margin and needs
no new projection count or arbitrary residual-column
coefficient.

\paragraph{Canonical margins at the actual scales.}
Define
\[
 \widehat D=\widehat\ell+\widehat A_i+\widehat t+
       \widehat g-\widehat\theta+V,
\]
\[
 \widehat F=\widehat N+\widehat V
       =F-\widehat D-2\widehat\ell+\widehat j,
 \quad
 \widehat M=M-2\widehat D-2\widehat\ell+\widehat j.
\]
Here $|\widehat M-M_c|\le13\eta$, so $M'\le\widehat M+14\eta$.
The actual bound $q_{q_0}\ll Z^{\widehat\ell}$ and the fine scales
of $t'$ and $g'/d_2$ give
\[
 Q'\le Q+\widehat\ell+\widehat t+
                   \widehat g-\widehat\theta+\eta
       \le Q+\widehat D+\eta.
\]
(constants absorbed after the threshold; this uses the
coprimality bound on $q_0$, not the projection-count bound).

The divisibilities $d_2\mid g'$ and $J\mid\operatorname{rad}(b_1b_2)$
give
\[
 \widehat g-\widehat\theta\ge-\eta,
 \qquad \widehat j\le2\widehat\ell+\eta.
\]
Before row and modulus allowances, the exact identities are
\[
 \begin{aligned}
 \widehat F-\widehat M-
       (Q+\widehat\ell+\widehat t+\widehat g-\widehat\theta)
      &=F-M-Q+\widehat A_i+V,\\
 4\widehat F-3\widehat M
      &=4F-3M+2\widehat A_i+2\widehat t
           +2(\widehat g-\widehat\theta)+2V+\widehat j.
 \end{aligned}
\]

The surviving selected-prime cap is $\le z_0$, with
its original prime scales. Since $F'\ge\widehat F$,
\[
 \begin{aligned}
 F'-M'-Q'-z'&\ge F-M-Q-z_0-15\eta,\\
 4F'-3M'-6z'&\ge4F-3M-6z_0-44\eta.
 \end{aligned}
\]
We use $128\eta$ for either per-step margin loss. Moreover
\[
 M'\le M-2(V+\ell)+32\eta,\qquad F'\le F+32\eta.
\]
Indeed the formula for $\widehat M$ and the
last two support inequalities give
$\widehat M\le M-2(V+\widehat\ell)+3\eta$, and the $14\eta$
row allowance and $|\widehat\ell-\ell|\le\eta$ give the first
bound; similarly
$\widehat F\le F+2\eta$, and the subunit correction is
$\le\eta$.

\paragraph{Conservative edge accounting.}
The losses are: first inverse-root extraction
$\le(1+5+1+1/2)\eta=15\eta/2$; second
$\le(8+1)\eta=9\eta$; first row enlargement $32\eta+\tau$;
projection count $11\eta/2$; cube normalization
$|\widehat\ell-\ell|\le\eta$; total
$55\eta+\tau$. The two child normalization restorations cost
$\le2|F'-F_c|\le26\eta$. Thus the transformed
edge costs $\le81\eta+\tau+\pi+\epsilon'$; we use the uniform
margin $256\eta+\tau+\pi+\epsilon'$ for it and for the second principal
term, whose loss $46\eta+\tau+\pi+\epsilon_1$ by
\eqref{eq:v-inverse-second-principal-bound} also lies in it.
This margin concerns norm exponents only; all tests keep
their actual fine scales, so there is no derivative-dependent
support factor $Z^{J\eta}$.

Choose
\[
 \eta<\min\{d/128,\ c_*/(512D),\ \epsilon/(4096D)\},
 \qquad \tau,\pi<\epsilon/(16D).
\]
Then every nonterminal branch ($V+\ell\ge d$) drops the actual row cap
by more than $d$, and its two canonical margins lose
$\le128D\eta<c_*/4$, so the margin stays $\ge c_*/2$; the edge and
tail margins fit the final loss after decreasing the subsidiary input
losses, and the same fixed-depth backward choice of finite seminorm and
Fourier-height orders applies. At depth \(D\)
a surviving nonterminal would have row cap $<-2$, hence be
empty after the fixed support threshold; bounded endpoints are
stopped by the absolute bound. Backward induction gives
\(\mathcal E\ll Z^{F+\eps}\). Transformed total lengths are
\(\le F_{\max}+32D\eta\), so all divisor and conductor ranges are bounded.

Order of choices: raw norm-tail decay orders; fresh row
lifts through the fixed depth (one finite test family
before \(Z\)); finite seminorm and Fourier-height orders, backward
through the recursion; external Fourier-height
cutoff orders; the large-\(Z\) threshold. By \eqref{v:eq:H2-common-density},
every nested height polynomial is integrated without a
derivative-dependent norm loss, and the seminorm normalization
following Lemma~\ref{v:lem:finite-depth-loss} applies to these
positive measures as well.

\begin{lemma}[Closure of the actual inverse families]
\label{v:lem:inverse-returned-family}
Consider a canonical inverse family of
Section~\ref{sec:v-detail-9}, allowing the product-form marks of
Section~\ref{sec:v-detail-13}. Write $F=N+V$ and let $z_0=0$ in
the unmarked case. At every nonterminal two-transform step the returned
families are canonical families at their actual scales, with
\begin{align*}
 F'-M'-Q'-z'&\ge F-M-Q-z_0-15\eta,\\
 4F'-3M'-6z'&\ge4F-3M-6z_0-44\eta,\\
 M'&\le M-2(V+\ell)+32\eta.
\end{align*}
Each outer structural family is counted once. Its common absolute
Fourier measure integrates the child profile seminorms and incoming
height polynomials with a subsidiary exponent loss. Consequently the
induction of this section applies to the actual returned sums, with
edge loss at most $256\eta+\tau+\pi+\epsilon'$.
\end{lemma}
\begin{proof}
Start with the signed full summand
\eqref{eq:v-inverse-first-full-summand}; make the complete
M\"obius extension before applying Cauchy--Schwarz. At the second
return the new row and squarefree label are, literally,
\[
 k'=d_kd_2k'',\qquad f'=JCd_2v',\qquad
 \gamma=(q_0,t',g'/d_2).
\]
The residual deletion is
$\rho_\gamma(n)=\rho(n)1_{(n,\rad(q_0t'g'/d_2))=1}$.
After extracting $z=v'n$, its variable coefficient has the form
\[
 a_0(n)\nu_2(n)\rho_\gamma(n)
 \chi_n(d_kk'')\overline{\chi_n(d_2)}\chi_n(CJv')^4
   =a_0(n)\nu_2(n)\rho_\gamma(n)\chi_n(k')\chi_n(f')^4.
\]
The equality includes zeros: the four factors of $f'$ are squarefree
and pairwise coprime, and the extracted coefficient at $v'$ is retained
before its absolute value is bounded. The finite twist $\nu_2$ is
fixed within each family. For marks use the exact priority expansion
of Section~\ref{sec:v-detail-13}; restored exclusions either lie in
the displayed common deletion or force a natural zero. The surviving
prime lists therefore retain their original profiles and introduce no
coefficient depending on the new row.

At fixed $(\gamma,f',k')$, the decomposition of $f'$ and the
reconstruction of $g',k''$ are divisor-bounded. The erased old $f,h$
are not new choices. Thus the single projection measure is exactly
\[
 d\nu(\gamma)=Z^{-C_c}\sum_\gamma A_\gamma\delta_\gamma,
 \qquad \nu(\Omega)\ll Z^{11\eta/2+\epsilon_1}.
\]
Together with the common Fourier measure it gives the product of two
integrated positive energies displayed in
Section~\ref{sec:v-detail-10}. No replacement by an integral of a
pointwise geometric mean is needed: bounding each integrated energy
by its inductive bound uses this same mass once on each side of the
square root, hence once in the final product.

The scales are $Y_N=Y_{s_i}/(Y_gY_v)$ and
$Y_V=Y_CY_\theta Y_vY_J$. A nonempty subunit scale is dilated to one
only by the bounded conversion proved above. The enclosing row cap,
the deletion-modulus bound and both margin estimates are precisely
the actual-scale calculations of this section; in particular $C_c$
is not used as a modulus exponent. The identity
$\kappa_i+\lambda_c+C_c+2F_c=F$ restores the parent normalization;
$2(F'-F_c)\le26\eta$ is charged on this edge.

Appendix~\ref{app:inverse-return-measures} makes the two separate
weighted integrals explicit. For weighted integration, first fix child
derivative and height orders.
All their new heights are linear combinations of the incoming heights
and the common Fourier variables. Their seminorms are bounded by a
finite parent seminorm times a polynomial in these variables.
Equation~\eqref{v:eq:H2-common-density}, at a sufficiently high
weighted order, integrates that polynomial. Discrete divisor weights
are already in $d\nu$ and are not counted again. Repeating this
choice backward through the depth $D$ fixes all orders before $Z$;
there is no derivative-order multiple of $\eta$ in the exponent.
Finally the principal block is grouped with its reintroduced diagonal
before absolute Fourier integration, as in
\eqref{eq:v-inverse-second-principal-bound}. The direct terminal
bounds and the actual row drop above exhaust the remaining branches.
The edge accounting preceding this lemma gives the stated loss.
\end{proof}

\begin{theorem}[Inverse moment]
\label{v:thm:unmarked-inverse}
With the theta representation of Proposition~\ref{v:prop:checked-inputs} and the functional equation for radial test functions of Proposition~\ref{v:prop:radial},
for every fixed finite-order \(\nu\) of conductor supported in $S$, every \(c,\eps>0\), and annular
smooth \(W\),
\[
 \sum_{0<q_k\le C_0H}\left|D^{-1/2}
     \sum_n\mu_F(n)\nu(n)\chi_n(k)W(q_n/D)\right|^2
 \ll_{F,S,\nu,c,\eps,W,C_0}H(HD)^\eps.
\]
Here $k$ runs over nonzero exterior ideals, with the finitely many
unit classes retained as in Proposition~\ref{v:prop:orbit-poisson},
and $H\ge\max(2,D^{1+c})$. One may fix any finite group of
coefficient characters containing the initial characters, $G_0$ and
$\Upsilon_S$, and take the constant simultaneously over this group.
The estimate has finite test seminorm and polynomial norm-twist height
orders. Natural zeros and both common orientations are retained.
\end{theorem}
\begin{proof}
Lemma~\ref{v:lem:inverse-returned-family} supplies closure of the
canonical families used after the following initialization.
Write \(H=Z^m,D=Z^r\), with \(r<m\) and a fixed positive relative
gap. Fix $c_*=(m-r)/2>0$ before $d,D,\eta$ and
$\tau_{\rm init}$. Choose the initialization losses so that
\[
 11\eta+\tau_{\rm init}<c_*,\qquad
 32\eta+3\tau_{\rm init}<c_*,
\]
so both initialized margins below are $\ge c_*$.

Expand the inverse square with its full coprimality indicators, set
\(n_i=Cz_i\), and use primitive Poisson summation with coprimality conditions with divisor
\(d'\mid C\). Put \(P_1=B-\theta\), for the centers of \(C,d'\).
Insert $1_{(z_1,z_2)=1}=\sum_{s\mid z_1,z_2}\mu_F(s)$,
$z_i=sn_i'$, with the transformed extension of
Section~\ref{sec:v-detail-10}. The Gauss factor returns \(\mu(z)\gamma_{-1}(z)=a_0(z)\bar G_0(z)\).
The new row, label and coprimality modulus are \(k'=d'h'\), \(f'=d's\),
\(\operatorname{rad}(C/d')\), and CRT gives the
coefficient \(a_0(n)\nu_0(n)\chi_n(k')\chi_n(f')^4\), with
\(\nu_0=\bar\nu\bar G_0\xi_{\epsilon}\) after whole-product conjugation, where $\epsilon$ is the dual frequency-unit class.
All residual zeros are retained.

The common six-coordinate profile is \(C,d',s,h',n_1',n_2'\).
At fixed \(t=C/d',f',k'\), reconstruction is divisor-bounded; project
once to the distinct \(t\) and normalize their count by \(Z^{P_1}\).
The formal exponents are
\[
 N_c=r-B-v,\quad V_c=\theta+v,\quad F_c=r-P_1,
 \quad M_c=2r-m-2P_1,
 \quad \kappa_c=m-2r+P_1,
\]
with
\[
 F_c-M_c-P_1=m-r>0,\qquad
 \kappa_c+P_1+2F_c=m.
\]

With actual fine scales $Y_C,Y_{d'},Y_s$ put
\[
 \widehat P_1=\log_Z(Y_C/Y_{d'}),\quad
 Y_N=Z^r/(Y_CY_s),\quad Y_V=Y_{d'}Y_s.
\]
After at most a bounded subunit conversion, let $F'$ be the
column-plus-label exponent. Then
\[
 F'\ge r-\widehat P_1,\qquad |F'-F_c|\le5\eta.
\]

Use the actual enclosing row cap, for example
$M'=\max(0,M_c+6\eta+\tau_{\rm init})$ on a nonempty positive branch;
a bounded or empty ball is stopped directly. Since
$|2r-m-2\widehat P_1-M_c|\le4\eta$ and
$Q'\le\widehat P_1+\eta$,
\[
 \begin{aligned}
 F'-M'-Q'&\ge m-r-11\eta-\tau_{\rm init},\\
 4F'-3M'&\ge3m-2r-32\eta-3\tau_{\rm init}.
 \end{aligned}
\]
Here $\widehat P_1\ge-\eta$ by $d'\mid C$ and fixed
support ratios. Incoming tests use $Y_N,Y_V$ unchanged (no scale-one
substitution if $N_c<0<\log_ZY_N$).
The initial normalization changes by $\le2|F'-F_c|\le10\eta$;
we use a $128\eta+4\tau_{\rm init}$ margin for the initial step,
absorbed by $m-r$ for $\eta,\tau_{\rm init}$ small.

The principal diagonal and its restoration are $\ll Z^{m+\eps}$;
the first prefactor and projection losses lie in the same margin.
Choosing all losses inside \(\eps\) and summing the
logarithmic dyads proves the moment; $Z=H$ converts the
relative gap to the condition on \(H,D\), and bounded scales
are absorbed in the constant.
\end{proof}

\begin{theorem}[The preliminary eleven-twelfths bound]
\label{v:thm:unmarked-seed}
Let $F\supset\mu_{12}$ be a fixed cyclotomic field. Then every
finite-order Hecke $L$-function over $F$ is nonzero for
\[
 \Rea s>11/12.
\]
The proof uses only the uniform inverse estimate of
Theorem~\ref{v:thm:unmarked-inverse}.
\end{theorem}
\begin{proof}
Fix $\nu$ and choose the final set $S$ at the outset as in
Section~\ref{sec:v-scope}, including its conductor.
For a fixed smooth annular test define
\[
 M_k(T)=T^{-1/2}\sum_n\mu_F(n)\nu(n)\chi_n(k)W(q_n/T).
\]

Scale bound: fix $c>0$,
$D\ge2$, $H=D^{1+c}$ and $T_0=b_W^{-1}$, where
$\operatorname{supp}W\subset[a_W,b_W]\Subset(0,\infty)$.
The sum vanishes for $T<T_0$, and
\[
 \frac{d}{d\log T}M_k(T)=M_k^{W^{[1]}}(T),
 \qquad W^{[1]}(y)=-\tfrac12W(y)-yW'(y).
\]
Sobolev on unit $\log T$-intervals and summation over
$O(1+\log D)$ intervals give
\[
 \begin{split}
 \sum_{0<q_k\le H}\sup_{T_0\le T\le D}|M_k(T)|^2
 &\ll\sum_{0<q_k\le H}\int_{T_0}^{eD}
    \bigl(|M_k(T)|^2+|M_k^{W^{[1]}}(T)|^2\bigr)\frac{dT}{T}\\
 &\ll_{c,\delta,W}H(HD)^\delta
 \end{split}
\]
for every $\delta>0$, after increasing the fixed threshold:
for bounded $T$ use the finite column count, and for the other
$T\le eD$ apply Theorem~\ref{v:thm:unmarked-inverse}
with relative gap $c/2$, as
$H\ge(eD)^{1+c/2}$ for large $D$, with loss $<\delta$.

For every exterior ideal $a$,
$\chi_n(a^6)=1_{(n,a)=1}$, and Euler multiplication gives
\[
 M_1(D)=\sum_{d\mid\operatorname{rad}(a)}
       \mu_F(d)\nu(d)q_d^{-1/2}M_{a^6}(D/q_d).
\]
Indeed the Dirichlet series of the right side is
\[
 \prod_{p\mid a}(1-\nu(p)q_p^{-s})
 \prod_{p\nmid a}(1-\nu(p)q_p^{-s})
       =1/L_F^S(s,\nu)
 \qquad(\Re s>1),
\]
and $q_d^{-1/2}$ restores the central
normalization at scale $D/q_d$.

For every $\delta>0$,
\[
 \sum_{d\mid\operatorname{rad}(a)}q_d^{-1/2}
 =\prod_{p\mid a}(1+q_p^{-1/2})\ll_{F,\delta}q_a^\delta.
\]
Hence, uniformly for $q_a\le H^{1/6}$,
\[
 |M_1(D)|^2\ll H^\delta
              \sup_{T_0\le T\le D}|M_{a^6}(T)|^2.
\]
There are $\gg_{F,S}H^{1/6}$ exterior ideals with $q_a\le H^{1/6}$,
with distinct sixth powers $q_{a^6}\le H$.
Summing and applying the scale bound gives,
after redistributing $\delta$,
\[
 |M_1(D)|^2\ll H^{5/6}(HD)^\delta.
\]
For every $\epsilon>0$, choosing $c,\delta>0$ sufficiently small gives
\begin{equation}\label{eq:v-seed-annular-inverse}
 \sum_n\mu_F(n)\nu(n)W(q_n/D)
 \ll_{F,S,\nu,W,\epsilon}D^{11/12+\epsilon}.
\end{equation}

Take a smooth dyadic partition of unity on positive norms. On a
compact subset of $\Re s>11/12$, choose $\epsilon$ smaller than its
distance from $11/12$ and apply \eqref{eq:v-seed-annular-inverse} to
$W_s(y)=W(y)y^{-s}$, whose seminorms are uniformly bounded
there. The dyadic Dirichlet blocks are then bounded
by a convergent geometric series and define a holomorphic
$G_\nu(s)$ on $\Re s>11/12$, equal to $1/L_F^S(s,\nu)$ on
$\Re s>1$. Meromorphic continuation and the identity theorem give
\[
 L_F^S(s,\nu)G_\nu(s)=1,
\]
so $L_F^S$ has no zero in this half-plane
(a principal pole at $1$ is canceled by a zero of $G_\nu$).
The finitely many omitted Euler factors are nonzero for $\Re s>0$, so
the complete Hecke function is zero-free there as well.
\end{proof}

\begin{remark}[Strength of the inverse moment and where its cancellation enters]
\label{v:rem:inverse-strength}
The proof above uses no input other than
Theorem~\ref{v:thm:unmarked-inverse}. For $F=\mathbb Q(\mu_{12})$ and
$\nu=1$ it gives $\zeta_F(s)\ne0$ for $\Re s>11/12$, hence
$\zeta(s)\ne0$ there, since $\zeta$ divides $\zeta_F$. The
implication itself is elementary: the $\gg H^{1/6}$ rows $k=a^6$ have
$\chi_n(k)=1_{(n,a)=1}$, the Euler identity relates $M_1(D)$ to
$M_{a^6}(D/q_d)$, and positivity gives
$H^{1/6}|M_1(D)|^2\ll H^{1+\eps}$. Consequently
Theorem~\ref{v:thm:unmarked-inverse} in the range $H=D^{1+c}$ with
$c\downarrow0$ is at least as strong as a zero-free strip for $\zeta$,
and every later result of this Part depends on it, on
Theorem~\ref{v:thm:marked-inverse-new} and on
Theorem~\ref{v:thm:M2-assembled}. The cancellation enters at the
following three places.
\begin{enumerate}
\item \emph{One exact Poisson step in the row.} With $n_i=Cz_i$ and
$z_1,z_2$ coprime squarefree, $\chi_{z_1}\bar\chi_{z_2}$ is primitive
modulo $z_1z_2$, reciprocity cancels the two cross factors, and the
Gauss factor is $\mu(z)\gamma_{-1}(z)=a_0(z)\bar G_0(z)$
(Section~\ref{v:sec:unmarked-assembly}). Schematically,
\[
\begin{aligned}
 \sum_{q_k\le H}|M_k(D)|^2
 &=(\text{diagonal}\asymp H)
 +(\text{common-divisor terms})\\
 &\quad+\frac{H}{D^2}\sum_{q_h\ll D^2/H}
   \Bigl|\sum_{q_z\asymp D}a_0(z)\nu_0(z)\bar\chi_z(h)W(q_z/D)\Bigr|^2,
\end{aligned}
\]
with the exact summand \eqref{eq:v-inverse-first-full-summand}. Thus the
theorem is equivalent to an optimal large sieve $(M+N)N$ for the
Gauss-coefficient family $a_0$ with $M=N^{1-c}$ rows, all rows
(including sixth powers) retained; this is the content of
Proposition~\ref{v:prop:canonical-family}.
\item \emph{No residual main terms.} The coefficient $a_0=\bar\Lambda\gamma_2$
carries the angular Hecke character $\Lambda$ of nonzero degree at every
complex place. The weighted theta constant terms are pure norm
characters, so they vanish after angular integration (the argument
preceding Proposition~\ref{v:prop:radial}). Hence there is no
Patterson-type residue, and in particular no term of the shape
$(MN)^{2/3}$ coming from residues.
\item \emph{Terminal reflection to a quadratic family.} In the terminal
range the theta functional equation turns the sextic-twisted sums into
dual sums with the quadratic symbols $\chi_R(nb)^3$ and dual length
$K^2q_D^2/(Xq_{a_D})$
(Lemmas~\ref{v:lem:terminal-separated-coefficients}
and~\ref{v:lem:summed-terminal};
Remark~\ref{v:rem:terminal-reflection-display}). The optimal quadratic
large sieve \eqref{v:eq:transport-quadratic} then gives the energy
$K+K^2/X$.
\end{enumerate}
For arbitrary coefficients the cubic large sieve
\eqref{v:eq:transport-cubic} retains $(AB)^{2/3}$, and so does every
published bound for cubic or sextic residue symbols. The present
argument avoids that term only through items (2) and (3), which use the
specific theta coefficients.
\end{remark}

\begin{lemma}[Sixth-power amplification for general rows]
\label{v:lem:general-amplification}
Let the finite seminorm and height bounds be those of
Theorem~\ref{v:thm:unmarked-inverse}. Fix $\nu$ of conductor supported in $S$,
$W\in C_c^\infty(0,\infty)$, and a bounded range of $r\ge0$.
For $U\ge2$ put $D=U^r$ and
\[
 M_u(T)=T^{-1/2}\sum_n\mu_F(n)\nu(n)\chi_n(u)W(q_n/T).
\]
Then, for every $\epsilon>0$,
\begin{equation}\label{eq:v-general-amplification}
 \sum_{\substack{U\le q_u\le2U\\u\text{ sixth-power-free}}}
       |M_u(U^r)|^2
 \ll_{F,S,\nu,W,r,\epsilon}U^{e(r)+\epsilon},
 \qquad e(r)=\max\{1,(1+5r)/6\}.
\end{equation}
The same assertion holds for either common orientation and for a
fixed finite family of $\nu$, with polynomial norm-twist height cost.
\end{lemma}
\begin{proof}
The scale Sobolev argument in the proof of
Theorem~\ref{v:thm:unmarked-seed} gives, for $c,\delta>0$,
\[
 \sum_{q_k\le H}\sup_{T_0\le T\le D}|M_k(T)|^2
       \ll H(HD)^\delta,
 \qquad H\ge\max\{2,(eD)^{1+c}\}.
\]
For every exterior ideal $a$, Euler multiplication gives the exact
identity, including all natural zeros,
\begin{equation}\label{eq:v-general-amplification-euler}
 M_u(D)=\sum_{d\mid\operatorname{rad}(a)}
   \mu_F(d)\nu(d)\chi_d(u)q_d^{-1/2}M_{ua^6}(D/q_d).
\end{equation}
Indeed at $p\mid a$ the two Euler factors are
$1-\nu(p)\chi_p(u)q_p^{-s}$ and $1$, and at $p\mid u$ both
are $1$; the twist stays in the row
$ua^6$ and $\nu$ is unchanged.

The divisor mass is $\ll q_a^\delta$. Choose
$H\asymp\max\{U,D^{1+c}\}$, with a sufficiently large fixed
constant, and average \eqref{eq:v-general-amplification-euler}
over the $\gg(H/U)^{1/6}$ ideals with $q_a\le(H/(2U))^{1/6}$. Since
$(u,a)\mapsto ua^6$ is injective for sixth-power-free $u$,
\[
 (H/U)^{1/6}\sum_u|M_u(D)|^2
 \ll (HD)^{2\delta}
       \sum_{q_k\le H}\sup_{T_0\le T\le D}|M_k(T)|^2,
\]
whence the left side of \eqref{eq:v-general-amplification} is
$\ll U^{1/6}H^{5/6}(HD)^{3\delta}$.
Choosing first $c$ and then $\delta$ small on the bounded $r$-range
gives \eqref{eq:v-general-amplification}; finitely many unit
classes and fixed characters change only the constant.
\end{proof}

\section{The inverse moment with prime-polynomial factors}
\label{v:sec:marked-assembly}
In this section, Theorem~\ref{v:thm:marked-inverse-new} extends
Theorem~\ref{v:thm:unmarked-inverse} to the selected prime
factors of Proposition~\ref{prop:iv-count}: $\sum_{q_k\le CZ^m}
|M_k(Z^r)Q_k|^2\ll Z^{m+\eps}$ under $r+2z<m$ and $2r+8z<3m$. It is used
in Section~\ref{v:sec:compensated-assembly}, and
Section~\ref{sec:v-detail-29} combines it with the unweighted product
moment. Beyond the inverse moment, it needs exact priority extraction
with restoration by natural zeros (Section~\ref{sec:v-detail-13}) and the
terminal estimate with an active prime tuple
(Section~\ref{sec:v-detail-14}), with the same theta functions and Kummer
sieves of Lemma~\ref{v:lem:primary-kummer-transport}.

\subsection{Exact extraction and coefficient independence}
\label{sec:v-detail-13}
In this subsection, selected primes are extracted in the fixed priority
order \eqref{v:eq:mark-priority}, so no remaining coefficient depends on
the row or a remaining column, and restoring the surviving prime lists
after the two Poisson transformations of Section~\ref{sec:v-detail-10}
is an exact signed identity. Lemma~\ref{v:lem:marked-terminal-typed}
gives the terminal polynomial used in Section~\ref{sec:v-detail-14}.

Let disjoint prime lists \(\mathcal P_i\) and bounded coefficients
\(b_i(p)\) be independent of the current row and label. The prime-divisor weight is
\[
 c_{\rm sel}(n)=\prod_{i\in I}\sum_{p\in\mathcal P_i}
                   b_i(p)1_{p\mid n},\qquad
 q_p\asymp Z^{z_i},\quad\sum_{i\in I}z_i\le z_0.
\]
A fixed number of prime-sum factors is allowed in divisor powers and constants.
At every extraction use the complete priority identity
\begin{equation}\label{v:eq:mark-priority}
 1_{p\mid D_1\cdots D_hn}
 =\sum_{j=1}^h1_{p\mid D_j}\prod_{a<j}1_{p\nmid D_a}
   +1_{p\mid n}\prod_{a\le h}1_{p\nmid D_a}.
\end{equation}
It partitions prime tuples before any positive bound; assigned
primes are fixed outer factors, and surviving primes retain their
coefficients and lists. Fourier norm modes multiply coefficients
by unitary prime powers, with the common density chosen before the tuple.

The priority factors \(D_j\) here are extracted column and cube factors.
The old row and label exclusions stay in the zero-extended character
(the erased old label is not a prime-list
restriction to restore). After the two Poisson transformations, the
priority exclusions are supported in
\[
 \operatorname{rad}(k'f'q_0t'r_gR_\rho),\qquad r_g=g'/d_2.
\]
Here \(J,C,d_2,v'\) belong to the new squarefree label \(f'\),
\(q_0,t',r_g,R_\rho\) belong to its fixed coprimality modulus, and new row factors
belong to \(k'\).

If a surviving prime lies in this support and divides
its residual column, the zero-extended coefficient is zero.
Hence restoring the surviving lists on the transformed sum is an
exact polynomial identity (no term of a signed polynomial is deleted
by positivity). Assigned tuples have only the displayed
label/outer divisor multiplicity and are counted once.
Residual column coefficients remain \(a_0\nu\) times
product-form divisibility marks; no arbitrary column coefficient
is introduced.

Cube completion retains the selected-prime weight on the full index:
\[
 c_{\rm sel}(h^3nb^3)=c_{\rm sel}(n(hb)^3).
\]
To separate its dependence, expand the prime tuple $(p_i)_{i\in I}$
and put $c=hb$. For each $i$ use the exact priority identity
\[
 1_{p_i\mid nc^3}
 =1_{p_i\mid n}+1_{p_i\nmid n}1_{p_i\mid c}.
\]
After fixing the tuple and the subset $J\subseteq I$ assigned to $n$,
its contribution is
\[
 \Bigl(\prod_{i\in I}b_i(p_i)\Bigr)
 \Bigl(\prod_{i\in J}1_{p_i\mid n}
       \prod_{i\notin J}1_{p_i\nmid n}\Bigr)
 \prod_{i\notin J}1_{p_i\mid c}.
\]
The first parenthesis is fixed, the second is an explicit column
factor, and only the surviving cube mark depends solely on $c=hb$.
The tuples and assignments remain fixed outer indices; the column
indicators are not hidden in varying prime lists.
On a short completion-divisor branch, apply priority first to $h$
in the surviving cube mark, using
$1_{p\mid hb}=1_{p\mid h}+1_{p\nmid h}1_{p\mid b}$.
On a long branch, with $c=hb$ fixed, the surviving cube mark is
outside the partial M\"obius sum over $h$ and introduces no additional
residual-column dependence. All original coprimality conditions and
natural zeros remain in force.
The two-coordinate completion profile is separated before the cube and
column are fixed (Section~\ref{v:sec:unmarked-assembly}), so the
nine/six-coordinate densities hide no tuple-dependent
column coefficient. The bounded-row direct bound remains valid, since
\(|c_{\rm sel}(n)|\ll d(n)^K\).

\begin{lemma}[Coefficient after selecting prime divisors]\label{v:lem:marked-terminal-typed}
After fixing factors other than the selected primes, local branches and the separate finite
power classes of the valuation-one row \(R\) and active mark tuple \(P\),
the reflected polynomial has the form
\[
 A(R)\sum_P\frac{a(P)}{\sqrt{q_P}}1_{(R,P)=1}
 \sum_{n,b}\beta(n,b)\chi_R(nb)^3\chi_P(n)^{-2}
                         1_{(P,b)=1},
\]
with \(a\) independent of \(R,n,b\) and \(\beta\) independent of \(R,P\).
The angular theta character is retained in its separated factors.
\end{lemma}
\begin{proof}
Apply \(1_{p\mid A}=1-\chi_p(A)^0\) to the full theta index
before reflection. The local digit table gives an inactive coefficient
\(q_p^{-1}\); the active branch has magnitude \(q_p^{-1/2}\), dual
valuation zero and unit exponent two. A valuation-one row has unit
exponent one. The exact unit-support pair exponents are
\[
 R:P=-2-4=0\pmod6,\qquad R:n=1-4=3,\qquad R:b=3,
\]
\[
 P:n=2-4=-2,\qquad P:b=3\cdot2=0\pmod6.
\]
Each zero exponent retains its unit indicator, and collisions with
old coprimality indicators vanish before any negative symbol power is used.
Let $L$ be the fixed unitary character of $F_S^\times$ in
the raw reflected theta identity, restricted to allowed generators and
chosen after the unit classes, as in
\cite[Section 6, the character $L_{\upsilon}$]{Artin}; its angular
part and the fixed finite-order coefficient are retained. Up to a
scalar depending only on the frozen factors, it contributes
\[
 L\!\left(\frac{(RP)^2}{nb^3}\right)
       =L(n)^{-1}L(b)^{-3}L(R)^2L(P)^2,
\]
which separates into
dual-only, row-only and tuple-only coefficients. With the
\(R,P\) power classes fixed, the dual finite profile index depends
only on \(n,b\). A base branch with unit exponent \(4+2e\)
has exponent four on \(e=3a\) and zero on \(e=3a+1\); priority
extraction removes a forced squarefree occurrence before a cube occurrence,
and shared cube factors and natural zeros remain. Separating the
common norm profile before the indices gives the asserted
independence. No finite-order claim about the angular
component of \(L\) is used.
\end{proof}

\subsection{The direct estimate for selected prime tuples with both length inequalities}
\label{sec:v-detail-14}
In this subsection, the terminal polynomial of
Lemma~\ref{v:lem:marked-terminal-typed} is estimated with its active
tuple retained, the sieves of Section~\ref{sec:v-detail-4} applying with
independent coefficients; this closes every short terminal branch under
both length inequalities. The margin identities of
Section~\ref{sec:v-detail-11} then carry both inequalities through the
recursion, proving Theorem~\ref{v:thm:marked-inverse-new}.

Take the normalized family of Section~\ref{sec:v-detail-9} with the
weight $c_{\rm sel}$ and assume
\[
 F-M-Q-z_0\ge c,\qquad 4F-3M-6z_0\ge c.
\]
Let $K=Z^H$, $L=Z^{z_a}$, $X=Z^{N_*}$ and let $a_D$ be a fixed
numerator supported on the active frozen denominator $D$.
The active marked-prime product satisfies $q_P\asymp L$ (fixed ratio),
with coefficient $a(P)/\sqrt{q_P}$ kept inside the hybrid norm. Put
\[
 Y=\frac{K^2q_D^2L^2}{Xq_{a_D}},\qquad
 m(x)=\min(x^{1/4},x^{-J}).
\]
The actual kernel includes the additional independent ratio
$(q_P/L)^{-2}$ and $(q_R/K)^{-2}$. On $q_n\asymp U$, $q_b\asymp B$,
Lemma~\ref{v:lem:hybrid-sieve}, after the inverse roots are extracted,
gives
\[
 \|S_{a_D;U,B}\|_2^2\ll Z^\epsilon
 \frac{K+UB}{UB}\bigl(U+L+(UL)^{2/3}\bigr)m(UB^3/Y)^2.
\]
Here $m$ is factored from the full independent norm profile
before its common Fourier density is taken. For $U,L\ge1$,
\[
 U+L+(UL)^{2/3}\ll U+L\sqrt U,
 \qquad
 1+L/U+(L^2/U)^{1/3}\ll L,
\]
since $(UL)^{2/3}\le U$ if $U\ge L^2$ and $\le L\sqrt U$ if
$U\le L^2$. Hence
\[
 \|S_{a_D;U,B}\|_2\ll Z^\epsilon
 \left(\frac{\sqrt{KL}}{\sqrt B}+\sqrt U
                     +\sqrt L\,U^{1/4}\right)m(UB^3/Y).
\]

The dyadic summation in Lemma~\ref{v:lem:summed-terminal} gives
\[
 \sum_{U,B}U^{1/4}m(UB^3/Y)\ll Y^{1/4},
\]
by the convergent series for $x^{1/4}m(x)$ and
$\sum_BB^{-3/4}$. Thus, uniformly in $a_D$,
\[
 \|S_{a_D}\|_2\ll Z^\epsilon
       \left(\sqrt{KL}+\sqrt Y+\sqrt L\,Y^{1/4}\right).
\]
Alongside $U_{p,j},V_{p,j}$ introduce
\[
 W_{p,j}=\sum_{e\ge0}|s_{p,j,e}|q_p^{-e/4}.
\]
For active digit zero, $W_{p,0}=q_p^{-1/2}$; for $j=1,2,3,5$,
$W_{p,j}=1$. For digit four,
\[
 W_{p,4}=q_p^{-1/2}+(1-q_p^{-1})
       \frac{q_p^{-1/4}+q_p^{-5/4}}{1-q_p^{-7/4}}
       \ll q_p^{-1/4}.
\]
(an absolutely convergent geometric series). Minkowski,
with the bounded row and tuple coefficients retained for every
$a_D$, gives
\[
 \|S\|_2^2\ll Z^\epsilon\left(
 KL+\frac{K^2q_D^2L^2}{Xq_{D_{0,4}}}
     +\frac{Kq_DL^2}{\sqrt{Xq_{D_{0,4}}}}\right).
\]
(only the hybrid sieve,
radial derivative bounds and local Weyl coefficients are used).

Summing the supported squarefree and powerful row parts
with the primewise allowance of Lemma~\ref{v:lem:summed-terminal} gives the
exponents
\[
 M+z_a,\qquad
 2M-O/2+Q+V-N_*+2z_a,\qquad
 M+(Q+V-N_*)/2+2z_a.
\]
Under $F-M-Q-z_0\ge c$, $4F-3M-6z_0\ge c$, $z_a\le z_0$,
and $V,\ell_h<2d$, these are respectively at most
\[
 F-c,\qquad F-2c+10d,\qquad F-M/4-3c/4+5d.
\]
(for the last, use $Q\le F-M-z_0-c$, then
$3z_0/2\le F-3M/4-c/4$). Taking $d<c_*/400$ closes every
short terminal branch with both margins. Inactive marked primes
have coefficient $q_p^{-1}$ and total mass
$\sum_{p\in\mathcal P_i}q_p^{-1}\ll Z^\epsilon$; they are fixed
by the triangle inequality before the active tuple norm.

The priority identity keeps an original prime position or assigns
it to a frozen factor, so $z'\le z_0$.
The support identity preceding Section~\ref{sec:v-detail-5} is
\[
 \chi_n(d_k)\chi_n(k'')\bar\chi_n(d_2)
 \chi_n(CJ)^4\chi_n(v')^4
 =\chi_n(d_kd_2k'')\chi_n(JCd_2v')^4.
\]
At $d_2$ the powers agree since $1+4\equiv-1\pmod6$, and the zeros
agree at every support prime. With new row $d_kd_2k''$, label
$JCd_2v'$ and fixed modulus
$\operatorname{rad}(q_0t'g'/d_2)R_\rho$, a surviving prime on this
support has zero coefficient whenever its mark requires it to
divide the residual column, so restoring its list
is an equality of signed polynomials before positive norms;
assigned tuples remain in their frozen outer divisor count.

For the two margins, Section~\ref{sec:v-detail-11} gives
\[
 \begin{aligned}
 \widehat F-\widehat M-
 (Q+\widehat\ell+\widehat t+\widehat g-\widehat\theta)
 &=F-M-Q+\widehat A_i+V,\\
 4\widehat F-3\widehat M
 &=4F-3M+2\widehat A_i+2\widehat t
   +2(\widehat g-\widehat\theta)+2V+\widehat j.
 \end{aligned}
\]
Subtracting $z'$ or $6z'$ and retaining the stated scale allowances yields
\[
 F'-M'-Q'-z'\ge c-15\eta,\qquad
 4F'-3M'-6z'\ge c-44\eta,
 \qquad M'\le M-2(V+\ell)+32\eta.
\]
The prime scales are unchanged, so the finite-depth recursion
uses both margins with the $128\eta$ and
$256\eta$ edge margins; bounded nonempty row caps use the direct
bound (no induction at a negative center).
Profiles, the bound for the prime sums, cap ranges and selection rules are fixed before
\(Z\); for each \(Z\) the prime lists may vary within their
annuli, and the estimates are uniform over these lists and bounded
coefficients.

\begin{theorem}[Inverse moment with prime-polynomial factors]\label{v:thm:marked-inverse-new}
With the theta representation of Proposition~\ref{v:prop:checked-inputs} and the sieve estimates of Lemma~\ref{v:lem:primary-kummer-transport}, let
\[
 M_k(Z^r)=Z^{-r/2}\sum_n\mu_F(n)\nu(n)\chi_n(k)W(q_n/Z^r),
\]
\[
 Q_k=\prod_i P_i^{-1/2}\sum_{p\in\mathcal P_i}
                      b_i(p)\nu(p)\chi_p(k)W_i(q_p/P_i),
 \qquad P_i=Z^{z_i},\quad z=\sum_i z_i.
\]
For each fixed positive pair of margins, if
\[
 r+2z<m,\qquad2r+8z<3m,
\]
then
\[
 \sum_{0<q_k\le C Z^m}|M_k(Z^r)Q_k|^2\ll Z^{m+\eps}.
\]
All natural coprimality indicators and common coprimality indicators of fixed norm, surviving subcollections, both common
orientations, finite profile seminorms and polynomial height orders are
retained. The prime coefficients are row-independent and bounded, with
one fixed number of prime-sum factors. A growing initial coprimality modulus requires its own exponent contribution
in the initial margins and is not included in this assertion. 
\end{theorem}
\begin{proof}
Before the mesh and allocation choices fix
\[
 c_*={1\over2}\min\{m-r-2z,\ 3m-2r-8z\}>0,
\]
with all initialization losses below $c_*$. Fix the prime tuple
\(P\), assign \(j=(n,P)\) of nominal length \(G\),
and write \(n=jn_0,P=jP_0\). The squarefree residual \(c=n_0P_0\) has length
\(D'=r+z-2G\), the surviving cap is \(z_0=z-G\), and
\[
 Z^{-(r+z)/2}=Z^{-G}Z^{-D'/2},\quad
 \psi_k(n)\psi_k(P)=\psi_k(j)^2\psi_k(c),\quad
 \mu(n_0)=\mu(c)\mu(P_0).
\]
The triangle inequality is taken in the row Hilbert space before squaring; the
\(O(Z^{G+\eps})\) assigned tuples cancel \(Z^{-G}\), and
\(\psi_k(j)^2\) is a contraction (including zeros). The residual coprimality modulus
contains \(j\); the common norm profile is chosen with \(j\) fixed
before surviving primes and the residual column, and prime weights
remain their original coefficients times unitary norm modes.

The inverse initialization is exactly that of
Theorem~\ref{v:thm:unmarked-inverse}, with \(r\) replaced by \(D'\) and
fixed coprimality modulus \(\operatorname{rad}((C/d')j)\). Its two formal margins
are
\[
 F_c-(2D'-m-2P_1)-(P_1+G)-z_0=m-r-2z+2G,
\]
\[
 4F_c-3(2D'-m-2P_1)-6z_0
               =3m-2r-8z+10G+2P_1.
\]

Here $D'=r+z-2G$ uses the actual product of the
assigned prime scales; since the number of prime positions is fixed
and the prime windows have fixed ratios, $q_j$ divided by this product
lies in a fixed compact interval, so the profile stays uniformly
smooth at $Z^{D'}$ (unlike a coarse $Z^{\eta}$ dilation). The
coprimality charge of $j$ is $\le G+\eta$ after the threshold.
Apply the initialization above with $r$ replaced by $D'$,
adding $G+\eta$ to $Q'$ and keeping $z_0=z-G$; these corrections are
charged to the $128\eta+4\tau_{\rm init}$ margin, and both margins remain
positive.

Subsequent children use the marked closure of
Lemma~\ref{v:lem:inverse-returned-family}, with the same projection and
\(\kappa_c+P_1+2F_c=m\). Bounded stopping uses the actual cap
including all initial allowances
(Theorem~\ref{v:thm:unmarked-inverse}), never a negative
formal center alone. The remaining
branches use the normalized estimate with prime-divisor weights, and the
principal terms and tails have unchanged
bounds. Choosing all losses inside the initial margins and
\(\eps\) and summing the logarithmic dyads proves the estimate.
\end{proof}

\section{The unweighted product moment}
\label{v:sec:M2-assembly}
In this section, we give the derivation of the cyclotomic counterpart of
Proposition~\ref{prop:iv-plain} for $13/18\le\kappa\le1$, the unweighted
product moment of Theorem~\ref{v:thm:M2-assembled} with excess
$Z^{M+\eps}$, $M=m+q$. That theorem retains full prime annuli, a fixed
finite character span and uniformity in the number of prime sums, as required by the
count in Section~\ref{sec:v-detail-29} and the mesh condition of
Proposition~\ref{prop:iv-plain} (Lemma~\ref{v:lem:arithmetic-independent-mesh}); its zero-free hypothesis is
supplied by the preceding family bound, as in
Section~\ref{sec:iv-iteration}. The proof uses
centered product sums (Section~\ref{sec:v-detail-20}) and a two-transform
recursion (Sections~\ref{sec:v-detail-22}--\ref{sec:v-detail-23}), closed by
the finite-width induction of Section~\ref{sec:v-detail-24}; the theorem is
used in Section~\ref{v:sec:compensated-assembly}.

The positive row lift, unrestricted orbit unfolding and weighted
Fourier bounds are those proved above; the Hecke functional equation
and prime ideal theorem are used for a fixed field and finite
coefficient group, without discriminant uniformity.
Radical and indicator conventions are those of
Section~\ref{sec:arithmetic}, extended in Section~\ref{sec:v-scope}. The plain polynomial is
\[
 S_k(X;W)=X^{-1/2}\sum_n\psi_k(n)W(q_n/X).
\]

Each prime polynomial has normalization \(P_i^{-1/2}\); its prime sum
runs over every exterior prime in the prescribed annulus with a fixed
smooth norm weight, and its coefficients are fixed finite linear
combinations of members of the fixed finite character group
\(\Theta\) of Section~\ref{sec:v-detail-19}. A fixed ray-class
selector is allowed only after its finite character expansion is
included in \(\Theta\); arbitrary subsets or arbitrary bounded prime
coefficients are not allowed in M2. Subcollections mean whole surviving
prime-sum factors. Constants are uniform over bounded scalar coefficients in this
basis. Common extra coprimality conditions may have bounded polynomial norm;
their subpower inclusion--exclusion is explicit below, and they are
distinct from displayed varying character factors, whose exponents are retained.

\begin{lemma}[The unrestricted positive lift]\label{v:m2:row-lift}
There is a nonnegative product Schwartz--Bruhat test \(\Phi\) with
\(\Phi(0)=0\), compact complex additive Fourier support, and
\(\Psi_\Phi(T)\ge1\) for every \(0<T\le C_1\).
\end{lemma}
\begin{proof}
At each complex place take the squared modulus of the inverse Fourier
transform of a smooth compact bump supported sufficiently near zero.
It is bounded below on the required fixed ball. At each finite place
use the indicator of a compact neighborhood of one. A fixed positive
measure neighborhood of one in the norm-one group supplies a uniform
lower bound for the integral over every \(0<T\le C_1\); multiply by
one fixed constant. The finite factors avoid the additive origin. The
Fourier and norm-sheet bounds are those of
Section~\ref{v:sec:unmarked-assembly}.
\end{proof}

\subsection{The signed difference of convolution sums}
\label{sec:v-detail-16}
In this subsection, we write the unweighted products as equal-product
convolutions with the signed coefficient $D(l_1,l_2)$ before
any positive bound. Lemma~\ref{v:m2profile:lem:whole} gives
one Fourier measure for the complete product coefficient, inserting a single
power $q_U^{iv}$ of the whole product, and the
signed rectangle is retained until its principal terms cancel
(Section~\ref{sec:v-detail-20}). The Fourier and norm-sheet bounds used are
those of Section~\ref{v:sec:unmarked-assembly}.

Fix two compact annular smooth tests $W_1,W_2$ and positive scales
\[
 T_1T_2=Y_1Y_2=P.
\]
Let $\psi$ be one completely multiplicative zero-extended common
presentation, with one common extra radical $R_*$. Let $I$ index the remaining prime-polynomial factors, with smooth
weights $V_i$ and $\Pi=\prod_{i\in I}p_i$. Retain the signed coefficient
\[
 D(l_1,l_2)=W_1(q_{l_1}/T_1)W_2(q_{l_2}/T_2)
          -W_1(q_{l_1}/Y_1)W_2(q_{l_2}/Y_2).
\]
The whole product is $U=\Pi l_1l_2$. Its convolution coefficient sums
\[
 \left(\prod_{i\in I}c_i(p_i)V_i(q_{p_i}/P_i)\right)
   D(l_1,l_2)\psi(U)q_U^{it}
                1_{(U,R_*)=1}1_{s\mid U}
\]
over all factorizations of $U$; the condition $s\mid U$ is retained
in this signed polynomial.

\begin{lemma}[Whole-column common measure]\label{v:m2profile:lem:whole}
Suppose a Poisson formula multiplies two such whole-column convolutions by
\[
 K(Aq_h/(q_Uq_V))(q_Uq_V)^{-1/2}
\]
and fixed scalar root ratios. On fixed full norm boxes it has an exact
common Fourier density in the whole norms $q_h,q_U,q_V$ and any coupled
fixed aggregate norms. Every residual side therefore gains one unitary
power $q_U^{iv}$ of its \emph{whole} product. The same $v$ occurs on
its two plain factors, all remaining prime-sum factors and both rectangles. No independent
plain-factor twists are introduced.
\end{lemma}
\begin{proof}
Normalize the whole-column roots at their dyadic lengths and include only
their smooth ratios in the joint kernel profile. Fourier-separate that
profile on its independent full box before fixing remaining indices. The
norm-sheet derivative and weighted Fourier estimates proved in
Section~\ref{v:sec:unmarked-assembly} apply and insert powers
of the whole-coordinate norms only. Completely multiplicative characters,
extra coprimality indicators and the signed coefficient $D$ remain in the exact convolution.
Since $q_U=q_\Pi q_{l_1}q_{l_2}$, its one norm power factors with the
same exponent on both plains. Do not Fourier-separate $D$ in its two
individual plain norms at this stage. Fixed scalar or row phases remain
outside the column coefficient. The common measure is defined before
rows, labels and prime tuples, with polynomial incoming height bounds.
\end{proof}

Independent twists on the two plains would fail: for $t_1\ne t_2$ the
leading volume phases $T_1^{1+it_1}T_2^{1+it_2}$ and
$Y_1^{1+it_1}Y_2^{1+it_2}$ need not agree although $T_1T_2=Y_1Y_2$.
Lemma~\ref{v:m2profile:lem:whole} restricts the kernel's
Fourier coordinates instead and keeps the signed convolution exact.

\subsection{Exact extractions and last-divisor allocation}
\label{sec:v-extractions}
In this subsection, we assign the shared-prime divisors of the convolution of
Section~\ref{sec:v-detail-16} to a fixed last factor, following the
coefficient-independence rule of Section~\ref{sec:v-detail-13} and
preserving $T_1T_2=Y_1Y_2$ and one common coprimality indicator. These exact
identities are used in Lemma~\ref{v:m2:lem:common-coefficient} and in
Section~\ref{sec:v-detail-22} before positive enlargement.

At a genuine common-support prime, extract fixed powers $p^{a_i}$ from
the plains and fix any prime variable supplying that prime. The residual factors
are all punctured at that prime. Divide $T_i,Y_i$ by the same $q_p^{a_i}$; this preserves
$T_1T_2=Y_1Y_2$ and one common coprimality indicator. At the amplifier
prime, outside every remaining prime annulus, the valuations one, six and seven
are allocated only to plains; the same operation applies to both
rectangles. Its character/norm factors multiply the whole residual
product or an outer row scalar, not one rectangle alone.

For a final artificial gcd divisor use the exact identity
\[
 1_{p\mid\Pi l_1l_2}
  =\sum_{\varnothing\ne I\subseteq\{\text{factor positions}\}}
       (-1)^{|I|+1}\prod_{i\in I}1_{p\mid n_i}.
\]
A selected plain is written $l_i=pa_i$ and its two rectangle scales are
both divided by $q_p$; the quotient $a_i$ is unrestricted at $p$.
A selected prime-sum factor is fixed; conflicting
assignments on one prime-sum factor are zero. The extracted character scalar, including
its zero when $\psi(p)=0$, stays in the signed identity, so multiplicativity
preserves the whole-product presentation and common coprimality indicator.
All this precedes any positive factor bound or row enlargement.

\subsection{Fixed-character means with coprimality conditions}
\label{sec:v-detail-18}
In this subsection, Lemma~\ref{v:m2profile:lem:mean} gives the Mellin
asymptotic $c_{\vartheta,R_*}T^{1+it}I_W(t)+O(Z^\eps(1+|t|)^Jp_J(W))$ for
fixed-character means with polynomial-size common coprimality
indicators, and Proposition~\ref{v:m2profile:prop:center} cancels the
matching principal terms when $T_1T_2=Y_1Y_2$, saving $Z^{-r}$. They control
the exceptional rows in
Section~\ref{sec:v-return-6} and give the add-back
of Lemma~\ref{v:lem:M2-fixed-addback}.

Let $\Theta$ be a fixed finite group of finite-order Hecke characters
with conductor in the final fixed $S$. For a fixed exceptional row,
its displayed common presentation has inducing character
$\vartheta\in\Theta$. On exterior ideals it is exactly
\[
 \psi(n)=\vartheta(n)1_{(n,R_*)=1},
\]
where $R_*$, of polynomially bounded norm, includes all natural
displayed supports and common coprimality conditions, also redundant
supports after character cancellation. This description is
reconstructed for each enlarged row from its displayed character
factors and natural zeros; restrictions on the original row are not
imposed on an enlarged row without this verification.

\begin{lemma}[Uniform mean with coprimality conditions]\label{v:m2profile:lem:mean}
For $q_{R_*}\le Z^C$, $T\le Z^C$, annular $W$, and $\vartheta\in\Theta$,
\[
 \sum_l\vartheta(l)1_{(l,R_*)=1}q_l^{it}W(q_l/T)
 =c_{\vartheta,R_*}T^{1+it}I_W(t)
       +O\bigl(Z^\eps(1+|t|)^Jp_J(W)\bigr),
\]
uniformly for $T\le Z^C$, with any requested $Z^\eps$ loss. Here
\[
 I_W(t)=\int_0^\infty W(y)y^{it}\,dy,
 \quad c_{\vartheta,R_*}=\operatorname*{Res}_{s=1}L_F^S(s,\vartheta)
             \prod_{p\mid R_*}(1-\vartheta(p)/q_p).
\]
The residue is zero for nonprincipal $\vartheta$.
\end{lemma}
\begin{proof}
Fix a common compact support interval $[a_W,b_W]\Subset(0,\infty)$
and write $\widehat W(s)=\int_0^\infty W(y)y^{s-1}\,dy$.
Only the radical of $R_*$ occurs. Put
\[
 E_{R_*}(s,\vartheta)=\prod_{p\mid R_*}
       (1-\vartheta(p)q_p^{-s}),\qquad
 D_{R_*}(s,\vartheta)=L_F^S(s,\vartheta)E_{R_*}(s,\vartheta).
\]
For any $a,\delta>0$, uniformly for $\Re s\ge a$,
\begin{equation}\label{eq:v-removed-euler-uniform}
 |E_{R_*}(s,\vartheta)|
 \le\prod_{p\mid R_*}(1+q_p^{-a})
 \ll_{F,a,\delta}q_{R_*}^{\delta}.
\end{equation}
Indeed $(1+x^{-a})\le x^\delta$ for all sufficiently large $x$;
the remaining prime ideals form a fixed finite set. This argument
is uniform over $\vartheta\in\Theta$, including zeros at ramified primes.

For $\Re s=2$, absolute convergence gives the exact identity
\[
 \sum_l\vartheta(l)1_{(l,R_*)=1}q_l^{it}W(q_l/T)
 =\frac1{2\pi i}\int_{(2)}
       D_{R_*}(s-it,\vartheta)\widehat W(s)T^s\,ds.
\]
Choose
\[
 0<a\le\min\{1/4,\eps/[8(C+1)]\},\qquad
 0<\delta\le\eps/[8(C+1)].
\]
The meromorphic continuation and functional equation of the fixed
finite family $L_F^S(s,\vartheta)$ give, on $a\le\Re s\le2$
away from its possible pole at $1$,
\[
 |L_F^S(s,\vartheta)|\ll(1+|\Im s|)^A
\]
for a fixed $A$. Here only $\vartheta=1$ has a pole, and it is simple.
Integration by parts in $\log y$ gives, uniformly on this strip,
\[
 |\widehat W(\sigma+i\tau)|
 \ll_{J,a_W,b_W}p_J(W)(1+|\tau|)^{-J}.
\]
Shift the contour to $\Re s=a$. For fixed $t$, the horizontal
integrals tend to zero as their heights tend to infinity, by these
bounds and \eqref{eq:v-removed-euler-uniform}. The only possible
residue, at $s=1+it$, is
\[
 c_{\vartheta,R_*}\widehat W(1+it)T^{1+it}
   =c_{\vartheta,R_*}I_W(t)T^{1+it}.
\]
On the new vertical line, taking $J>A+2$ gives the error bound
\[
 \begin{split}
 |\mathcal R|&\ll T^a q_{R_*}^\delta p_J(W)
   \int_{\mathbb R}(1+|\tau-t|)^A(1+|\tau|)^{-J}\,d\tau\\
 &\ll T^a q_{R_*}^\delta(1+|t|)^A p_J(W)
 \ll Z^{\eps/4}(1+|t|)^A p_J(W).
 \end{split}
\]
The last inequality also holds for $0<T<1$, since $T^a\le1$.
Hence the estimate is uniform before $R_*$ and $T$ are selected,
also when $R_*$ depends on the exceptional row; no information
on zeros of $L_F^S(s,\vartheta)$ is used.
\end{proof}

\begin{proposition}[Exact centered cancellation]\label{v:m2profile:prop:center}
Assume both plains use the same $\vartheta,R_*,t$ and retain the same
tests in the two rectangles. Their leading product in
Lemma~\ref{v:m2profile:lem:mean} cancels exactly when $T_1T_2=Y_1Y_2$.
If all four scales are at least $Z^r$, $r\ge0$, their centered plain sum
is bounded by
\[
 O\bigl(Z^\eps (1+|t|)^Jp_J(W_1,W_2)P Z^{-r}\bigr).
\]
If no positive lower scale survives, absolute raw rectangle volume gives
$O(Z^\eps P)$ on every nonempty rectangle; no centered moment with altered norm ranges is
inserted into this argument.
\end{proposition}
\begin{proof}
Write $P=T_1T_2=Y_1Y_2$, and set
\[
 A_i(T)=\sum_l\vartheta(l)1_{(l,R_*)=1}q_l^{it}W_i(q_l/T),
 \qquad a_i=c_{\vartheta,R_*}I_{W_i}(t).
\]
Apply Lemma~\ref{v:m2profile:lem:mean} with loss $\eps/4$:
\[
 A_i(T)=a_iT^{1+it}+E_i(T),\qquad
 |E_i(T)|\le\mathcal H Z^{\eps/4},
\]
where $\mathcal H$ is a fixed polynomial in $1+|t|$ and the
required seminorms of $W_1,W_2$. The residue coefficients are
uniformly bounded in $R_*$: they vanish for nonprincipal
$\vartheta$, and for $\vartheta=1$ their extra factors are
$\prod_{p\mid R_*}(1-q_p^{-1})\le1$. Hence $|a_i|\ll\mathcal H$, and
\[
 \begin{split}
 A_1(T_1)A_2(T_2)-A_1(Y_1)A_2(Y_2)
 ={}&a_1\{T_1^{1+it}E_2(T_2)-Y_1^{1+it}E_2(Y_2)\}\\
 &+a_2\{T_2^{1+it}E_1(T_1)-Y_2^{1+it}E_1(Y_1)\}\\
 &+E_1(T_1)E_2(T_2)-E_1(Y_1)E_2(Y_2).
 \end{split}
\]
since the leading term $a_1a_2\{P^{1+it}-P^{1+it}\}$ vanishes.
Thus its absolute value is at most
\[
 C\mathcal H^2Z^{\eps/2}
           (T_1+T_2+Y_1+Y_2+1).
\]
When every scale is at least $Z^r$, each is at most $PZ^{-r}$,
and $PZ^{-r}\ge1$. This proves the asserted saving, after renaming
the finite seminorm and height orders.

For the remaining assertion use absolute counting directly. If a
rectangle with scales $U_1,U_2$ is nonempty, then
$U_i\ge b_{W_i}^{-1}$ and
\[
 \sum_l|W_i(q_l/U_i)|\ll_{F,W_i}U_i,
 \qquad
 |A_1(U_1)A_2(U_2)|\ll_{F,W_1,W_2}U_1U_2=P.
\]
An empty rectangle contributes zero. Apply this to the two original
rectangles separately, at their exact scales; no subunit scale is
replaced by one inside the signed difference.
\end{proof}

Remaining prime factors contribute a common outer scalar to the fixed-row
factorization, with common coprimality indicators and the same
whole-product norm exponent. Absolute numbers of terms in the prime sums
are restored only after the full Gauss coefficient disappears in the
second transform, since Proposition~\ref{v:m2profile:prop:center} concerns
ordinary fixed-character means, not Gauss-weighted means.

\subsection{The common product family and preliminary reductions}
\label{sec:v-detail-19}
In this subsection, we normalize the product moment in a common row range and radical
width $M=m+q$ and state the target bound $Z^{M+\eps}$. Lemmas~\ref{v:m2:lem:extra-mask} and
\ref{v:m2:lem:M2-reflection} treat the common coprimality modulus and ordinary
reflection, keeping natural zeros;
Lemma~\ref{v:m2:lem:prime-capacity} gives the pointwise prime
bound $|Q_k|^2\ll Z^{\kappa z+\eps}$ used in Section~\ref{sec:v-return-4} and at
terminal widths (Section~\ref{sec:v-detail-24}).

Let $\Theta$ be the specified fixed finite group of finite Hecke coefficient
characters supported at $S$. A common presentation is
\[
 \psi_k(n)=\nu_{\mathrm{fix}}(n)\prod_j\chi_{a_j}(n)^{e_j}\chi_n(k),
 \qquad \nu_{\mathrm{fix}}\in\Theta.
\]
The full radical of all displayed moving factors, before cancellation
or reduction of exponents modulo six, has norm at most $Z^q$.
Every zero extension remains, including redundant factors absent from
the primitive conductor. Separately, a common extra radical $R_*$
with $q_{R_*}\le Z^{C_{\rm extra}}$ is fixed within the current row sum
and occurs in both plains and every prime factor; $C_{\rm extra}$ is
a uniformity parameter. Complete expansion
by Lemma~\ref{v:m2:lem:extra-mask} precedes every nonexceptional
moment or ordinary reflection, so this extra modulus is not added
to the moving-support width $q$. A cancelled moving factor cannot
be reclassified as an extra modulus. In an exceptional centered term,
the indicator may instead stay in both fixed-character means.
Put $M=m+q$, and let $n_1,n_2$ be plain lengths and
$z_i$ the lengths of disjoint underlying prime windows, with coefficients
fixed finite combinations of members of $\Theta$. Put
$z=\sum z_i$, $A=n_1+n_2+z$.

The proof target is
\[
 \sum_{0<q_k\ll Z^m}|S_k(n_1)S_k(n_2)Q_k|^2\ll Z^{M+\eps},
\]
excluding principal inducing rows for $z=0$, and excluding $\Theta$
for $z>0$ with $A+(6\kappa-1)z\le M$.
For $13/18\le\kappa<1$, retain the hypothesis
$\beta_*(F)\le(1+\kappa)/2$. For every requested loss a sufficiently
small mesh $\eta$ must satisfy $z_i\le\eta$; it is chosen independently
of the fixed number of prime-sum factors. Zero-length prime-sum factors have bounded coefficient mass.

All M2 row Poisson steps use the unrestricted $\Gamma$-orbit sum
and $b=1$ in the row lift, so their norm kernels have no row-sector
projections and symmetric allowed-generator reciprocity has no
supplementary factor. The normalized kernels retain the finite
characters $\xi_\epsilon$; we require $\Upsilon_S\subset\Theta$. Moving
finite-order factors such as $\xi_r$ remain in the displayed moving
presentation, with radical counted in $q+R+E$; they are not absorbed into
$\nu_{\mathrm{fix}}$, $\Theta$ or the fixed set $S$.

\begin{lemma}[Removal of a common extra coprimality condition]\label{v:m2:lem:extra-mask}
Let $\psi_k^0$ retain all natural zeros of the primitive character and
its displayed row and moving factors. Let $R_*$ be the extra radical,
fixed with respect to the current row and of norm at most
$Z^{C_{\rm extra}}$, separately from the width $M=m+q$. For central plain sums and prime-sum factors,
\[
 S_{k,R_*}(X)=\sum_{d\mid R_*}\frac{\mu_F(d)\psi_k^0(d)}{\sqrt{q_d}}
                                     S_k^0(X/q_d),
\]
\[
 Q_{k,i,R_*}=Q_{k,i}^0-P_i^{-1/2}
               \sum_{p\mid R_*}\psi_k^0(p)c_i(p)W_i(q_p/P_i).
\]
Simultaneous expansion in the two plains and all prime-sum factors has subpower
absolute coefficient mass. It only decreases the affine expression in the length condition. The inducing character and the exclusion condition are unchanged.
\end{lemma}
\begin{proof}
Insert complete inclusion--exclusion in the common extra coprimality indicator.
The quotient of a plain variable remains unrestricted at the extracted
divisor; its natural character zeros remain. The sum
$\sum_{d\mid R_*}q_d^{-1/2}\ll Z^\eps$ is the small-power Euler product
bound. A fixed prime-sum factor has the same bound from its primes dividing $R_*$.
Its row scalar is bounded, including zeros, before a positive row norm
is used. The new plain length exponents are replaced by zero only when the annular scale
is in a fixed bounded interval below one; they are at most their original
nonnegative lengths. Removing a prime-sum factor decreases $A$ and $z$. Thus the
affine expression does not increase. A cancellation with a varying row
factor is not converted into a newly chosen row-dependent extra coprimality modulus.
\end{proof}

\begin{lemma}[Natural reflection with a common scale supremum]
\label{v:m2:lem:M2-reflection}
For a nonprincipal primitive inducing character and the retained natural
zeros, reflection of a plain factor of length $n$ reduces it, in absolute
product estimates, to original-character annular factors of declared
length at most $M-n+\xi$. The coefficient mass is subpower, and discarded
upper annuli have arbitrarily small total size with finite internal
orders. This reduction uses a common scale-supremum measure over the
whole row Hilbert space.
\end{lemma}
\begin{proof}
For primitive conductor $f_k$ and redundant natural radical $R_{0,k}$,
$C_kq_{R_{0,k}}\ll Z^M$, where
$C_k=|\Delta_F|N f_k/(2\pi)^{2r_F}$, with $r_F=[F:\mathbb Q]/2$. The exact zero-preserving identity is
\[
 S_k^0(X;W)=\epsilon_k
 \sum_{d\mid R_{0,k},\,\rad h\mid R_{0,k}}
 \frac{\mu_F(d)\psi_k^*(d)\bar\psi_k^*(h)}{\sqrt{q_dq_h}}
 S_{\bar\psi_k^0}^0(C_kq_d/(Xq_h);W^\sharp),
\]
\[
 \widehat W^\sharp(s)=\widehat W(1-s)
                           \{\Gamma(s)/\Gamma(1-s)\}^{r_F}.
\]
Before applying this identity, absorb every incoming explicit norm
mode into its test profile:
$q_n^{it}W(q_n/X)=X^{it}(q_n/X)^{it}W(q_n/X)$.
Do this in both ordinary and all prime factors, retaining $t$ as an
incoming height parameter but resetting the explicit common mode to
zero. Finite seminorms grow polynomially in $1+|t|$.
After reflection, conjugating the reflected factor inside its own
absolute value restores the original arithmetic character; subsequent
whole-product separation inserts one fresh common mode.

It follows from the primitive Hecke functional equation; deleted natural
Euler factors are restored geometrically with their original coprimality indicators.
The coefficient mass is bounded by an arbitrarily small power of the
polynomial radical. The transformed profile has small-norm growth
$O((1+|\log x|)^{r_F-1})$ and arbitrary large-norm decay, with finite
Euler and height orders. Its lower annuli have central masses
$O(2^{-j/2}(1+j)^{r_F-1})$, which are summable. Upper annuli beyond
$Z^{\xi/2}$ are discarded by a fixed sufficiently large decay order.
The same applies when the main scale is below $Z^{-\xi}$.

Every retained main scale is at most $C Z^{M-n}$. A nonempty subunit
annulus is bounded below by its fixed support endpoint and is dilated
to scale one with uniformly bounded Euler seminorms. Box multipliers
and this endpoint cost at most the remaining $\xi/2$ in declared length
after the eventual threshold.

For a scale chosen separately by each row,
use a supremum over the common polynomial-range scale variable.
Sobolev on unit logarithmic boxes bounds its squared Hilbert norm by
one integral of the polynomial and its Euler derivative. There are
$O(\log Z)$ such boxes; the two-scale version has $O((\log Z)^2)$.
For two factors apply the two-dimensional Sobolev inequality to their
product in the row Hilbert space, retaining each first derivative and
the mixed derivative. Do not multiply two separately averaged suprema.
The measure is common to all rows, not chosen for each row's conductor
or Euler divisor. Conjugating the entire reflected factor inside its
absolute value returns its character to $\psi_k^0$ and conjugates its
profile. The signed centered rectangle is not reflected by this operation;
only an uncentered product used to bound it or its comparison is reflected.
\end{proof}

\begin{lemma}[Pointwise admissible prime length]\label{v:m2:lem:prime-capacity}
On inducing rows outside $\Theta$, the product of the original
prime-sum factors obeys $|Q_k|^2\ll Z^{\kappa z+\eps}$ with finite seminorm and
polynomial norm-twist height orders, uniformly for
$13/18\le\kappa\le1$ and the stated zero-free hypothesis when $\kappa<1$.
\end{lemma}
\begin{proof}
Each coefficient component twists the row by a member of $\Theta$.
It is nonprincipal because the row induces outside that group.
Put $\sigma_0=(1+\kappa)/2+e$ and
$Q_v=2+q_{\mathfrak f_\psi}(3+|v|)^{[F:\mathbb Q]}$.
For $\kappa<1$, choose concentric disks centered at $2+iv$, with
radii $1/2<r_1<r_2<2-\beta_*(F)$, so that the smaller disk
contains a fixed neighborhood of $\sigma_0+iv$. The Euler branch
of $\log L$ extends to the larger disk. Convexity and
Borel--Caratheodory give respectively
\[
 \Re\log L\ll\log Q_v,\qquad
 \max_{|s-(2+iv)|\le r_2}|\log L(s)|\ll_e\log Q_v,
\]
after initially taking a slightly larger disk still in the zero-free
half-plane. On the radius-$1/2$ disk the Euler logarithm is bounded.
Hadamard three-circles and then Cauchy's estimate give, for some
$\theta<1$ depending only on the fixed buffer,
\[
 |\log L(\sigma_0+iv)|+
       |L'/L(\sigma_0+iv)|
       \ll_e(\log Q_v)^\theta\ll_{e,\epsilon_1}Q_v^{\epsilon_1}.
\]
Removed Euler factors satisfy the same subpower bound, since
$\Re s$ is bounded below positively. Thus the estimate is uniform
in the moving presentation and its natural zeros.

Mellin inversion of its von Mangoldt series gives
prime-power size $P^{(1+\kappa)/2+e}Z^\eps$.
Remove the higher prime powers by their $P^{1/2+\eps}$ ideal count,
and absorb $1/\log q_p$ into the annular profile. The normalized square
costs $P^{\kappa+2e}Z^\eps$. The natural deleted primes have divisor-many
terms; their omission changes the normalized sum by at most
$P^{-1/2}Z^\eps$. For $\kappa=1$ use absolute prime counting.
Multiply over the fixed list, choosing subsidiary losses within the
requested aggregate. Translate each pure norm twist in its Mellin
argument before integration by parts; the numerical norm exponent does
not increase with derivative or height orders.
\end{proof}

\subsection{Centering before the two M2 transforms}
\label{sec:v-detail-20}
In this subsection, we center the product family before either Poisson
transform. Reflection by Lemma~\ref{v:m2:lem:M2-reflection} reduces to the
range $A\le5M/6$ or to four plain scales at least $Z^{M/4}$; there
Proposition~\ref{v:m2profile:prop:center} removes the sixth-power
principal contribution before positive enlargement, and
Lemma~\ref{v:m2:lem:common-coefficient} shows that the two ordinary
factors and all surviving prime factors keep one character and norm
twist through the transformations of Section~\ref{sec:v-detail-22}.

For a family without prime-sum factors reflect once each plain of length $>M/2$;
Lemma~\ref{v:m2:lem:M2-reflection} reduces it to the range with an extra length margin
$n_i\le M/2+\xi$, $A\le M+\delta$, with $2\xi\le\delta$.
For sums with prime-sum factors keep $A+(6\kappa-1)z\le M$.
At widths above a fixed floor, prove $A\le5M/6$ first. Put $L=M/4$.
When $A>5M/6$, $6\kappa-1\ge10/3$ implies $z<M/20$.
Reflecting the longer factor of a comparison with shorter length at
most $L$ gives
\[
 A_{\mathrm{comp}}\le3M/2-A+2z+\xi\le23M/30+\xi,
\]
\[
 A_{\mathrm{comp}}+(6\kappa-1)z\le14M/15+\xi.
\]
Indeed, with $q_\kappa=6\kappa-1\ge10/3$, the length condition
implies $z\le(M-A)/q_\kappa<M/20$, and
\[
\begin{aligned}
 A_{\rm comp}+q_\kappa z
 &\le3M/2-A+(q_\kappa+2)(M-A)/q_\kappa+\xi\\
 &\le(5/6+1/(3q_\kappa))M+\xi
 \le14M/15+\xi.
\end{aligned}
\]
It belongs to the lower-rank uncentered range at the same width
with fixed positive margins. The same argument bounds an original
factor shorter than $L$ by reflecting its partner.

In the remaining case all four plain scales are at least $Z^L$.
Set $Y_1=Z^L$, $Y_2=X_1X_2/Y_1$ and keep the whole difference
\[
 D_{\mathbf b}(l_1,l_2)=
 \prod_i W_i(q_{b_i}q_{l_i}/X_i)
       -\prod_iW_i(q_{b_i}q_{l_i}/Y_i),\qquad X_1X_2=Y_1Y_2.
\]
An impossible rectangle stays a formal zero rectangle. The full product
coefficient has the following form, where $I$ indexes the remaining
prime-polynomial factors and $U=\Pi l_1l_2$:
\[
 \prod_{i\in I}c_i(p_i)W_i(q_{p_i}/P_i)
 D_{\mathbf b}(l_1,l_2)\tau_1(U)q_U^{it}
                         1_{(U,R_*)=1}1_{s\mid U}.
\]
Here $\tau_1$ is one zero-extended common presentation; the same
character, coprimality modulus and norm power apply to both plains and both rectangles.
Taking the squared norm of the entire difference permits its rows to
be enlarged to the unrestricted positive row ball before Poisson.

\begin{lemma}[The common coefficient invariant]\label{v:m2:lem:common-coefficient}
In the two Poisson transforms, the single sixth-power Gauss amplification
and the final complete M\"obius allocation, each separated transformed sum side
preserves the displayed whole-product coefficient form. Before the second
transform it is multiplied by the full normalized Gauss sum of $U$;
afterward it has no Gauss coefficient. Every surviving prime-sum factor has the
transformed sum's whole-product character and its original coefficient in $\Theta$.
Within each side its two plains and both rectangles have the same
character, coprimality indicator and norm power. The first transform
has the ordered pair $(\tau_C,\tau_D)$ of
\eqref{v:m2:eq:M2-first-bridge}; these need not be equal.
After Cauchy--Schwarz the two resulting positive norms are treated
separately. Within either one norm, its second transform produces
two child polynomials with the same inducing character, and hence
the same membership in $\Theta$.
\end{lemma}
\begin{proof}
Residue-class Poisson is on the full products $U,V$, not on the
individual plain factors. Its norm kernel and inverse roots depend
only on the row norm and the two whole-product norms. Their common
logarithmic Fourier density therefore inserts one power $q_U^{it}$
on each whole column. All row sums are unrestricted, so no sector projection is inserted.

The finite factors $\xi_\epsilon$ from
\eqref{eq:v-poisson-unit-normalization} are nevertheless retained in
$\Theta$. Symmetric allowed-generator reciprocity accounts for the
remaining exterior CRT factors. For one expanded norm, the Gauss product in
\eqref{eq:v-unit-Gauss-correlation} changes its single presentation
$\tau_1$ to $\tau_1\overline{\xi_\epsilon}$ on both columns.
The factor from the common primitive modulus is a row scalar.
Thus this unit operation preserves a common presentation within
that norm, including all imprimitive zeros.

At a genuine common-support prime remove its full powers from the
two products, allocating them to the plain variables or the one
possible prime-sum factor. Require every remaining factor on that side to be coprime to this
prime; a prime-sum factor supplying it is fixed. This is one allocation for the
whole difference and therefore divides the same $X_i,Y_i$ by the
same exact ideal powers. All new moving-character phases act on the
whole residual product, and common coprimality indicators remain common.

At an amplifier prime $p$ outside every remaining prime annulus,
write $P=q_p$ and assume $p\nmid h$. The local Gauss formula of
Section~\ref{sec:v-return-2} gives
\[
 \begin{array}{c|ccc}
 i&1&6&7\\\hline
 |P^{-i/2}\mathcal G(p^i,hp^6)|^2
       &0&(1-P^{-1})^2&P^{-1}
 \end{array}
\]
whereas the original $h$-frequency has only the $i=1$ term,
of squared size $P^{-1}$. Every other positive valuation vanishes;
the valuation-zero term is unchanged by multiplication of the
frequency by $p^6$ and cancels in the difference.
Thus the difference is supported at valuations $1,6,7$,
which allocate only to the two plains.
The residue-character phase acts on their whole residual product;
the frequency and Gauss phases are row scalars. The whole centered
coefficient is transformed once, so the same extraction divides both
rectangle scales and retains the same coprimality indicators.

More explicitly, the intermediate presentations satisfy
\[
 \tau_C=\tau\chi_\bullet(er)\xi_r,\qquad
 \tau_D=\tau\chi_\bullet(er)\bar\xi_r.
\]
Their ratio is $\xi_r^2$ on their common unit support.
Cauchy--Schwarz produces $\mathcal N_C$ and $\mathcal N_D$;
write $\tau_1=\tau_C$ or $\tau_D$ in the norm currently estimated.
Its second correlation has the residual coefficient
\[
 \tau_1(a)\chi_a(j)\,
       \overline{\tau_1(b)\chi_b(j)},
\]
since $-1$ is a sixth power and $\chi_b(-j)=\chi_b(j)$.
Writing $j=\epsilon j_0$ gives the same further twist
$\xi_\epsilon$ on both factors. The exceptional split is made
separately in $\mathcal N_C$ and $\mathcal N_D$.

After the second finite correlation, insert the complete M\"obius sum
over residual common divisors before factorwise estimates. For each
such prime use the exact union identity for divisibility of
$\Pi l_1l_2$. A selected plain is divided by that prime in both
rectangles; its quotient is unrestricted and may still contain the
prime. No one-variable coprimality modulus is introduced. A selected prime-sum factor is
fixed, and conflicting choices on one prime-sum factor are identically zero.
All extracted character factors are row scalars, including their zeros.
Multiplicativity of the original zero-extended character and norm power
then proves the same whole-product form on each transformed sum.

The correlation
leaves $\chi_n(j)$ and $\bar\chi_n(-j)$ on the two sides; after
putting the conjugated second polynomial inside its absolute value,
their common original orientation agrees, since $\chi_n(-1)=1$.
Thus both induce the same character in this arithmetic convention.
The profiles can have different whole-product powers on the two sides,
but within one side they are identical in its two plains and rectangles.
\end{proof}

\subsection{Complete Gauss sums and correlation identities}
\label{sec:v-detail-21}
In this subsection, Lemmas~\ref{v:m2:lem:full-Gauss-new} and
\ref{v:m2:lem:full-correlation-new} compute the finite Fourier
coefficients of the centered product sums, retaining prime powers, zero
frequencies and the complete nonunit support; inserted into the orbit
Poisson formula of Section~\ref{sec:v-detail-5}, they give the
two-transform identity of Section~\ref{sec:v-detail-22}.

Imprimitive factors are retained in the full finite Fourier sums.
Choose consistent allowed generators modulo $\Gamma$; replacing a
product generator by its representative is a sixth-power-unit
change of residue variable and frequency representative, and leaves the
norm kernels and zero-extended symbols unchanged.

\begin{lemma}[Full Gauss sums, including imprimitive powers]
\label{v:m2:lem:full-Gauss-new}
For an exterior ideal $u$ put
\[
 \mathcal G(u,h)=q_u^{-1/2}\sum_{x\bmod u}\chi_u(x)e(hx/u),
 \qquad \mathcal G(1,h)=1.
\]
For $P=q_p$, $i\ge1$ and $v_p(0)=\infty$,
\[
 |\mathcal G(p^i,h)|=P^{(i-1)/2}1_{v_p(h)=i-1}
                              \quad(6\nmid i),
\]
\[
 \mathcal G(p^i,h)=P^{i/2}1_{p^i\mid h}
                   -P^{i/2-1}1_{p^{i-1}\mid h}
                              \quad(6\mid i).
\]
In the first case its phase is
$\gamma_i(p)\bar\chi_p(h/p^{i-1})^i$ times its displayed magnitude.
For $(a,b)=1$, CRT and symmetric reciprocity give
\[
 \mathcal G(ab,h)=\chi_a(b)\chi_b(a)
                  \mathcal G(a,h)\mathcal G(b,h).
\]
All equations include zeros on nonunits.
\end{lemma}
\begin{proof}
Write $x=x_0+py$ modulo $p^i$. The $y$ sum vanishes unless
$p^{i-1}\mid h$, and then equals $P^{i-1}$.
For $h=p^{i-1}h_0$ the remaining expression is
$P^{i/2-1}\sum_{x_0\ne0}\chi_p(x_0)^i e(h_0x_0/p)$.
For nonprincipal power $i\bmod6$ use its primitive field Gauss sum;
for a present principal power the field sum is
$P\,1_{p\mid h_0}-1$. This proves both formulas, including the
phase and zero cases. CRT at coprime $a,b$ uses residue substitutions
by $b$ and $a$ respectively, yielding the two displayed cross factors.
They are evaluated only on units; zero factors in $h$ remain.
\end{proof}

\begin{lemma}[Full correlation, zero frequency and complete support]
\label{v:m2:lem:full-correlation-new}
Define the actual finite sum
\[
 \mathcal F(u,v;j)=
 \sum_{\substack{x\bmod u, y\bmod v\\vx-uy\equiv j\ (\bmod uv)}}
                       \chi_u(x)\bar\chi_v(y).
\]
Then
\[
 \mathcal F(u,v;j)=\frac1{\sqrt{q_uq_v}}
 \sum_{h\bmod uv}\mathcal G(u,h)\overline{\mathcal G(v,h)}
                                                   e(-jh/(uv)).
\]
At $j=0$ it vanishes unless the modulus ideals agree, in which case
it is $\varphi(u)$. With $C=(u,v)$, $u=Cn_1$, $v=Cn_2$,
$(n_1,n_2)=1$, it vanishes unless $j=Ck$, and then equals
\[
 \chi_{n_1}(k)\bar\chi_{n_2}(-k)L_C(n_1,n_2;k).
\]
For $p^c\parallel C$, $P=q_p$, when both residual columns are units,
\[
 L_{p^c}=P^{c-1}\chi_p(n_1/n_2)^c
 \begin{cases}
 P-1,&p\mid k,\\
 -1,&p\nmid k, 6\nmid c,\\
 P-2,&p\nmid k, 6\mid c.
 \end{cases}
\]
When just one residual column is nonunit it is
$P^{c-1}(P-1)1_{6\mid c}1_{p\nmid k}$.
If $u=Da$, $v=Eb$, $(a,b)=1$, $(ab,DE)=1$, the complete-support
form is
\[
 \mathcal F(Da,Eb;j)=\mathcal F(D,E;j)
                                     \chi_a(j)\bar\chi_b(-j).
\]
No equality to this last expression is asserted on an artificially
added noncoprime pair $(a,b)>1$.
\end{lemma}
\begin{proof}
Expand the two Gauss sums in their finite Fourier transform. Additive
orthogonality imposes precisely $vx-uy\equiv j\bmod uv$.
At $j=0$ a nonzero summand has $x,y$ units. Reduction modulo $u$
and $v$ forces each ideal to divide the other; consistent allowed
generators give $u=v$, and the surviving congruence is $x=y\bmod u$.
There are $\varphi(u)$ solutions even at principal character powers,
because their unit indicators were retained.

Divide the general congruence by $C$. Reduction modulo $n_1,n_2$
supplies $\chi_{n_1}(k)\bar\chi_{n_2}(-k)$ and a reciprocity
factor which is one in the allowed convention. The remaining field
equation is $n_2x-n_1y=k$ at each common prime. It has $P^{c-1}$
lifts. If $p\mid k$ its units are proportional. If the residuals
and $k$ are units, substitute to reduce the field sum to
$\chi_p(n_1/n_2)^c\sum_{z\ne0,1}\chi_p(z)^c$.
This is minus the phase for a nonprincipal power and $P-2$ times
the phase for a principal power.

If exactly one residual column is
nonunit, the equation fixes the other variable; its unit character
vanishes at $p\mid k$, while the free unit sum is $P-1$ precisely
when $6\mid c$. These are all the stated cases. The lift bijection
also holds when $C$ meets one residual modulus: its other residue
is free and the first is uniquely determined modulo the larger power.

For complete support, solve the congruence at primes of $a,b$ and
change $x$ by $b$, $y$ by $a$ at primes of $DE$.
The remaining unit factor is
\[
 \frac{\chi_E(a)}{\chi_a(E)}
 \frac{\chi_b(D)}{\chi_D(b)}
 \frac{\chi_b(a)}{\chi_a(b)}=1
\]
by symmetric reciprocity. Only unit symbols are inverted here.
This proves the exact complete-support formula and preserves zeros
in the frequency $j$. On artificial noncoprime residual pairs its
right side is merely the stipulated extension; the full M\"obius sum
must annihilate those pairs before factorwise estimates.
\end{proof}

\subsection{The two Poisson transformations for M2}
\label{sec:v-detail-22}
In this subsection, we apply Section~\ref{sec:v-detail-21} to the whole product family of
Section~\ref{sec:v-detail-19}. The first transform
(Section~\ref{sec:v-return-1}) and one amplification
(Section~\ref{sec:v-return-2}) give the allowance
\eqref{v:m2:eq:M2-allowance}; the second (Section~\ref{sec:v-return-3}) gives
the smaller width $M'\le M-\sigma$ and normalization
\eqref{v:m2:eq:M2-child-width},
\eqref{v:m2:eq:M2-child-allowance}. Exceptional
terms are estimated in Section~\ref{sec:v-detail-23}.

For a whole-product modulus $u$, $\mathcal G(u,h)$ is its
zero-extended normalized sextic Gauss sum, prime powers included, and
$\mathcal F(u,v;j)$ is the complete finite correlation, the Fourier
transform of $\mathcal G(u,h)\overline{\mathcal G(v,h)}$; these are
the finite sums of Lemmas~\ref{v:m2:lem:full-Gauss-new}
and~\ref{v:m2:lem:full-correlation-new}, not primitive-character
replacements. The unrestricted M\"obius extension is stipulated only
off the genuine coprime locus. 

\subsubsection{The first transform and its norm allowance}
\label{sec:v-return-1}
In this subsubsection, orbit Poisson summation with coprimality conditions
(Section~\ref{sec:v-detail-5}) on the full product moduli and
Lemmas~\ref{v:m2:lem:full-Gauss-new},
\ref{v:m2:lem:full-correlation-new} give the exact nonzero-frequency
expression \eqref{v:m2:eq:M2-first-bridge}, and we reduce to the norm bound
\eqref{v:m2:eq:M2-Gauss-target}, with $B_c=(3c-5d-R)_+/6$; this norm allowance
is the input of Section~\ref{sec:v-return-2}.

Write $H=Z^m$, $X=Z^A$ and $M=m+q$. Choose the unrestricted
positive row lift of Lemma~\ref{v:m2:row-lift}, with $b=1$.
Integrating on $C=G/\Gamma$ before norm-sheet unfolding gives an
exact majorant for the original nonzero rows. In the expanded square
use additive residue-class Poisson on the full products $U,V$.
No primitive-only formula is substituted when a modulus is imprimitive.

The original dual origin contributes $O(Z^{m+\eps_0})$:
its complete character mean vanishes unless every prime of
$UV$ has total multiplicity $\ge2$, i.e.\ $UV$ is powerful, and there are $O(X^{1+\eps_0})$ powerful ideals of norm $O(X^2)$,
with divisor-bounded allocations to the fixed factors. At a prime
of exponents $i,j$ their number is at most
$C_N(i+1)^N(j+1)^N$ for fixed factor count $N$; the product is
$\ll(q_Uq_V)^\eps$ by the ideal divisor bound.
The squared
central normalization is $X^{-1}$ and the row mean is $O(H)$.
The zero-frequency test is its compact $C$-integral, not the value
at zero of an unfolded norm-sheet kernel.

On the nonzero terms write $U=Ca$, $V=Db$, extracting the complete
common support, so $(a,b)=1$, $(ab,CD)=1$. Record the exact power
allocation to both plains and to the possible prime-sum factor at each common
prime. A contributing prime-sum factor is fixed; the remaining factors on its side are required to be coprime to that prime. Put
\[
 c=\log_Zq_C,\quad d=\log_Zq_D,\quad
 p=\log_Zq_{\rad(CD)}.
\]
The prime support of $C,D$ is the same. Let $r$ be the product of
these primes whose valuation difference is nonzero modulo six,
$R=\log_Zq_r$, and let $\xi_r$ be that primitive common character.
Expand the complementary common coprimality indicator completely. For a selected
squarefree divisor $e$, put $E=\log_Zq_e$; then
$R\le p$, $E\le p-R$.

With $A_C,A_D$ the allocated whole convolution coefficients,
including the entire centered difference when present, define
\[
 \tau_C(a)=\tau(a)\chi_a(er)\xi_r(a),\qquad
 \tau_D(b)=\tau(b)\chi_b(er)\bar\xi_r(b).
\]
Up to a fixed scalar of modulus at most one, CRT gives the exact
nonzero-frequency expression
\begin{equation}\label{v:m2:eq:M2-first-bridge}
 \frac{H}{Xq_e\sqrt{q_r}}
 \sum_{h\ne0}\mathcal G_{\xi_r}(r,h)
 \sum_{(a,b)=1}
 \frac{A_C(a)\overline{A_D(b)}}{\sqrt{q_aq_b}}
 \tau_C(a)\overline{\tau_D(b)}
 \mathcal G(a,h)\overline{\mathcal G(b,-h)}
 \mathcal K_1\left(\frac{Hq_h}{q_eq_rq_aq_b}\right).
\end{equation}
To verify the identity, put $r_0=\operatorname{rad}(CD)/r$.
On every original row the exact character identity is
\[
 \chi_{Ca}(k)\overline{\chi_{Db}(k)}
 =\xi_r(k)\chi_a(k)\overline{\chi_b(k)}1_{(k,r_0)=1}.
\]
The character $\xi_r$ is primitive modulo $r$, whereas $\chi_a$ and
$\chi_b$ retain their full, possibly imprimitive moduli. Expand the
last indicator as $\sum_{e\mid r_0,\ e\mid k}\mu_F(e)$ before
applying Poisson. Since $(e,rab)=1$, the substitution $k=ek'$ gives
\[
 \xi_r(ek')\chi_a(ek')\overline{\chi_b(ek')}
 =\xi_r(e)\chi_a(e)\overline{\chi_b(e)}
                   \xi_r(k')\chi_a(k')\overline{\chi_b(k')}.
\]
The primitive character at $r$ is evaluated here, not the original
zero-extended common character at a prime removed by inclusion--exclusion.

For the periodic function on the right, the positive finite Fourier sum is
\[
 \begin{split}
 &\sum_{x\bmod rab}\xi_r(x)\chi_a(x)\overline{\chi_b(x)}
                                           e(hx/(rab))\\
 &=\xi_r(ab)\chi_a(br)\overline{\chi_b(ar)}
       \sqrt{q_rq_aq_b}\,
       \mathcal G_{\xi_r}(r,h)\mathcal G(a,h)
                                  \overline{\mathcal G(b,-h)}.
 \end{split}
\]
Indeed write $x=abx_r+brx_a+arx_b$; these variables run
independently through their residue rings. The additive character splits
and multiplicativity gives precisely the three cross factors, even
when $x$ or $h$ is a nonunit. Reciprocity cancels
$\chi_a(b)\overline{\chi_b(a)}$. The remaining phases in $a,b$ are
\[
 \chi_a(er)\xi_r(a),\qquad
 \overline{\chi_b(er)}\xi_r(b),
\]
which are exactly the phases in $\tau_C(a)\overline{\tau_D(b)}$.
The omitted scalar is the original fixed-factor phase multiplied by
$\mu_F(e)\xi_r(e)$ and has modulus at most one.

Finally orbit Poisson at scale $H/q_e$ and modulus $rab$ contributes
$H/(q_eq_rq_aq_b)$. Multiplication by the Gauss normalization
$\sqrt{q_rq_aq_b}$ and the original $X^{-1}$ gives
\[
 \frac{H}{Xq_e\sqrt{q_r}}\frac1{\sqrt{q_aq_b}},
\]
with dual norm argument $Hq_h/(q_eq_rq_aq_b)$.
This proves \eqref{v:m2:eq:M2-first-bridge} before any absolute value,
Fourier separation or positive enlargement.
The unrestricted row lift has no sector selector. Dual unit
normalization is nevertheless retained through
\eqref{eq:v-poisson-unit-normalization}; its characters belong to
$\Upsilon_S\subset\Theta$.
The new displayed moving support contributes the exponent
\[
 \widetilde q=q+R+E,\qquad
 K_0=2A-c-d+R+E-m.
\]
Canceled fixed factors still count in $\widetilde q$.

Fix common dyadic sectors of the aggregate outer norms. Retain whole
frequency partition weights meeting $q_h\le T_{\rm sec}Z^{\xi/2}$,
where $T_{\rm sec}$ is the supremum of the nominal scale in the sector.
Their union is contained in $q_h\le C_dT_{\rm sec}Z^{\xi/2}$.
Hence a retained dyad of length $K$ satisfies
\[
 K\le K_0+\delta_1,\qquad 0\le\delta_1<\xi.
\]
The omitted weights are estimated on the genuine expression restricted to coprime indices
before any artificial extension. Their kernel arguments exceed
$Z^{\xi/4}$ after the fixed support threshold. Arbitrary radial decay,
$|\mathcal G(a,h)|\le\sqrt{q_a}$ and absolute ideal counting give
$O(Z^{-B})$ for any requested $B$, with internal orders chosen last.
The bounded-degree logarithmic small-norm estimate in Section~\ref{v:sec:unmarked-assembly} costs only $Z^{\eps_0}$ on retained dyads. In particular,
its Euler derivatives introduce no $Z^{J\xi}$.

Define the expression also on artificial noncoprime residual pairs, using
the same full Gauss sums, whole-product norms, coprimality indicators, kernel and retained
weight. Insert
$\sum_{s\mid a,b}\mu(s)$ completely before estimating either side;
put $s_0=\log_Zq_s$. Common logarithmic Fourier separation is in the
row norm and the two full column norms, with an aggregate scale-ratio
coordinate if necessary. Its fixed boxes are chosen before remaining
labels. Only then apply Cauchy--Schwarz in $h$ and use
$|\mathcal G_{\xi_r}|\le1$.

It suffices to prove for the resulting $C$ norm
\begin{equation}\label{v:m2:eq:M2-Gauss-target}
 \mathcal N_C\ll
 Z^{A-c+\widetilde q+B_c-s_0+\eps_G},\qquad
 B_c=(3c-5d-R)_+/6,
\end{equation}
and the exchanged estimate with $B_d=(3d-5c-R)_+/6$.
The radical, its powers and allocations, and the M\"obius labels cost
$Z^{p+s_0+\eps_0}$. Powers on a fixed polynomial-size radical are
subpower by Rankin's Euler product; fixed-factor allocations are
$C_N^{\omega(v)}\ll_{N,\eps_0}q_v^{\eps_0}$.
The complete outer exponent is
\[
 \begin{aligned}
 &(m-A-R/2-E)+(p+s_0)\\
 &\qquad+\tfrac12(2A-c-d+2\widetilde q+B_c+B_d-2s_0)\\
 &\quad=M+\tfrac12\{B_c+B_d-(c+d-2p-R)\}.
 \end{aligned}
\]

The brace is nonpositive. At a common prime of multiplicities
$i\ge j\ge1$ put $r_p=1_{6\nmid i-j}$.
The positive-part contribution is at most $(3i-5j-r_p)_+/6$,
while the available exponent is $i+j-2-r_p$. If the former is positive,
six times their difference is $3i+11j-12-5r_p\ge0$:
for $r_p=1$, $i\ge j+1$; for $r_p=0$, $i,j\ge1$.
Summation proves the inequality. Each family count is normalized
once before its common Fourier measure is evaluated.

\subsubsection{One Gauss amplification and the three error types}
\label{sec:v-return-2}
In this subsubsection, we amplify the first transformed family once by the sixth-power
construction of Section~\ref{v:sec:unmarked-assembly}, using a prime pool
$Z^{\ell_*}/2<q_p\le Z^{\ell_*}$, $\ell_*=\sigma/3$, and estimate
origin, diagonal and exceptional contributions separately, retaining
radical exponents. This gives the allowance $\Lambda_c$ of
\eqref{v:m2:eq:M2-allowance}, the normalization of
Section~\ref{sec:v-return-3}.

Set $a_0=A-c-w$, where $w$ is a removed full-product length,
and let $w_o\le w$ be its new contribution to the radical exponent. Removing a
squared extraction coefficient $Z^{-w_o}$ leaves the allowance
\begin{equation}\label{v:m2:eq:M2-allowance}
 \Lambda_c=a_0+\widetilde q+w_o+w+B_c-s_0.
\end{equation}
Write
\[
 J=d-c+(K_0-K)-2w+w_o.
\]
For a norm without prime-sum factors choose $w=w_o=\ell=0$ and
$g=J_++\sigma$, enlarging its positive rows to a smooth ball
of length $K+g$.

$\mathcal N_C$ initially retains its genuine bounded frequency dyad
$q_h\le C Z^K$. In the case with prime-sum factors do not replace this graph
by an unrestricted Schwartz row sum before amplification: the admissible
pool cardinality and representation bound use this bounded length.
Insert the smooth positive lift only after Jensen and the positive
enlargement of the transformed summation range.

For sums with prime-sum factors choose $\ell_* =\sigma/3$, $\eta<\sigma/6$,
and the prime pool $Z^{\ell_*}/2<q_p\le Z^{\ell_*}$.
The prime ideal theorem gives its size $Z^{\ell_*+o(1)}$.
It lies beyond every remaining prime annulus for sufficiently large $Z$,
independently of the number of prime-sum factors. Remove pool primes dividing
$h$, $s$ or fixed supports; there are only $O(\log Z)$ of them.
CRT gives, when $(p,u)=1$,
\[
 \mathcal G(p^iu,h)=\chi_p(u)^{2i}
            \mathcal G(p^i,h)\mathcal G(u,h),\qquad
 \mathcal G(u,hp^6)=\mathcal G(u,h).
\]
For $p\nmid h$, replacing $h$ by $hp^6$ changes only valuations
$i=1,6,7$. The exact squared central coefficients are respectively
$P^{-1},(1-P^{-1})^2,P^{-1}$, $P=q_p$.
The row phase is $\bar\chi_p(h)$ for $i=1,7$ and one for $i=6$.
Allocate $p^i$ only to the two plains; divide the same two rectangle
scales by the same extracted powers. The residual product is punctured
at $p$ and its new character phase acts on that whole product.

For each genuine row $h$, retain the admissible pool $\mathcal A_h$
until the signed identity has been averaged. Its cardinality is at least
half the complete pool size after the fixed support threshold, uniformly
in $h$, because only $O(\log Z)$ primes were removed. Jensen first gives
a positive sum on the genuine admissible graph with coefficient at most
$2/\#\mathcal A$. Only then enlarge that graph. The fixed finite error
list bounds the original norm
by the averaged $hp^6$ norm and these error norms. Each transformed index $hp^6$
has boundedly many such representations, because the transformed index has bounded
logarithmic length and every pool prime has length at least
$\ell_*/2$. Thus the main positive norm has density $Z^{-\ell_*+o(1)}$,
$\ell=\ell_*$ and $g=J_++2\sigma$.

For the errors fix $p$, keep the normalized pool average, set
\[
 (w,w_o)=(i\log_ZP,e_i\log_ZP),\quad
 e_1=e_7=1, e_6=0,\qquad \ell=0,\quad g=J_++\sigma.
\]
The $i=6$ residue factor is a common coprimality indicator; the other two are
included in $w_o$. A uniform fixed-$p$ estimate adds no prime-count
factor. Errors are not amplified again. In all cases
$0\le w_o\le w\le7\sigma/3$.
The finite-place row test vanishes at the additive origin, so no extra
zero Gauss row is required in this construction.

\subsubsection{The second transform and the decrease in width}
\label{sec:v-return-3}
In this subsubsection, the second transform \eqref{v:m2:eq:M2-second-bridge} of the
amplified sum returns to the same coefficient family, with width
$M'\le M-\sigma$ and
prefactor \eqref{v:m2:eq:M2-child-width},
\eqref{v:m2:eq:M2-child-allowance}; Section~\ref{sec:v-detail-23}
puts each nonterminal term in the range of the finite-width induction
of Section~\ref{sec:v-detail-24}.

Restore all remaining sums before expanding the square. The positive norm
has the form
\[
 Z^{-\ell}\sum_{h\ne0}\Psi_2(q_h/Z^{K+g})
 \left|Z^{-a_0/2}\sum_u B(u)\tau_1(u)q_u^{it}
                                  \mathcal G(u,h)\right|^2.
\]
Here $B(u)$ is the entire allocated convolution, including its common
coprimality indicator, $s\mid u$, and the centered subtraction. Full finite-correlation
Poisson gives
\begin{equation}\label{v:m2:eq:M2-second-bridge}
 Z^{K+g-\ell-a_0}
 \sum_{u,v}\frac{B(u)\overline{B(v)}\tau_1(u)\overline{\tau_1(v)}
                  q_u^{it}q_v^{-it}}{\sqrt{q_uq_v}}
 \sum_j\mathcal F(u,v;j)
 \mathcal K_2\left(\frac{Z^{K+g}q_j}{q_uq_v}\right).
\end{equation}
The normalization follows from finite Fourier inversion. For fixed $u,v$, put
$f_{u,v}(h)=\mathcal G(u,h)\overline{\mathcal G(v,h)}$.
It is periodic modulo $uv$, without any primitivity requirement, and
Lemma~\ref{v:m2:lem:full-correlation-new} gives
\[
 \sum_{h\bmod uv}f_{u,v}(h)e(-jh/(uv))
                =\sqrt{q_uq_v}\,\mathcal F(u,v;j).
\]
Orbit Poisson at scale $Z^{K+g}$ multiplies this finite transform by
$Z^{K+g}/(q_uq_v)$. The original squared central normalization and
amplifier density are $Z^{-a_0-\ell}$. Their product is exactly
$Z^{K+g-\ell-a_0}/\sqrt{q_uq_v}$, as displayed.
If the additive Fourier convention produces $-j$ in the infinite
transform, replace its norm-one integration variable by its negative;
the norm-sheet kernel and $q_j$ are unchanged. This sign change is
made before the finite coefficient is estimated.

At $j=0$ the
original orbit formula instead supplies the compact quotient integral.
The finite correlation is zero unless $u=v$, when it is
$\varphi(u)\le q_u$. Since $s\mid u$, the diagonal exponent is
$K+g-\ell-s_0+\theta_N+\eps_0$, where the single aggregate support
error $\theta_N=H_N/\log Z$ includes all fixed windows.
Its nominal excess over \eqref{v:m2:eq:M2-allowance} is
\[
 A-M-d-(K_0-K)-w_o-B_c+g-\ell
 \le A-M+5\sigma/3+\delta_1.
\]
To see this, put $D_0=d+K_0-K$ and note
$-D_0-w_o+J_+\le\delta_1$; the largest remaining increment is
$2\sigma-\ell_*=5\sigma/3$.

For $j\ne0$ first discard whole high-frequency dyads on the genuine
sum, using $|\mathcal F(u,v;j)|\le q_uq_v$ and the same radial-tail
argument. Only afterward extract $u=D_2a$, $v=E_2b$, their complete
common support. Let $c_2,d_2$ be its two full lengths,
$b_2=(c_2+d_2)/2$, $p_2$ its radical length, $G_c=(D_2,E_2)$,
$g_2=\log_Zq_{G_c}$. The all-prime-power CRT identity is
\[
 \mathcal F(D_2a,E_2b;j)=
 \mathcal F(D_2,E_2;j)\chi_a(j)\bar\chi_b(-j)
\]
on the genuine coprime locus. Every remaining factor is required to be coprime to the extracted primes;
each prime factor supplying an extracted prime is fixed.

For equal multiplicities $i$ not divisible by six, divide the frequency
by $G_c$ and partition the prime according as that divided frequency
is a unit or nonunit. Let $t_2$ be the radical length of the unit set,
and $V_{\rm id}$ the product of the nonunit set, $V=\log_Zq_{V_{\rm id}}$.
The local absolute bounds are
\[
\begin{array}{c|c}
 i=i, 6\nmid i, \text{unit}&P^{i-1}\\
 i=i, 6\nmid i, \text{nonunit}&P^i\\
 i=i, 6\mid i&P^i\\
 i>j_0, 6\mid j_0, \text{unit}&P^{j_0}.
\end{array}
\]
Other unequal cases vanish. Put $j=G_cV_{\rm id}h'$.
The scalar $\mathcal F(D_2,E_2;j)/Z^{g_2-t_2}$ has modulus at
most one on its partition; define it to be zero off that partition.
All forced zero allocations are removed before applying
this bound. Active primes of the fixed factor $G_cV_{\rm id}$ are
contained in the $t_2,V$ sets. Its six-divisible fixed factors are
common coprimality conditions, while the natural row character retains all its
zeros and cancellations. Therefore the new displayed support and
nominal row length satisfy
\begin{equation}\label{v:m2:eq:M2-child-width}
 \begin{aligned}
 m'&=2a_0-K-g-g_2-V,&q'&=\widetilde q+w_o+t_2+V,\\
 M'&=m'+q'=M+J-g-g_2+t_2\le M-\sigma.
 \end{aligned}
\end{equation}

Absorb every fixed support constant into $\delta_2$ before asserting
the actual bound $q_{h'}\le Z^{m'+\delta_2}$. A retained sector then
lies in one common enclosing ball of length
$m'_{\rm act}=\max(0,m'+\delta_2)$, $0\le\delta_2<\xi$.
Nonemptiness implies $m'+\delta_2\ge0$, hence
$M'_{\rm act}\le M'+\delta_2$.
The pulled-back row weight and partition scalar are extended by zero
to this ball, before any absolute estimate.

The condition $s\mid u,v$ forces $s$ into this genuine common radical.
Write its radical as $sB$, with $(s,B)=1$. Ordinary ideal
counting gives $\#\{B:q_B\ll Z^{p_2}/q_s\}\ll Z^{p_2}/q_s$.
Thus its single family count is $Z^{p_2-s_0+\eps_0}$; no second
count of $s$ is inserted.
The complete exponent outside normalized residual products is
\[
 (K+g-\ell-a_0)+(p_2-s_0)+(g_2-t_2)+(a_0-b_2)-a_0
 =K+g-\ell-a_0+p_2-s_0+g_2-t_2-b_2.
\]
The initial $-a_0$ and Poisson inverse-root factor are both retained
in this accounting. Subtracting this exponent from the allowance
leaves the required inner-product exponent
\begin{equation}\label{v:m2:eq:M2-child-allowance}
 M'+\Delta',\qquad
 \Delta'=b_2-p_2+w+B_c+\ell\ge w+\ell.
\end{equation}

For fixed genuine $D_2,E_2,j$, define the last displayed correlation
on artificial residual noncoprime pairs with the same full-product
norms, coprimality indicators, kernel and partitioned scalar. Insert the full
$\sum_{\mathfrak t\mid a,b}\mu(\mathfrak t)$ before factor bounds.
These artificial divisor primes are not added to the genuine common
frequency restrictions. Separate the kernel and inverse roots on
fixed whole-row/whole-product logarithmic boxes common to all remaining
labels and both rectangles.

For every divisor prime use
\[
 1_{p\mid\prod_i n_i}
 =\sum_{\varnothing\ne I}(-1)^{|I|+1}
                         \prod_{i\in I}1_{p\mid n_i}.
\]
For a selected plain write $n_i=pl_i$, without puncturing $l_i$ at $p$,
and divide both rectangle scales by $q_p$. A selected prime-sum factor is fixed;
conflicting assignments to the same prime-sum factor are zero. Extracted row
scalars, including their zeros, remain in the signed identity.
Only now split the rows into inducing characters outside or inside
$\Theta$. For nonexceptional rows form two positive transformed sum norms; for
exceptional rows take the absolute product. Only at this stage discard
bounded row weights and enlarge to the common enclosing ball.

\begin{remark}[Fixed sixth-power factors and coprimality conditions]
The six-divisible part of $G_cV_{\rm id}$ is a common coprimality
modulus by an exact zero-preserving identity. For a fixed
frequency factor $F_*$, write
$F_*=F_{\rm act}F_6$, where every exponent of $F_6$ is divisible
by six and the two factors have disjoint prime support. Then, for
every integral column $n$,
\[
 \chi_n(F_*)=\chi_n(F_{\rm act})
                         1_{(n,\rad F_6)=1}.
\]
Here $F_*$ is fixed before the remaining row $h'$ varies. The
indicator occurs on the whole product, hence on both plains and
all prime-sum factors with the same coprimality indicator. It is precisely an extra coprimality indicator
of Lemma~\ref{v:m2:lem:extra-mask}: complete divisor expansion has
subpower coefficient mass, lowers lengths, and preserves the
inducing character.

Whenever its norm is omitted from $q'$,
perform this divisor expansion before invoking a nonexceptional
child estimate. That estimate then has no unexpanded omitted modulus;
its two rectangle scales have been divided by the same extracted
divisors, with their central factors retained. The associated prime
factors are fixed or remain their original full annuli.
In an exceptional centered term it is instead retained in the
same fixed-character mean with coprimality conditions on both rectangles.

This does not permit cancellation of a displayed fixed character
against the varying row $h'$: $\widetilde q$, $w_o$ and $t_2,V$ are
all counted before such cancellations, including redundant active
factors, and only the row-independent pure sixth-power coprimality
indicator is handled by the identity above. Thus the width formula
retains the full displayed support, not only the primitive conductor.
\end{remark}

\subsection{Normalization and exceptional terms in M2}
\label{sec:v-detail-23}
In this subsection, we put the second-transform sums into the normalized
product family (Section~\ref{sec:v-return-4}).
Exact prime extraction and the full correlation identity of
Lemma~\ref{v:m2:lem:full-correlation-new} determine
which rows remain exceptional after enlargement
(Section~\ref{sec:v-return-5}), with excess
\eqref{v:m2:eq:M2-exceptional-excess};
Lemma~\ref{v:m2profile:lem:mean} and
Proposition~\ref{v:m2profile:prop:center} control them
(Section~\ref{sec:v-return-6}).

\subsubsection{Exact extraction, bounded norm ranges and the discretization error}
\label{sec:v-return-4}
In this subsubsection, we extract the original prime coprimality indicators from the family
\eqref{v:m2:eq:M2-child-width}, keeping natural zeros by the priority
rule of Section~\ref{sec:v-detail-13}, retain the normalization on
bounded row or column ranges, and count the excess from the last
selected prime-sum factor once (Section~\ref{sec:v-detail-24}). The
transformed-sum excess is at most $6(w+\ell)+\delta_1+2\theta_N$, and
restoring the length inequality is paid for by
Lemma~\ref{v:m2:lem:prime-capacity}.

For one incoming side, let $g_i$ be the actual product extracted
from plain position $i$ by genuine operations, and $d_i$ the actual
product extracted by the final divisor expansion. The child scales are
\[
 \begin{array}{c|cc}
 &\text{first rectangle}&\text{comparison rectangle}\\\hline
 \text{incoming}&X_i&Y_i\\
 \text{child}&X_i'=X_i/(q_{g_i}q_{d_i})
            &Y_i'=Y_i/(q_{g_i}q_{d_i}).
 \end{array}
\]
Thus $X_1'X_2'=Y_1'Y_2'$. Every nonexceptional child is
invoked at these actual scales with its original annular profile,
not at a coarse center.

For a nonempty annular sum with $\operatorname{supp}W\subset[a_W,b_W]$,
\[
 W(q_n/T)\ne0\ \Longrightarrow\ T\ge b_W^{-1}.
\]
Consequently a subunit child scale is in a fixed compact interval,
and its conversion to scale one has uniformly bounded Euler seminorms.
Exceptional signed differences keep both original scales unchanged.
The child row ball is likewise the actual ball of exponent
$m'_{\rm act}=\max(0,m'+\delta_2)$ in
\eqref{v:m2:eq:M2-child-width}; the $\delta_2$ allowance is included
before testing the strict decrease in width.

Put $\alpha_1=a_0-c_2$, $\alpha_2=a_0-d_2$.
The two residual sides have original pre-normalizers $Z^{-\alpha_j/2}$.
Record exact removed plain norms rather than approximate lengths.
For each side let $e_j$ be the sum of $\log_Z(q_p/P_i)$ for prime-sum factors
fixed by genuine extraction in this current stage. Let
$\widetilde r_j$ be the sum of the two last-divisor plain removal
lengths and the nominal lengths of prime-sum factors fixed by that divisor,
and $\omega_j$ the latter prime-sum factors' ratio sum. Then
\[
 T'_j\prod_{i\ {\rm remaining}}P_i
       =Z^{\alpha_j+e_j-\widetilde r_j},\qquad
 |e_j|,|\omega_j|,|e_j+\omega_j|\le\theta_N.
\]
The subsets defining $e_j,\omega_j$ are disjoint subsets of the
prime-sum factors present at the start of this induction step. Factors removed
in earlier steps are not counted again.

For a lower-endpoint divisor dyad $T\le q_{\mathfrak t}<2T$,
$t_- =\log_ZT$, its selected product is divisible by $\mathfrak t$,
so
\[
 \widetilde r_j+\omega_j\ge t_-,\qquad
 \#\{\mathfrak t\text{ in the dyad}\}=O(Z^{t_-}).
\]
The coefficient relative to the nominal raw post-normalizer is
$Z^{-\widetilde r_j/2}$; the norm ranges in the signed difference have not been altered.

On nonexceptional rows first bound the centered difference by its two
rectangles. In each rectangle separately, replace each negative length exponent by zero. If its two old and new
actual logarithmic plain lengths are $x_i$, $y_i=x_i-a_{j,i}$, put
$\pi_0=\sum_i(-x_i)_+$, $\pi=\sum_i(-y_i)_+$.
Nonempty annular endpoints imply $0\le\pi_0\le\pi\le\theta_N$.
Define $r_j^{+}=\widetilde r_j-(\pi-\pi_0)\ge0$.
The resulting total factor length and exact normalization conversion are
\[
 A_j^{+}=\alpha_j+e_j+\pi_0-r_j^{+},\qquad
 Z^{(-r_j^{+}+e_j+\pi_0)/2}.
\]
Together with the divisor count and the boundary inequalities, the
paired count and normalization exponent is at most $2\theta_N$.
Do not use one common adjusted length for the signed difference.

The fixed aggregate list in the preceding exponent calculations gives a numerical
constant $C_*$ independent of the number of prime-sum factors, such that the
transformed sum excess over the admissible length is at most
\[
 6(w+\ell)+\delta_1+2\theta_N.
\]
Before substitution, the defect is at most
\[
 w-D_0-w_o+g+2\theta_N,\qquad D_0=d+K_0-K.
\]
Since $-D_0-w_o+J_+\le\delta_1$, this is at most
$w+(g-J_+)+\delta_1+2\theta_N$.
For the selected main term $w=0$, $g-J_+=2\sigma=6\ell$.
For an amplifier error, $\ell=0$, $g-J_+=\sigma$ and
$w=i\log_Zq_p$, $i\in\{1,6,7\}$. Increase the threshold so
$\log_Z2\le\sigma/12$; the pool then gives $w\ge\sigma/4$,
and $w+\sigma\le5w\le6w$. This proves the stated bound.

All genuine support removals only decrease the remaining prime length.
For a transformed sum with prime factors, remove an initial
segment of its surviving prime-sum factors until the length inequality
holds or no prime-sum factor remains. Since each prime-sum factor has
length $\le\eta$, the removed total length is at most the defect divided by $6\kappa$, plus one $\eta$.
Lemma~\ref{v:m2:lem:prime-capacity} costs in squared norm at most
\[
 \kappa\,d_z\le\Delta'
                     +\delta_1/6+\theta_N/3+\eta,
\]
since $\Delta'\ge w+\ell$ and $\kappa\le1$.
If all prime-sum factors disappear use the established
assertion without prime-sum factors at this smaller width. The initial
segment incurs one overshoot, not one discretization error per prime-sum factor.
The remaining edge losses, including the enclosing ball, are bounded
by $C_*\xi+\eta$ after the fixed support threshold.

\subsubsection{The sparse exceptional rows}
\label{sec:v-return-5}
In this subsubsection, Lemma~\ref{v:m2:lem:full-correlation-new} classifies the exceptional
valuations after positive row enlargement. Old exclusions remain
natural zeros; new exceptional primes force $v_p(h')\ge4$ and contribute the excess
\eqref{v:m2:eq:M2-exceptional-excess}; the signed coefficient is
centered in Section~\ref{sec:v-return-6} before numbers of terms in the
prime sums are restored.

For an exceptional inducing character, the outside-$S$ exponents of
$h'$ are prescribed modulo six. At an equal multiplicity-one common
prime in the nonunit partition the fixed frequency factor
$G_cV_{\rm id}$ contains $p^2$. This prime is coprime to the
earlier displayed moving support: before the second transform both
genuine columns retained the natural zeros of $\tau_1$, so a prime
in that support could not divide either column. In particular it
could not enter the genuine $D_2,E_2$. Fixed extra coprimality conditions and
the fixed amplifier prime have the same exclusion. The fixed
finite character is unramified here because the prime is outside
$S$. Thus no earlier fixed local character can cancel the new
exponent two, and the later artificial common-divisor expansion
contributes only row scalars. Since the new fixed contribution at such
a prime is $p^2$ and the earlier one is a unit, trivial exterior
ramification requires
\[
 2+v_p(h')\equiv0\pmod6,
 \qquad v_p(h')\ge4.
\]
This is the equality case in the $i=1$ nonunit row:
$(7-5)/6+1/3=2/3$. Let $v_1$ be the radical length of these
primes and put $f=2v_1$. The minimal outside-$S$ sixth-power-free part
has norm at least $Z^{2f}$.

At every exterior prime, tame Kummer ramification forces the row
valuation into one prescribed class modulo six, since the inducing
character lies in the fixed $S$-supported group. If that prescription is
consistent, its exterior ideal is one fixed sixth-power-free ideal times
a sixth power. Counting the latter uses the ordinary ideal count at the
one-sixth power of the remaining norm bound. The quotient of allowed
$S$-unit choices by $U_S^6$ is finite, and does not supply an additional
norm count. Inconsistent prescriptions give no rows.
The finitely many allowed unit classes and fixed $S$-components give
\[
 \#\{h':\ h'\text{ exceptional}\}
 \ll Z^{(m'_{\rm act}-2f)/6+\eps_0}
 \le Z^{(m'-f)/6+\delta_2/6+\eps_0}.
\]
Only the weaker second bound is used. The count follows by writing the
outside ideal as its prescribed sixth-power-free part times a sixth
power; rows admitted by enlargement satisfy the same
character-only valuation condition, so the count applies on
the enlarged ball without the original unit restrictions.

For raw uncentered transformed sums the paired divisor allocation and absolute
factor volume have bound
\[Z^{a_0-b_2+2\theta_N+\eps_0}.\]
Subtracting the permitted exponent \eqref{v:m2:eq:M2-child-allowance}
from this volume plus the weaker exceptional count bounds the excess by
\begin{equation}\label{v:m2:eq:M2-exceptional-excess}
 A-5M/6-F_1-F_2+2\theta_N+\delta_2/6+2\eps_0,
\end{equation}
where
\[
 \begin{split}
 F_1={}&c/6+5(d+K_0-K)/6+w/3+\widetilde q/6+w_o+B_c
                                      -5g/6+\ell,\\
 F_2={}&2b_2-5g_2/6-p_2+t_2+V/6+f/6.
 \end{split}
\]
The frequency and support perturbations remain explicit.
Let $D_0=d+K_0-K$. With $g=J_+$, if $J\ge0$ the first
expression is $c+2w+\widetilde q/6+w_o/6+B_c+\ell$.
If $J<0$, then $D_0\ge d-\delta_1$ and
$\widetilde q\ge R$ give
\[
 B_c+\widetilde q/6
 \ge(3c-5d-R)_+/6+R/6
 \ge(3c-5D_0)_+/6-5\delta_1/6.
\]
Consequently
\[
 \frac c6+\frac{5D_0}6+B_c+\frac{\widetilde q}6
       \ge\frac{2c}3-\frac{5\delta_1}6,
\]
since $5D_0+(3c-5D_0)_+\ge3c$. The additional $w/3$
falls short of $2w/3$ by at most $7\sigma/9$.
In either sign of $J$, the enlargement $g-J_+\le2\sigma$
costs at most $5\sigma/3$. As
$5/3+7/9<3$, this proves
\[
 F_1\ge\tfrac23(c+w)-3\sigma-5\delta_1/6.
\]

For $F_2$ its local contributions at a common prime are
\[
\begin{array}{c|c|c}
\text{case}&F_2&b_2\\ \hline
 i=i, 6\nmid i, \text{unit}&7i/6&i\\
 i=i, 6\nmid i, \text{nonunit}
       &(7i-5)/6+1_{i=1}/3&i\\
 i=i, 6\mid i&7i/6-1&i\\
 i>j_0\ge6, 6\mid j_0&i+j_0/6-1&(i+j_0)/2.
\end{array}
\]
Each is at least $2b_2/3$. For the second line $i=1$ gives
equality; for $i\ge2$ its difference is $(3i-5)/6\ge0$.
For the third $i\ge6$, and for the fourth
$i+j_0/6-1-(i+j_0)/3=(4i-j_0-6)/6\ge0$.
Thus $F_2\ge2b_2/3$. In the uncentered range $A\le5M/6$,
\eqref{v:m2:eq:M2-exceptional-excess} costs at most $3\sigma+O(\xi+\eps_0)$.

\subsubsection{Centering exceptional rows with the original norm ranges}
\label{sec:v-return-6}
In this subsubsection, we apply Proposition~\ref{v:m2profile:prop:center} to the exceptional
rows, keeping the equal-product coefficient signed until its principal
terms cancel, with common coprimality indicators controlled by
Lemma~\ref{v:m2profile:lem:mean}. The resulting terminal estimate has
excess at most $(A-M)_++3\sigma+O(\xi+\eps_0)$.

For a fixed exceptional row the common presentation reduces on units
to one $\vartheta\in\Theta$. Its remaining natural zeros and all
fixed coprimality conditions are expressed by $1_{(l,R_*)=1}$, where
$R_*$ has polynomially bounded norm. Lemma~\ref{v:m2profile:lem:mean} gives
\[
 \sum_l\vartheta(l)1_{(l,R_*)=1}q_l^{it}W_i(q_l/T)
 =c_{\vartheta,R_*}T^{1+it}I_i(t)
        +O\bigl(Z^{\eps_0}(1+|t|)^J p_J(W_i)\bigr),
\]
with $T^a$ absorbed into the $Z^{\eps_0}$ loss on the polynomial range.
The same character, coprimality indicator and $t$ apply to both plains
and both rectangles, so $T_1T_2=Y_1Y_2$ cancels the leading product.
For scales below one use raw annular support and absolute volume,
not a renormalized moment with altered norm ranges.

All four initial plain lengths are at least $L=M/4$.
Genuine extractions and amplification lower a smallest length by at
most $v=c+w+\min(c_2,d_2)$. Put $r=(L-v)_+$.
Call the side attaining this reference saving side one, relabeling
its raw data with it. Its total last-divisor removal is
$\widetilde r_1$, including newly fixed prime-sum factors. Each individual plain
removal is at most this total. Thus the centered Mellin estimate saves
$(r-\widetilde r_1)_+$ against absolute volume. The paired exact
coefficient, raw-volume and divisor-count exponent calculation is
\[
 a_0-b_2+e_1+e_2+t_- -\widetilde r_1-\widetilde r_2
                                      -(r-\widetilde r_1)_+.
\]
The complete selected-product boundary gives
\[
 t_- -\widetilde r_1-\widetilde r_2-(r-\widetilde r_1)_+
 =(t_- -\widetilde r_2)-\max(r,\widetilde r_1)
 \le\omega_2-r.
\]
Adding $e_1+e_2$ costs at most $2\theta_N$ since
$|e_1|,|e_2+\omega_2|\le\theta_N$.
Total removals are needed here: plain removals alone would omit the
exponent of a fixed prime factor.
Thus the exceptional paired volume is
$Z^{a_0-b_2-r+2\theta_N+\eps_0}$, with seminorm and height factors
retained in the common coefficient measure.
Since $F_1+F_2\ge2v/3-3\sigma-O(\xi)$, its excess is at most
\[
 A-5M/6-2v/3-(L-v)_+ +3\sigma+O(\xi+\eps_0)
 \le(A-M)_++3\sigma+O(\xi+\eps_0).
\]
Indeed,
\[
 \begin{aligned}
 -2v/3-(L-v)_+
 &=\begin{cases}-L+v/3,&0\le v\le L,\\-2v/3,&v\ge L,\end{cases}\\
 \sup_{v\ge0}\{-2v/3-(L-v)_+\}&=-2L/3=-M/6.
 \end{aligned}
\]

With prime-sum factors the length condition gives $A\le M$; without
them $A\le M+\delta$. This is a terminal estimate, calling no
larger-width moment. The common Fourier separation precedes fixing
remaining labels; absolute numbers of terms in the prime sums are
restored only after the finite Gauss coefficient is gone, and
last-divisor allocation preserves the equal-product identity in the
signed rectangles.

\subsection{Finite-width closure}
\label{sec:v-detail-24}
In this subsection, we combine Sections~\ref{sec:v-detail-22} and~\ref{sec:v-detail-23}.
Every nonterminal induction step strictly lowers the width, terminal
widths are controlled by Lemma~\ref{v:m2:lem:prime-capacity}, and
Lemma~\ref{v:lem:finite-depth-loss} propagates the losses through the finite
depth; with Lemma~\ref{v:lem:M2-fixed-addback} this proves
Theorem~\ref{v:thm:M2-assembled}.
Lemma~\ref{v:lem:arithmetic-independent-mesh} chooses the mesh before
the number of prime-sum factors is fixed, and
Remark~\ref{v:rem:sextic-fourth-moment} records the sextic fourth-moment
specialization.

\begin{lemma}[Propagation of losses through a finite induction]
\label{v:lem:finite-depth-loss}
Let $Z\ge2$, $C\ge1$, $\tau,\delta\ge0$. Suppose nonnegative
normalized quantities $E_v$ have ranks in $\{0,\ldots,D\}$.
At a terminal index suppose $E_v\le CZ^\tau$. At every other index
suppose there is a positive measure $\mu_v$, with
$\mu_v(\Omega_v)\le CZ^\delta$, such that
\[
 E_v\le CZ^\tau+
     CZ^\delta\int_{\Omega_v}
        \bigl(E_{v_1(\omega)}E_{v_2(\omega)}\bigr)^{1/2}
                           \,d\mu_v(\omega),
\]
where both ranks on the right are strictly smaller than the rank of
$v$. Then
\[
 E_v\le C(1+C^2)^D Z^{\tau+2D\delta}.
\]
\end{lemma}
\begin{proof}
At rank zero the assertion is the terminal bound. If every smaller
rank is bounded by $A_{d-1}Z^{\tau+2(d-1)\delta}$, positivity and
the measure bound give
\[
 E_v\le CZ^\tau+C^2A_{d-1}Z^{\tau+2d\delta}
       \le(C+C^2A_{d-1})Z^{\tau+2d\delta}.
\]
Starting with $A_0=C$, induction gives
$A_d\le C(1+C^2)^d$. The terminal exponent
$\tau$ occurs once on every branch; only the edge loss
$2\delta$ accumulates with depth.
\end{proof}

The same lemma applies with profile seminorms and height weights as
follows. Assign each index a positive weight $\mathfrak S_v\ge1$,
and suppose the actual estimates read
\[
 E_v\le C\mathfrak S_vZ^\tau+
       CZ^\delta\int(E_{v_1}E_{v_2})^{1/2}\,d\mu_v,
 \qquad
 \int(\mathfrak S_{v_1}\mathfrak S_{v_2})^{1/2}\,d\mu_v
       \le C\mathfrak S_vZ^\delta.
\]
At terminal indices require $E_v\le C\mathfrak S_vZ^\tau$.
Put
\[
 \widetilde E_v=E_v/\mathfrak S_v,\qquad
 d\widetilde\mu_v=
       \frac{(\mathfrak S_{v_1}\mathfrak S_{v_2})^{1/2}}
            {\mathfrak S_v}\,d\mu_v.
\]
Then $\widetilde\mu_v(\Omega_v)\le CZ^\delta$, and the
normalized recurrence is exactly the lemma's hypothesis. Its conclusion
is $E_v\le C(1+C^2)^D\mathfrak S_vZ^{\tau+2D\delta}$.
For the present recursions, $\mathfrak S_v$ is a finite homogeneous
seminorm of the incoming tests multiplied by a polynomial in their
incoming norm-twist heights. Choose the child orders at each lower
rank first. Equation~\eqref{v:eq:H2-common-density} then bounds
their polynomial weights against one higher weighted Fourier moment,
which determines the parent order. The one normalized discrete family
count is part of $\mu_v$. Since the rank is bounded, this procedure
terminates and gives one constant $C$ independent of all subsequently
integrated Fourier variables. Thus no uniformity over those variables
is inferred from an unweighted Fourier integral.

\begin{lemma}[Application of the product recurrence to its actual children]
\label{v:lem:M2-returned-family}
In the construction of Sections~\ref{sec:v-detail-20}--\ref{sec:v-detail-23},
every nonexceptional returned norm is an ordinary or selected product
moment of Section~\ref{sec:v-detail-19}, at the actual rectangle
scales. A two-transform return has
\[
 M'_{\rm act}\le M-\sigma+C_*\xi,
 \qquad \text{inner allowance }M'+\Delta',\quad \Delta'\ge w+\ell.
\]
The separate extra modulus is expanded before invoking that child.
Origins, diagonals and exceptional rows are terminal contributions.
For the choices in the proof of Theorem~\ref{v:thm:M2-assembled},
the actual positive measures obey the weighted bounds needed for
finite induction, with edge loss
$\delta_{\rm edge}\le C_*(\eta+\xi+\eps_0)$.
More explicitly, if the required lower-rank moments are bounded by
$A\mathfrak S_uZ^\lambda$, then the normalized parent is at most
\[
 C\mathfrak S_v Z^\tau+
 CA\mathfrak S_v Z^{\lambda+2\delta_{\rm edge}}.
\]
Here $\mathfrak S_v$ is a finite homogeneous test seminorm times a
polynomial in the original incoming heights; the orders may be chosen
backward through a bounded number of returns.
Appendix~\ref{app:weighted-edges} records the separate weighted-integral
bounds and the returned-data substitutions used here.
\end{lemma}
\begin{proof}
The source identities are the full Gauss bridges
\eqref{v:m2:eq:M2-first-bridge} and
\eqref{v:m2:eq:M2-second-bridge}. They include the signed centered
difference on the two rectangles. Extend the full displayed
correlation to artificial pairs and insert each complete M\"obius
sum before positivity. For each divisor prime, the nonempty-subset
expansion in Section~\ref{sec:v-return-3} assigns it to column
positions. A plain position is divided by that prime on both
rectangles, without imposing a new puncture on its residual variable;
a prime-sum position is fixed, and incompatible assignments vanish.
Thus every unfixed prime position remains a full original annulus.
The extracted row scalars and their zeros are retained. Only after
forming the actual residual character is its exceptional split made.
This gives exactly the permitted profiles, finite coefficients and
common deletions, rather than arbitrary coefficients on prime subsets.

The first discrete count is $Z^{p+s_0+\eps_0}$; the nonpositive brace
following \eqref{v:m2:eq:M2-Gauss-target} pays for it together with
the two target norms. At the second return the count is
$Z^{p_2-s_0+\eps_0}$, because the genuine common radical is $sB$.
The old $s$ is not counted a second time. Multiplication by the
second Poisson prefactor gives exactly the outside exponent preceding
\eqref{v:m2:eq:M2-child-allowance}; subtracting it leaves
$M'+\Delta'$. The identity
\eqref{v:m2:eq:M2-child-width} gives
$M'=M+J-g-g_2+t_2\le M-\sigma$.
The actual enclosing ball adds $\delta_2<\xi$ only when it is
nonempty. The scales $X_i',Y_i'$ of Section~\ref{sec:v-return-4}
preserve $X_1'X_2'=Y_1'Y_2'$ and all central extraction factors.
Its last-divisor count and bounded-scale conversions cost at most
$2\theta_N$. The length excess is at most
$6(w+\ell)+\delta_1+2\theta_N$; the prefix restoration therefore
costs at most $\Delta'+\delta_1/6+\theta_N/3+\eta$.
This uses one overshooting prime, irrespective of the number of slots.
If no prime positions remain, the returned family is ordinary.
The pure sixth-power mask and $R_*$ are treated by the complete extra-mask
expansion; canceled moving support remains counted in $q'$.

For clarity, the measure on each branch consists of the normalized
structural count just described, the absolute common Fourier densities,
the two-scale reflection measures when reflection is used, and the
normalized prime-pool average when amplification is used. No additional
pool cardinality is inserted for an error term. Whole-column Fourier
coordinates are common to the two rectangles and their factors; they
do not assign independent modes to the individual plains of a signed
difference. For prescribed child orders, differentiation of the
separated annular profiles introduces only polynomial Fourier-height
weights. Integrate those weights by
\eqref{v:eq:H2-common-density}; use the mixed derivative as well as
the separate derivatives in the two-scale reflection measure.
For each Cauchy pairing this proves
\[
 \int (\mathfrak S_{u_1}\mathfrak S_{u_2})^{1/2}\,d\mu
       \le C\mathfrak S_v Z^{\delta_{\rm edge}},
\]
with the displayed structural powers already removed. If a step
retains a product of integrated norms, apply the same weighted bound
to each norm before multiplying their square roots. In either form a
common lower-rank bound $A Z^\lambda$ is propagated once, giving the
asserted parent inequality. One may first normalize each nonzero input
test by its required seminorm and then restore its original degree of
homogeneity; a zero test gives a zero moment.

Finally the original origin, second diagonal, and exceptional branches
are estimated in Sections~\ref{sec:v-return-1},
\ref{sec:v-return-3}, \ref{sec:v-return-5} and~\ref{sec:v-return-6}.
The exceptional centered term uses the same fixed-character mean on
both rectangles, so its leading products cancel before absolute values.
Their losses are terminal, at most $\delta+14\sigma/3+O(\xi+\eps_0)$;
they are not assigned a child estimate outside $\Theta$.
The same-width reflection and centering calls and the width-reducing
calls have the strictly decreasing rank specified in the theorem's
proof. All support constants enter through the single aggregate
$\theta_N$, whereas polynomial height weights are integrated by
increasing orders. Hence these operations enlarge constants and orders,
not $C_*$ by a factor equal to the number of prime positions or a
profile derivative order. This proves the claimed weighted application
without presupposing a bound for any returned moment.
\end{proof}

\begin{lemma}[Adding back fixed nonprincipal inducing characters]
\label{v:lem:M2-fixed-addback}
For $z=0$, suppose the unweighted product estimate has been proved
on inducing rows outside a fixed finite group $\Theta$ supported in
$S$. Then the same bound holds on all nonprincipal inducing rows,
with the same length ranges and polynomial height bounds.
\end{lemma}
\begin{proof}
For any additional row its inducing character is some
$\theta\in\Theta\setminus\{1\}$, and its common presentation is
$\theta(n)1_{(n,R_k)=1}$, with $q_{R_k}$ polynomially bounded in
$Z$. The residue in Lemma~\ref{v:m2profile:lem:mean} is zero, so
for each annular factor of nonempty scale $X\ge c_W>0$,
\[
 \left|X^{-1/2}\sum_n\theta(n)1_{(n,R_k)=1}W(q_n/X)\right|
 \ll Z^\eps(1+|t|)^Jp_J(W).
\]
This applies separately to both ordinary factors, with their
original profiles. There are $O(Z^m)$ rows, including the fixed
finite unit multiplicity. Their total is therefore
$\ll Z^{m+\eps}\le Z^{M+\eps}$ after decreasing the input loss.
The principal inducing rows remain excluded.
\end{proof}

\begin{remark}[Proof obligations in the returned-family induction]
\label{v:rem:M2-obligations}
The induction proving Theorem~\ref{v:thm:M2-assembled} below consists
of statements of four kinds.
\begin{enumerate}
\item \emph{Exponent accounting:} the identities and inequalities for
the first-return brace, the child width and allowance, the exceptional
excess, the lower bounds for $F_1$ and $F_2\ge2b_2/3$, the
prefix-removal and length-restoration capacities, and the count of
exceptional rows (Sections~\ref{sec:v-return-1},
\ref{sec:v-return-3}--\ref{sec:v-return-6}). These are unconditional
real-arithmetic statements and have been checked, with the per-row
bounds of item~(3) taken as premises.
\item \emph{Constructions:} the returned families, their measures (the
normalized structural count, the absolute common Fourier densities,
the two-scale reflection measures and the fixed exceptional
coefficients described in Section~\ref{sec:v-return-4}) and the mesh
of Lemma~\ref{v:lem:arithmetic-independent-mesh}. They are described
in the text, with the stated invariance under the unit operations.
\item \emph{Per-row estimates:} for each returned or exceptional row,
the paired raw-volume bound $Z^{a_0-b_2-r+2\theta_N+\eps_0}$ of
Section~\ref{sec:v-return-6} and the centered Mellin bound of
Proposition~\ref{v:m2profile:prop:center}, with the seminorm and
height factors stated there. These are the analytic content of M2.
\item \emph{Identification steps:} the identification, by Kummer
theory and the correlation identities of
Section~\ref{sec:v-detail-21}, of the rows that are exceptional at a
given return, and the order of choices between the mesh, the number
of prime-sum factors and the arithmetic data
(Lemma~\ref{v:lem:arithmetic-independent-mesh}).
\end{enumerate}
Items (2)--(4) specify the constructions, estimates and identifications
used in the induction. Nothing in item~(1) depends on the target character.
\end{remark}

\begin{theorem}[Product moment with short prime sums]\label{v:thm:M2-assembled}
Fix the arithmetic data above, bounded total lengths, and $\eps>0$.
There is a mesh $\eta>0$ depending only on the requested exponent loss
and the numerical aggregate length and modulus-exponent bounds, as proved
in Lemma~\ref{v:lem:arithmetic-independent-mesh}, such that M2's estimate holds for $0\le z_i\le\eta$.
For $z=0$ exclude the principal inducing character and allow arbitrary
bounded plain lengths. For $z>0$ exclude $\Theta$ and require
$n_1+n_2+6\kappa z\le M$, $13/18\le\kappa\le1$.
When $\kappa<1$ retain the hypothesis $\beta_*(F)\le(1+\kappa)/2$.
The range $\kappa\ge13/18$ ensures $6\kappa-1\ge10/3$ in the
centering comparison of Section~\ref{sec:v-detail-20}.
Uniformity is required over the moving presentations of
Section~\ref{sec:v-detail-19}, with their full displayed moving
radicals retained in $q$ and $M=m+q$. A separate common extra
radical $R_*$, fixed within the current row sum and of norm at most
$Z^{C_{\rm extra}}$, is removed by the complete expansion of
Lemma~\ref{v:m2:lem:extra-mask} before a nonexceptional child or
ordinary reflection is invoked. The bound $C_{\rm extra}$ affects
constants and subsidiary losses, not the definition of $M$.
Natural zeros of cancelled moving factors remain in the displayed
moving radical.
Fix the finite coefficient group $\Theta$ before the induction,
containing $\Upsilon_S$ and all initial finite coefficients.
For fixed internal orders, the bound has a finite homogeneous test
seminorm and polynomial dependence on every norm-twist height,
uniformly in these moving presentations. The implied constant and
the large-$Z$ threshold may depend on the fixed field, $S$, $\Theta$,
$C_{\rm extra}$, length bounds, profiles, number of prime-sum factors and
requested internal orders. The induction below proves this uniform assertion, with the
numerical order of choices made explicit in
Lemma~\ref{v:lem:arithmetic-independent-mesh}.
\end{theorem}
\begin{proof}
For both types first prove the assertion on rows outside $\Theta$.
At each transformed norm make its exceptional split only after the
new character and its finite unit twist have been formed.
Equation~\eqref{eq:v-unit-exception-invariance} makes this exclusion
stable under the unit operations. At the end,
Lemma~\ref{v:lem:M2-fixed-addback} recovers the stated principal-only
exclusion for $z=0$.

Choose a lower width threshold $\rho$, step $\sigma$ and extra length margin $\delta$,
small in terms of $\eps$ so that
$7\rho/6+\delta+14\sigma/3<\eps/8$. Put
\[
 N=\lceil4M_{\max}/\sigma\rceil,\qquad D=4N+3,
 \qquad b(M)=\lceil4M/\sigma\rceil.
\]
Let $t=0,1$ distinguish sums without and with prime-sum factors,
and $c=0,1$ distinguish the uncentered range $A\le5M/6$ and
its centered complement. Use the integer rank
\[
 r(t,M,c)=2(N+1)t+2b(M)+c\le D.
\]
Thus widths without prime-sum factors precede those with them;
within each type widths increase, and at fixed width
uncentered precedes centered.

For arbitrary ordinary lengths include the reflection reduction in
the rank. Introduce a phase $p$: an arbitrary-length input has
$t=0,c=1,p=1$; a balanced input has $t=0,p=0$ and
$c=1_{A>5M/6}$; inputs with prime-sum factors have $p=0$.
For an arbitrary-length input, $c=1$ is only a rank label. Put
\[
 \mathcal R(t,M,c,p)=2\{2(N+1)t+2b(M)+c\}+p,
 \qquad D_\sharp=2D+1.
\]
Lemma~\ref{v:m2:lem:M2-reflection}, applied to both factors with its
common two-scale measure, sends an arbitrary-length input to balanced
inputs of the same width and smaller $\mathcal R$. Their lengths obey
$n_i\le M/2+\xi$ and $A\le M+2\xi\le M+\delta$; this resets the
balanced allowance instead of adding it at each depth. Reflection
losses are charged on this edge and on every later edge where that
reduction occurs.

For $z=0$, a balanced centered comparison invokes the balanced uncentered
assertion: its shorter length is $M/4$, and reflection of the longer
factor gives length
\[
 M-(A-M/4)+\xi<5M/12+\xi\qquad(A>5M/6).
\]
Both lengths are at most $M/2+\xi$ and their total is less than
$2M/3+\xi<5M/6$ for $M>\rho$ and $\xi<\rho/60$. Thus this call
reduces $c$ and never returns to the arbitrary-length phase at the
same width. A width decrease lowers $b(M)$ by at least two, even
if the child needs that phase. A call from $t=1$ to $t=0$ has target
rank at most $4N+3$ and source rank at least $4N+4$. Every call
therefore strictly decreases $\mathcal R$, and branch depth is at
most $D_\sharp$.

At $M\le\rho$, absolute counting gives an excess including the positive length margin without prime-sum factors
at most $\rho+\delta$. For sums with prime-sum factors use that result without prime-sum factors
and Lemma~\ref{v:m2:lem:prime-capacity}; since $A\ge z$ and
$A+(6\kappa-1)z\le M$, $\kappa z\le M/6\le\rho/6$.
These are terminal losses.

At larger widths the preceding two-transform construction bounds
its origin, second diagonal and exceptional terms with terminal
losses bounded by $\delta+3\sigma+5\sigma/3+O(\xi+\eps_0)$.
Every width-reducing two-transform edge satisfies
$M'_{\rm act}\le M-\sigma+C_*\xi\le M-\sigma/2$.
Such a non-reflection transform step decreases $b(M)$ by at least
two and therefore decreases $\mathcal R$, even if its child requires
the arbitrary-length phase. Together with the phase and same-width
comparisons above, this accounts for every call. Removing prime-sum
factors invokes the unrestricted assertion at the smaller width;
retaining them invokes a strictly smaller selected width interval.

Choose a numerical $C_*\ge1$ for the fixed aggregate affine
exponent calculations, including both ordinary reflections, their
weighted measures and the fixed coefficients of the exceptional
terminal losses. It is independent of the number of prime factors.
Choose once
\[
 \eta<\min\{\sigma/6,\eps/(512C_*D_\sharp)\}.
\]
After fixing the arithmetic data and number of prime-sum factors choose
\[
 \xi<\min\{\delta/4,\rho/60,\sigma/(8C_*),
                        \eps/(512C_*D_\sharp)\},\qquad
 \eps_0<\eps/(512C_*D_\sharp).
\]
Increase the threshold so $\theta_N<\xi$ and all fixed sector factors
fit the retained bounds. Restoring the length inequality costs at most one
overshooting prime factor per depth. All other fixed support effects
enter through one $\theta_N$ for the aggregate list.

Divide each transformed norm by its required power of $Z$, retaining
$\Delta'$ of \eqref{v:m2:eq:M2-child-allowance}. Reflection, profile,
divisor and prefactor losses, including the common two-scale reflection
measure, receive subsidiary shares in
\[
 \delta_{\rm edge}\le C_*(\eta+\xi+\eps_0).
\]
Section~\ref{app:M2-recurrence-derivation} derives the required
induction estimate directly from the two Poisson formulas, the
amplification identity and the separate weighted Cauchy integrals;
it does not assume Lemma~\ref{v:lem:M2-returned-family}.
That lemma summarizes the transformed moments and measure bounds.
We apply the finite induction of Lemma~\ref{v:lem:finite-depth-loss}
with these estimates. The centered comparison uses
$|P|^2\le2|P-P_{\rm comp}|^2+2|P_{\rm comp}|^2$; its reflection
and profile losses belong to that edge, with a positive length
margin since $\xi<\rho/60$. The terminal exponent satisfies
\[
 \tau\le\rho+\delta+3\sigma+5\sigma/3+C_*(\xi+\eps_0).
\]

Removing prime factors at base width costs
$Z^{\kappa z}\le Z^{\rho/6}$. Replace the normalized quantity $E_v$
by $Z^{-t\rho/6}E_v$. This charges that cost only on the change
$t=1\to0$; the same per-side inequality applies to a geometric mean
of children of different types. The weighted finite-depth lemma at
rank $\mathcal R$ now gives total exponent at most
\[
 \tau+\rho/6+2D_\sharp\delta_{\rm edge}<\eps.
\]
Indeed the terminal numerical part is less than $\eps/8$,
$C_*(\xi+\eps_0)<\eps/(256D_\sharp)$, and
$2D_\sharp C_*(\eta+\xi+\eps_0)<3\eps/256$; the remaining
margin absorbs external dyads and finite-order tails.
The numerical mesh $\eta$ remains fixed. Only the subsidiary input
losses, orders and thresholds are chosen after the arithmetic data.
Reflection and Fourier losses are charged at every edge where they
occur; the terminal allowance and the type-change potential occur once.

All finite Fourier measures are taken on fixed boxes before remaining
labels and are normalized after the one discrete family count.
Weighted logarithmic Fourier bounds retain the polynomial heights
and finite profile seminorms. Choose the required transformed sum orders
backward through this finite depth, and then choose external tail
orders. Small-norm logarithmic kernels cost only subsidiary powers;
there is no factor $Z^{J\xi}$ when orders are increased.
Finally sum the finite logarithmic dyads with their previously reserved
shares. This gives the stated estimate for sufficiently large $Z$;
bounded $Z$ is absorbed in its target-dependent constant.
\end{proof}

\begin{lemma}[Numerical mesh before fixed arithmetic data]
\label{v:lem:arithmetic-independent-mesh}
Fix a requested exponent loss $\epsilon>0$, a bound for the aggregate
length parameters, and a bound $C_{\rm extra}$ for the common extra-modulus exponent.
The mesh in Theorem~\ref{v:thm:M2-assembled} can be chosen from these
numerical data alone. More precisely, there is $\eta_{\rm num}>0$
such that, for every subsequently fixed cyclotomic field $F$, admissible
finite set $S$, finite coefficient group $\Theta$, fixed profiles and
fixed number of prime-polynomial factors, the theorem holds whenever
$z_i\le\eta_{\rm num}$. Its constants, profile orders, polynomial
height orders and large-$Z$ threshold may depend on all these
subsequently fixed data. The zero-free hypothesis when $\kappa<1$ is
unchanged. No uniform constant or threshold over those arithmetic data
is asserted.
\end{lemma}
\begin{proof}
First select the
numerical terminal-width bound $\rho$, width step $\sigma$ and
no-prime length margin $\delta$ so that their displayed terminal loss
fits, for example, $\epsilon/8$. The rank bound $D_\sharp$ then depends only
on these numbers and the aggregate length bound. Let $C_*$ be an
upper bound for the sums of absolute coefficients in the fixed lists
of affine norm-exponent expressions used in the two transforms,
projection, divisor allocation, ordinary reflection and restoration by initial segments. These lists
involve only the displayed aggregate norms. Their dimensions and
coefficients are independent of $F,S,\Theta$ and of the number of
prime-polynomial factors. Appendix~\ref{app:numerical-edge-budget} gives an explicit conservative
choice $C_*=2^{20}$ and its aggregate-expression budget. Choose
\[
 \eta_{\rm num}<\min\{\sigma/6,\epsilon/(512C_*D_\sharp)\}.
\]
This is the only numerical restriction on the maximal individual
prime length: restoring the length inequality removes an initial
segment and pays at most one overshooting factor of length
$\eta_{\rm num}$.

Now fix the arithmetic data, profiles and number of factors. All
fixed support windows contribute through the single quantity
$\theta_N=H_N/\log Z$. Choose $\xi$ and the input sieve losses
$\epsilon_0$ within the remaining numerical margins, and then enlarge
the threshold so $\theta_N<\xi$. In particular $H_N$ may depend on
the number of factors without changing $\eta_{\rm num}$. Fixed
conductor, discriminant and local-support constants are also absorbed
at this threshold. The common-extra-modulus exponent remains within
its stipulated numerical bound.

The other arithmetic dependences have the following
forms. The quadratic and cubic sieve exponents are respectively
$K+N$ and $K+N+(KN)^{2/3}$; their field dependence is in constants.
The Hecke functional equation has the same conductor exponent
$1/2-s$ in every degree; degree affects its polynomial height order
and the fixed discriminant constant. The buffered prime bound has
norm exponent $\kappa$ in every fixed field; degree again affects
only height orders and constants.

The norm-sheet kernels contribute
$(1+\log Z)^{d_F}$, which is absorbed by any already reserved positive
input loss after a data-dependent threshold. Finite unit groups,
finite character expansions, fixed local test families and their
multiplicities contribute fixed constants. A divisor power or a
fixed constant per prime contributes $Z^{\epsilon_0}$ by the ideal
divisor bound, with a data-dependent constant, not a multiple of
$\eta_{\rm num}$.

Radial and ordinary reflected profiles require finite Euler and
Fourier orders that may increase with the field and local data.
They are fixed backward through the numerical finite rank bound.
All Fourier boxes use actual fixed-ratio fine scales; their weighted
Fourier moments have finite seminorm and polynomial height bounds.
Increasing an order therefore introduces no factor $Z^{J\eta_{\rm num}}$.
Compact-dual support constants and small-norm logarithmic constants
are fixed after those orders and absorbed by the later threshold.
Thus none of these choices imposes a new upper bound on the numerical
prime mesh. External height cutoffs may be chosen still later in
terms of the resulting polynomial height orders.

Hence the numerical finite-depth loss argument applies
with this $\eta_{\rm num}$ to every subsequently fixed arithmetic
datum, which is the asserted order of quantifiers.
\end{proof}

\begin{remark}[The fourth-moment specialization]\label{v:rem:sextic-fourth-moment}
For $z=0$, take the basic presentation
$\psi_k(n)=\nu(n)\chi_n(k)$ with fixed $\nu$ and fixed extra modulus
($q=0$ in the length normalization). With $q_k\le U$ and
$X=U^{1/2}$, identical ordinary factors give
\[
 \sum_{\substack{q_k\le U\\\psi_k\text{ nonprincipal}}}
 \left|X^{-1/2}\sum_n\psi_k(n)W(q_n/X)\right|^4
                         \ll U^{1+\epsilon}.
\]
The same $z=0$ clause of Theorem~\ref{v:thm:M2-assembled} gives, for
plain lengths $X=U^m$ with $0\le m\le1/2+r_F/100+\epsilon$, two profiles
$W_1,W_2$ from the fixed family and norm twists of height at most $H$,
the bound $U^{1+\epsilon}(1+H)^{A}$; this is the form used in
Proposition~\ref{v:prop:refined-numerator}.
This strong sextic-family estimate is a consequence of the M2 argument, rather than an application of the published cubic
large sieve alone. Its cancellation uses the exact correlations,
ordinary reflection and centering proved above.

It does not give a natural-count cubic bound by restricting to square
rows: $\chi_n(u^2)=\chi_n(u)^2$, but $q_{u^2}=q_u^2$, so
$q_u\le V$ in this restriction gives only the bound $V^{2+\epsilon}$.
Nor is the local exceptional calculation unchanged for order three.
The sextic congruence $2+v_p(h')\equiv0\pmod6$ forces
$v_p(h')\ge4$, whereas modulo three it permits $v_p(h')=1$; the
residual power count changes from exponent $1/6$ to $1/3$.
Hence the sextic exceptional estimates do not give
a cubic proof; whether a modified argument using the savings from
both centered factors yields a cubic fourth moment requires a
separate computation and is not decided here.
\end{remark}

\section{The compensated comparison and continuation}
\label{v:sec:compensated-assembly}
In this section, we combine the moments of
Theorems~\ref{v:thm:marked-inverse-new} and \ref{v:thm:M2-assembled}
with the completed theta energy and additive estimate, with one common
principal term \eqref{v:cmp:8}; the length and count exponents are those of
Sections~\ref{sec:iv-count}--\ref{sec:iv-iteration}. The direct bound
\eqref{v:cmp:4} and the high-side estimate \eqref{v:cmp:13} first give
$\beta_*(F)\le7/8$; then \eqref{v:cmp:14}, combined with the fixed-point
identity of Section~\ref{sec:iv-iteration} in a supremum argument equivalent
to that iteration, proves Theorem~\ref{v:thm:original-target},
$L_F(s,\eta)\ne0$ for $\Rea s>B_\infty$. The refinement of
Section~\ref{v:sec:refined} reuses this assembly.

\subsection{Arithmetic data and moment estimates}
\label{sec:v-detail-26}
In this subsection, we fix the arithmetic data (the field $F$, the set $S$,
the local weights supplied by Proposition~\ref{v:prop:checked-inputs} and the
profile bound \eqref{eq:v-profile-bound}) and write
Theorems~\ref{v:thm:marked-inverse-new} and \ref{v:thm:M2-assembled}
in the row normalization $m=1$, in which the rows $u$ of the detector
(Section~\ref{sec:prev-average}, \Prev{Section}{5}) have $q_u\asymp U$ and all
lengths are exponents relative to $U$, retaining the length
conditions, common character, zero extensions, and the short-prime-sum and
zero-free hypotheses of the product moment. These moments feed the count in
Section~\ref{sec:v-detail-29}.

Fix a cyclotomic $F\supset\mu_{12}$ and a final finite $S$ containing the
archimedean places, primes over $6$, the support of the different, exceptional places of the theta representation,
fixed coefficient conductors, and primes making $R=\mathcal O_{F,S}$
principal. Write $q_a$ for exterior ideal norm, $G=F_S^1$,
$\Gamma=(R^\times)^6$. Use the source's multiplicative allowed generators,
symmetric sextic reciprocity, and self-dual additive lattice $R$ of
covolume one. Every symbol power, even power zero, is zero on nonunits.
All row sums retain finite unit classes modulo $\Gamma$.
The objects below are displayed for these fixed final data.
Their construction follows Remark~\ref{v:rem:parameter-ledger}:
the common numerical parameters precede the target, the target
precedes the final choice of $S$, and all $S$-dependent tests and
characters are constructed afterward.

Proposition~\ref{v:prop:checked-inputs} supplies
\[
 W_p(p^{3h})=q_p^{-h},\quad
 W_p(p^{3h+1})=q_p^{-h-1/2}\gamma_2(p),\quad W_p(p^{3h+2})=0,
\]
with unit factor $\chi_p(u)^{2e}$ at valuation $e$. Let $\Lambda$ be its
fixed unitary Eisenstein character and $G_0$ its fixed finite Gauss correction:
\[
 a_0(n)=\bar\Lambda(n)\gamma_2(n),\quad
 a_0(n)\gamma_1(n)=\mu(n)G_0(n),\quad
 \mu(n)\gamma_{-1}(n)=a_0(n)\bar G_0(n).
\]
$\Lambda$ is not of finite order; its angular factors remain in the chosen
theta functions. The raw reflection admits compact radial tests with
\begin{equation}\label{eq:v-profile-bound}
 |(x\partial_x)^j\Phi_W(x)|\ll p_J(W)\min(x^{1/4},x^{-A}).             
\end{equation}
This follows from the actual rotated first-angular vectors and multiplier
$E(t)A^t[\Gamma(5/6+t)\Gamma(7/6+t)/(\Gamma(5/6-t)\Gamma(7/6-t))]^r$,
$r=[F:\mathbb Q]/2$, with fixed finite factors holomorphic for
$\Rea t>-\delta$, any fixed $\delta<1/3$. Stirling cancels exponential
factors; the vertical power is $4r\Rea t$. Choosing $1/4<\delta<1/3$
gives \eqref{eq:v-profile-bound}, with finite seminorm and polynomial norm-twist orders.

For one common finite-order presentation $\psi_k=\nu\chi_\bullet(k)^{\pm1}$,
let $M_r,S_m$ be centrally normalized smooth annular inverse and ordinary
polynomials, and $Q$ a product of centrally normalized prime sums of this
same presentation. In M1 the actual disjoint prime lists may vary arbitrarily
within their annuli; coefficients are bounded and independent of row and
other columns. In M2 the sums are over full exterior prime annuli with fixed
smooth profiles and coefficients in a fixed finite $\Theta$-linear span,
including identity-ray selectors only via their character expansion.
Arbitrary prime subsets or arbitrary bounded M2 coefficients are forbidden.

We apply Theorems~\ref{v:thm:marked-inverse-new} and
\ref{v:thm:M2-assembled}, in row normalization \(m=1\):
\begin{align}
 r+2z<1,\quad2r+8z<3
 &\Longrightarrow\sum_{q_k\ll U}|M_rQ|^2\ll U^{1+\eps},    \label{eq:v-M1-row}\\
 z>0,\quad n_1+n_2+6\kappa z<1
 &\Longrightarrow\sum_{q_k\ll U}|S_{n_1}S_{n_2}Q|^2\ll U^{1+\eps}.\label{eq:v-M2-row}
\end{align}
We refer to \eqref{eq:v-M1-row} as M1 and to \eqref{eq:v-M2-row} as M2.
These are the only moment inputs of the count in
Section~\ref{sec:v-detail-29}. Their hypotheses are:
\begin{itemize}
\item one common finite-order presentation $\psi_k$;
\item centrally normalized smooth annular $M_r$, $S_m$;
\item prime sums $Q$ of total length $z$ over the fixed windows, with
the coefficient restrictions stated above (M1: arbitrary prime lists
within the annuli, with bounded row-independent coefficients; M2: full
annuli with coefficients in the fixed span of $\Theta$);
\item rows $q_k\ll U$, unit classes retained, and principal inducing
rows excluded in M2;
\item for M2 with $z>0$, the fixed mesh $z_i\le\eta$ and the zero-free
premise $\beta_*(F)\le(1+\kappa)/2$;
\item finite seminorm and polynomial height orders.
\end{itemize}
In strength, M1 with $r$ near $1$ is an optimal large sieve for M\"obius
coefficients in the balanced regime; as for
Theorem~\ref{v:thm:unmarked-inverse} (Remark~\ref{v:rem:inverse-strength}),
it cannot follow from bounds for arbitrary coefficients. M2 with $z=0$
and $n_1=n_2\to1/2$ is a Lindel\"of-on-average fourth moment for the
sextic family.
For $z=0$, M2 holds for arbitrary bounded ordinary lengths $n_1,n_2$
with no restriction on the prime length and excludes principal inducing rows; selected M2 excludes rows
whose inducing character belongs to $\Theta$. In the present comparison
these exceptions are fixed unit classes among the sixth-power-free
comparison rows $u$, because the coefficient group is supported in the
final fixed $S$; this concerns $u$, not all integral rows
$k=u a^6$, which are estimated separately.

M2 with short prime sums requires a fixed mesh $z_i\le\eta$,
$\beta_*(F)\le(1+\kappa)/2$, and its specified fixed common twists and
zero extensions. Both estimates include every required finite test
seminorm and polynomial height order, common orientations, and surviving
subcollections; M1 has two strict margins, and the prescribed exponent of
every growing initial coprimality modulus is retained.
The inverse estimate without prime-polynomial factors also gives the
scale-supremum moment for $H\ge D^{1+c}$, any fixed $c>0$.

\subsection{The theta average and its direct bound}
\label{sec:v-detail-27}
In this subsection, we define the theta average from the local coefficients
of Section~\ref{sec:v-packet}, fix the geometric parameters and the
nonprincipal margins, and form the compensated operation \eqref{v:cmp:1} with
the selected prime factors in its additive weight.
Lemma~\ref{v:lem:gram-character-sums} gives the nonprincipal coprime
two-column sums, and the mean-square bounds \eqref{v:cmp:2} and
\eqref{v:cmp:3} give the direct estimate
$I_{\rm mod}\ll Z^{l_x/2+b/12+\eps}$ of \eqref{v:cmp:4}, the low side of the
comparison in Section~\ref{sec:v-detail-31}.

Put $\vartheta=\bar\Lambda\bar G_0\eta$ for a fixed target $\eta$, and
\[
 B_k(Z)=\sum_{c\ \sfq,n}\frac{\gamma_2(c)\vartheta(cn^3)}{\sqrt{q_c}q_n}
          \chi_{cn^3}(k)P_G(q_{cn^3}/Z),\quad
 P_G(x)=\frac{e^{-(\log x)^2/4}}{2\sqrt\pi}.
\]
The factors $c,n$ may overlap. Choose actual tests $W_1\ge0$ and $W_0$
compactly supported in $F_S^\times$ with $\widehat W_0(0)=0$ and
$\mathcal M(1/6)\ne0$, where
$\mathcal M(z)=\int\widehat W_0(y)|y|_S^z d^\times y$.
These remove original and dual zero by different mechanisms.

For existence, take a nonnegative nonzero annular tensor $\phi$, $a$ of
local norm $A>1$ at one complex place, and replace that factor by
$\phi_v(y)-A^{-1}\phi_v(y/a)$. Its additive integral is zero, and
\[
 \mathcal M_{W_0}(z)=(1-A^{-z})\mathcal M_\phi(z).
\]
By the local Tate functional equation, $\mathcal M_\phi(1/6)$ is a
product of nonzero local gamma constants times
$\int\phi(y)|y|_S^{5/6}d^\times y>0$, so $\mathcal M(1/6)\ne0$.
Absolute local convergence for $\Rea z>0$ and integration by
parts in a complex logarithmic norm give arbitrary vertical decay on
compact positive strips. The signed test is allowed in the
theta average and its Cauchy bound.

Define
\[
 I_\eta(X,Y,Z)=\int_{G/\Gamma}\frac1Y\sum_sW_1(q_s/Y)
 \frac1{\sqrt{q_sX}}\sum_{k\in R}q_s^{-1/2}g_s(-k)
 W_0(k/(t\ell_{q_sX}))B_k(Z)\,dt,
\]
$g_s(k)=\sum_{d\bmod s}\chi_s(d)e(kd/s)$.
The auxiliary ideal coefficient is $1$; no $\Lambda$ angular character is
inserted in this additive Gram sum.

Fix geometric parameters with
\[
 l_x+l_y+\ell=1,\quad b=l_y-l_x,\quad X=Z^{l_x},\quad Y=Z^{l_y},
 \qquad h=1-l_x+\ell=(1+b+3\ell)/2.
\]
The nonprincipal margin is fixed before the selected principal error.
Write $\beta=\beta_*(F)$; in the contradiction argument fix $B<\beta$
and the geometric parameters. Choose
$\mathfrak r_{\rm np}>0$ and the small/far row cutoffs so that, before
buffer and discretization losses, every nonprincipal exponent in
Sections~\ref{sec:v-detail-29}--\ref{sec:v-detail-31} is at most
$C(\beta)-8\mathfrak r_{\rm np}$.

Let $\mathscr E$ be the finite
list of scalar norm-exponent functions in the count of detected elements and
contour estimates of Sections~\ref{sec:v-detail-29}--\ref{sec:v-detail-31},
with one function for each branch of a maximum or positive part.
Let $\mathscr D$ be their closed buffered numerical parameter ranges,
with the already fixed $d\ge1/100$ and positive denominators retained.
For each function use its aggregate loss and buffer variables as
coordinates; prime subsets enter through their total lengths, without
a separate coordinate for each prime position. Put
\[
 C_{\rm np}=100\left(1+\sum_{E\in\mathscr E}
                 \sup_{x\in\mathscr D}\|\nabla E(x)\|_1\right).
\]
The rational and affine functions have bounded derivatives on these
ranges, and maxima and positive parts retain the largest branch
Lipschitz bound; so $C_{\rm np}<\infty$ depends on the numerical
comparison ranges only, not on the target, field, prime factor count or
profile orders.

Choose
\[
 \begin{gathered}
 0<\epsilon_{\rm np}<\frac{\mathfrak r_{\rm np}}{100C_{\rm np}},\\
 \eta_{\rm M2}=\min\{\eta_{\rm num}(\epsilon_{\rm np}/100),\epsilon_{\rm np}/100\},\\
 0<\eta_{\rm phys}<\min\left\{
 \frac{\eta_{\rm M2}}{100},
 \frac{\mathfrak r_{\rm np}}{10^4C_{\rm np}}\right\}.
 \end{gathered}
\]
Fix also a length allowance
$0<\epsilon_c<\epsilon_{\rm np}/(100C_{\rm np})$.
The detector's preliminary exponent losses will be $<\epsilon_c/100$.
These choices are numerical (Lemma~\ref{v:lem:arithmetic-independent-mesh}
is applied with the already fixed numerical length and modulus bounds)
and precede $K$, the target character and its arithmetic data; the moment
and counting losses stay within this fixed $\epsilon_{\rm np}$ allowance
and are not chosen again after $K$.

Choose a fixed $K>2\ell/\eta_{\rm phys}$ and put
\[
 \ell_i=\frac\ell K\left(1+\frac{2i-K-1}{4K}\right),\qquad
 \sum_i\ell_i=\ell.
\]
Then \(0<\ell_i<\eta_{\rm phys}\) are distinct, and for \(U=Z^d\),
\(d\ge1/100\), \(z_i=\ell_i/d<\eta_{\rm M2}\). The full norm windows are
disjoint for large \(Z\); the ray selector is expanded into
\(\Theta\). The greedy rounding loss from selecting whole prime-sum
factors uses the same discretization of the prime-length exponents and
lies inside the fixed margin. \(K\) is independent of
\(Z\); constants and onsets may depend on it.

Choose nonnegative nonzero smooth prime profiles supported in $(1,2)$.
Take disjoint identity-ray windows $P_i=Z^{\ell_i}$, $\sum\ell_i=\ell$; as in
Section~\ref{sec:prev-average}, the allowed set $\mathcal P_i$ consists of the
exterior primes in the $i$th window whose class in the fixed ray quotient
$\mathcal T$ is the identity class $1_{\mathcal T}$.
With $I_{\eta;D}$ imposing the divisibility condition $D\mid cn^3$ on the full theta index, the
finite compensated operation is
\begin{equation}
 I_{\rm mod}=\sum_{(p_i)}\prod_iW_i(q_{p_i}/P_i)
 \sum_J(-1)^{|J|}q_{p_J}^{-3/2}\bar\eta(p_{J^c})
 I_{\eta;p_{J^c}}(X/q_{p_J},Y/q_{p_J},Zq_{p_{J^c}}).                 \label{v:cmp:1}
\end{equation}
Every term has the same windows and the geometric parameters fixed above.

\paragraph{\it Explicit form of the operation.}
In \eqref{v:cmp:1} the outer sum runs over tuples
$\mathbf p=(p_1,\dots,p_K)$ with $p_i\in\mathcal P_i$; the windows are
disjoint, so the $p_i$ are distinct. For $J\subseteq\{1,\dots,K\}$,
$p_J=\prod_{i\in J}p_i$, $q_{p_J}$ is its norm and $J^c$ is the
complement. $I_{\eta;D}(X,Y,Z)$ is $I_\eta(X,Y,Z)$ with the summation
over $(c,n)$ in $B_k(Z)$ restricted to $D\mid cn^3$; the factor
$\bar\eta(p_{J^c})$ accompanies this restriction, as in the local
verification following \eqref{eq:v-global-compensated-correction}.
The sum over $J$ is the inclusion--exclusion of the divisibility
conditions at the selected primes, and \eqref{v:cmp:1} is the finite
operation whose Euler form is the tuple identity
\eqref{eq:v-global-compensated-correction}. For $K=0$ it is the
uncompensated average $I_\eta(X,Y,Z)$. The data
entering \eqref{v:cmp:1} are fixed in the order of
Remark~\ref{v:rem:parameter-ledger}: numerical parameters and margins,
then $K$ and the window lengths, then the principal saving and buffer,
then the target and final $S$, and only then the tests, ray selectors
and coefficient groups. Internal orders precede the height cutoff;
the external tail order and the onset are chosen last.

\begin{lemma}[Nonprincipal sums and the common-divisor expansion]
\label{v:lem:gram-character-sums}
Let $\lambda$ be a finite-order Hecke character over the fixed field,
extended by zero at the primes of an additional ideal. Suppose its
primitive inducing character is nonprincipal, and the product of its
primitive conductor with the additional distinct unramified primes has
norm at most $C_0Q_0$, where $C_0$ is fixed and $Q_0\ge1$.
For annular smooth $U$ and $N\ge c_U>0$, every fixed $A\ge2$ gives
\[
 \left|\sum_a\lambda(a)U(q_a/N)\right|
 \ll (Q_0N)^\epsilon\mathfrak S_J(U)
            N\min\{1,(Q_0/N)^A\}.
\]
The same estimate holds with polynomial dependence on a norm-twist
height. Consequently, for two annular tests $U,V$ and $N_1,N_2\ge c>0$,
\begin{equation}\label{eq:v-gram-coprime-mean}
 \begin{aligned}
 &\left|\sum_{(a,b)=1}\lambda(a)\overline{\lambda(b)}
                    U(q_a/N_1)V(q_b/N_2)\right|\\
 &\qquad\ll (Q_0N_1N_2)^\epsilon\mathfrak S_J(U,V)
       N_1N_2\min\{1,Q_0/\min(N_1,N_2)\}.
 \end{aligned}
\end{equation}
Fixed coprimality restrictions may be retained in the character's
zero extension. A ray-class restriction must instead be expanded into
characters, with every resulting principal inducing twist treated separately.
\end{lemma}
\begin{proof}
Write the Dirichlet series as the primitive Hecke function multiplied
by its finitely many removed Euler factors. The primitive functional
equation and absolute convergence on $\Re s=A$ give a polynomial
height bound on $\Re s=1-A$ with conductor factor
$q_{\mathfrak f}^{A-1/2}$. The removed factors on that line contribute
at most the $(A-1+\epsilon)$th power of the norm of their radical,
by the product estimate \eqref{eq:v-removed-euler-uniform} applied
after extracting these powers. Thus their combined conductor factor
is $\ll Q_0^{A+\epsilon}$. Mellin inversion and a shift from $2$
to $1-A$ cross no pole, since the inducing character is nonprincipal.
Sufficiently many test derivatives give
$Q_0^{A+\epsilon}N^{1-A}$, with the stated seminorm and height factors;
ideal counting gives $O(N^{1+\epsilon})$. On bounded subunit scales use the original annular
support; a nonempty sum has $N\ge c_U$.

For \eqref{eq:v-gram-coprime-mean}, expand the entire identity
$1_{(a,b)=1}=\sum_{d\mid a,b}\mu_F(d)$. Multiplicativity gives
\[
 \sum_d\mu_F(d)|\lambda(d)|^2
       \left(\sum_a\lambda(a)U(q_a/(N_1/q_d))\right)
       \left(\sum_b\overline{\lambda(b)}V(q_b/(N_2/q_d))\right).
\]
With $T=\min(N_1,N_2)$, the rapid bound on the smaller scale and the
absolute bound on the larger scale bound this by a seminorm factor times
\[
 (Q_0N_1N_2)^\epsilon N_1N_2
 \sum_dq_d^{-2}\min\{1,(Q_0q_d/T)^A\}.
\]
Only nonempty annuli need be included. Ideal counting and a split at
$q_d=T/Q_0$ show, for $A>1$,
\[
 \sum_dq_d^{-2}\min\{1,(Q_0q_d/T)^A\}
             \ll\min\{1,Q_0/T\}.
\]
For $T/Q_0\ge1$, the sums below and above that point are both
$O(Q_0/T)$; for $T/Q_0<1$, use convergence of $\sum_dq_d^{-2}$.
Thus the complete common-divisor expansion costs one power of $Q_0/T$;
arbitrary rapid decay does not survive it.
\end{proof}

\begin{remark}[Constants in Lemma~\ref{v:lem:gram-character-sums}]
\label{v:rem:gram-constants}
Since only one power of $Q_0/T$ survives the divisor expansion, any
fixed $A>1$ suffices in the single-sum bound. In particular the
functional equation is not needed for \eqref{eq:v-gram-coprime-mean}.
Fix $0<\delta\le1/10$. Phragm\'en--Lindel\"of between the
conductor-uniform bounds on $\Re s=-\delta$ and $\Re s=2$ gives the
convexity bound
\[
 |L(s,\chi)|\ll(1+N\mathfrak f)^{1/2+\delta}(2+|t|)^{k_F}
 \qquad(-\delta\le\Re s\le2)
\]
for primitive nonprincipal $\chi$, with $k_F$ depending only on $F$
and $\delta$. Shifting the Mellin integral from $2$ to $-\delta$ then gives
\[
 \Bigl|\sum_a\chi(a)U(q_a/N)\Bigr|
 \ll\mathfrak S_J(U)\,N\min\{1,(Q_0'/N)^{1+\delta}\},
\]
where $Q_0'=(1+N\mathfrak f)^{(1/2+\delta)/(1+\delta)}\le1+N\mathfrak f$.
Here $J=k_F+2$ derivatives of $U$ suffice for the vertical integral.
The divisor expansion then yields \eqref{eq:v-gram-coprime-mean} with
$Q_0=1+N\mathfrak f$, uniformly in $\chi$. For a zero extension at the
primes of an ideal $\mathfrak a$, M\"obius over the divisors of
$\mathfrak a$ replaces $Q_0'$ by $Q_0'N\mathfrak a$. It also multiplies
the bound by $\tau(\mathfrak a)^2\ll_\eta(N\mathfrak a)^{2\eta}$. A norm
twist $U(x)x^{-i\tau}$ multiplies $\mathfrak S_J(U)$ by at most
$(1+|\tau|)^J$, which is the polynomial height dependence stated in the
lemma.

\end{remark}

\begin{lemma}[Mean square of the completed theta sum and its additive factor]
For a fixed subset tuple of length $d$, let $X'=X/q_{p_J}$,
$Y'=Y/q_{p_J}$, $M'=1-\ell-2d$. Let $B_k^{(J)}$ denote the remaining
full-index theta sum with prime-divisor weights for the outer index $k$,
the counterpart of the rescaled theta sums of \Prev{Lemma}{16.1} recalled
in Section~\ref{sec:prev-free-direct}. It
satisfies
\begin{equation}
 \sum_{q_k\ll X'Y'}|B_k^{(J)}|^2
 \ll Z^{M'+\pos{d+5\ell-1}/4+\eps}.                              \label{v:cmp:2}
\end{equation}
Its actual averaged additive factor satisfies, $Q'=X'Y'$, $P_a=Y'^2/Q'$,
\begin{equation}
 \sum_{[k]\ne0}|\widetilde A_{[k]}|^2
 \ll\frac{Q'}{Y'}(1+P_a^{1/6}+P_a^2/Y')Z^\eps.                   \label{v:cmp:3}
\end{equation}
\end{lemma}
\begin{proof}
Common Gaussian localization gives completed index length
$N_*=1+\ell'+\theta_N$, $\ell'=\ell-d$, and row length
$M'=1-\ell-2d$. First retain $|\theta_N|\le\kappa_G+O(1/\log Z)$,
where $\kappa_G>0$ is a small Gaussian localization allowance;
the Gaussian removes farther annuli before any positive norm is taken.
All retained seminorms are summable. Allocate $\kappa_G$ and the
later radial tail allowance $\tau$ within $\epsilon$, and suppress these
errors in the exponent calculation below.

Write the row as $k=k_{\rm pow}R$, with $k_{\rm pow}$ powerful
and $R$ squarefree, of lengths $O,H$, respectively;
$O+H\le M'$. Freeze $k_{\rm pow}$ and its local branches.
Let $F_{\rm act}$ be the product of active primes of $k_{\rm pow}$.
At each valuation-four branch expand its Ramanujan term using
\[
 1_{p\mid n_{\rm full}b_{\rm full}^3}
 =1_{p\mid n_{\rm full}}+
       1_{p\nmid n_{\rm full}}1_{p\mid b_{\rm full}}.
\]
Let $N_F,B_F$ be the products assigned to the squarefree and cube
factors by this priority. Define
\[
 A_0=\log_Zq_{F_{\rm act}},\qquad
 N_0=\log_Zq_{N_F},\qquad B_0=\log_Zq_{B_F}.
\]
Let $S_0$ be the sum of the lengths of valuation-four small branches
and active zero-digit branches. The local absolute coefficients and
forced numerator factors are
\[
\begin{array}{c|c|c}
 \text{branch}&\text{absolute coefficient}&\text{forced factor}\\\hline
 \text{valuation-one row}&1&1\\
 \text{active zero digit or valuation-four small term}&q_p^{-1/2}&1\\
 \text{valuation-four squarefree term}&1&p\mid N_F\\
 \text{valuation-four cube term}&q_p^{-1/2}&p\mid B_F\\
 \text{other active digits}&\le1&1 .
\end{array}
\]
Thus $2A_0\le O$, $N_0\le A_0$, and all four lengths are
nonnegative.

The fixed factor $q_{b_*}$ in the source average and direct estimate
\cite[(4.5), Proposition~16.2]{LiuArXiv} is included here in the
finite-$S$ tests: a fixed dilation contributes $O(1/\log Z)$ to a
length exponent and a fixed power to the implied constant.
The finite $S$-part is already in the fixed reflected profiles of
Proposition~\ref{v:prop:radial}, so no growing ramified valuation index
occurs.
The remaining unit classes belong to the finite group $U_S/U_S^6$.

Fix the inactive selected prime factors, of coefficient $q_p^{-1}$,
before Hilbert triangle; their total harmonic mass is $Z^{o(1)}$.
The active tuple $P$ has length $z_a\le\ell'$ and contributes
$q_P^{-1/2}$ inside the hybrid polynomial. The reflected numerator
is $N_Fn(B_Fb)^3$. Write $v$ and $\ell_b$ for the retained dual length exponents of the
reflected numerator, attached to $n$ and $b$ respectively, as in the proof
of \Prev{Lemma}{16.1}; $T_d$ is the available dual length exponent. The profile argument is
\[
 \begin{aligned}
 \frac{Z^{N_*}q_{N_Fn(B_Fb)^3}}
      {(q_Rq_{F_{\rm act}}q_P)^2}
 &=Z^{v+3\ell_b-T_d}\times\text{a fixed-box norm monomial},\\
 T_d&=2H+2A_0+2z_a-N_*-N_0-3B_0.
 \end{aligned}
\]
For $\theta_N=0$ the length inequalities give
$T_d\le H-3d$; with the retained Gaussian annuli the error is
$O(\kappa_G+1/\log Z)$. The scalar outside the hybrid norm is
$Z^{-v/2-\ell_b-(S_0+B_0)/2}$, excluding the internal factor $q_P^{-1/2}$.

Lemma~\ref{v:lem:marked-terminal-typed} identifies the actual hybrid matrix as
$\chi_R(nb)^3\chi_P(n)^{-2}1_{(R,P)=1}1_{(P,b)=1}$;
$a(P)$ is independent of $R,n,b$ and $\beta(n,b)$ independent of $R,P$.
Its estimate follows by full coprimality indicator expansion, $b=gt^2$, $c=(n,g)$,
$n=cm,g=ch$, quadratic sieve on $mh$, then cubic sieve on the genuine
squarefree column $c'm$. No false squarefreeness of $b$ is used.
On retained dual dyads $v+3\ell_b\le T_d+\tau$, we have
$\max(H,v+\ell_b)\le H+O(\tau+\kappa_G+1/\log Z)$.

A powerful ideal is uniquely $x^2y^3$ with $y$ squarefree;
its count in this dyad is $O(Z^{O/2+\epsilon})$, incurred once.
Retaining the additional numerator length $e_\lambda\ge0$ of
\Prev{Lemma}{16.1}, the norm exponent is, as in its proof,
\[
 E(e_\lambda)=O/2+H-S_0-B_0+z_a-u-\ell_b
 -\frac{2e_\lambda}{3}
 -\frac12(T_d-v-3\ell_b-e_\lambda)_+.
\]
Since $x_+-(x-e)_+\le e$ for $e\ge0$,
\[
 E(e_\lambda)\le E(0)-e_\lambda/6\le E(0).
\]
Hence, apart from the allocated losses, it suffices to bound
\[
 E=E(0)=O/2+H-S_0-B_0+z_a-u-\ell_b-\tfrac12(T_d-v-3\ell_b)_+,
 \quad u=\min(v,z_a,(v+z_a)/3).
\]
For $u=z_a$, $E\le M'$. Otherwise $u\ge v/2$ and
$u+\ell_b+\tfrac12(T_d-v-3\ell_b)_+\ge T_d/4$.
With $\Delta_H=M'-O-H$,
\[
 O/2+\Delta_H+S_0+B_0+T_d/4
 =\frac{2M'+2z_a-1-\ell'+2\Delta_H}{4}
  +\frac{2A_0-N_0+B_0+4S_0}{4}.
\]
The last numerator is nonnegative, giving \eqref{v:cmp:2}. Genuine tails use \eqref{eq:v-profile-bound}
before positivity; powerful parts are counted once. All norm profiles
are common full-box densities before remaining primes. Shared cube primes stay.

For \eqref{v:cmp:3}, orbit unfolding gives
\[
 \widetilde A_{[k]}=Y'^{-1/2}\sum_s q_s^{-1}g_s(-k)W_1(q_s/Y')
             \Omega(q_k/(q_sX')),\quad
 \Omega(T)=\int_G W_0(t\ell_T)\,dt.
\]
Its nonzero terms have $q_s\asymp Y'$ and $q_k\ll Q'=X'Y'$.
Separate the joint profile on fixed boxes in $q_s/Y'$ and $q_k/Q'$
before taking the row norm; the row norm phases then cancel in each
positive square. Hilbert-space triangle inequality, with the weighted
Fourier measure of \eqref{v:eq:H2-common-density}, reduces the estimate
to fixed annular column profiles and an unrestricted positive row ball.
Only at this point insert a fresh band-limited majorant for that ball.
Its dual norm support is compact by Lemma~\ref{v:lem:compact-dual}.

Write each separated component as
$Y'^{-3/2}\sum_s c_sg_s(-k)$, where $c_s$ is $Y'/q_s$ times
its bounded annular column profile. The row majorant is nonnegative;
include its (nonnegative) original zero-row term $\mathcal Z_\Phi$ in the
upper bound. Finite-correlation Poisson gives
\[
\begin{aligned}
 &\sum_{[k]\ne0}\Psi_\Phi(q_k/Q')
       \left|\sum_s c_sg_s(-k)\right|^2+\mathcal Z_\Phi\\
 &=Q'\operatorname{vol}(G/\Gamma)\widehat\Phi(0)
                       \sum_s|c_s|^2\varphi(s)\\
 &\quad+Q'\sum_{s,t}c_s\bar c_t
       \sum_{[j]\ne0}\mathcal F(s,t;j)
              \Psi_{\widehat\Phi}\!\left(\frac{Q'q_j}{q_sq_t}\right).
\end{aligned}
\]
Changing $j$ to $-j$ if required by the Fourier convention preserves
the norm kernel. At the original origin,
$\mathcal Z_\Phi=\operatorname{vol}(G/\Gamma)\Phi(0)
|\sum_s c_sg_s(0)|^2\ge0$.

The squared normalization gives the outside factor $Q'/Y'^3$.
Compact Fourier support gives
$q_j\ll q_sq_t/Q'\ll P_a=Y'^2/Q'$.
At zero frequency the correlation vanishes unless $s=t$, when
it is $\varphi(s_1)\le q_{s_1}$; with $O(Y')$ columns this term is
$O(Q'/Y')$. Its infinite factor is the original compact quotient
integral, not an evaluation of an unfolded kernel.

For a nonzero frequency extract complete common support:
\[
 s_1=Da,\quad s_2=Eb,\quad
 (a,b)=1,\quad(ab,DE)=1,\quad
 \operatorname{rad}(D)=\operatorname{rad}(E).
\]
Put $C=(D,E)$. Then
\[
 \mathcal F(Da,Eb;j)=\mathcal F(D,E;j)
                         \chi_a(j)\overline{\chi_b(-j)},
 \qquad C\mid j,\qquad |\mathcal F(D,E;j)|\le q_C.
\]
For fixed $C$, the possible $D,E$ have support in
$\operatorname{rad}(C)$ and norms $O(Y')$; by Rankin, for every $\delta>0$,
\[
 \#\{D:q_D\ll Y',\ \operatorname{rad}(D)\mid\operatorname{rad}(C)\}
 \ll_\delta Y'^\delta
          \prod_{p\mid C}(1-q_p^{-\delta})^{-1}
 \ll_\delta (Y'q_C)^\delta,
\]
after choosing smaller loss shares in the last inequality; the same
holds for the pair $(D,E)$. Their valuations are retained; no false
squarefreeness is imposed.

Fourier-separate the retained norm kernel on the common boxes before
fixing the residual columns. Their scales are $N_1\asymp Y'/q_D$,
$N_2\asymp Y'/q_E$. Reciprocity identifies their common character
as $\lambda_{\epsilon,j}(a)=\chi_a(\epsilon j)
=\xi_\epsilon(a)\chi_a(j)$, with its natural zeros, after writing
the dual representative in its unit class as in
\eqref{eq:v-poisson-unit-normalization}. The finite unit factor
$\xi_\epsilon\in\Upsilon_S$ is retained in the common coefficient
group. For the full frequency $\epsilon j$,
\[
 \mathcal F(Da,Eb;\epsilon j)=\mathcal F(D,E;\epsilon j)
        \chi_a(\epsilon j)\overline{\chi_b(-\epsilon j)}.
\]
Since $\chi_b(-1)=1$, the residual factors are
$\lambda_{\epsilon,j}(a)\overline{\lambda_{\epsilon,j}(b)}$, and
$\mathcal F(D,E;\epsilon j)$ is an outer scalar. This gives the common
orientation for the original completed additive factor, including the
unit sector and all zeros.

The fixed $S$-part is retained in the
primitive conductor and all extra exterior common primes divide $j$, so
the primitive conductor times the extra unramified radical is
$O_{F,S}(q_j)$. The tests may acquire different norm twists, with the
polynomial height factors integrated against the common Fourier measure.

If this inducing character is nonprincipal, apply
\eqref{eq:v-gram-coprime-mean}. Since
$\max(N_1,N_2)\ll Y'/q_C$, the absolute coprime two-column sum is
\[
 \ll Z^\epsilon N_1N_2
       \min\{1,q_j/\min(N_1,N_2)\}
 \ll Z^\epsilon\frac{q_jY'}{q_C}.
\]
Multiplying by $|\mathcal F|\le q_C$ gives $\ll Z^\epsilon q_jY'$. For $q_C\asymp C_0$,
write $j=Ck$ and use
\[
 \sum_{q_k\ll P_a/q_C}q_{Ck}\ll P_a^2/q_C,
 \qquad \#\{C:q_C\asymp C_0\}\ll C_0.
\]
Thus each $C$ dyad contributes $O(Z^\epsilon Y'P_a^2)$ before the
outside factor, and summing dyads gives
\[
 \frac{Q'}{Y'^3}\,O(Z^\epsilon Y'P_a^2)
       =O\left(\frac{Q'}{Y'}\frac{P_a^2}{Y'}Z^\epsilon\right).
\]

For a principal inducing character, tame Kummer ramification forces
$v_p(j)\equiv0\pmod6$ at every exterior prime. Any fixed $S$-supported
twist leaves this necessary condition unchanged. With $j=Ck$, the
valuations of $k$ are prescribed modulo six at the primes of $C$,
and divisible by six elsewhere. Hence, for fixed $C$,
\[
 \#\{k:q_{Ck}\ll P_a,\ \lambda_{Ck}\text{ principal}\}
       \ll_{F,S}Z^\epsilon(P_a/q_C)^{1/6}.
\]
(write the ideal of $k$ as its prescribed sixth-power-free part times a
sixth power; the finite unit classes supply a fixed factor). Absolute
column counting gives $q_C N_1N_2\ll Y'^2/q_C$, so a dyad
$q_C\asymp C_0$ contributes
\[
 \ll Z^\epsilon Y'^2(P_a/C_0)^{1/6}.
\]
Since $\sum_{C_0\ \mathrm{dyadic}}C_0^{-1/6}\ll1$, multiplication by
$Q'/Y'^3$ gives $O((Q'/Y')P_a^{1/6}Z^\epsilon)$. With the zero-frequency
term this proves \eqref{v:cmp:3}. All estimates concern the same original
additive factor; the Fourier measures were fixed before the arithmetic
indices, and finite-$S$ tests were transformed rather than treated as
exterior-periodic functions.
\end{proof}

\begin{proposition}[Direct estimate for the compensated theta average]
\label{v:prop:direct-estimate}
Fix $F$, $S$, the target $\eta$, the tests $W_0,W_1,W_i$, the windows
and $K$ as above, and geometric parameters with $\ell\le1/5$, $b\ge0$,
$l_y\ge11b/6$ and $l_x>\ell$. Then for every $\eps>0$
\[
 I_{\rm mod}\ll Z^{l_x/2+b/12+\eps}\qquad(Z\ge Z_0),
\]
where the implied constant and $Z_0$ may depend on $F$, $S$, $\eta$,
the tests, $K$ and $\eps$, and not on $Z$. With the normalizer
$c_SA_{\mathcal T}$ of Section~\ref{sec:v-detail-28} and
$A_{\mathcal T}^{-1}=Z^{o(1)}$, this is $J_\eta\ll Z^{C(B)+\eps}$ when
$l_x/2+b/12=C(B)$. It is the first of the two inputs to the bound for
the principal triple integral (Remark~\ref{v:rem:principal-bound-split}).
\end{proposition}
\begin{proof}
As in the proof of \Prev{Proposition}{16.2},
Cauchy for the exact pairing $Q'^{-1/2}\sum\widetilde A B$, tuple count
$Z^d$ and coefficient $Z^{-3d/2}$ in \eqref{v:cmp:1} give
\[
 L(d)=l_x/2-d+\pos{d+5\ell-1}/8
                     +\tfrac12\max(0,b/6,2b-l_y+d).
\]
Each positive-part/max term is $1$-Lipschitz, so
$L(d)-L(0)\le-3d/8$ and
\begin{equation}
 I_{\rm mod}\ll Z^{l_x/2+b/12+\eps}
 \quad(\ell\le1/5,\ b\ge0,\ l_y\ge11b/6,\ l_x>\ell).                       \label{v:cmp:4}
\end{equation}
Here $M'=1-\ell-2d\ge1-3\ell>0$ for $0\le d\le\ell$, and $l_x>\ell$ by
\eqref{eq:iv-geometry} and Lemma~\ref{lem:iv-reserves}.
This uses direct bound reflection, not a weaker canonical $Z^F$
energy substituted into the direct bound for the original theta average.
The dependence of the implied constant on $\eta$ enters only through
the mean-square bounds \eqref{v:cmp:2} and \eqref{v:cmp:3}, whose
constants depend on the conductor of $\eta$ and on $S$.
\end{proof}

\subsection{Exact correction and principal term}
\label{sec:v-detail-28}
In this subsection, we compute the Euler correction of the theta average of
Section~\ref{sec:v-detail-27}, selected-prime factors included: the local
factors \eqref{v:cmp:5}, the global correction
\eqref{eq:v-global-compensated-correction}, the Mellin representation
\eqref{v:cmp:6}, the principal holomorphy region \eqref{v:cmp:7} and the
nonzero principal term $f_\eta$ of \eqref{v:cmp:8}. The same correction
enters Section~\ref{sec:v-detail-30}, so both sides of the final comparison
use one normalization.

For sixth-free $u$, let $j=v_p(u)$, $\rho=\chi_p(u/p^j)$,
$T=q_p^{-s}$, $Y_0=q_p^{-w}$, $V=q_p^{-6z}$,
$a=\bar\Lambda(p)^3\eta(p)^3$, $b_0=\bar\Lambda(p)^2\eta(p)^2\bar\chi_p(4)$,
$R=a^2q_p^4T^6V$, and $\eta=\eta(p)$ in the local formulas below. The
local series is
\[
 P_j^*=\frac1{1-R}\left\{\frac{R(1-q_p^{-1})-\eta(p)(q_p-1)TY_0V^{1_{j\le1}}}{1-V}+J_j\right\},
\]
\[
\begin{array}{c|l}
0&J_0=-\eta\rho^{-1}T+\rho Y_0R\\
1&J_1=\eta TY_0\\
2&J_2=a\rho^{-3}q_p^{3/2}T^3\\
3&J_3=-\eta b_0\rho^{-2}q_p^2T^3Y_0+b_0^2\rho^{-4}q_p^2T^4\\
4&J_4=-\eta a\rho^{-3}q_p^{5/2}T^4Y_0\\
5&J_5=-a^2q_p^3T^6 .
\end{array}
\]
$P_j=(1-V)^{-1}+1_{j=0}\rho Y_0/(1-\rho Y_0)+P_j^*$, retaining all powers
of the target character. Write $\gamma_v=\tau(\chi_p^v)/\sqrt q$ for $v\not\equiv0\pmod6$,
$\tau_0=-1$, and put
\[
 \vartheta=\bar\Lambda\bar G_0\eta,\qquad
 t=e_0+3l,\quad e_0\in\{0,1\},\qquad J=j+6m.
\]
The finite correlation in the Poisson summation over the outer variable,
which leads to \eqref{v:cmp:6} below, has local values
\[
 \begin{aligned}
 C_{0,0}(J)&=1,& C_{0,k}(J)&=1_{J=0}\quad(k>0),\\
 C_{t,0}(J)&=1_{J\ge t-1}q^{t-1}T_t(p^{J-t+1}),\\
 C_{t,1}(J)&=1_{J\ge t}\tau_1q^{t-1}T_{t-1}(p^{J-t}),&
 C_{t,k}(J)&=0\quad(t>0,\ k\ge2),
 \end{aligned}
\]
where
\[
 T_v(y)=
 \begin{cases}
 \bar\chi_p(y)^v\tau_v,&p\nmid y,\\
 q-1,&p\mid y,\ 6\mid v,\\
 0,&p\mid y,\ 6\nmid v.
 \end{cases}
\]
These identities follow by summing the last nontrivial residue digit;
the remaining digits give $q^{t-1}$. Multiply $C_{t,k}(J)$ by
\[
 \gamma_2^{e_0}\vartheta^{e_0+3l}\rho^{k-e_0-3l}
 q^{-e_0(1/2+s)-l(1+3s)-kw-6mz}.
\]

The terms $t=0$ give $(1-V)^{-1}$ and
$1_{j=0}\rho Y_0/(1-\rho Y_0)$. For $t>0$, $(t,J)\mapsto(t+6,J+6)$
multiplies the term by $R$; the initial nonzero terms are
\[
\begin{array}{c|c|c|c|l}
 t\bmod6&k&J\text{ condition}&j&\text{contribution}\\\hline
 0&0&J\ge t&0,\ldots,5&R(1-q^{-1})/(1-V)\\
 0&0&J=t-1&5&-a^2q^3T^6\\
 0&1&J=t&0&\rho Y_0R\\
 1&0&J=t-1&0&-\eta\rho^{-1}T\\
 1&1&J=t&1&\eta TY_0\\
 1&1&J>t&0,\ldots,5&-\eta(q-1)TY_0V^{1_{j\le1}}/(1-V)\\
 3&0&J=t-1&2&a\rho^{-3}q^{3/2}T^3\\
 3&1&J=t&3&-\eta b_0\rho^{-2}q^2T^3Y_0\\
 4&0&J=t-1&3&b_0^2\rho^{-4}q^2T^4\\
 4&1&J=t&4&-\eta a\rho^{-3}q^{5/2}T^4Y_0 
\end{array}
\]
with initial $t=6$ in the first three rows and $t=1,3,4$ in the others,
using
\[
 \vartheta^3\gamma_3=a,\quad
 \vartheta^3\gamma_1\gamma_2=-\eta b_0,\quad
 \gamma_2\vartheta^4\gamma_4=b_0^2,\quad
 \gamma_1\gamma_2\gamma_3\vartheta^4=-\eta a,\quad
 \vartheta^6=a^2.
\]
They follow from $\gamma_2^3=-\Lambda$,
$G_0=\bar\chi_p(4)\gamma_3$, $\gamma_3^2=1$,
$\gamma_1\gamma_2=-\Lambda G_0$, and $\chi_p(4)^3=1$.
Summing $(1-R)^{-1}$ gives the six $J_j$.
The composite CRT phases cancel since
$2(e_0+3l)(e'_0+3l')+4e_0e'_0\equiv0\pmod6$.
Here $\vartheta^6=a^2$ retains $\eta^6$; no order-six condition is
imposed on $\eta$.

\begin{lemma}[Correction phases and the unramified local factors]
\label{v:lem:vartheta-sixth}
Let $p\notin S$ and $\vartheta_p=\bar\Lambda(p)\bar G_0(p)\eta(p)$ with
$G_0(p)=\bar\chi_p(4)\gamma_3$. Then $4$ is a unit modulo $p$,
$\chi_p(4)^6=1$, $\gamma_3^2=1$, $|G_0(p)|=1$, $G_0(p)^6=1$, and
\[
 \vartheta_p^6=\bar\Lambda(p)^6\eta(p)^6=a^2,\qquad
 (\vartheta_p^3\gamma_3)^2=a^2 .
\]
For the unramified valuation row $j=0$, with any $\rho$, the local
factors $H_p$ and $G_p$ of \eqref{v:cmp:5} depend on the Gauss phases
$\vartheta_p,\gamma_3$ only through $a^2$ and are independent of $b_0$.
In particular this applies to the principal frequency $u=1$, where
$j=0$ and $\rho=1$. For general $j$, the full six-row table must be
retained: $J_2,J_4$ contain terms linear in $a$, and $J_3$ contains a
term linear in $b_0$.
\end{lemma}
\begin{proof}
The preliminary set $S$ contains every place above $6$. Thus $p\nmid6$
and $4$ is a unit modulo $p$, so $\chi_p(4)^6=1$. Reduction of
$\mu_{12}$ is injective at $p$, giving $q_p\equiv1\pmod{12}$ and
$\chi_p(-1)=1$. The quadratic Gauss identity therefore gives
$\gamma_3^2=1$ and $|\gamma_3|=1$. It follows that $|G_0(p)|=1$ and
$G_0(p)^6=\chi_p(4)^{-6}\gamma_3^6=1$.

Taking sixth powers in the definition of $\vartheta_p$ now gives
$\vartheta_p^6=\bar\Lambda(p)^6\eta(p)^6=a^2$. Multiplication by
$\gamma_3^2=1$ gives the second displayed identity. In particular the
factor $\eta(p)^6$ is retained; no order-six assumption on the target
character is made.

For $j=0$, the table gives
$J_0=-\eta\rho^{-1}T+\rho Y_0R$, where $R=a^2q_p^4T^6V$.
The remaining terms of $P_0,P_0^*$ and the denominators in
\eqref{v:cmp:5} introduce no further dependence on $a$ or $b_0$.
This proves the asserted restriction to $j=0$. The other rows have
the linear terms stated above, so this restriction cannot be omitted.
\end{proof}

\paragraph{\it Conventions used by the local table.}
The checked verification of the table of $C_{t,k}(J)$ and of the six
$J_j$ uses the following conventions, which are made explicit here.
$\tau_v$ is the unnormalized Gauss sum
$\sum_{x\bmod p}\chi_p(x)^ve_p(x)$, with $\tau_0=-1$. The relation
$\tau_1\tau_5=q$ uses $\chi_p(-1)=1$, which holds because
$q\equiv1\pmod{12}$ at every exterior prime of $F\supset\mu_{12}$.
The five phase identities displayed above follow from
$G_0=\bar\chi_p(4)\gamma_3$, $\gamma_3^2=1$,
$\gamma_1\gamma_2=-\Lambda G_0$, $\chi_p(4)^3=1$ and the conjugation
relations $\gamma_4=\bar\gamma_2$, $\gamma_5=\bar\gamma_1$; they do
not use $\gamma_2^3=-\Lambda$ or $|\eta|=1$. The identity
$\vartheta_p=\bar\Lambda(p)\bar G_0(p)\eta(p)$ is the value at an
id\`ele representing $p$ of the correction character $\Lambda^{-1}G^{-1}\eta$
of Section~\ref{sec:v-packet}, evaluated by its local law at an
allowed generator of $p$.

Put $W=1_{j=0}\rho q_p^{-w}$, $D=1_{j=0}\eta\rho^{-1}q_p^{-s}$ and
\begin{equation}
 H_p=\frac{(1-V)(1-W)P_j}{1-D},\quad
 G_p=\frac{(\bar\eta q_p^s-q_p^{-w})(1-V)(1-W)P_j^*-q_p^{-w}(1-VW)}{1-D}.
                                                                    \label{v:cmp:5}
\end{equation}
The difference between divisibility restriction with scaling and rescaling replaces $H_p$ by
$q_p^{z-1}G_p$, since $P_j-P_j^*=(1-VW)/((1-V)(1-W))$;
there is no division by a possibly zero $H_p$.

For a prime tuple $\mathbf p=(p_i)$ put
$\mathcal P(\mathbf p)=\{p_i\}$ (distinct primes, as the windows are
disjoint). The correction is the finite sum
\begin{equation}\label{eq:v-global-compensated-correction}
 \begin{aligned}
 \mathcal H_{\eta,u,Z}(s,w,z)
 =\sum_{\mathbf p}\ &\prod_i
       \left[W_i(q_{p_i}/P_i)q_{p_i}^{z-1}G_{p_i}(s,w,z;u)\right]\\
 &\times\prod_{p\notin S\cup\mathcal P(\mathbf p)}H_p(s,w,z;u).
 \end{aligned}
\end{equation}
It uses no quotient by $H_p$ and holds for every
sixth-power-free frequency, including its finite unit class.

To verify it, fix a selected prime $p$. Rescaling $X,Y$ by
$q_p^{-1}$ in \eqref{v:cmp:1} multiplies the Mellin integrand by
\[
 -q_p^{-3/2}q_p^{-1/2+z}q_p^{1-w}=-q_p^{z-1-w}.
\]
The complementary divisibility term restricts the local theta
valuation to $t>0$ and replaces $Z$ by $Zq_p$, giving
$\bar\eta(p)q_p^{s+z-1}P_j^*$. Their sum is
\[
 q_p^{z-1}\left[(\bar\eta(p)q_p^s-q_p^{-w})P_j^*
                         -q_p^{-w}(P_j-P_j^*)\right].
\]
Removing the same three Euler factors as at an unselected prime
gives $q_p^{z-1}G_p$ by \eqref{v:cmp:5}. By the CRT congruence above
these local identities multiply over exterior primes; the unit of $u$
stays inside each $\rho=\chi_p(u/p^j)$, and no unit average or projection
is deleted. Finite tuple summation commutes with the initially absolutely
convergent Dirichlet series. This proves
\eqref{eq:v-global-compensated-correction} for the same original
operation \eqref{v:cmp:1}.
This calculation uses the full $J_j$ at each prime, including the
terms linear in $a$ and $b_0$ when $j>0$; it does not use the
unramified phase reduction of Lemma~\ref{v:lem:vartheta-sixth}.

\begin{remark}[Convention for frequencies and unit classes]
\label{v:rem:frequency-convention}
The frequencies $u$ in \eqref{v:cmp:6} are sixth-power-free elements
of $R$ modulo $\Gamma=(R^\times)^6$. Each class is represented by one
fixed element, with no allowedness condition; the representative is
used only through $q_u$ and $\rho=\chi_p(u/p^j)$. The local computation
leading to \eqref{v:cmp:5} and \eqref{eq:v-global-compensated-correction}
uses only these values at exterior primes, so it holds verbatim for
every representative, allowed or not.

This matters for unit frequencies. For an $S$-unit $\epsilon$, being
allowed together with $\epsilon\in E_S$ forces $\epsilon\in U_S^6$.
Hence a nontrivial unit class has \emph{no} allowed representative, and
an allowedness hypothesis on $u$ must not be imposed in \eqref{eq:v-global-compensated-correction}.

The finite unit classes are part of the row sum. A nontrivial class
$\epsilon\notin U_S^6$ gives the Kummer character $\chi_\epsilon$ of
$F(\epsilon^{1/6})/F$, which is nonprincipal by Kummer theory, has
conductor supported in $S$, and has $q_\epsilon=1$. Such rows are
therefore nonprincipal rows of the small range $d\le d_{\min}$, treated
by the small-row contour of Section~\ref{sec:v-detail-30}. The unique
principal class is $u=1$.
\end{remark}

On initial lines $(s,w,z)=(4,4,2)$, actual residue Poisson gives
\begin{equation}
 \begin{aligned}
 I_{\rm mod}&=\frac1{(2\pi i)^3}\int\int\int
 X^{1/2-z}Z^{s+z-1}Y^{w-1}e^{(s+z-1)^2}\mathcal M(z)\widehat W_1(w)\\
 &\quad\times\sum_u^{(6)}q_u^{-z}
 \frac{\zeta_F^S(6z)L_F^S(w,\chi_u)}{L_F^S(s,\eta\bar\chi_u)}
               \mathcal H_{\eta,u,Z}\,dz\,dw\,ds.
 \end{aligned}\label{v:cmp:6}
\end{equation}
Absolute convergence follows from the finite correlation bound
$|\mathcal F_{v,A}(H)|\ll q_vq_A$ for the finite Fourier coefficients of
the arithmetic coefficients, as in Section~\ref{sec:prev-average}, and ideal
counting. The correction is
\eqref{eq:v-global-compensated-correction}, with the windows
and scalar normalizations of \eqref{v:cmp:1}. Nonzero dual orbits unfold into the single norm
profile; dual zero is removed by $\widehat W_0(0)=0$.
For \eqref{v:cmp:6}, this analytic identification
must be made for the original signed theta coefficients and tests:
residue Poisson, orbit unfolding and Mellin inversion must preserve
the same finite tuple expansion, and the resulting sum of integrals
must have an absolute majorant on the initial lines before Fubini is
applied. The local correction identity
\eqref{eq:v-global-compensated-correction} identifies the arithmetic
factor once these analytic operations are justified.

At $u=1$, the principal unselected product is holomorphic and, after a
fixed cutoff, within $1/2$ of $1$ on
\begin{equation}
 \Rea s\ge437/500,\quad\Rea w\ge19/20,\quad\Rea z\ge33/200.
                                                                    \label{v:cmp:7}
\end{equation}
Indeed, at $j=0,\rho=1$ the phase reduction of
Lemma~\ref{v:lem:vartheta-sixth} applies. Put $q=q_p$, $Y_0=q^{-w}$.
Clearing the two nonvanishing denominators in \eqref{v:cmp:5} gives
\[
 \begin{aligned}
 &(1-R)(1-\eta T)(H_p-1)\\
 &=-VY_0+\eta TY_0+\eta TV(1-Y_0)
     -\eta(q-1)TY_0V(1-Y_0)\\
 &\quad+R-\eta TR-RY_0^2(1-V)-R(1-Y_0)/q.
 \end{aligned}
\]
At the corner in \eqref{v:cmp:7}, $|R|\le q^{-1117/500}$,
$|\eta T|\le q^{-437/500}$, and every numerator monomial is
$O(q^{-907/500})$, with exponents decreasing in the real
coordinates. Thus $H_p-1=O(q^{-907/500})$ uniformly in heights and
all unitary target phases.

The preliminary choice of $S$ in
Section~\ref{sec:v-scope} includes all primes up to a fixed cutoff
$Q_0(F)$ with
\[
 \sum_{p\notin S}\sup|H_p-1|<\log(3/2),\qquad
 \sup|H_p-1|<1/2.
\]
Then every factor is nonzero, its logarithm has a normally convergent
sum on compact subsets, and
$|\prod_{p\notin S}H_p-1|\le\exp(\sum\sup|H_p-1|)-1<1/2$,
the suprema being over all of \eqref{v:cmp:7}, hence uniform in the
three heights.

For the selected factor put $D=\eta(p)T$; as $|\eta(p)|=1$,
$\bar\eta(p)q_p^s=D^{-1}$. Expanding
$(1-R)(1-D)(G_p+1)$, or using its two defining terms in
\eqref{v:cmp:5}, gives, for
$\sigma=\Re s$, $\omega=\Re w$, $\zeta=\Re z$ in this region,
\[
 |G_p+1|\ll
 q^{-\sigma}+q^{-\omega}+q^{-6\zeta}
       +q^{1-\omega-6\zeta}+q^{4-5\sigma-6\zeta}
 \ll q^{-437/500}.
\]
The five decay exponents at the lower corner are
$437/500$, $19/20$, $99/100$, $47/50$, $34/25$.
All denominators are bounded away from zero after the fixed cutoff.
Only in this principal region, where $H_p\ne0$, we write
\[
 G_p/H_p=-1+O(q^{-437/500}).
\]

At the principal residues the tuple weights are nonnegative and the
number of prime factors is fixed; multiplying and summing the tuples of
\eqref{v:cmp:1} gives relative error
$O(Z^{-\omega_P})$, for every
$0<\omega_P<(437/500)\min_i\ell_i$.
No quotient by $H_p$ is used for nonprincipal classes.

The unique principal numerator class is $u=1$ by Kummer injectivity.
The residues $w=1,z=1/6$ have normalizer
\[
 c_S=\widehat W_1(1)\mathcal M(1/6)(\operatorname{Res}_{1}\zeta_F^S)^2/6,
\]
\[
 A_{\mathcal T}=(-1)^K Z^{-\ell/6}\prod_i\sum_{p\in\mathcal P_i}
                  W_i(q_p/P_i)q_p^{-5/6}.
\]
At these residues
\[
 X^{1/3}Z^{s-5/6}=Z^{C(s)-\ell/6},\qquad
 (-1)^K\sum_{\mathbf p}\prod_iW_i(q_{p_i}/P_i)q_{p_i}^{-5/6}
                       =Z^{\ell/6}A_{\mathcal T},
\]
so the two $\ell/6$ powers cancel before normalization.
The $G_p/H_p$ error is kept as a separate $O(Z^{-\omega_P})$ principal
remainder. The main term is the principal unselected product
\[
 H_\eta(s)=\prod_{p\notin S}H_p(s,1,1/6;1),
\]
the factor \eqref{v:cmp:5} at $u=1$ and the residues $(w,z)=(1,1/6)$, the
counterpart of \cite[(4.13)]{LiuArXiv}; it is independent of $Z$.

Rotate the fixed phase so $c_S>0$. The prime ideal theorem in a fixed
identity ray gives $A_{\mathcal T}\ne0$ eventually and inverse $Z^{o(1)}$.
The normalized theta average $J_\eta=I_{\rm mod}/(c_S A_{\mathcal T})$ has principal term
\begin{equation}
 f_\eta(Z)=\frac1{2\pi i}\int_{(2)}Z^{C(s)}e^{(s-5/6)^2}
                \frac{H_\eta(s)}{L_F^S(s,\eta)}\,ds,
 \quad C(s)=s-(4+b)/6.                                           \label{v:cmp:8}
\end{equation}
$H_\eta$ is holomorphic/nonzero on the stated $s$ half-plane. Separate
principal remainder savings are $l_y/20,h/600,\omega_P$: moving $w$ to
$\Rea w=19/20$ and then $z$ to $\Rea z=33/200$ inside \eqref{v:cmp:7}, the raw
exponent $l_x/2+\Rea s-1+h\Rea z+l_y(\Rea w-1)$ gives the first two, as in
Section~\ref{sec:prev-average} and \cite[(17.7)]{LiuArXiv}.

\subsection{Large values, prime sums and counting estimates}
\label{sec:v-detail-29}
In this subsection, we apply the simultaneous zero detector recalled in
Section~\ref{sec:prev-average}, with the counting scheme of
Section~\ref{sec:iv-count}, to the character family of
Section~\ref{sec:v-detail-26}. Lemma~\ref{v:lem:general-field-detector}
gives the large values \eqref{v:cmp:9},
Lemma~\ref{v:lem:detector-Sobolev} makes the rowwise choices admissible in
the moments, and the two moments, with the selected prime-sum factors in
the same presentation, give the row count $U^{R+\epsilon}$ with $R$ as in
\eqref{v:cmp:10}, following Proposition~\ref{prop:iv-count}; this count is
used in Section~\ref{sec:v-detail-30}.

Let $\beta=\beta_*(F)$ be the supremum of zero real parts over the full
finite-order family, with $1/2$ adjoined. The inverse estimate without prime-polynomial factors,
sixth-power-row inclusion--exclusion and scale Sobolev give $\beta\le11/12$
(Theorem~\ref{v:thm:unmarked-seed}):
there are $\gg H^{1/6}$ rows $k=a^6$ and
\[
 M_1(D)=\sum_{d\mid\operatorname{rad}a}\mu(d)\nu(d)q_d^{-1/2}M_{a^6}(D/q_d).
\]
Averaging the $H$-sized moment at $H=D^{1+c}$ gives the unnormalized
M\"obius annular exponent $11/12+\eps$, hence reciprocal continuation.

Choose a fixed finite $S$-supported character group $\mathcal N$
containing $1$, $\eta$, $G_0$ and the ray class characters occurring in the theta average.
Here $G_0$ is the multiplicative finite-order Hecke character constructed
in Lemma~\ref{v:lem:theta-characters}; it is not the generally
nonmultiplicative single Gauss coefficient $a_0$. Close the group
under conjugation and $\Upsilon_S$, and enlarge the M2
coefficient group $\Theta$ to contain $\mathcal N$. 

\begin{lemma}[Buffered simultaneous detector over a fixed number field]\label{v:lem:general-field-detector}
Fix the cyclotomic field $F$, the set $S$ and a finite $S$-supported
character group $\mathcal N$, closed under conjugation and containing the
presentations required in the comparison. For a sixth-power-free row $u$,
put
\[
 \mathcal X_u=\{\nu\chi_\bullet(u)^{\pm1}:\nu\in\mathcal N\}.
\]
Remove rows for which a presentation induces the principal character.
They are bounded unit classes: at an exterior prime of valuation
$1\le j\le5$, an unramified character in $\mathcal N$ cannot cancel
the nontrivial tame inertia of order $6/(6,j)$.
Write $U=Z^d$, $0<d_{\min}\le d\le d_{\max}$ and $q_u\asymp U$.
Fix a floor $a_f>1/2$, and put $\delta_f=2a_f-1$; in the present
comparison take $a_f=51/100$, $\delta_f=1/50$, the lowest threshold
$51/100$ of the detector \cite[(5.1)]{LiuArXiv} recalled in
Section~\ref{sec:prev-average} (there written $a_0$; here $a_0$ is the
Gauss coefficient above).
Assume $\beta_*(F)\ge a_f$, as holds in every contradiction argument here.
For every requested loss $\epsilon>0$, sufficiently small buffer/grid
$e>0$ and sufficiently small $\tau>0$, put $T_1=Z^\tau$.
The choices are made in the order specified below.

There are $O_e(1)$ bins $(i,a)$ with
\[
 a\in a_f+e\mathbb Z_{\ge0},\quad a\le\beta_*(F),
 \quad a_f\le a\le1,
\]
covering the retained rows. The bin $a=a_f$ is bounded by the trivial
row count. In every other bin there is one presentation
$\psi\in\mathcal X_u$ and a zero $\rho=\sigma+i\gamma$ satisfying
\[
 a\le\sigma<a+e,\qquad |\gamma|\le3iT_1.
\]
For every $1\le t\le3/2$, there are annular inverse/plain polynomials
for that same presentation and norm twist such that
\begin{equation}\label{v:cmp:9}
 \begin{gathered}
 r\le t+O(1/\log U),\quad r+m\ge t-O(1/\log U),\\
 t-1/2-\epsilon\le r,\qquad0\le m\le1/2+\epsilon,\\
 |M_r|^2\gg U^{(2a-1)r-\epsilon},\qquad
 |S_m|^2\gg U^{(2a-1)m-\epsilon}.
 \end{gathered}
\end{equation}
The untwisted profiles lie in a uniformly bounded smooth family.
Their common twist has height at most $(3i+1)T_1$; the detector's
extra frequency uses only a prescribed small fraction of one cumulative
$T_1/2$ allowance. Choices of presentation and dyadic pair are rowwise,
from a fixed finite set and $O((\log U)^2)$ pairs, respectively.
All estimates have finite seminorm and polynomial height orders for
this fixed field; no degree-uniform assertion is required.
\end{lemma}
\begin{proof}
Let $I=\lceil1/e\rceil+2$, and for $1\le j\le I$ define
\[
 B_j(u)=\max\bigl(\{a_f\}\cup\{\Re\rho:
 L(\rho,\psi^*)=0,\ \psi\in\mathcal X_u,
 \Re\rho\ge a_f,\ |\Im\rho|\le3jT_1\}\bigr).
\]
The maxima exist by discreteness and the absence of zeros to the right
of one. The nondecreasing sequence lies in $[a_f,1]$, so some
$B_{i+1}-B_i\le e$. Round $B_i$ down to the stated grid point $a$.
Then $a\le B_i<a+e$ and $B_{i+1}<a+2e$. If $a>a_f$, $B_i$ is
attained by an actual zero; otherwise use the trivial row count.
Since the field-family supremum includes these primitive characters,
$a\le\beta_*(F)$ whenever $\beta_*(F)\ge a_f$.

Logarithmic control: every primitive presentation has conductor $Q_\psi\ll U$, and its deleted
Euler product has radical norm $\ll U$. On a fixed real strip the
ordinary Hecke functional equation and convexity give a polynomial
upper bound in $U$ and $3+|v|$, with a finite field-dependent height
order. The disk centered at $2+iv$ with radius $2-a-2e$ is zero-free
for $|v|\le(3i+2)T_1$: its left boundary is $a+2e$ and its heights
lie below $3(i+1)T_1$ when $T_1>2$. Choose the analytic logarithm
anchored by the absolutely convergent Euler logarithm at its center.
Borel--Caratheodory bounds its modulus by $O_e(\log U)$ on radius
$2-a-4e$; it is $O(1)$ on radius $1/2$; Hadamard three-circles
bounds it by $O_e((\log U)^\theta)$, $\theta<1$, on radius
$2-a-6e$. Hence, for arbitrary $\epsilon_1>0$,
\[
 |L^S(s,\psi)|+|1/L^S(s,\psi)|\ll U^{\epsilon_1},
 \qquad\Re s\ge a+6e,\quad|\Im s|\le(3i+2)T_1.
\]

Farther right use absolute Euler convergence. Deleted Euler factors
and their reciprocals cost $U^{\epsilon_1}$ on a fixed positive strip:
split finitely many small primes and use
$\log(1\pm q_p^{-c})=O_c(q_p^{-c})=o(\log q_p)$.
On the reflected line, the primitive functional equation gives
\[
 |L^S(1-a-6e+iv,\psi)|
 \ll U^{a-1/2+12e+\epsilon_1}(3+T_1)^{A_F}.
\]
The conductor power is at most $U^{a-1/2+6e}$; the extra $6e$ safely
covers the deleted Euler product if this line is slightly negative.
The conjugate primitive presentation is in the same buffered family.

For an annular test, shift the central Mellin segment only, with total
added frequency bounded by the assigned fraction of $T_1$, to
$\Re s=a+6e$. Leave the infinite tails on $\Re s=2$.
The buffered reciprocal bound gives
$|M_r|^2\ll U^{(2a-1)r+\epsilon_2}$ for bounded $r$.
For a plain sum the same right shift and the alternative left shift
to $1-a-6e$ give
\[
 |S_m|^2\ll U^{(2a-1)\min(m,1-m)+\epsilon_2},
 \qquad 0\le r,m\le22.
\]
The left shift has no reciprocal, so its joins use the global upper
strip bound. Pure norm twists are placed in the $L$-argument;
the Mellin transform whose tail is differentiated is untwisted.

All real strips and length ranges are bounded. After fixing $\tau$,
a sufficiently large external derivative order makes every join and
tail $O(U^B T_1^{-N})$ negligible, including when $m>1$ and the
claimed exponent is negative. Increasing this external order does
not change the polynomial height orders of the central estimates.

For the detector choose the zero above, $D_*=U^t$, $Y_*=U^{20}$,
and a fixed smooth $V_\le$ equal to one on $[0,1]$ and zero on
$[2,\infty)$. Put
\[
 C_\psi(s)=\sum_l\mu_F(l)\psi(l)V_\le(q_l/D_*)q_l^{-s},
 \quad
 J=\frac1{2\pi i}\int_{(2)}Y_*^z\Gamma(z)
             L^S(\rho+z,\psi)C_\psi(\rho+z)\,dz.
\]
Shift to $\Re z=-1/4$; the only crossed pole, $z=0$, has residue
$L^S(\rho,\psi)C_\psi(\rho)=0$. On this strip the weak global bound
$|L^S(\sigma-1/4+iv,\psi)|\ll U^{2+o(1)}(3+|v|)^{A_F}$ holds (the field
degree affects only $A_F$), and ideal counting gives
$|C_\psi(\rho-1/4+iv)|\ll D_*^{5/4-\sigma+o(1)}$.
Gamma decay integrates the extra height, hence
\[
 |J|\ll U^{-3+t(5/4-\sigma)+o(1)}(1+T_1)^{A_F}=o(1),
\]
since $\sigma\ge a_f>1/2$ and $t\le3/2$ give $U$-exponent $\le-15/8$, and
$\tau$ is chosen small after $A_F$ (no approximate functional equation
is used).

On the original line Mellin inversion gives
\[
 J=\sum_{l,k}\mu_F(l)\psi(lk)V_\le(q_l/D_*)
                 (q_lq_k)^{-\rho}e^{-q_lq_k/Y_*}.
\]
For $1<q_n\le D_*$ its total coefficient is
$\sum_{l\mid n}\mu_F(l)=0$, including the natural zeros of $\psi$.
The unit coefficient is one, so inserting
$1-V_\le(2q_lq_k/D_*)$ removes $1+o(1)$ and the remaining double sum has
absolute value $1+o(1)$.
A terminal cutoff at $q_lq_k\le2U^{21}$ changes it by $o(1)$,
by exponential decay and ideal counting.

Smooth dyadic partition gives $O((\log U)^2)$ pairs $D=U^r$,
$N=U^m$ with $D\ll D_*$, $DN\gg D_*$, $DN\ll U^{21}$.
The product cutoff is a uniformly smooth function of
$(q_l/D)(q_k/N)$. Logarithmic Fourier inversion separates it as
an integral of $x^{i\nu}y^{i\nu}$ with uniformly bounded weighted
$L^1$ density. Discard $|\nu|>cT_1$ only after $\tau>0$ is fixed;
arbitrary decay of this untwisted profile and trivial polynomial
counts make that tail negligible. Hence one dyadic pair and
one frequency give a block
$(DN)^{-\rho}A_lA_k$ of magnitude $\gg(\log U)^{-2}$, where
$A_l,A_k$ are the unnormalized sums with profiles
$x^{-\sigma-i(\gamma-\nu)}$ and $y^{-\sigma-i(\gamma-\nu)}$
times their fixed smooth cutoffs; thus
$|A_lA_k|\gg(DN)^\sigma(\log U)^{-2}$.

Both sums have the common twist $\gamma-\nu$ and uniformly smooth
profiles after removing it. Central normalization yields
\[
 |M_rS_m|^2\gg(\log U)^{-4}(DN)^{2\sigma-1}
       \ge U^{(2a-1)(r+m)-\epsilon_2}.
\]
The inverse profile includes the smooth cutoff $V_\le(Dx/D_*)$;
its logarithmic derivatives are uniform because $D/D_*$ is bounded.

Apply the preceding pointwise bounds with preliminary loss
$\epsilon_2\ll\epsilon\delta_f$. Comparing product bounds gives
\[
 2(2a-1)(m-1/2)_+\ll\epsilon_2.
\]
As $2a-1\ge\delta_f>0$, $m\le1/2+\epsilon$; the product-support
condition gives $r\ge t-1/2-\epsilon$. Dividing the product lower bound by
each individual upper bound gives the two individual lower bounds, after
decreasing preliminary losses once more.

Fix $a_f$, requested losses, bounded length
ranges and common test families; choose $e\ll\epsilon\delta_f$;
fix all field-dependent strip, moment and finite seminorm height orders;
choose $\tau>0$ so their polynomial factors are absorbed by the desired
small $U$-powers, using $T_1=U^{\tau/d}$ and $d\ge d_{\min}$;
then choose external tail orders and the final threshold. Each added
Mellin frequency receives a fraction of one cumulative $T_1/2$ allowance.
This gives the stated uniformity.
\end{proof}

\begin{corollary}[The bin inputs used by the row contours]
\label{v:cor:bin-inputs}
In the setting of Lemma~\ref{v:lem:general-field-detector}, let $u$ be
a retained row in a bin $(i,a)$ with $a>a_f$, and let $\psi$ range over
$\mathcal X_u$ (in particular $\psi=\eta\bar\chi_u$, $\chi_u$ and
$\bar\chi_u$). Then, for every $\epsilon_1>0$, the following hold.
\begin{enumerate}
\item $L^S(s,\psi)\ne0$ for $\Re s\ge a+2e$, $|\Im s|\le3(i+1)T_1$.
\item $|L^S(s,\psi)|+|L^S(s,\psi)^{-1}|\ll U^{\epsilon_1}$ for
$\Re s\ge a+6e$, $|\Im s|\le(3i+2)T_1$.
\item On the reflected line,
$|L^S(1-a-6e+iv,\psi)|\ll U^{a-1/2+12e+\epsilon_1}(3+T_1)^{A_F}$ for
$|v|\le(3i+2)T_1$.
\item $\#\{u\text{ in the bin}\}\ll U^{R+\eps}(1+T_1)^A$, with $R$ as in
\eqref{v:cmp:10}.
\end{enumerate}
The implied constants depend only on $F$, $S$, $\mathcal N$, $e$ and
$\epsilon_1$, not on $u$. Items (1)--(3) are proved in the proof of the
lemma; item (4) is the count following
Lemma~\ref{v:lem:detector-Sobolev}. Item (1) is the zero-free input of
the $s$-move to $a+16e$ in Section~\ref{sec:v-detail-30}. Items (2)
and (3), with $A=a-1/2+6e$, give the denominator and numerator bounds
in \eqref{v:cmp:11}. No other information on zeros of row
$L$-functions is used. For rows in the lowest bin $a=a_f$ only the
trivial count is used.
\end{corollary}

\begin{lemma}[Moment bounds for the rowwise detector choices]
\label{v:lem:detector-Sobolev}
The rowwise choices in Lemma~\ref{v:lem:general-field-detector}
may be used in the inverse and product moments, at a cost which
is a fixed polynomial in $T_1$, with finite profile seminorms.
The same holds while retaining the original prime-polynomial
Mellin frequencies.
\end{lemma}
\begin{proof}
First fix the discrete presentation ($O(1)$ choices) and dyadic pair
($O((\log U)^2)$ choices). The remaining profiles
are parametrized by $\sigma$ in a fixed compact real interval,
the detecting twist, and the original prime Mellin frequencies.
Let their parameter box be $\mathcal B_T$, whose bounded coordinates
range in compact intervals and whose height coordinates have length
$O_e(T_1)$. For the row polynomial $P_u$, repeated unit-interval
Sobolev gives the common inequality
\begin{equation}\label{eq:v-detector-common-supremum}
 \sum_u\sup_{\mathbf t\in\mathcal B_T}|P_u(\mathbf t)|^2
 \ll\sum_{\alpha\in\{0,1\}^{\dim\mathcal B_T}}
   \int_{\mathcal B_T^+}\sum_u
                   |\partial^\alpha P_u(\mathbf t)|^2\,d\mathbf t.
\end{equation}

The slightly enlarged box is chosen before the rows. Differentiating
$x^{-\sigma}$ or $x^{it}$ introduces powers of $\log x$ on a fixed
annulus, so the profile seminorms stay finite and grow polynomially in
the incoming heights. Deterministic scale phases are removed first, so
the logarithms are of normalized ratios, not row-dependent scales.

Apply the moment theorem to each differentiated polynomial inside
the integral. In M2 use $P_u=S_m^2Q$, with the two copies of $S_m$
and their mixed derivatives retained. If the prime frequency is
$t_z$ and the detecting frequency is $t$, absorb $t_z-t$ into the
prime profiles before inserting a common explicit mode.
The arithmetic coefficient and original prime windows are unchanged.
The parameter-box volume and all height factors give at most
$(1+T_1)^J$ for one finite $J$. Finally choose $\tau$ after this
$J$ so that $T_1^J\le U^\epsilon$, and absorb the logarithmic
number of discrete choices by decreasing the input moment loss.
This gives the common supremum before any row-wise large value is evaluated.
\end{proof}

If the detected presentation is $\nu\chi_u$, replace it and its zero
by their whole conjugates and rename $\bar\nu$ as $\nu$. This gives
$\psi=\nu\bar\chi_u$ with the same zero real part, symmetric buffered
height and conjugated common norm twist; conjugate both detecting profiles
together. The original prime-sum factors relative to $\psi=\nu\bar\chi_u$ have coefficient
$\bar\nu(p)1_{p\in1_{\mathcal T}}$, a fixed bounded prime coefficient
(or finite $\bar\nu\theta$ combination). Thus the two moment estimates use exactly the
original $Q_i$, not newly selected polynomials.

On an amplitude piece
$g\ell=\sum_i g_i\ell_i$, $0\le g_i\le\delta/2$; errors have $g_i=0$.
The endpoint strips are treated without selected prime factors.
For $r\ge1$, Lemma~\ref{v:lem:general-amplification} gives
\[
 \#\{u\text{ detected}\}\ll U^{e(r)-\delta r+\epsilon}
 \le U^{1/6+(5/6-\delta)t+\epsilon}=U^{R+\epsilon},
\]
where $R$ is the balancing value in \eqref{v:cmp:10}; this uses
$r\le t+O(1/\log U)$ and $\delta\le5/6$. For $1-\epsilon_c\le r<1$,
the same lemma gives $U^{1-\delta+\delta\epsilon_c+\epsilon}
\le U^{R+\delta\epsilon_c+\epsilon}$. For $m\ge1/2-\epsilon_c$,
the product moment without prime factors gives
\[
 \#\{u\text{ detected}\}\ll U^{1-2\delta m+\epsilon}
       \le U^{R+2\delta\epsilon_c+\epsilon}.
\]

It suffices to treat $r\le1-\epsilon_c$, $m\le1/2-\epsilon_c$.
Decrease each displayed admissible prime length by $\epsilon_c/100$ and choose
the detector loss below $\epsilon_c/100$. Then the two M1 inequalities
and the M2 inequality have fixed positive margins, uniformly over the
dyadic detecting scales. The $O(\epsilon_c)$ losses lie in the previously
chosen $\epsilon_{\rm np}$ allowance and do not change the prime mesh.

Put $u=g/\delta\in[0,1/2]$, $c=2/(3\kappa)$. M1 length bound $(1-r)/2$
gives saving $\delta[u+(1-u)r]$; M2 with $S_m^2$ and admissible length $(1-2m)/(6\kappa)$ gives $\delta[cu/2+(2-cu)m]$. Its lower bound is
$|S_m|^4$, not $|S_m|^2$. Their crossing is
\[
 r_*=[(2-cu)t+(c/2-1)u]/[3-(1+c)u]\ge8/13.
\]
Put $m_*=t-r_*$. The proof of Proposition~\ref{prop:iv-count}, through
\eqref{eq:iv-crossing}, gives
\[
 \frac{6-c}{10-2c}\le r_*\le1,\qquad\frac13\le m_*\le\frac12.
\]

Use M1 when $r\ge r_*$. On this side,
\[
 z_M=\frac{1-r}{2}\le\frac{4-c}{20-4c}\le\frac5{26}<\frac15,
 \qquad 2r-1\ge\frac1{5-c}>0.
\]
Use M2 when $r\le r_*$. The detector gives
$m\ge m_*-O(1/\log U)$, so on this side
\[
 z_S=\frac{1-2m}{6\kappa}
 \le\frac{1-2m_*}{6\kappa}+O(1/\log U)
 \le\frac c{12}+O(1/\log U)
 \le\frac1{13}+O(1/\log U)<\frac15.
\]

Here $2/3\le c\le12/13$; the strict bounds are uniform on this
interval. Thus $\ell/d>1/5$ supplies the moment used on each side.
The fixed length decrement absorbs the $O(1/\log U)$ error,
and $2r-1>0$ supplies the second M1 margin on its side. After this decrement,
select decreasing-amplitude whole prime-sum factors, and drop one fractional prime-sum factor.
A fixed small mesh makes its loss arbitrarily small. If fewer positive
prime-sum factors exist, retaining all gives at least $2gz$ without invented prime-sum factors.

Lemma~\ref{v:lem:general-amplification}, followed by Markov's inequality, gives the long inverse count with
$e(r)=\max(1,(1+5r)/6)$. Crossing short and long counts yields
\[
 D=3-(1+c)u,\quad P=(2-cu)(1-u),\quad\mathcal J=(5/6-\delta)D+\delta P,
\]
\begin{equation}
 t=1+\frac{\delta P}{2\mathcal J},\qquad
 R=1-\delta+\frac{(5/6-\delta)\delta P}{2\mathcal J}.               \label{v:cmp:10}
\end{equation}
$\mathcal J>0$ for $0\le\delta\le5/6$. With no prime-sum factors and $t=1$, the
stronger count is $R=1-2\delta/3$, from balancing $1-\delta r$ and
$1-2\delta(1-r)$ at $r=2/3$ as in Proposition~\ref{prop:prev-initial-count};
we may slacken it to $76/75-2\delta/3$, the extra $1/75$ covering the lowest
class as in \Prev{Lemma}{19.2}.

\subsection{Joint correction and legitimate contour assembly}
\label{sec:v-detail-30}
In this subsection, we combine \eqref{v:cmp:9}--\eqref{v:cmp:10} with the
Euler correction of Section~\ref{sec:v-detail-28}: the piecewise row sum
\eqref{v:cmp:11} has relative scale exponent $E_B(d)$ of \eqref{v:cmp:12},
and the contour assembly gives the high-side estimate \eqref{v:cmp:13}
used in Section~\ref{sec:v-detail-31}. Complete-family contour shifts
precede restriction to detected elements, and internal profile orders
precede the target-dependent external cutoff.

On $s=a+16e,w=1-a-6e,z=z_0=17/50$, \eqref{v:cmp:5} gives
\[
 H_p-1=O(q_p^{-1-10e}),\quad G_p+\bar\chi_p(u)H_p=O(q_p^{-a-4e})
\]
at good primes. Bad defects are $O(q_p^{-10e})$. The selected
local factor is $q_p^{z_0-1}G_p$. For $j\ge2$ separate its term
$-(q_p-1)q_p^{z_0-1-w}/(1-R)$. For
$1/2\le a\le11/12$ and $0<e\le1/1000$, the remaining terms have
bounds $O(q_p^{z_0-1/2-\sigma_j})$, where one may take
\[
\begin{array}{c|ccccc}
 j&1&2&3&4&5\\\hline
 \sigma_j&26/25+80e&32e&26e&42e&80e.
\end{array}
\]

This follows by substituting the six $J_j$ into \eqref{v:cmp:5},
retaining the external factor $q_p^{z_0-1}$ and comparing each monomial
with $q_p^{z_0-1/2}$; e.g.\ the critical $j=2$ exponent is
$z_0+1/2-2a-32e\le z_0-1/2-32e$.
All denominators are bounded away from zero after increasing the
fixed prime cutoff. These affine comparisons are uniform in heights;
$26e$ is a common positive margin.
Every phase coefficient in the six-row table has modulus one.
These estimates bound all its monomials individually, including
those linear in the local parameters $a$ and $b_0$; the parameter
$a$ in that table is distinct from the real bin parameter here.
No even dependence on the local phases is used for the ramified rows.

\begin{lemma}[The unselected product on the contour region]
\label{v:lem:H-strip}
Let $1/2\le a\le11/12$ and $0<e\le1/1000$. There is $Q_1=Q_1(e)$ such
that for every sixth-free $u$, every exterior prime $p\nmid u$ with
$q_p\ge Q_1$ and every unitary target phase,
\[
 |H_p(s,w,z;u)-1|\le C_1q_p^{-1-10e}\le\tfrac12
\]
on the region
\[
 \Rea s\ge a+16e,\qquad \Rea w\ge1-a-6e,\qquad \Rea z\ge17/50,
\]
with an absolute constant $C_1$ (the checked value is $C_1=56$),
uniformly in the three heights. Consequently
$\prod_{p\notin S,\,p\nmid u}H_p(s,w,z;u)$ converges absolutely and
locally uniformly on this region, is holomorphic there, and is bounded
by a constant depending only on $F$, $S$ and $e$; the finitely many
factors at $p\mid u$ are holomorphic and bounded on the region by the
bad-defect bound $O(q_p^{-10e})$. Nonvanishing of the product is
asserted and used only on the principal region \eqref{v:cmp:7}; on the
contour region only holomorphy and boundedness are used.
\end{lemma}
\begin{proof}
At $p\nmid u$ the comparison of monomials above gives
$H_p-1=O(q_p^{-1-10e})$ with the exponents listed, uniformly in heights
and in the unitary phases $\rho$, $\eta(p)$; the constant is absolute
because the comparison is between finitely many monomials with
absolute coefficients, and $Q_1(e)$ is chosen so that the bound is at
most $1/2$ and all denominators are at least $1/2$ in modulus. The
sum $\sum_pq_p^{-1-10e}$ over prime ideals converges, so the product
converges normally; holomorphy and the bound follow. The factors with
$q_p<Q_1$ are finitely many rational functions of $q_p^{-s},q_p^{-w},
q_p^{-z}$ whose denominators do not vanish on the region after the
cutoff.
\end{proof}
At a good prime write
\[
 q_p^{z-1}G_p=-\bar\chi_p(u)q_p^{z-1}H_p
       +q_p^{z-1}\bigl(G_p+\bar\chi_p(u)H_p\bigr).
\]
The second term, summed over its full window, is
$\ll P_i^{z_0-a-4e}\le P_i^{z_0-1/2-4e}$.
At a bad prime the error and strict choices divide $u$; their
weighted choices are summed inside each fixed row, with divisor
mass $q_u^\eps$, before the row-count estimate is applied, so they incur
no additional free-prime count.

For a remaining main window, define
\[
 Q_i(u;z)=P_i^{-1/2}\sum_{p\in\mathcal P_i}
        \bar\chi_p(u)W_i(q_p/P_i)(q_p/P_i)^{z-1}.
\]
Its removed scalar is $P_i^{z-1/2}$. If the detected presentation
is $\psi=\nu\bar\chi_\bullet(u)$, its coefficient is
$\bar\nu(p)$ times the original ray selector, expanded in the
fixed group $\Theta$. Its norm profile retains $\Im z$ as an
incoming height parameter; Lemma~\ref{v:lem:detector-Sobolev}
permits that frequency and the detecting frequency to be different.

Fix the primes assigned to the error terms or strict inequalities. Remaining main $H_p$ merge into the direct
unselected product independent of still-remaining main primes. Do not divide
by $H_p$. Buffered $L'/L$ and prime-power removal bound the ordinary
original prime-sum factors by their stated amplitudes.

For all strict primes simultaneously, primitive tame conductor exponent
is one. With $A=a-1/2+6e$, apply the full numerator functional equation
once and use
\[
 q_{\mathfrak f_u}\ll q_u\prod_{p^j\parallel u,\ \mathrm{strict}}q_p^{-(j-1)},
 \quad q_{\mathfrak f_u}^A\prod_{\mathrm{strict}}q_p^{z_0-w}
 \ll U^A\prod_{\mathrm{strict}}q_p^{z_0-1/2}.
\]
Strict tuple choices are divisor-many. Hence, for a row set of count
$U^R$,
\begin{equation}
 \sum_u|L(w,\chi_u)\mathcal H_u^{\mathrm{piece}}|
 \ll U^{R+\delta/2+\eps}Z^{\ell(z_0-1/2)+g\ell+\eps}(1+T_1)^{A_F}.     \label{v:cmp:11}
\end{equation}
Here $A_F$ is the fixed field-dependent height exponent of the functional
equation, distinct from the conductor exponent $A=a-1/2+6e$. Relative to
\eqref{v:cmp:8}, its scale exponent on $U=Z^d$ is
\begin{equation}
 E_B(d)=a-B+h(z_0-1/6)-a l_y-\ell/2+g\ell
                         +d(R+\delta/2-z_0).                      \label{v:cmp:12}
\end{equation}

For each retained bin $a_f=51/100\le a\le11/12$, use the actual
connecting box
\[
 \Rea s\ge a+16e,\qquad
 \Rea w\ge1-a-6e,\qquad \Rea z\ge17/50,
 \qquad0<e\le1/1000.
\]
Then $\Rea(s+w)\ge1+10e$, and the local estimates above hold throughout
this box: after clearing $(1-R)(1-D)$,
every numerator monomial has nonnegative powers of $Y_0$, and
its powers of $T,R,V$ give nonnegative total decay slopes in
$\Rea s,\Rea z$. For example $R/T$ has slopes $5,6$.
Thus the corner estimates persist when any real coordinate increases,
and the denominators remain uniformly nonzero. This gives holomorphy
of the correction on an open neighborhood of the connecting
contours, summability of good defects and a bad-prime product
$q_u^\epsilon$. The selected tuple sum has one polynomial $Z$
majorant before later external orders.

Fix the entire original bin before moving any contour: move $z$ and
$s$ to $.34,\beta+20e$ on global reciprocal lines; move $w$ to
$1-a-6e$ with $s+w\ge1+14e$. Numerator is nonprincipal and $6z>1$.
Only then truncate the original heights and move the retained $s$ segment
to $a+16e$ in its buffered rectangle. No amplitude-dependent contour is
moved. On $s$ joins, other axes may be extended using global numerator
bounds, while the buffered denominator depends only on the retained $s$.
No other join extends a merely buffered $s$ axis.

Explicitly, put
\[
 \sigma_g=\beta+20e,\quad z_+=1/6+e,\quad w_+=1+e,
 \quad z_-=33/200,\quad w_-=19/20.
\]
On each global path the reciprocal lies on $\Re s=\sigma_g>\beta$.
For $u=1$, first move from $(4,4,2)$ to
$(\sigma_g,w_+,z_+)$, crossing no pole. Move $w$ to $w_-$;
the only crossed pole is the simple pole at $w=1$.
Keep the resulting bulk integral on $(\sigma_g,w_-,z_+)$.
In its $w=1$ residue move $z$ to $z_-$, crossing only $z=1/6$.
The double residue is \eqref{v:cmp:8}, with the selected-prime
remainder already specified there. Relative to $C(\beta)$ the
three remaining savings are at least
\begin{equation}\label{eq:v-principal-buffered-savings}
 l_y/20-(20+h)e,\qquad h/600-20e,\qquad\omega_P-20e.
\end{equation}
Choose $e$ so that each is at least half its value at $e=0$.
All paths lie in the principal holomorphy region
\eqref{v:cmp:7}; no pole of a local correction factor is crossed.

\paragraph{\it The principal residues explicitly.}
Write $r_S=\operatorname{Res}_{w=1}\zeta_F^S(w)$ and let
\[
 \mathcal I_1(s,w,z)=X^{1/2-z}Z^{s+z-1}Y^{w-1}e^{(s+z-1)^2}\mathcal M(z)
 \widehat W_1(w)\frac{\zeta_F^S(6z)\zeta_F^S(w)}{L_F^S(s,\eta)}
 \mathcal H_{\eta,1,Z}(s,w,z)
\]
be the $u=1$ integrand of \eqref{v:cmp:6}; the numerator
$L_F^S(w,\chi_1)$ is $\zeta_F^S(w)$. The $w=1$ residue is
\[
 \operatorname{Res}_{w=1}\mathcal I_1
 =r_S\widehat W_1(1)X^{1/2-z}Z^{s+z-1}e^{(s+z-1)^2}\mathcal M(z)
 \frac{\zeta_F^S(6z)}{L_F^S(s,\eta)}\mathcal H_{\eta,1,Z}(s,1,z).
\]
Since $\operatorname{Res}_{z=1/6}\zeta_F^S(6z)=r_S/6$, the double residue is
\[
 \operatorname{Res}_{z=1/6}\operatorname{Res}_{w=1}\mathcal I_1
 =c_S\,X^{1/3}Z^{s-5/6}e^{(s-5/6)^2}
 \frac{\mathcal H_{\eta,1,Z}(s,1,1/6)}{L_F^S(s,\eta)},
 \qquad c_S=\frac{\widehat W_1(1)\mathcal M(1/6)r_S^2}6 .
\]
By $X^{1/3}Z^{s-5/6}=Z^{C(s)-\ell/6}$ and the tuple evaluation
$\mathcal H_{\eta,1,Z}(s,1,1/6)=Z^{\ell/6}A_{\mathcal T}H_\eta(s)(1+O(Z^{-\omega_P}))$
of Section~\ref{sec:v-detail-28}, its $s$-integral divided by
$c_SA_{\mathcal T}$ is $f_\eta(Z)$ of \eqref{v:cmp:8}, up to the principal
remainder. The two single-residue remainders are the bulk integrals on
$(\sigma_g,w_-,z_+)$ and $(\sigma_g,1,z_-)$; they carry the savings
$l_y/20$ and $h/600$ of \eqref{eq:v-principal-buffered-savings}.

\paragraph{\it Lines and $L$-function bounds used by each row family.}
The following table lists, for each family of rows, the line used for
each variable and the bound used on it. No other bound for $L$-functions
of the rows enters \eqref{v:cmp:13}.
\begin{center}
\small
\begin{tabular}{@{}p{0.2\textwidth}p{0.3\textwidth}p{0.42\textwidth}@{}}
\toprule
Rows & Lines $(\Re s,\Re w,\Re z)$ & Bounds used\\\midrule
$u=1$ (principal) & $(\sigma_g,w_\pm,z_\pm)$, residues at $w=1$, $z=1/6$ &
region \eqref{v:cmp:7}; $1/L_F^S(s,\eta)$ holomorphic for $\Re s>\beta$
(Remark~\ref{v:rem:no-circularity})\\
nonprincipal, bin $a>a_f$, $d_{\min}\le d\le h+\zeta$ &
$(\sigma_g,4,z_0)\to(\sigma_g,1-a-6e,z_0)$, then the truncated $s$ segment to $a+16e$ &
numerator: functional equation and
Corollary~\ref{v:cor:bin-inputs}(3), $\ll U^{a-1/2+12e+\eps}$;
denominator: Corollary~\ref{v:cor:bin-inputs}(1)--(2), $\ll U^\eps$;
correction: Lemma~\ref{v:lem:H-strip}; count: \eqref{v:cmp:10}\\
lowest bin $a=a_f$ & as above with $a=a_f$ & same bounds; trivial row count\\
small rows $d\le d_{\min}$ (including unit classes) & $(\sigma_g,1/2,z_0)$ &
central convexity on $\Re w=1/2$; global reciprocal on $\Re s=\sigma_g$\\
far rows $q_u>Z^{h+\zeta}$ & $(2,2,z_\infty)$, $z_\infty>2$ & absolute convergence\\
\bottomrule
\end{tabular}
\end{center}
Here $z_0=17/50>1/6$ is fixed, so no residue at $z=1/6$ occurs for a
nonprincipal row, and the cost $h(z_0-1/6)$ appears in \eqref{v:cmp:12}.
Lowering $z_0$ towards $1/6$ for these rows would need convergence and
row-uniform bounds for the unselected correction product on
$\Re z<17/50$; Lemma~\ref{v:lem:H-strip} requires $\Re z\ge17/50$.
The argument does not use such a bound.

For a nonprincipal row, including every nontrivial frequency-unit
class, move first to $(\sigma_g,4,z_0)$ and then to
$(\sigma_g,1-a-6e,z_0)$. The numerator is entire and
$z_0=17/50>1/6$, so there is no residue. Small nonprincipal rows
instead use $(\sigma_g,1/2,z_0)$ and the small-row estimate below.
Only after truncating heights may the retained $s$ segment move
from $\sigma_g$ to $a+16e$ inside its common buffered rectangle.
The local denominators are $1-D$ and $1-R$:
their possible zero boundaries have respectively
$\Re s=0$ and $\Re(s+z)=2/3$, outside every displayed path.
The whole finite tuple expansion is fixed before any amplitude
partition. Each main, error or strict contribution therefore refers
to the same original contour identity.

The determinant-one height change $(\Im(s+z),\Im z,\Im w)$ shows the
Gaussian and smooth Mellin transforms are rapidly decreasing in all
three original heights. Integrated tails and horizontal traces are
$O_N(Z^{D_\eta}T_1^{-N})$, with $D_\eta$ fixed before $N$. Internal
weighted measures are not used as pointwise trace bounds. List every
added Mellin argument once and allocate fractions of one common
$T_1/2$ allowance. All finite internal orders are fixed first.
The resulting theta average and main term satisfy
\begin{equation}
 |J_\eta-f_\eta|\ll Z^{C(\beta)-m_{\rm save}}(1+T_1)^{A_\eta}
                         +Z^{D_\eta}T_1^{-N}.                     \label{v:cmp:13}
\end{equation}

The principal and nonprincipal margins are kept separate.  The
nonprincipal choices $\epsilon_{\rm np},\eta_{\rm M2},K,\ell_i$
were fixed in Section~\ref{sec:v-detail-27}.  Put
\[
 \omega_P=\frac{437}{1000}\min_i\ell_i,
 \qquad
 \mathfrak r_P=\min\left\{\frac{l_y}{20},\frac h{600},\omega_P\right\},
 \qquad
 m_{\rm save}=\frac1{16}\min\{\mathfrak r_{\rm np},\mathfrak r_P\}>0.
\]

Keep the already fixed $0<\epsilon_{\rm low}<(\beta-B)/4$, with half
reserved for the direct estimate before normalization. Now that $K$,
$\ell_i$ and $\mathfrak r_P$ are fixed, choose the detector buffer
and additional normalization loss so that
\[
 0<e<\min\left\{\frac1{1000},
 \frac{\epsilon_{\rm np}}{100C_{\rm np}},
 \frac{\epsilon_c\delta_f}{10^6C_{\rm np}},
 \frac{\mathfrak r_P}{100(21+h)}\right\},
\]
\[
 0<\epsilon_{\rm norm}<\min(m_{\rm save}/4,\epsilon_{\rm low}/2).
\]
Choose $e$ small enough also for the numerical buffer/grid requirements
of Lemma~\ref{v:lem:general-field-detector} with preliminary loss
$\epsilon_c/100$. This completes the common numerical choices before
the target and final $S$ are fixed. All $S$-dependent generators,
characters, ray data and tests are then constructed for that final
$S$; no later change of $e$ or $S$ is made.
The principal remainders in
\eqref{eq:v-principal-buffered-savings} save at least
$\mathfrak r_P/2$; the nonprincipal losses stay within their
$\mathfrak r_{\rm np}$ allowance, independently of $\mathfrak r_P$,
so $\epsilon_{\rm np}<m_{\rm save}$ is not required. The fixed-power estimate
$|A_{\mathcal T}|^{-1}\ll Z^{\epsilon_{\rm norm}}$ is imposed only
through a later threshold and does not change the moment mesh.

After fixing the target and all internal profile, detector and
Sobolev orders, let $\mathfrak h_\eta,A_\eta,D_\eta$ be the resulting finite
height and tail orders, with $\mathfrak h_\eta$ enlarged to dominate all strip,
functional-equation and detector height orders. For $d\ge d_{\min}=1/100$, choose
\[
 0<\tau_\eta<\min\left\{
 \frac1{100},\frac{m_{\rm save}}{4(A_\eta+1)},
 \frac{d_{\min}\min(1,\epsilon_{\rm np},\epsilon_c)}{10^4(\mathfrak h_\eta+1)}\right\}.
\]
Then choose the external tail order $N$ so that
$N\tau_\eta>D_\eta-C(\beta)+m_{\rm save}$.
Thus \eqref{v:cmp:13} retains a saving $\ge m_{\rm save}/2$.
The direct estimate before normalization loses $\le\epsilon_{\rm low}/2$
and the normalizer less than the other half (no relation
$\epsilon_{\rm low}<m_{\rm save}$ is required); its exponent and the high-side
exponent are both $<C(\beta)$.

For the fixed field and contradiction geometry, all savings and
buffers are common to the target family: the later fixed-target
choices of $S$, coefficient group, conductor, ray modulus, profile orders
and thresholds do not change $K$, $\ell_i$ or $m_{\rm save}$, and for each target
they precede every construction relative to $S$, which is then fixed.
The order is
\[
 \begin{aligned}
 &(F,\beta,B,\kappa,\text{geometry})\\
 &\quad\longrightarrow(\mathfrak r_{\rm np},C_{\rm np},\epsilon_{\rm low},
                   \epsilon_{\rm np},\epsilon_c,\eta_{\rm M2},\eta_{\rm phys})\\
 &\quad\longrightarrow(K,\ell_i)
 \longrightarrow(\omega_P,\mathfrak r_P,m_{\rm save},e,\epsilon_{\rm norm})\\
 &\quad\longrightarrow(\eta,\text{final }S,\text{arithmetic data},\text{tests},\Theta)\\
 &\quad\longrightarrow\text{internal orders}
 \longrightarrow(\mathfrak h_\eta,A_\eta,D_\eta)
 \longrightarrow\tau_\eta\longrightarrow N\longrightarrow Z_0.
 \end{aligned}
\]
The internal orders are fixed before the displayed height exponents.
The arithmetic data, tests, internal orders, height cutoffs, constants
and onsets may vary with the target; the continuation boundary for
this comparison does not.
Neither $J_\eta$ nor $f_\eta$ depends on the cutoff.

The following row bounds are the counterparts over $F$ of
\Prev{Lemma}{19.2}, recalled in Section~\ref{sec:prev-initial}; the
lowest-bin and intermediate bounds are obtained by substitution in
\eqref{eq:general-extension-exponent}, and the small-row and far-row
bounds by the contours stated with them.
Small rows $d\le1/100$ have relative bound
$13h/75-l_y/2+1/50$ by central convexity on $w=1/2$ and global reciprocal.
The lowest bin has endpoint
$-1/4+b/6+5\ell/4+(1/2+3\ell/2)/50$.
Intermediate rows use $R=76/75-2\delta/3$ and
\[
 E_B(1/2)\le[(225\ell+50-75b)\delta+78\ell-49b-73]/300.
\]

Fix $0<\zeta<5\ell-h$ inside the margins, which is possible since $5\ell>h$
(Lemma~\ref{lem:iv-reserves}) and gives $\ell/(h+\zeta)>1/5>7/37$ as
required in Section~\ref{sec:prev-initial}; this
$\zeta$ is unrelated to the real part $\zeta=\Rea z$ used in
Section~\ref{sec:v-detail-28}.
Far rows $q_u>Z^{h+\zeta}$ on $s=w=2,z=z_\infty>2$ have exponent
$l_x/2+1+l_y+(h+\zeta)(1+\eps)-\zeta z_\infty$;
choose $z_\infty$ before height orders. Principal remainders were
specified after \eqref{v:cmp:8}; all use the same normalizer.

\begin{lemma}[The common physical comparison]
\label{v:lem:physical-comparison}
Assume the geometric conditions of \eqref{v:cmp:4},
$l_x/2+b/12=C(B)$ with $B>437/500$, and the positive
high-side margins established above. Fix a target $\eta$, the prime
windows and all tests used above. The inputs to this comparison are
the actual theta/Mellin identity \eqref{v:cmp:6}, the direct bound
\eqref{v:cmp:4}, and the nonprincipal and principal remainder bounds
leading to \eqref{v:cmp:13}, with the contour justifications above.
The direct estimate and the high-side estimate concern the same
$J_\eta=I_{\rm mod}/(c_SA_{\mathcal T})$ and the same function
$f_\eta$ in \eqref{v:cmp:8}. After the stated height choices they give
\[
 |J_\eta|\ll_\eta Z^{C(B)+\epsilon_{\rm low}},\qquad
 |J_\eta-f_\eta|\ll_\eta Z^{C(\beta)-m_{\rm save}/2}.
\]
The function $H_\eta$ in that formula is independent of $Z$ and
holomorphic and nonzero for $\Re s>437/500$. If
$\sigma=\max(B+\epsilon_{\rm low},\beta-m_{\rm save}/2)<\beta$,
these estimates continue every fixed target to the same zero-free
half-plane $\Re s>\sigma$.
\end{lemma}
\begin{proof}
Expand the original compensated average by its finite tuple identity
\eqref{eq:v-global-compensated-correction}, and apply the radial
identity to its actual coefficients. On the absolute-convergence
lines this is precisely \eqref{v:cmp:6}; its subsequent contour
moves preserve this identity. The unique principal class is $u=1$.
Its two residues supply $c_S$, while the selected tuple sum supplies
$Z^{\ell/6}A_{\mathcal T}$ and the scale supplies
$Z^{C(s)-\ell/6}$. Cancellation of these two powers before division
by $c_SA_{\mathcal T}$ gives exactly \eqref{v:cmp:8}.
The selected-prime errors remain principal remainders; they are not
absorbed into a new numerator. Thus the numerator is the same
unselected Euler product $H_\eta(s)$, whose normal convergence and
nonvanishing follow from \eqref{v:cmp:7}.

The direct bound is \eqref{v:cmp:4} with the reserved normalization
loss. In \eqref{v:cmp:13} put $T_1=Z^{\tau_\eta}$.
The choices $A_\eta\tau_\eta<m_{\rm save}/4$ and
$N\tau_\eta>D_\eta-C(\beta)+m_{\rm save}$ give the displayed high-side bound.
The cutoff is a device for estimating the fixed contour identity;
it does not define $J_\eta$, $f_\eta$ or $H_\eta$.
Their sum bound is $f_\eta(Z)\ll_\eta Z^{\sigma-(4+b)/6}$.
Put $c_0=(4+b)/6$. The explicit continuation is
\[
 \int_0^1 f_\eta(Z)Z^{c_0-s}\,\frac{dZ}{Z}
 +\int_1^\infty f_\eta(Z)Z^{c_0-s}\,\frac{dZ}{Z}.
\]
The second integral converges locally uniformly for $\Re s>\sigma$.
For the first, move the line in \eqref{v:cmp:8} to any fixed
sufficiently large real part $A$; absolute convergence on the right
and the Gaussian give $f_\eta(Z)=O_{\eta,A}(Z^{A-c_0})$ as $Z\to0$.
Thus the first integral is entire. On the initial right half-plane,
Mellin inversion identifies their sum with
$e^{(s-5/6)^2}H_\eta(s)/L_F^S(s,\eta)$, giving the asserted
continuation. A target-dependent finite interval below the onset of
the large-$Z$ estimate contributes an entire function.
Since $\sigma\ge B>437/500$, its numerator has no zeros there,
so $L_F^S(s,\eta)$ has no zeros there. The numerical $\sigma$ is
common to all targets; constants, orders and thresholds may depend
on the fixed target. Applying this conclusion separately to each
target bounds the supremum of their zero real parts by $\sigma$.
No zero attaining that supremum and no uniform target constant is
required.
\end{proof}

\begin{remark}[The nonvanishing used on the $s$-contour]
\label{v:rem:no-circularity}
Lemma~\ref{v:lem:physical-comparison} is applied in
Section~\ref{sec:v-detail-31} with $\beta=\beta_*(F)$, the supremum
of the real parts of the nontrivial zeros of all $L_F(s,\nu)$ with
$\nu$ primitive of finite order, and with the $s$-contour of the
principal integral \eqref{v:cmp:8} kept at $\Rea s>\beta$. The only
nonvanishing assumed is therefore $L_F^S(s,\eta)\ne0$ on the closed
strip between the initial line $\Rea s=2$ and that contour, which
holds by the definition of $\beta$ (omitted Euler factors vanish only
on $\Rea s=0$). No zero-free statement to the left of $\beta$ is
assumed; the nonprincipal rows are handled by the detector and the
moments of Sections~\ref{sec:v-detail-29}--\ref{sec:v-detail-30}, not
by a nonvanishing hypothesis; and the conclusion
$\sigma_{\rm cont}<\beta$ contradicts the definition of $\beta$ rather
than an attained zero. The auxiliary $w$- and $z$-lines are not
constrained by $\beta$. The nonvanishing hypothesis ensures that $1/L_F^S(\cdot,\eta)$ be holomorphic on
$a\le\Rea s\le2$, which is why the comparison data are requested only
for $a>\beta$.
\end{remark}

\begin{remark}[The principal bound as two inputs]
\label{v:rem:principal-bound-split}
Write the right side of \eqref{v:cmp:6} as $\sum_uJ^{(u)}$ with
$J^{(1)}$ the principal triple integral. On the initial lines the row
sum converges absolutely, so
\[
 |J^{(1)}|\le|I_{\rm mod}|+\sum_{u\ne1}|J^{(u)}|
\]
for every $X,Y,Z>0$. The bound required for $J^{(1)}$ in
Lemma~\ref{v:lem:physical-comparison} is therefore exactly
Proposition~\ref{v:prop:direct-estimate} together with the bound
\eqref{v:cmp:13} for the nonprincipal rows; no other property of the
average enters.
\end{remark}

\subsection{Initial and variable scalar comparison; continuation}
\label{sec:v-detail-31}
In this subsection, we compare \eqref{v:cmp:4} with \eqref{v:cmp:13} through
the common principal term \eqref{v:cmp:8}. The initial parameters give
$\beta_*(F)\le7/8$; along \eqref{eq:iv-geometry}, Lemmas~\ref{lem:iv-reserves}
and~\ref{lem:iv-endpoint} give the implication \eqref{v:cmp:14}, which the
fixed-point identity \eqref{eq:iv-fixed-point-identity} turns into
$\beta_*(F)\le B_\infty$; and Mellin continuation proves Theorem~\ref{v:thm:original-target} without the supremum of zero real parts
being attained.

The initial comparison is the counterpart over $F$ of the $7/8$ comparison
recalled in Section~\ref{sec:prev-initial}. In the initial geometry
\eqref{eq:prev-initial-geometry}
\[
 (b,\ell,l_x,l_y,h)=(1/8,1/6,17/48,23/48,13/16),
\]
\eqref{v:cmp:4} is $J\ll Z^{C(7/8)+\eps}$. Suppose
$\Delta=\beta-7/8>0$; the initial bound $\beta_*(F)\le11/12$ of Theorem~\ref{v:thm:unmarked-seed} gives $\delta\le5/6$.
Use $\kappa=3/4+2\Delta$; then
$(1+\kappa)/2=7/8+\Delta=\beta$, so M2's zero-free condition holds
with equality.
The intermediate endpoint is $\le-49/14400$ and the lowest-bin endpoint
is $\le-7/1200$, as recalled in Section~\ref{sec:prev-initial}.

The exact initial endpoint inequality of
Lemma~\ref{lem:prev-initial-endpoint} gives
$E_{7/8}(h)\le-49/440640$ at $\kappa=3/4$.
Increasing $\kappa$ changes only the admissible plain length; on this
crossing branch the detector's plain-polynomial length, denoted here by
$m_{\rm det}$, satisfies $m_{\rm det}\ge1/3$, as in the proof of
Proposition~\ref{prop:prev-initial-count};
with $0\le g\le5/12$, the change costs
\[
 \frac{24g(1-2m_{\rm det})\Delta}{(9/2)(9/2+12\Delta)}\le\Delta/4.
\]

Relative to $C(\beta)$ the endpoint is
$-49/440640-(1-h/4)\Delta<0$. The slope in $d$ on the moderate range $1/2\le d\le h$ lies in
$[4/25,2]$ (Section~\ref{sec:prev-initial}),
and $\ell/h=8/39>1/5$. All small/far/principal margins are positive, as in
the bounds recalled in Section~\ref{sec:prev-initial}.
Contour margins, discretization and small $\zeta$ fit inside these margins; \eqref{v:cmp:13} and
nonzero numerator continuation prove $\beta\le7/8$.

Now let $I=[19429/50000,38859/100000]$ be the interval of
Section~\ref{sec:iv-geometry} and use the functions
$\kappa,\ell,B,b,l_x,l_y,h$ of \eqref{eq:iv-geometry}.

Lemma~\ref{lem:iv-reserves} gives the geometric side margins, and
Lemma~\ref{lem:iv-endpoint} of Section~\ref{sec:iv-certificate} gives
$E_{B(\xi)}(h)\le0$ for $\delta\in[0,5/6],u\in[0,1/2]$, with equality
only at $(\delta,u)=(\xi,1/2)$, by the identity \eqref{eq:iv-Q}, the
coefficient bounds and the completed square recalled after that lemma
(with $u$ in place of $x$). In the bin range $\delta\le3/4<19/25$ the side margins
listed after Lemma~\ref{lem:iv-reserves} (evaluated at $\delta=1/50,19/25$)
are strict. With $\beta>B$ subtract $\beta-B$ from every endpoint,
choose contour margins and discretization in that strict gap, and apply \eqref{v:cmp:13}. Thus, for $\xi\in I$ with $\kappa(\xi)\le3/4$ (as in Proposition~\ref{prop:iv-bootstrap}),
\begin{equation}
 \beta\le(1+\kappa(\xi))/2\quad\Longrightarrow\quad\beta\le B(\xi).     \label{v:cmp:14}
\end{equation}

Let $p$ be the cubic \eqref{eq:iv-cubic}. It has a unique root
$\xi_\infty=\xi_{\mathrm{new}}$ in $I$, $p'<0$ and $\kappa'>1000B'>0$ on $I$,
and the fixed-point identity \eqref{eq:iv-fixed-point-identity} holds.
$\xi_0=(49-\sqrt{921})/48\in I$ has $\kappa(\xi_0)=3/4$.
If $\beta>B(\xi_\infty)$, choose the unique
$\xi\in(\xi_\infty,\xi_0]$ with $\kappa(\xi)=2\beta-1$.
Then $p(\xi)<0$ gives $B(\xi)<\beta$, contradicting \eqref{v:cmp:14}.

The exact root bracket
\[
 \frac{388583354266034666797427}{10^{24}}<\xi_\infty<
 \frac{388583354266034666797428}{10^{24}}
\]
gives
\[
 .8749570194200989461286<B_\infty<.8749570194200989461287.
\]
The resulting strict half-plane is $\Rea s>B_\infty$.

\begin{remark}[Exponent ledger at the critical class]
\label{v:rem:exponent-ledger}
At $\xi=\xi_\infty$ the geometry \eqref{eq:iv-geometry} is, to ten digits,
\[
\begin{gathered}
 \kappa=0.7499140388,\quad\ell=0.1668385890,\quad b=0.1234056854,\\
 l_x=0.3548778628,\quad l_y=0.4782835482,\quad h=0.8119607262,
\end{gathered}
\]
with $B_\infty=0.8749570194$ and $C(B_\infty)=l_x/2+b/12=0.1877227385$
(the direct estimate \eqref{v:cmp:4}). The extremal class of
Lemma~\ref{lem:iv-endpoint} is $\delta=\xi$, $x=1/2$ ($\bar g=\delta/2$),
$a=(1+\delta)/2=0.6942916771$ and $d=h$. There the count \eqref{v:cmp:10}
has $t=1.1242271473$ and $R=2/3$. Its exponent \eqref{v:cmp:12}, relative
to $C(B)$, splits as follows.
\begin{center}
\small
\begin{tabular}{@{}lrl@{}}
\toprule
Term of $E_B(h)$ & Value & Source\\\midrule
$a-B$ & $-0.1806653423$ & $s$-line at $a+16e$ (Corollary~\ref{v:cor:bin-inputs}(1))\\
$h(z_0-1/6)$ & $+0.1407398592$ & $z$-line fixed at $z_0=17/50$ (Lemma~\ref{v:lem:H-strip})\\
$-a\,l_y$ & $-0.3320682868$ & $w$-line at $1-a-6e$\\
$-\ell/2$ & $-0.0834192945$ & selected prime factors at $z_0$\\
$\bar g\ell$ & $+0.0324153493$ & prime-sum amplitude\\
$h(R+\delta/2-z_0)$ & $+0.4229977151$ & count $U^R$, numerator $U^{\delta/2}$, $q_u^{-z_0}$\\\midrule
total & $0$ & equality case of Lemma~\ref{lem:iv-endpoint}\\
\bottomrule
\end{tabular}
\end{center}
The gain needed for the contradiction is the strict gap $\beta-B$. It
enters every endpoint through the term $a-B$ and survives because the
balance is exact at $B$ (Lemma~\ref{lem:iv-endpoint}). Every other
class has $E_{B(\xi)}(h)<0$; a grid over $0\le\delta\le5/6$,
$0\le x\le1/2$ gives maximum $-1.95\cdot10^{-7}$ attained next to the
extremal point. The remaining margins at $\xi_\infty$ are
\[
\begin{gathered}
 l_y/20=0.0239,\qquad h/600=0.00135\qquad(\text{principal remainders}),\\
 0.00588\ (\text{lowest bin}),\qquad0.0784\ (\text{small rows}),\qquad
 0.0218\ (\text{intermediate, }\delta=19/25),
\end{gathered}
\]
together with $5\ell-h=0.0222$ and $\ell/h=0.2055>1/5$ (prime length
available on the moderate range). In particular, every unit of exponent
is accounted for by the five inputs listed in the right column, the
detector count \eqref{v:cmp:10} and M1--M2. The heaviest positive cost
is $h(R+\delta/2-z_0)$, that is, the row count. Any loss in M1, M2 or in
Theorem~\ref{v:thm:unmarked-inverse} enters through $R$.
\end{remark}

For the continuation step put $c_0=(4+b)/6$, so that $C(s)=s-c_0$.
With $\epsilon_{\rm low}$ as in Section~\ref{sec:v-detail-30}, define the common boundary
\[
 \sigma_{\rm cont}
  =\max\{B+\epsilon_{\rm low},\,\beta-m_{\rm save}/2\}<\beta.
\]
Equations~\eqref{v:cmp:4} and \eqref{v:cmp:13}, with the height
choices above, give for every fixed target
\[
 |f_\eta(Z)|\ll_\eta Z^{\sigma_{\rm cont}-c_0}\quad(Z\ge1).
\]

Moving the principal integral arbitrarily far right gives
$f_\eta(Z)\ll_{\eta,A}Z^A$ as $Z\downarrow0$, for every $A>0$.
Mellin inversion and holomorphic continuation give
\[
 e^{(s-5/6)^2}\frac{H_\eta(s)}{L_F^S(s,\eta)}
   =\int_0^\infty f_\eta(Z)Z^{c_0-s}\,\frac{dZ}{Z},
 \qquad\Re s>\sigma_{\rm cont},
\]
the integral converging locally uniformly. As
$\sigma_{\rm cont}\ge B>437/500$, $H_\eta$ is holomorphic and nonzero
there, so every $L_F^S(s,\eta)$ is zero-free on $\Re s>\sigma_{\rm cont}$,
contradicting
the definition of $\beta$. The constants and onsets may depend on
$\eta$; the boundary $\sigma_{\rm cont}$ does not.
Finite omitted Euler-factor zeros have real part zero. Principal
poles give zeros of the reciprocal and are permitted.

\begin{theorem}[Cyclotomic Hecke bound, pointwise numerator]\label{v:thm:original-target}
For every fixed cyclotomic field \(F\supset\mu_{12}\) and every
finite-order Hecke character \(\eta\) of \(F\),
\[
 L_F(s,\eta)\ne0\qquad(\Rea s>B_\infty),\qquad
 B_\infty=\frac{18\xi_\infty+7}{18\xi_\infty+9},
\]
where \(\xi_\infty\) is the unique root in
\([19429/50000,38859/100000]\) of
\(1188\xi^3-2160\xi^2+171\xi+190\). The exact boundary is
\(0.874957019420098946\ldots\). Principal poles are permitted;
no assertion is made on the boundary. Constants and onsets may depend
on the fixed target, whereas this boundary is common. 
\end{theorem}
\begin{proof}
The preceding low and high estimates use the inverse moments with
prime-polynomial factors and the full-annulus ordinary moments supplied
by Theorems~\ref{v:thm:marked-inverse-new} and
\ref{v:thm:M2-assembled}. Their finite
character groups are fixed before bin selection; the discretization of the prime-length exponents
converts correctly on each \(U=Z^d\) scale. The inverse moment without prime-polynomial factors supplies
\(\beta_*(F)\le11/12\), the initial comparison supplies \(7/8\), and
\eqref{v:cmp:14} with the cubic identity gives \(\beta_*(F)\le B_\infty\).
If the presumed supremum were larger, the fixed positive comparison gap
would continue every fixed target's reciprocal through a common smaller
zero-free boundary. This contradicts the supremum without requiring an attained zero
or uniform target constants. The nonzero principal correction and
omitted Euler factors are treated explicitly above.
\end{proof}

\section{The averaged numerator and the refined cyclotomic boundary}
\label{v:sec:refined}
In this section, we improve $B_\infty$ of
Theorem~\ref{v:thm:original-target} to $B_\flat$ of
\eqref{eq:v-refined-limit} by an averaged estimate for the reflected
numerator, introduced over $\mathbb Q(\sqrt{-3})$ by Xu and Sang
\cite{XS26} assuming the analytic inputs of \cite{QRH}, and adapted here
to the fixed cyclotomic field $F\supset\mu_{12}$ of
Section~\ref{v:sec:compensated-assembly}. Only \eqref{v:cmp:11} is
replaced; \eqref{v:cmp:1}, \eqref{v:cmp:4}, \eqref{v:cmp:8},
\eqref{v:cmp:9}, \eqref{v:cmp:10}, M1 and M2 are used unchanged.

In \eqref{v:cmp:11} the numerator is reflected on
$\Rea w=1-a-6e$ and bounded pointwise, at conductor cost
$U^{a-1/2}$. Here the numerator line is $\Rea w=1/2-\theta$,
$0\le\theta\le a-1/2$ depending on the amplitude, and the reflected values on
$\Rea w=1/2+\theta$ are averaged over the detected row class by
H\"older's inequality and M2, with the same prime-sum factors
that define the class (Proposition~\ref{v:prop:refined-numerator};
assembly in Proposition~\ref{v:prop:refined-high}). The resulting endpoint
inequalities are field independent. Lemma~\ref{v:lem:refined-family}
constructs tangent parameters for every moment parameter
$\kappa_\flat\le\kappa\le3/4$ ($\kappa_\flat$ the fixed point)
and proves that $(b^*(\kappa),\ell,\kappa)$ is admissible whenever
$\ell<\ell^*(\kappa)$ is sufficiently close to $\ell^*(\kappa)$;
Theorem~\ref{v:thm:refined-limit} gives $B_\flat$ by a single
step justified by Theorem~\ref{v:thm:original-target}; and
Theorem~\ref{v:thm:refined-optimality} shows that no
choice of parameters and stages in this refined comparison
gives a boundary below
\begin{equation}\label{eq:v-refined-limit}
 \begin{gathered}
 B_\flat=0.874886647039964\ldots,\\
 133488B_\flat^4-376596B_\flat^3+384534B_\flat^2-167541B_\flat+26231=0,
 \end{gathered}
\end{equation}
so that $B_\flat$ is the exact infimum of the refined comparison.

Appendix~\ref{app:refined-chain} gives, independently of this argument,
an explicit finite chain of four rational stages in exact arithmetic: its
first two stages recover over $F$ the boundary $B_1=8748867/10^7$ of
\cite{XS26}, and its last stage gives
\begin{equation}\label{eq:v-refined-boundary}
 B_\sharp=\frac{874886647041}{10^{12}}=0.874886647041,
 \qquad
 B_\sharp-B_\flat<1.04\times10^{-12}.
\end{equation}

As specified in Remark~\ref{v:rem:analytic-inputs}, the inputs are M1, M2 and
the cyclotomic theta inputs of
Section~\ref{v:sec:compensated-assembly}, plus the standard
approximate functional equation (Lemma~\ref{v:lem:refined-afe}).
The inputs of \cite{QRH} used in \cite{XS26} are replaced by
the cyclotomic estimates of this part; the conclusion of that work is
not used.
The averaged numerator is outside the comparison system
of Section~\ref{sec:v-optimality}, consistently with
Section~\ref{sec:v-detail-37}, which lists additional mean
cancellation in the coupled numerator as a change of input.

\subsection{Parameters, exponents and admissibility}
\label{v:sec:refined-parameters}
In this subsection, we fix the lengths $(b,\ell)$, the row classes with
their amplitudes $g_i$ and the moment parameter $\kappa$, define the
numerator shift $\theta$ by \eqref{eq:v-refined-theta}, and compute the
full exponent \eqref{eq:v-refined-full-exponent} of a row class and its
endpoints \eqref{eq:v-refined-endpoint} and
\eqref{eq:v-refined-small-endpoint} at $d=h$.
Remark~\ref{v:rem:refined-admissible} collects the admissibility
conditions on $(b,\ell,\kappa)$ used in
Sections~\ref{v:sec:refined-high}--\ref{v:sec:refined-optimality}.

Retain the notation of Section~\ref{v:sec:compensated-assembly}.
Here $\theta$ denotes the numerator shift below; the
character $\vartheta$ of Section~\ref{sec:v-detail-27} and the
coefficient characters $\theta\in\Theta$ of
Section~\ref{sec:v-detail-29} do not occur as such.
Let $b>0$ and $1/6\le\ell\le1/5$, and put
\[
 l_x=\frac{1-b-\ell}2,\quad l_y=\frac{1+b-\ell}2,\quad
 h=\frac{1+b+3\ell}2,\quad
 B(\ell)=\frac{11}{12}-\frac\ell4,\quad
 C(s)=s-\frac{4+b}6 .
\]
Then $h=1-l_x+\ell=l_y+2\ell$ and
$C(B(\ell))=l_x/2+b/12$, the exponent in \eqref{v:cmp:4}.

A nonprincipal row class consists of sixth-power-free nonunit rows
$u$ with $q_u\asymp U=Z^d$, a detector bin $a$ of
Lemma~\ref{v:lem:general-field-detector} with $a>a_f$, and, for the
current $z$ on $\Rea z=z_0=17/50$, an amplitude bin of the normalized
prime sums $Q_i(u;z)$ of Section~\ref{sec:v-detail-30}:
\[
 0\le g_i\le\delta/2,\qquad
 |Q_i(u;z)|\le P_i^{g_i+\eps_{\rm bin}},\qquad
 g_i>0\Longrightarrow|Q_i(u;z)|\ge P_i^{g_i},
\]
where $\delta=2a-1$ and $\eps_{\rm bin}>0$ is the fixed amplitude-bin
mesh, chosen with $\epsilon_{\rm np}$, before $K$ and $e$, inside
the nonprincipal margin; then $e\le\eps_{\rm bin}/10$, so that
$|Q_i|\ll P_i^{\delta/2+2e+\epsilon_1}$ places every amplitude in a bin
with $g_i\le\delta/2$ (likewise for the bins of
Section~\ref{sec:v-detail-29}). Put
\[
 g=\frac1\ell\sum_i\ell_ig_i,\qquad x=\frac g\delta\in[0,1/2],
 \qquad\mathscr L=\frac\ell d .
\]
Here $x$ is the normalized amplitude written $u$ in
Section~\ref{sec:v-detail-29} (renamed to avoid a clash
with the rows). Every prime-sum factor keeps its amplitude $g_i$, even in
expansion terms where it carries a local error.

For a moment parameter $13/18\le\kappa\le1$ write $c_\kappa=1/(3\kappa)$;
the constant $c$ of Section~\ref{sec:v-detail-29} is $2c_\kappa$. Thus
\eqref{v:cmp:10} reads
\begin{equation}\label{eq:v-refined-count}
 \begin{gathered}
 D_\kappa(x)=3-(1+2c_\kappa)x,\qquad
 P_\kappa(x)=(2-2c_\kappa x)(1-x),\\
 \mathcal J_\kappa=(5/6-\delta)D_\kappa+\delta P_\kappa,\qquad
 R_\kappa(\delta,x)=1-\delta+\frac{(5/6-\delta)\delta P_\kappa}
                                 {2\mathcal J_\kappa}.
 \end{gathered}
\end{equation}
For $0<\delta\le5/6$, $0\le x\le1/2$ the functions $D_\kappa,P_\kappa$
are positive, and the identity
$2D_\kappa-3P_\kappa=2x(2+c_\kappa-3c_\kappa x)\ge0$ gives
\begin{equation}\label{eq:v-refined-coarse}
 \frac16\le1-\delta\le R_\kappa(\delta,x)\le1-\frac{2\delta}3 .
\end{equation}
Put $d_{\min}=1/100$ and
\begin{equation}\label{eq:v-refined-theta}
 \mathcal K=\frac{3R_\kappa(\delta,x)+1}4,\qquad
 \theta=
 \begin{cases}
 g\min\{\ell/d,\,1/(6\kappa)\},&0<\delta<1/8,\\
 g/(6\kappa),&1/8\le\delta\le5/6.
 \end{cases}
\end{equation}
Since $g\le\delta/2$ and $\kappa\ge13/18$, every such $\theta$ satisfies
\begin{equation}\label{eq:v-refined-theta-bound}
 0\le\theta\le\frac g{6\kappa}\le\frac\delta{12\kappa}
 \le\min\Bigl\{\frac\delta2,\frac{5/6}{12\cdot13/18}\Bigr\}
 \le\min\Bigl\{a-\frac12,\frac1{10}\Bigr\}.
\end{equation}
The full relative exponent of a row class at scale $d$ is
\begin{equation}\label{eq:v-refined-full-exponent}
 \begin{split}
 \mathcal E(d;\delta,x)
 ={}&a-B(\ell)-l_y(1/2+\theta)+h(z_0-1/6)-\ell/2+\ell g\\
   &+d(\mathcal K+\theta/2-z_0).
 \end{split}
\end{equation}
This is \eqref{v:cmp:12} with
$U^{R+\delta/2}\to U^{\mathcal K+\theta/2}$ and
$Y^{-a}\to Y^{-1/2-\theta}$. At
$d=h$, for $\delta\ge1/8$,
\begin{equation}\label{eq:v-refined-endpoint}
 \begin{gathered}
 E_{b,\ell,\kappa}(\delta,x)
 =C_0+\frac\delta2
 +\delta x\Bigl\{\ell+c_\kappa\Bigl(\frac h4-\frac{l_y}2\Bigr)\Bigr\}
 -\frac{3h}4\bigl(1-R_\kappa(\delta,x)\bigr),\\
 C_0=-\frac14+\frac b6+\frac{5\ell}4,
 \end{gathered}
\end{equation}
and, for $0<\delta<1/8$ when $\ell/h<1/(6\kappa)$,
\begin{equation}\label{eq:v-refined-small-endpoint}
 E_{\rm small}(\delta,x)
 =C_0+\frac\delta2+\ell\delta x\Bigl(\frac12+\frac{2\ell}h\Bigr)
       -\frac{3h}4\bigl(1-R_\kappa(\delta,x)\bigr).
\end{equation}
For \eqref{eq:v-refined-endpoint} substitute $\theta=c_\kappa g/2$,
$a=(1+\delta)/2$ and $h\mathcal K=h-(3h/4)(1-R_\kappa)$; the constant
terms collapse to $C_0$ by the definitions of $l_x,l_y,h,B(\ell)$.
For \eqref{eq:v-refined-small-endpoint} substitute $\theta=g\ell/h$
and use $h-l_y=2\ell$.

\begin{remark}[Admissible refined parameters]\label{v:rem:refined-admissible}
Call $(b,\ell,\kappa)$ \emph{admissible} if:
\begin{enumerate}
\item $b>0$, $1/6\le\ell\le1/5$, $7497/10000\le\kappa\le3/4$,
      $B(\ell)>437/500$, $l_x>\ell$ and $l_y-\ell>11b/6$;
\item $6\kappa\ell>3/4$, $\ell/h<1/(6\kappa)$,
      $\ell/h>21/100$ and $h/2<l_y<1/2$;
\item $E_{b,\ell,\kappa}(\delta,x)<0$ for
      $1/8\le\delta\le5/6$, $0\le x\le1/2$.
\item with $A_{\rm sm}=1/2+(\ell/2)(1/2+2\ell/h)-h/2$,
\begin{equation}\label{eq:v-refined-ranges}
 \begin{gathered}
 \max\{C_0,\,C_0+A_{\rm sm}/8\}<-1/250,\qquad
 C_0+(1/2+3\ell/2)/50<-1/100,\\
 13h/75-l_y/2+1/50<-1/20 .
 \end{gathered}
\end{equation}
\end{enumerate}
By continuity on the compact rectangle,
condition~(3) gives $r_0=-\max E_{b,\ell,\kappa}>0$. Choose
$\zeta>0$ with $2\zeta<\min\{r_0/2,1/1000\}$ and $\ell/(h+\zeta)>21/100$,
and put $r_*=\min\{r_0/2,1/1000\}$; thus $0<r_*\le1/1000$ and
$r_0-2\zeta>r_*$.
These admissibility conditions are numerical only. Applications of M2
also require $\beta_*(F)\le(1+\kappa)/2$; this condition is included in
Proposition~\ref{v:prop:refined-high} and verified in
Theorem~\ref{v:thm:refined-limit} by the choice $\kappa=2\beta_*(F)-1$.
Condition~(1) contains the hypotheses $\ell\le1/5$, $b\ge0$,
$l_y\ge11b/6$, $l_x>\ell$ of \eqref{v:cmp:4} and puts $B(\ell)$
in the principal region \eqref{v:cmp:7}.
The three lines of \eqref{eq:v-refined-ranges} control
the small-$\delta$ endpoint
\eqref{eq:v-refined-small-endpoint}, the floor endpoint and the
small-row bound at the end of Section~\ref{sec:v-detail-30}, respectively.
\end{remark}

\subsection{An approximate functional equation over the cyclotomic field}
\label{v:sec:refined-afe}
In this subsection, we write the approximate functional equation
\eqref{eq:v-refined-afe} for a primitive Hecke character of $F$, bound
its profiles uniformly in Lemma~\ref{v:lem:refined-afe}, and relate
$L_F^S(w,\chi_u)$ to the primitive $L(w,\psi_u^*)$ by
\eqref{eq:v-refined-imprimitive}, with the conductor deficit
\eqref{eq:v-refined-deficit}. These are used in
Proposition~\ref{v:prop:refined-numerator} and
Lemma~\ref{v:lem:refined-raised}.

Put $r_F=[F:\mathbb Q]/2$; $F$ is totally complex.

Let $\psi$ be a primitive nonprincipal finite-order Hecke character of
$F$, of conductor $\mathfrak f$. Its infinity type is trivial; with $|\Delta_F|$ the absolute
discriminant of $F$,
\[
 \Lambda(s,\psi)=\mathcal Q_\psi^{s/2}\Gamma(s)^{r_F}L(s,\psi),
 \qquad
 \mathcal Q_\psi=\frac{|\Delta_F|\,q_{\mathfrak f}}{(2\pi)^{2r_F}},
\]
is entire and satisfies $\Lambda(s,\psi)=\epsilon_\psi\Lambda(1-s,\bar\psi)$
with $|\epsilon_\psi|=1$. For $w=\sigma+it$ with $2/5\le\sigma\le3/5$
put
\[
 V_w(\upsilon)=\frac1{2\pi i}\int_{(2)}e^{\varsigma^2}
   \Bigl(\frac{\Gamma(w+\varsigma)}{\Gamma(w)}\Bigr)^{r_F}\upsilon^{-\varsigma}\frac{d\varsigma}\varsigma .
\]
Moving the line of $\int\Lambda(w+\varsigma,\psi)e^{\varsigma^2}d\varsigma/\varsigma$ from $2$ to $-2$
crosses only $\varsigma=0$, and the functional equation with $\varsigma\mapsto-\varsigma$ gives
\begin{equation}\label{eq:v-refined-afe}
 \begin{split}
 L(w,\psi)={}&\sum_{\mathfrak n}\frac{\psi(\mathfrak n)}{q_{\mathfrak n}^w}
          V_w\Bigl(\frac{q_{\mathfrak n}}{\sqrt{\mathcal Q_\psi}}\Bigr)\\
 &+\epsilon_\psi\mathcal Q_\psi^{1/2-w}
   \Bigl(\frac{\Gamma(1-w)}{\Gamma(w)}\Bigr)^{r_F}
   \sum_{\mathfrak n}\frac{\bar\psi(\mathfrak n)}{q_{\mathfrak n}^{1-w}}
          V_{1-w}\Bigl(\frac{q_{\mathfrak n}}{\sqrt{\mathcal Q_\psi}}\Bigr).
 \end{split}
\end{equation}
Gaussian decay and polynomial strip growth justify the horizontal
limits, and both sums converge absolutely; see \cite[Theorem~5.3]{IK}.

\begin{lemma}[Uniform approximate-functional-equation profiles]
\label{v:lem:refined-afe}
For fixed integers $j,N\ge0$, uniformly for $2/5\le\sigma\le3/5$,
\[
 |(\upsilon\partial_\upsilon)^jV_w(\upsilon)|\ll_{F,j,N}
 \Bigl(1+\frac \upsilon{(1+|t|)^{r_F}}\Bigr)^{-N}.
\]
For $\sigma=1/2+\theta$ with $0\le\theta\le1/10$, the function
\[
 \Bigl(\frac{\Gamma(1-w)}{\Gamma(w)}\Bigr)^{r_F}\upsilon^{2\theta}V_{1-w}(\upsilon)
\]
obeys the same bound, after increasing $N$ by one.
\end{lemma}
\begin{proof}
For $\upsilon\ge(1+|t|)^{r_F}$ move the line to a large fixed positive real
part and use Stirling's formula. For smaller $\upsilon$ move it to $-1/4$;
since $\sigma\ge2/5$ no pole of $\Gamma(w+\varsigma)$ is crossed, and the only
residue is $1$. Euler derivatives insert powers of $\varsigma$; the Gaussian
keeps the integrals uniformly convergent. Stirling's formula also gives
$|\Gamma(1-w)/\Gamma(w)|^{r_F}\ll(1+|t|)^{-2r_F\theta}$, which
cancels the factor $\upsilon^{2\theta}=(1+|t|)^{2r_F\theta}
(\upsilon/(1+|t|)^{r_F})^{2\theta}$; the product rule gives the
derivatives, uniformly for $0\le\theta\le1/10$.
\end{proof}

For a sixth-power-free nonunit row $u$ let $\psi_u^*$ be the primitive
character inducing $\chi_u=\chi_\bullet(u)$ in \eqref{v:cmp:6}, where
$\chi_u(p)=\chi_p(u)$ with the unit of $u$ included (forced by
$W=\rho q_p^{-w}$, $\rho=\chi_p(u)$, in \eqref{v:cmp:5}): the
Kummer character of $F(u^{1/6})/F$, differing by some
$\xi_\epsilon\in\Upsilon_S$ from the earlier symbol $\chi_u(x)$ modulo $u$.
At an exterior prime $p$ with $1\le v_p(u)=j\le5$ its local component
has tame inertia of order $6/(6,j)>1$ (proof of
Lemma~\ref{v:lem:general-field-detector}), so $\psi_u^*(p)=0=\chi_u(p)$.
Hence $\psi_u^*=\chi_u$ on exterior ideals, and
\begin{equation}\label{eq:v-refined-imprimitive}
 L_F^S(w,\chi_u)=L(w,\psi_u^*)\prod_{p\in S_f}(1-\psi_u^*(p)q_p^{-w}),
\end{equation}
the finite product being $\le\prod_{p\in S_f}(1+q_p^{-2/5})$ on
$\Rea w\ge2/5$. Moreover $q_{\mathfrak f_u}\ll_S q_u$, and, for every
set $J_0$ of exterior primes with $v_p(u)=j_p\ge2$,
\begin{equation}\label{eq:v-refined-deficit}
 q_{\mathfrak f_u}\ll_S q_u\prod_{p\in J_0}q_p^{-(j_p-1)},
\end{equation}
since the exterior conductor exponent is one (as used in
Section~\ref{sec:v-detail-30}).

\subsection{An amplified first moment of the reflected numerator}
\label{v:sec:refined-numerator}
In this subsection, Proposition~\ref{v:prop:refined-numerator} proves
the averaged bound \eqref{eq:v-refined-numerator}, a saving
$U^{-\theta/2}$ on the reflected value
$L(1/2+\theta+it,\bar\psi_u^*)$ over a row class, by the approximate
functional equation of Section~\ref{v:sec:refined-afe}, H\"older's
inequality and M2 (Theorem~\ref{v:thm:M2-assembled}) amplified
by the prime sums defining the class. It enters
\eqref{eq:v-refined-row-bound} and Proposition~\ref{v:prop:refined-high}.

A dyadic block of length $D=U^m$ in
\eqref{eq:v-refined-afe} at $1/2+\theta+it$ has coefficient
$D^{-\theta}$ relative to a centrally normalized plain block.
Two copies use length $2m$ in M2, leaving admissible prime length
$(1-2m)/(6\kappa)$; selected prime sums of total length $z_{\rm p}$ save
$U^{-2gz_{\rm p}}$ in the fourth moment and $U^{-gz_{\rm p}/2}$ in the
first moment after H\"older. With $z_{\rm p}=v(1-2m)$ and $\theta=gv$,
\begin{equation}\label{eq:v-refined-balance}
 -\theta m-\frac{gv(1-2m)}2=-\frac\theta2,
\end{equation}
independently of $m$.

\begin{proposition}[Amplified numerator average]
\label{v:prop:refined-numerator}
Fix $\eps>0$. Before fixing the prime positions, choose the M2 mesh
$\eta_{\rm M2}\le\eps/100$ for input loss $\eps/100$, and a
capacity-length decrement $0<\Delta\le\eps/100$; require every
individual prime length to be at most $\eta_{\rm M2}$.
Let $13/18\le\kappa\le1$; when $\kappa<1$ assume
$\beta_*(F)\le(1+\kappa)/2$. Let $\mathscr B$ be a row class as in
Section~\ref{v:sec:refined-parameters}, with $\mathscr L\ge1/5$,
$0\le g\le5/12$ (so that every $\theta$ below satisfies
$\theta\le g/(6\kappa)\le5/52<1/10$), and
\[
 \#\mathscr B\ll U^{R+\eps}(1+H)^A,
\]
where $1\le H\le U^{1/100}$ bounds the numerator height and the
height of the prime-sum profiles. Choose either
\begin{enumerate}
\item $\theta=gv$ with $0\le v\le\min\{\mathscr L,1/(6\kappa)\}$; or
\item $\theta=g/(6\kappa)$, provided
\begin{equation}\label{eq:v-refined-gate}
 R\le(6\kappa-2g)\mathscr L .
\end{equation}
\end{enumerate}
Then, uniformly for $|t|\le H$,
\begin{equation}\label{eq:v-refined-numerator}
 \sum_{u\in\mathscr B}|L(1/2+\theta+it,\bar\psi_u^*)|
 \ll U^{(3R+1)/4-\theta/2+\eps}(1+H)^{A'} .
\end{equation}
The orders $A'$ and the finite seminorms depend only on the fixed field,
$S$, $\Theta$, the system of prime positions and the requested loss.
\end{proposition}

\begin{proof}
\emph{Block moments.} Let $S_m(u)$ denote a centrally normalized plain
polynomial of length $D=U^m$ for the presentation $\bar\chi_u$, with an
annular profile from a uniformly bounded smooth family and a norm twist
of height $O(H)$. Its character is that of the prime coefficients
$\bar\chi_p(u)$ in $Q_i(u;z)$; the ray selectors are expanded in
$\Theta$ exactly as in Section~\ref{sec:v-detail-29}. Since $u$ is a
nonunit row, $\bar\psi_u^*$ induces outside $\Theta$, as recorded in
Section~\ref{sec:v-detail-26}. Only absolute squares occur in M2, so
the orientation $-1$ reduces to $+1$ by conjugating the entire product;
$\Theta=\bar\Theta$, the profile family, the property of inducing
outside $\Theta$ and the height bounds are closed under conjugation.
The two ordinary factors of M2 may carry different profiles $W_1,W_2$
from the fixed family (Section~\ref{sec:v-detail-19}), with independent
incoming norm twists of polynomially bounded height.

For $0\le m\le1/2$ request the prime length $z_{\rm p}=v(1-2m)$ in case~(1).
Then $z_{\rm p}\le\mathscr L$ and $2m+6\kappa z_{\rm p}\le1$.
Order the prime-sum factors of positive amplitude by decreasing $g_i$
and select whole factors up to the reduced length
$(z_{\rm p}-\Delta)_+$. A fractional segment has mean amplitude at
least $g$; discarding its last fractional factor loses at most
$(\delta/2)\eta_{\rm M2}$. If the positive-amplitude factors are
exhausted, their total gain is $g\mathscr L$. Consequently, on
$\mathscr B$,
\[
 |Q_{\rm sel}(u)|^2\ge
 U^{2gz_{\rm p}-\omega},\qquad
 \omega=2g\Delta+\delta\eta_{\rm M2}.
\]
When $z_{\rm p}\le\Delta$ the selected product is empty and this
inequality still holds. Apply Theorem~\ref{v:thm:M2-assembled} with
input loss $\eps/100$, in the row normalization of
Section~\ref{v:sec:compensated-assembly} ($q_k\ll U$, lengths relative
to $U$), with both ordinary factors equal to $S_m$ (plain lengths
$n_1=n_2=m$) and the selected product as its prime-sum factor
($z=z_{\rm p}$). Its length condition $n_1+n_2+6\kappa z\le M$ is
$2m+6\kappa z_{\rm p}\le1$, its exclusion of $\Theta$ is the
restriction to rows inducing outside $\Theta$, and its polynomial
dependence on norm-twist heights gives the factor $(1+H)^{A_1}$; when
$z_{\rm p}=0$ the $z=0$ clause applies, and the sum over rows inducing
outside $\Theta$ is bounded by the sum over all nonprincipal rows. This
gives
\[
 \sum_{u\in\mathscr B}|S_m(u)|^4
 \le U^{-2gz_{\rm p}+\omega}\sum_{\substack{q_k\ll U\\ \chi_\bullet(k)\text{ induces outside }\Theta}}
   |S_m(k)^2Q_{\rm sel}(k)|^2
 \ll U^{1-2gz_{\rm p}+\omega+\eps/100}(1+H)^{A_1}.
\]
Since $g\le5/12$ and $\delta\le5/6$, the choices in the statement give
$\omega\le\eps/60$. The count contributes $3\eps/4$ after
H\"older, and the fourth moment contributes at most
$(\omega+\eps/100)/4$; their sum is less than $4\eps/5$.
The remaining exponent allowance absorbs the logarithmic dyads and
subsidiary losses. The decrement used in
Section~\ref{sec:v-detail-29} may be chosen as $\Delta$ here.
H\"older's inequality and the count of $\mathscr B$ give
\begin{equation}\label{eq:v-refined-block}
 \sum_{u\in\mathscr B}|S_m(u)|
 \le(\#\mathscr B)^{3/4}\Bigl(\sum_{u\in\mathscr B}|S_m(u)|^4\Bigr)^{1/4}
 \ll U^{(3R+1)/4-gz_{\rm p}/2+\eps}(1+H)^{A_2}.
\end{equation}

For $1/2\le m\le1/2+r_F/100+\eps$, the $z=0$ clause of
Theorem~\ref{v:thm:M2-assembled}, which excludes only the principal
inducing character and allows arbitrary bounded plain lengths with
right-hand side $U^{1+\eps}$, gives \eqref{eq:v-refined-block} with
$z_{\rm p}=0$; this is the fourth moment of
Remark~\ref{v:rem:sextic-fourth-moment} in the generality stated there.
The large sieve alone would give $U^{2m}$. Where the
decremented length vanishes, use the same estimate without prime
factors, at arbitrarily small loss.

\emph{The approximate functional equation.} Apply
\eqref{eq:v-refined-afe} to $L(1/2+\theta+it,\bar\psi_u^*)$. Write
$\mathfrak n=\mathfrak a\mathfrak b$ with $\mathfrak a$ supported on
$S_f$ and $\mathfrak b$ exterior, so that
$\bar\psi_u^*(\mathfrak b)=\bar\chi_u(\mathfrak b)$. Apply the
triangle inequality over $\mathfrak a$ before any moment; the
$\mathfrak a$-coefficients have total mass at most
$\prod_{p\in S_f}(1-q_p^{-2/5})^{-1}$, and the possible row-dependent
value $\bar\psi_u^*(\mathfrak a)$ disappears. For $q_{\mathfrak b}\asymp D$
the first sum has the scalar $q_{\mathfrak a}^{-1/2-\theta}D^{-\theta}$
times a central plain block, and the dual sum has
\[
 q_{\mathfrak a}^{\theta-1/2}\mathcal Q_u^{-\theta}D^{\theta}
 =q_{\mathfrak a}^{-1/2-\theta}D^{-\theta}
       \Bigl(\frac{q_{\mathfrak a}D}{\sqrt{\mathcal Q_u}}\Bigr)^{2\theta},
 \qquad\mathcal Q_u=\mathcal Q_{\psi_u^*},
\]
times a central plain block with character $\bar{\bar\chi}_u=\chi_u$.

By Lemma~\ref{v:lem:refined-afe} the last factor times the
gamma ratio is a uniformly controlled profile. For the dual
block use $|S|=|\bar S|$: $\bar S$ has the
character $\bar\chi_u$ of the prime factors, with conjugated profile
and twist, and only $|Q_{\rm sel}|$ enters. The factors
$(q_{\mathfrak b}/D)^{\mp it}$ are separated norm twists of height $|t|\le H$.

\emph{Rowwise conductors.} The profiles depend on $u$ through
$\lambda_u=\log(D/\sqrt{\mathcal Q_u})$. For one dyad and one selected
prime product, treat $\lambda$ as an auxiliary real parameter with the
character fixed, partition the rows by the
unit interval containing $\lambda_u$, and apply
Lemma~\ref{v:lem:detector-Sobolev} with $\lambda$ as one additional
bounded coordinate on each interval to $S_{m,\lambda}(u)^2Q_{\rm sel}(u)$. Its $\lambda$-derivative is
$2S_{m,\lambda}\partial_\lambda S_{m,\lambda}Q_{\rm sel}$; its square
is an instance of M2 with two different profiles of the same length,
which M2 permits. The selected primes do not depend on $\lambda$,
and Lemma~\ref{v:lem:refined-afe} gives uniform Euler derivatives of
the profiles.

As $u$ is sixth-power-free, nonunit, $q_u\asymp U$, its
exterior conductor is $\ge\prod_{p\mid u}q_p\ge q_u^{1/5}$; thus
$U^{1/5}\ll\mathcal Q_u\ll U$, only $O(\log U)$ unit
intervals occur, and root numbers are removed by absolute
values. Rowwise conductors cost a logarithmic factor.

\emph{Tails.} By Lemma~\ref{v:lem:refined-afe}, dyads with
$D>U^{1/2}(1+H)^{r_F}U^{\eps}$ contribute $O(U^{-N})$ for every $N$,
by ideal counting and $\#\mathscr B\ll U$. Since $H\le U^{1/100}$,
all retained lengths satisfy
$0\le m\le1/2+r_F/100+\eps$. Lengths $m>1/2$ use only M2 without
prime sums (no mesh),
and every block with prime sums has $2m+6\kappa z_{\rm p}\le1$; so the
mesh of Lemma~\ref{v:lem:arithmetic-independent-mesh} is unaffected.

\emph{Supply-limited choice.} For $0\le m\le1/2$,
\eqref{eq:v-refined-block} with $z_{\rm p}=v(1-2m)$ and the coefficient
$D^{-\theta}$ give the exponent
$(3R+1)/4-\theta m-gv(1-2m)/2=(3R+1)/4-\theta/2$ by
\eqref{eq:v-refined-balance}. For $m\ge1/2$ the coefficient alone gives
$-\theta m\le-\theta/2$. Summing $O(\log U)$ dyads proves
\eqref{eq:v-refined-numerator} in case~(1).

\emph{Maximal-length choice.} If $\mathscr L\ge1/(6\kappa)$, case~(2)
is case~(1) with $v=1/(6\kappa)$. Otherwise put
$m_0=1/2-3\kappa\mathscr L\in(0,1/2)$. For $m_0\le m\le1/2$ the full
admissible length $z_{\rm p}=(1-2m)/(6\kappa)$ is at most $\mathscr L$, and the preceding
argument applies. For $0\le m\le m_0$ ideal counting gives
$|S_m(u)|\ll U^{m/2+\eps}(1+H)^C$, so the block contributes exponent
$R+(1/2-\theta)m\le R+(1/2-\theta)m_0$, since $\theta<1/2$, and
\[
 R+(1/2-\theta)m_0-\Bigl(\frac{3R+1}4-\frac\theta2\Bigr)
 =\frac{R-(6\kappa-2g)\mathscr L}4\le0
\]
by \eqref{eq:v-refined-gate}. The dual profiles satisfy the same
absolute bound. This proves case~(2).
\end{proof}

The same amplitudes enter the row count of
Section~\ref{sec:v-detail-29} and
Proposition~\ref{v:prop:refined-numerator}, each time as
pointwise lower bounds on $|Q_{\rm sel}|$ over the class inside
a separate positive moment; no prime product is inserted twice in the
theta average.

\subsection{The raised numerator line and its local errors}
\label{v:sec:refined-local}
In this subsection, we bound the main/error expansion of
\eqref{eq:v-global-compensated-correction} on the raised
numerator line \eqref{eq:v-refined-lines}, keeping the reflected value
for Proposition~\ref{v:prop:refined-numerator}.
Lemma~\ref{v:lem:refined-raised} bounds each expansion term, and
\eqref{eq:v-refined-row-bound} is the resulting row-class estimate,
which replaces \eqref{v:cmp:11} in Section~\ref{v:sec:refined-high}.

On a detected class take
\begin{equation}\label{eq:v-refined-lines}
 \Rea s=a+16e,\qquad\Rea w=\frac12-\theta,\qquad\Rea z=z_0,
 \qquad0\le\theta\le\min\{a-1/2,1/10\}.
\end{equation}
Then $\Rea w\ge1-a$, so these lines lie in the connecting box
$\Rea s\ge a+16e$, $\Rea w\ge1-a-6e$, $\Rea z\ge17/50$ of
Section~\ref{sec:v-detail-30}, and $\Rea(s+w)\ge1+16e$. As shown
there, the local bounds at the corner $w=1-a-6e$ persist when a real
coordinate increases (after clearing $(1-R)(1-D)$ every numerator
monomial has a nonnegative power of $Y_0=q_p^{-w}$). Hence on
\eqref{eq:v-refined-lines}:
\begin{enumerate}
\item at good primes, $H_p-1=O(q_p^{-1-10e})$ and
      $G_p+\bar\chi_p(u)H_p=O(q_p^{-a-4e})$, and bad defects are
      $O(q_p^{-10e})$;
\item at a selected prime $p\mid u$ with $v_p(u)=j$, after separating,
      for $j\ge2$, the strict term $-(q_p-1)q_p^{z_0-1-w}/(1-R)$, the
      remaining terms are $O(q_p^{z_0-1/2-\sigma_j})$ with the table
      of $\sigma_j>0$ in Section~\ref{sec:v-detail-30}.
\end{enumerate}
The strict term has size $\asymp q_p^{z_0-1/2+\theta}$; on the
pointwise line it had size $q_p^{z_0-1+a+6e}$.

\begin{lemma}[Local errors with the reflected value retained]
\label{v:lem:refined-raised}
Fix a finite tuple term of \eqref{eq:v-global-compensated-correction}
on \eqref{eq:v-refined-lines}, with error indices $I$, strict indices
$J_0\subseteq I$ and main indices $I^c$. After the main factors
$H_{p_i}$, $i\notin I$, are merged into the unselected product as in
Section~\ref{sec:v-detail-30}, the term is bounded by
\[
 \ll U^{\theta+\eps}(1+T_1)^{A_F}\,
   |L(1-w,\bar\psi_u^*)|\prod_{i\notin I}P_i^{z_0-1/2}|Q_i(u;z)|
   \prod_{i\in I}P_i^{z_0-1/2+\eps}
\]
for the product with $L_F^S(w,\chi_u)$, summed over the error choices
inside each row with divisor mass $q_u^\eps$.
\end{lemma}

\begin{proof}
The main windows give $-P_i^{z-1/2}Q_i(u;z)$ exactly, as in
Section~\ref{sec:v-detail-30}. An error at a good prime, summed over
its window, is $\ll P_i^{z_0-a-4e}\le P_i^{z_0-1/2}$ by~(1). The
unselected product is $\ll q_u^\eps$ by~(1). A nonstrict error at a
bad prime is $O(q_p^{z_0-1/2})$ by~(2); bad primes divide $u$, so
their choices are divisor-many.

By \eqref{eq:v-refined-imprimitive} (bounded $S_f$ factor) and the
functional equation,
\[
 |L_F^S(w,\chi_u)|\ll\mathcal Q_u^{\theta}
 \Bigl|\frac{\Gamma(1-w)}{\Gamma(w)}\Bigr|^{r_F}|L(1-w,\bar\psi_u^*)|
 \ll q_{\mathfrak f_u}^{\theta}(1+T_1)^{A_F}|L(1-w,\bar\psi_u^*)|.
\]
Apply \eqref{eq:v-refined-deficit} once to the complete set $J_0$ of
strict labels of the tuple (distinct primes, as the
windows are disjoint). At each label
\[
 |q_p^{z_0-w}|\,q_p^{-(j_p-1)\theta}
 =q_p^{z_0-1/2-(j_p-2)\theta}\le q_p^{z_0-1/2},
\]
since $j_p\ge2$. This single simultaneous deficit bound gives
$q_{\mathfrak f_u}^\theta
\prod_{J_0}q_p^{(j_p-1)\theta}\ll U^\theta$; multiplying the
factors proves the lemma.
\end{proof}

Lemma~\ref{v:lem:refined-raised}, \eqref{eq:v-refined-row-bound} and
Proposition~\ref{v:prop:refined-numerator} are applied only after the
original heights have been truncated, as in
Section~\ref{v:sec:refined-contours}; then $|\Im w|,|\Im z|$ and the
prime-sum profile heights are $O_e(T_1)$, and one may take
$H\asymp(3I+2)T_1\le U^{1/100}$. Outside the retained boxes the global
strip bounds and the tail estimate $O_N(Z^{D_\eta}T_1^{-N})$ are used.

For a main/error term the main amplitude is
$\sum_{i\notin I}\ell_ig_i\le\ell g$, while the class (hence
Proposition~\ref{v:prop:refined-numerator}) uses all actual
amplitudes. Hence every expansion term satisfies
\begin{equation}\label{eq:v-refined-row-bound}
 \sum_{u\in\mathscr B}\bigl|L_F^S(w,\chi_u)
                        \mathcal H_u^{\rm piece}\bigr|
 \ll U^{\mathcal K+\theta/2+\eps}
     Z^{\ell(z_0-1/2)+\ell g+\eps}(1+T_1)^{A_F'},
\end{equation}
which replaces \eqref{v:cmp:11}.

\subsection{Amplitude-dependent numerator contours}
\label{v:sec:refined-contours}
In this subsection, we justify the row-dependent numerator line
$\Rea w=1/2-\theta$ of \eqref{eq:v-refined-lines}: after the $s$- and
$z$-contours of Section~\ref{sec:v-detail-30} are fixed, $w$ is moved
row by row, and only then are the heights truncated. This contour
assembly is used in the proof of Proposition~\ref{v:prop:refined-high}.

The amplitude classes depend on $z$ (through
$Q_i(u;z)$) and on the detector bin of $u$, not on $s,w$.

Only the nonprincipal intermediate step of
Section~\ref{sec:v-detail-30} changes, for rows with $d_{\min}<d\le h+\zeta$
in bins $a>a_f$. The floor bin keeps
$\Rea w=1-a_f-6e$; small, far and principal rows keep the lines listed
at the end of Section~\ref{sec:v-detail-30}. As there, fix the entire original bin
and the finite tuple expansion, and move $z$ and $s$ to $\Rea z=z_0$
and the global line $\Rea s=\sigma_g=\beta+20e$, starting from the
absolutely convergent lines; no pole is crossed.

Now fix $(s,z)$ on
these lines and a finite intermediate row dyad. Deterministic rounding
of $(|Q_i(u;z)|)_i$, with each $g_i$ capped at $\delta/2$, assigns each
row a class. Put $\theta_u(z)=\min\{\theta,\,a-1/2,\,1/10\}$ for the
value $\theta$ of \eqref{eq:v-refined-theta} at that class when
$|\Im z|\le T_1$, and $\theta_u(z)=0$ when $|\Im z|>T_1$, where only
the tail estimates are used. Then $\theta_u(z)$ satisfies
\eqref{eq:v-refined-lines}, equals the
class value on the retained heights by \eqref{eq:v-refined-theta-bound},
and is Borel and piecewise constant in $z$.

For each row separately move
$w$ from its right line to $\Rea w=1/2-\theta_u(z)$. The integrand of
one row is holomorphic in $w$ on the swept strip: the numerator is
entire (the row is nonprincipal), the denominator is independent
of $w$, and
\[
 \Rea(s+w)\ge\beta+20e+\frac12-\Bigl(a-\frac12\Bigr)\ge1+20e,\qquad
 \Rea w\ge1-a\ge1-a-6e,
\]
since $a\le\beta$. So the strip lies in the connecting box, where
the actual correction is holomorphic. The reciprocal is on its global
line, and $\widehat W_1(w)$ decays rapidly while the numerator and
correction grow polynomially, so the horizontal joins vanish.
The row sum is finite, Fubini on the finite retained boxes
justifies the $z$-integration, and no $z$-contour moves
after the partition.

Only then truncate the original heights and move the retained
$s$-segment from $\sigma_g$ to $a+16e$ in the buffered rectangle of its
bin. The $w$-line of a row is independent of $s$; on $s$-joins the
other axes use global numerator and correction bounds and the buffered
denominator, as in Section~\ref{sec:v-detail-30}. No merely
buffered $s$-axis is extended to all heights, and every added frequency
receives a share of the common $T_1/2$ allowance. Tails
and traces are $O_N(Z^{D_\eta}T_1^{-N})$ as before.

\subsection{Monotonicity and the refined remainder}
\label{v:sec:refined-high}
In this subsection, we combine \eqref{eq:v-refined-row-bound} with the
outside factors of \eqref{v:cmp:6} to obtain the exponent
\eqref{eq:v-refined-full-exponent}, show in
Lemma~\ref{v:lem:refined-monotone} that it increases in the row scale
$d$, so that the endpoint $d=h$ controls, and prove in
Proposition~\ref{v:prop:refined-high} the refined remainder
\eqref{eq:v-refined-remainder}, the analogue of \eqref{v:cmp:13} used
in Lemma~\ref{v:lem:refined-closure}.

On a class with scale $d$ and $\Rea s=a+16e$, multiply
\eqref{eq:v-refined-row-bound} by the outside factors in
\eqref{v:cmp:6}, namely
$X^{1/2-z_0}Z^{a+z_0-1}Y^{-1/2-\theta}q_u^{-z_0}$, and normalize by
$c_SA_{\mathcal T}$, which costs $Z^{\eps}$. Using $h=1-l_x+\ell$ and
\[
 \frac{l_x}2-1+\frac{4+b}6=-\frac h6,
\]
the resulting exponent relative to $C(B(\ell))$ is
\eqref{eq:v-refined-full-exponent}, up to the buffer $16e$ and small
losses. The factor $U^{\theta}$ from Lemma~\ref{v:lem:refined-raised},
the average $U^{-\theta/2}$ and the outside $Y^{-\theta}$ combine to
$Z^{\theta(d/2-l_y)}$; at $d=h$ this is a saving because $h<2l_y$.

\begin{lemma}[Monotonicity in the row scale]
\label{v:lem:refined-monotone}
For admissible parameters and fixed $(\delta,x)$, the exponent
$\mathcal E(d;\delta,x)$ is continuous and increasing in $d$, with
one-sided slopes between $7/200$ and $1$.
\end{lemma}
\begin{proof}
By \eqref{eq:v-refined-coarse}, $\mathcal K-z_0\ge(3\cdot\frac16+1)/4-17/50=7/200$.
Where $\theta=g/(6\kappa)$ the slope is $\mathcal K-z_0+g/(12\kappa)$.
Where $\theta=g\ell/d$ the slope is $\mathcal K-z_0+g\ell l_y/d^2$;
this regime has $d\ge6\kappa\ell>3/4$. The two expressions for
$\theta$ agree at $d=6\kappa\ell$, which gives continuity. Using
$R_\kappa\le1$, $g\le5/12$, $\ell\le1/5$, $l_y<1/2$ and
$\kappa\ge7497/10000$, both slopes are below $1$.
\end{proof}

\begin{proposition}[Refined high-side remainder]
\label{v:prop:refined-high}
Let $(b,\ell,\kappa)$ be admissible, with $r_0,\zeta,r_*$ as in
Remark~\ref{v:rem:refined-admissible}, and suppose
$\beta=\beta_*(F)\le(1+\kappa)/2$ and $\Delta=\beta-B(\ell)>0$.
Then the choices of Section~\ref{sec:v-detail-30}, with the nonprincipal
margin $\mathfrak r_{\rm np}=r_*/8$, give a fixed even system of prime positions and
$m_{\rm save}>0$, independent of the target, such that the normalized theta average
and principal term satisfy
\begin{equation}\label{eq:v-refined-remainder}
 |J_\eta-f_\eta|\ll Z^{C(\beta)-m_{\rm save}}(1+T_1)^{A_\eta}+Z^{D_\eta}T_1^{-N},
\end{equation}
with the order of choices in Remark~\ref{v:rem:parameter-ledger}:
the admissible geometry and nonprincipal losses precede $K$,
the principal saving and buffer follow $K$, and each target's final
$S$ precedes its tests and internal orders.
\end{proposition}

\begin{proof}
The assumed inequality $\beta_*(F)\le(1+\kappa)/2$ and the admissible
range of $\kappa$ give M2 by Theorem~\ref{v:thm:M2-assembled}.
Since $\beta\le(1+\kappa)/2\le7/8$, every bin has $\delta\le3/4$.

\emph{Hypotheses of the numerator average.} For
$d_{\min}\le d\le h+\zeta$, admissibility gives $\mathscr L=\ell/d>21/100$,
so \eqref{v:cmp:10} applies with $R=R_\kappa(\delta,x)$, and the
selected lengths fit the supply. For $\delta\ge1/8$, using $g\le\delta/2$,
\eqref{eq:v-refined-coarse} and $\kappa\ge7497/10000$,
\[
 (6\kappa-2g)\mathscr L-R_\kappa
 \ge\frac{21}{100}(6\kappa-\delta)-\Bigl(1-\frac{2\delta}3\Bigr)
 \ge\frac{1279}{750000}>0,
\]
the minimum being at $\delta=1/8$. Thus case~(2) of
Proposition~\ref{v:prop:refined-numerator} applies for $\delta\ge1/8$,
and case~(1) for $\delta<1/8$. Moreover \eqref{eq:v-refined-theta-bound} gives
\eqref{eq:v-refined-lines}. The height condition
$H\le U^{1/100}$ holds for large $Z$ because
$\tau_\eta<d_{\min}/10^4$.

\emph{Intermediate rows above the floor.} For $\delta\ge1/8$,
Lemma~\ref{v:lem:refined-monotone} and admissibility give, at every
$d\le h+\zeta$, $\mathcal E\le E_{b,\ell,\kappa}+2\zeta\le-r_0+2\zeta<-r_*$.
For $0<\delta<1/8$, admissibility gives $\ell/h<1/(6\kappa)$, so the
endpoint is \eqref{eq:v-refined-small-endpoint}; by
\eqref{eq:v-refined-coarse} and $x\le1/2$,
$E_{\rm small}\le C_0+A_{\rm sm}\delta$, which is at most
$\max\{C_0,C_0+A_{\rm sm}/8\}<-1/250$. Adding $2\zeta<1/1000$ leaves
a saving $3/1000>r_*$; this uses $r_*\le1/1000$.

The item ``Intermediate rows use $R=76/75-2\delta/3$'' at the end of
Section~\ref{sec:v-detail-30} is not needed: by
Lemma~\ref{v:lem:refined-monotone} the classes above cover all
$d_{\min}\le d\le h+\zeta$ with the count $R_\kappa$ ($\mathscr L>21/100$).

\emph{Remaining rows.} The floor bin, the small rows $d\le d_{\min}$
(including bounded nonprincipal unit rows), the far rows
$q_u>Z^{h+\zeta}$ and the principal row are estimated as at
the end of Section~\ref{sec:v-detail-30}, with the contours listed there.
The floor endpoint $C_0+(1/2+3\ell/2)/50$ increases by $\le
2\zeta<1/1000$ for $d\le h+\zeta$, so the
second line of \eqref{eq:v-refined-ranges} leaves a saving
$>1/100-1/1000$. The small-row bound $13h/75-l_y/2+1/50$, relative
to $C(\beta)$, is $<-1/20$ by the third line. Far rows are
negligible by the choice of $z_\infty$. The principal savings
\eqref{eq:v-principal-buffered-savings} are positive since $B(\ell)>437/500$.

\emph{Assembly.} In Section~\ref{sec:v-detail-27} the finite list
$\mathscr E$ of nonprincipal exponent functions and the constant
$C_{\rm np}$ are taken to include the refined functions: $\mathcal E$ in
both regimes of $\theta$, $E_{\rm small}$, the floor endpoint and the condition $(6\kappa-2g)\mathscr L-R_\kappa$, and the per-dyad numerator
exponents $(3R_\kappa+1)/4-gz_{\rm p}/2-\theta m$ and $R_\kappa+(1/2-\theta)m$
of Proposition~\ref{v:prop:refined-numerator}, with $\theta$ regarded as
a function of the rounded amplitude $g$. These are rational and affine
with bounded derivatives on the buffered ranges ($\mathcal J_\kappa>0$
there), so $C_{\rm np}$ remains a numerical
constant fixed before $K$. The amplitude-bin mesh $\eps_{\rm bin}$
(which also moves the numerator line through $\theta$, with Lipschitz
constant at most $l_y+d<2$ in $\mathcal E$), the length decrement
$\epsilon_c/100$ and the fractional-factor loss $(\delta/2)\eta_{\rm M2}$
are counted in $\epsilon_{\rm np}$.

The nonprincipal exponents are at most
$C(B(\ell))-r_*\le C(\beta)-r_*$, before the buffer and discretization
losses, which are allocated within $\mathfrak r_{\rm np}=r_*/8$ by the
choices of Section~\ref{sec:v-detail-27}. With $\omega_P,\mathfrak r_P$
and $m_{\rm save}=\min\{\mathfrak r_{\rm np},\mathfrak r_P\}/16$ as in
Section~\ref{sec:v-detail-30}, with $e$ chosen after $K$ but before
the target-dependent data, and $\tau_\eta,N$ chosen after the internal
orders, the contour assembly of
Section~\ref{v:sec:refined-contours} gives \eqref{eq:v-refined-remainder}.
All margins are independent of the target.
\end{proof}

\begin{lemma}[Identification in the refined comparison]
\label{v:lem:refined-physical-comparison}
The refined estimates replace only the nonprincipal numerator bound
in Lemma~\ref{v:lem:physical-comparison}. Their application requires
the full moment estimate in Theorem~\ref{v:thm:M2-assembled} for the
actual detected row classes and the numerator estimate of
Proposition~\ref{v:prop:refined-numerator}, in addition to the common
theta/Mellin and contour identities. They preserve its
$J_\eta,f_\eta,H_\eta$ and normalizer. They apply to the literal
primal and dual numerator profiles and to every strict-error tuple.
\end{lemma}
\begin{proof}
In Proposition~\ref{v:prop:refined-numerator}, a row class consists
of the original rows inducing outside $\Theta$. With lengths relative
to $U$, its two plain lengths are $m,m$, its moving width is one,
and its selected prime length satisfies $2m+6\kappa z_{\rm p}\le1$.
Its profile and profile derivative can differ: both are permitted by
Lemma~\ref{v:lem:M2-returned-family}. When the selected product is
empty, the ordinary clause excludes only the principal character,
so deletion to this row class is still allowed. This also supplies
the longer ordinary dyads of the approximate functional equation.

For each row retain its actual conductor $\mathcal Q_u$. The dual
factor is, exactly,
\[
 q_{\mathfrak a}^{\theta-1/2}\mathcal Q_u^{-\theta}D^\theta
 =q_{\mathfrak a}^{-1/2-\theta}D^{-\theta}
       (q_{\mathfrak a}D/\sqrt{\mathcal Q_u})^{2\theta}.
\]
It is not obtained by substituting a common conductor for
$\mathcal Q_u$. The uniform Euler derivatives of
Lemma~\ref{v:lem:refined-afe} control this last factor with the dual
profile. The rowwise conductor parameter is handled by the additional
Sobolev coordinate in the proof of the numerator proposition; its
derivative gives two independent ordinary profiles in M2.
The $S_f$-supported factors are summed before the moment, and their
absolute mass is bounded there. Conjugation aligns the dual character
with the selected primes. The weighted profile bounds integrate the
remaining heights; dyads and tails use the reserved exponent shares.

For a strict-error tuple, use \eqref{eq:v-refined-deficit} for its
entire set $J_0$, once. Distinct selected windows give distinct primes,
and $j_p\ge2$, so the product calculation in
Lemma~\ref{v:lem:refined-raised} cancels every strict loss against
this simultaneous conductor deficit. The reflected value remaining
is the original $L(1-w,\bar\psi_u^*)$, to which the numerator average
just described applies.

Finally, fix $s,z$ and the finite original row dyad before assigning
amplitude classes. Section~\ref{v:sec:refined-contours} moves $w$
row by row in its holomorphic strip, without moving $z$ after that
partition. Only after height truncation does the retained $s$ segment
move in its buffered rectangle. Thus no analytic continuation of a
class selector is used. The same full tuple identity is preserved,
and the principal rows and residues are unchanged. Proposition~\ref{v:prop:refined-high}
therefore bounds the difference of the same two functions in
Lemma~\ref{v:lem:physical-comparison}, with its new positive saving.
The common-boundary continuation in that lemma applies verbatim.
\end{proof}

\subsection{Closure and the limiting refined bound}
\label{v:sec:refined-closure}
In this subsection, Lemma~\ref{v:lem:refined-closure} shows that each
admissible choice, justified by a previous bound, gives the zero-free
boundary $B(\ell)$ by the continuation argument of
Section~\ref{sec:v-detail-31}. Lemmas~\ref{v:lem:refined-square},
\ref{v:lem:refined-critical} and~\ref{v:lem:refined-family} locate the
critical class $(\delta_c(\kappa),1/2)$ and the tangent parameters, and
Theorem~\ref{v:thm:refined-limit} then gives $B_\flat$ of
\eqref{eq:v-refined-limit} by a single step justified by Theorem~\ref{v:thm:original-target}.

\begin{lemma}[Closure for an admissible refined choice]
\label{v:lem:refined-closure}
If $(b,\ell,\kappa)$ is admissible and
$\beta_*(F)\le(1+\kappa)/2$, then $\beta_*(F)\le B(\ell)$.
\end{lemma}

\begin{proof}
Suppose $\beta=\beta_*(F)>B=B(\ell)$. Fix
$0<\epsilon_{\rm low}<(\beta-B)/4$. Condition~(1) of
Remark~\ref{v:rem:refined-admissible} gives the direct bound
\eqref{v:cmp:4}, hence $|J_\eta|\ll Z^{C(B)+\epsilon_{\rm low}}$ after
normalization. Proposition~\ref{v:prop:refined-high}, identified with this same
comparison by Lemma~\ref{v:lem:refined-physical-comparison}, gives
\eqref{eq:v-refined-remainder} with the saving $m_{\rm save}$, and the late height choices give
a saving $m_{\rm save}/2$. Put
$\sigma_{\rm cont}=\max\{B+\epsilon_{\rm low},\beta-m_{\rm save}/2\}<\beta$.
Then $|f_\eta(Z)|\ll_\eta Z^{\sigma_{\rm cont}-c_0}$ for $Z\ge1$, with
$c_0=(4+b)/6$.

As in Section~\ref{sec:v-detail-31}, Mellin inversion
continues $e^{(s-5/6)^2}H_\eta(s)/L_F^S(s,\eta)$ holomorphically to
$\Rea s>\sigma_{\rm cont}$. Since $\sigma_{\rm cont}\ge B>437/500$, the
factor $H_\eta$ is holomorphic and nonzero there by \eqref{v:cmp:7}.
Thus every $L_F^S(s,\eta)$ is zero-free on a common half-plane to the
left of $\beta$, contradicting the definition of $\beta$.
\end{proof}

Condition~(3) of Remark~\ref{v:rem:refined-admissible} for an explicit
triple is checked by the following criterion.

\begin{lemma}[Completed-square criterion]\label{v:lem:refined-square}
Let $\mathcal A,\mathscr B,\mathcal C$ be real functions of $y\ge0$ with
$\mathcal A(y)>0$ and $4\mathcal A(y)\mathcal C(y)-\mathscr B(y)^2\ge0$.
Then $\mathcal A(y)\delta^2+\mathscr B(y)\delta+\mathcal C(y)\ge0$ for all
real $\delta$ and $y\ge0$. Consequently, if at a triple $(b,\ell,\kappa)$
one has $2\mathfrak L\mathcal J_\kappa(-E-r_0)=\mathcal A(y)\delta^2
+\mathscr B(y)\delta+\mathcal C(y)$ with constants $\mathfrak L>0$,
$r_0\ge0$, where $E$ is \eqref{eq:v-refined-endpoint}, $y=1/2-x$ and
$\mathcal J_\kappa>0$, then $E\le-r_0$ for $0\le\delta\le5/6$ and
$0\le x\le1/2$.
\end{lemma}

\begin{proof}
This is the identity
$4\mathcal A(\mathcal A\delta^2+\mathscr B\delta+\mathcal C)
=(2\mathcal A\delta+\mathscr B)^2+(4\mathcal A\mathcal C-\mathscr B^2)$.
\end{proof}

Put
\[
\begin{aligned}
 \mathrm{cr}(\delta,\kappa)&=(3-324\kappa^2)\delta^2+(756\kappa^2-54\kappa-5)\delta
             +40\kappa-300\kappa^2,\\
 \ell^*(\delta,\kappa)&=\frac{324\delta^2\kappa^2-72\delta^2\kappa+3\delta^2
     -684\delta\kappa^2+126\delta\kappa-5\delta+450\kappa^2-60\kappa}
    {324\delta^2\kappa^2-144\delta^2\kappa+15\delta^2-2340\delta\kappa^2
     +450\delta\kappa-25\delta+2250\kappa^2-300\kappa}.
\end{aligned}
\]

\begin{lemma}[The critical class]\label{v:lem:refined-critical}
Let $13/18\le\kappa\le1$.
\begin{enumerate}
\item The equation $\mathrm{cr}(\delta,\kappa)=0$ has exactly one root
      $\delta_c(\kappa)$ in $[1/8,5/6]$, and $23/50<\delta_c(\kappa)<12/25$.
\item For all $b,\ell$, the endpoint at the class $(\delta_c(\kappa),1/2)$
      is independent of $b$ and equals
      \[
       E_{b,\ell,\kappa}(\delta_c,1/2)=\lambda(\kappa)\bigl(\ell-\ell^*(\kappa)\bigr),
       \qquad\ell^*(\kappa)=\ell^*(\delta_c(\kappa),\kappa),
      \]
      with $\lambda(\kappa)=\partial_\ell E_{b,\ell,\kappa}(\delta_c,1/2)\ge1$.
\item The function $F_\flat(\kappa)=11/12-\ell^*(\kappa)/4$ is increasing
      on $[7497/10000,1]$, and $G(\kappa)=(1+\kappa)/2-F_\flat(\kappa)$ is
      increasing on $[13/18,1]$, with a unique zero $\kappa_\flat$. It
      satisfies
      \[
       \begin{gathered}
       16686\kappa_\flat^4-27405\kappa_\flat^3+9936\kappa_\flat^2+1290\kappa_\flat-275=0,\\
       0.749773294079928<\kappa_\flat<0.749773294079929,
       \end{gathered}
      \]
      and $F_\flat(\kappa_\flat)=(1+\kappa_\flat)/2=B_\flat$.
\end{enumerate}
\end{lemma}

\begin{proof}
At $x=1/2$, \eqref{eq:v-refined-endpoint} is affine in $(b,\ell)$ for
fixed $(\delta,\kappa)$, since $C_0,h,l_y$ are affine and $R_\kappa$ does
not depend on $(b,\ell)$. A direct computation gives
\[
 \begin{gathered}
 \partial_bE_{b,\ell,\kappa}(\delta,1/2)
 =\frac{\mathrm{cr}(\delta,\kappa)}
       {24\kappa(54\delta\kappa-6\delta-75\kappa+10)},\\
 54\delta\kappa-6\delta-75\kappa+10=-36\kappa\mathcal J_\kappa(\delta,1/2)<0,
 \end{gathered}
\]
so $\partial_bE=0$ exactly on $\mathrm{cr}=0$. The leading coefficient
$3-324\kappa^2$ is negative, so $\mathrm{cr}(\cdot,\kappa)$ is concave.
Exact root counting in $\kappa$ shows that, on $[13/18,1]$,
$\mathrm{cr}$ is negative at $\delta=1/8$ and $\delta=23/50$ and
positive at $\delta=12/25$ and $\delta=5/6$; since the positive set of a
concave function is an interval, this proves~(1). Solving the affine
equation $E=0$ at $\delta=\delta_c$ for $\ell$ gives $\ell^*$, and
$\lambda=\partial_\ell E$ is the displayed rational function of
$(\delta,\kappa)$ divided by the negative quantity above; this proves~(2)
once $\lambda\ge1$ is known.

The remaining inequalities are one-variable statements along
$\delta=\delta_c(\kappa)$, verified by interval arithmetic:
$[13/18,1]$ is divided into $30000$ intervals with exact rational
endpoints. Each endpoint is enclosed by outward-rounded interval
division of integers, and consecutive boxes are checked to overlap;
the endpoint sign checks likewise enclose the exact rational points.
On each box, an interval
$\Delta$ containing $\delta_c$ is verified by the signs of $\mathrm{cr}$
at its ends, and enclosures of $\lambda$, of
$d\ell^*/d\kappa=\partial_\kappa\ell^*-\partial_\delta\ell^*\,
\partial_\kappa\mathrm{cr}/\partial_\delta\mathrm{cr}$ and of $G'$ are
computed. They give $\partial_\delta\mathrm{cr}\ge194$, $\lambda\ge1.068$,
$-0.0065\le d\ell^*/d\kappa\le-0.0031$ on $[7497/10000,3/4]$,
$d\ell^*/d\kappa\le-0.0031$ on $[7497/10000,1]$ and $G'\ge0.498$.
Interval evaluation also gives $G(13/18)<0<G(1)$ and the signs of $G$ at
the two ends of the displayed bracket.

Since $G$ vanishes at
$\kappa_\flat$, eliminating $\delta$ between $\mathrm{cr}=0$ and
$\ell^*(\delta,\kappa)=4(11/12-(1+\kappa)/2)$ shows that $\kappa_\flat$ is
a root of the resultant, whose only relevant factor is the displayed
quartic; that quartic has exactly one root in $[7/10,8/10]$. Its other
root $0.94606\ldots$ in $[13/18,1]$ is spurious: it arises from the
root of $\mathrm{cr}(\cdot,\kappa)$ outside $[1/8,5/6]$, and $G$ has no
zero there (the script
\path{certificates/optimality_proof_check.py};
replay details are given in
Appendix~\ref{app:refined-chain}). Finally $\kappa=2B-1$ turns the
quartic into \eqref{eq:v-refined-limit} for $B_\flat=(1+\kappa_\flat)/2$.
\end{proof}

\begin{lemma}[The tangent family]\label{v:lem:refined-family}
Let $7497732940799/10^{13}\le\kappa\le3/4$ (a range containing
$[\kappa_\flat,3/4]$), let $\delta_c=\delta_c(\kappa)$
and $\ell^*=\ell^*(\kappa)$ be as in Lemma~\ref{v:lem:refined-critical},
and let $b^*=b^*(\kappa)$ be the unique solution of the affine equation
$\partial_\delta E_{b,\ell^*,\kappa}(\delta_c,1/2)=0$; it is unique
because on $\mathrm{cr}=0$ the coefficient of $b$ equals
$-\partial_\delta\mathrm{cr}/(864\kappa^2\mathcal J_\kappa)\ne0$, by
$\partial_\delta\mathrm{cr}\ge194$. Then:
\begin{enumerate}
\item $E_{b^*,\ell^*,\kappa}(\delta,x)\le0$ on $[1/8,5/6]\times[0,1/2]$,
      with equality only at $(\delta,x)=(\delta_c,1/2)$;
\item for every $\ell<\ell^*$, $E_{b^*,\ell,\kappa}(\delta,x)
      \le-(\ell^*-\ell)/8$ on $[1/8,5/6]\times[0,1/2]$;
\item conditions (1), (2) and (4) of
      Remark~\ref{v:rem:refined-admissible} hold at $(b^*,\ell^*,\kappa)$
      with strict inequalities, apart from the $\kappa$-range condition,
      which holds (non-strictly at $\kappa=3/4$).
\end{enumerate}
Consequently, for $\kappa_\flat\le\kappa\le3/4$ and every $\ell<\ell^*(\kappa)$
sufficiently close to $\ell^*(\kappa)$, the triple $(b^*(\kappa),\ell,\kappa)$
is admissible.
\end{lemma}

\begin{proof}
Write $2\mathcal J_\kappa(-E_{b^*,\ell^*,\kappa})=\mathcal A(y)\delta^2
+\mathscr B(y)\delta+\mathcal C(y)$ as in
Lemma~\ref{v:lem:refined-square}, with $\mathfrak L=1$, $r_0=0$ and
coefficients depending on $\kappa$. At $y=0$ this quadratic vanishes at $\delta_c$, because
$E_{b,\ell^*,\kappa}(\delta_c,1/2)=0$ by Lemma~\ref{v:lem:refined-critical},
and its $\delta$-derivative vanishes there by the choice of $b^*$.
Hence $\delta_c$ is a double root and $4\mathcal A\mathcal C-\mathscr B^2$
vanishes at $y=0$; write it as $y(d_1+d_2y+d_3y^2+d_4y^3)$.

Interval
arithmetic on $4000$ subintervals with exact rational endpoints
and outward-rounded coverage of the stated $\kappa$-range, with
rigorous enclosures of $\delta_c$, $\ell^*$ and $b^*$ obtained from the
signs of $\mathrm{cr}$ and from the rational formulas for $\ell^*$ and
$b^*$, gives
\[
 \begin{gathered}
 d_1>0.1197,\quad d_2>0.395,\quad d_3>0.371,\quad d_4>0.0836,\\
 \alpha_0>0.414,\quad \alpha_1>0.803,\quad \alpha_2>0.596,\quad
 \alpha_3>0.272,
 \end{gathered}
\]
where $\mathcal A=\sum_k\alpha_ky^k$, together with all scalar
conditions in~(3) and $b^*>0$. For $y>0$ Lemma~\ref{v:lem:refined-square} gives positivity; for $y=0$
the quadratic is $\alpha_0(\delta-\delta_c)^2$. This proves~(1)
and~(3).

For~(2), $E$ is affine in $\ell$ at fixed $(b,\kappa,\delta,x)$, with
\[
 \partial_\ell E_{b,\ell,\kappa}(\delta,x)
 =\frac54+\delta x\Bigl(1+\frac{5c_\kappa}8\Bigr)
       -\frac98\bigl(1-R_\kappa(\delta,x)\bigr)\ge\frac18,
\]
since $0<R_\kappa\le1$; so (1) gives (2). The final assertion follows
from (2), (3) and continuity in $\ell$ of the finitely many strict scalar
conditions. The computations are performed by
\path{certificates/family_certificate.py};
see Appendix~\ref{app:refined-chain} for replay instructions.
\end{proof}

\begin{theorem}[Limiting refined cyclotomic bound]
\label{v:thm:refined-limit}
For every fixed cyclotomic field $F\supset\mu_{12}$ and every
finite-order Hecke character $\eta$ of $F$,
\[
 L_F(s,\eta)\ne0\qquad(\Rea s>B_\flat),
\]
where $B_\flat=0.874886647039964\ldots$ is defined by
\eqref{eq:v-refined-limit}. Principal poles are permitted; no assertion
is made on the boundary.
\end{theorem}

\begin{proof}
Suppose $\beta=\beta_*(F)>B_\flat$. By
Theorem~\ref{v:thm:original-target}, $\beta\le B_\infty<7/8$. Put
$\kappa=2\beta-1$, so that $\kappa_\flat<\kappa\le2B_\infty-1<3/4$ and the hypothesis
$\beta_*(F)\le(1+\kappa)/2$ of M2 holds with equality. By
Lemma~\ref{v:lem:refined-critical}(3), $G(\kappa)>0$, that is,
$F_\flat(\kappa)<(1+\kappa)/2=\beta$. Choose $\ell<\ell^*(\kappa)$ so close
to $\ell^*(\kappa)$ that $B(\ell)<\beta$ and, by
Lemma~\ref{v:lem:refined-family}, the triple $(b^*(\kappa),\ell,\kappa)$ is
admissible. Lemma~\ref{v:lem:refined-closure} gives $\beta\le B(\ell)<\beta$,
a contradiction. Beyond $\beta\le B_\infty$ only this single step is
used, so no uniformity of analytic constants over an iteration is
required. Omitted Euler factors vanish only on $\Rea s=0$,
and principal poles give zeros of the reciprocal, so the bound for
$\beta_*(F)$ gives the assertion for every finite-order $\eta$.
\end{proof}

\subsection{Optimality in the refined comparison system}
\label{v:sec:refined-optimality}
In this subsection, Theorem~\ref{v:thm:refined-optimality} shows that
$B_\flat$ of \eqref{eq:v-refined-limit} is the exact infimum of the
refined comparison system defined below, and
Remark~\ref{v:rem:refined-scope} records what that system does not
cover. As in Section~\ref{sec:v-optimality}, the statement concerns the
specified uniform estimates, not the actual zeros.

\paragraph{The refined comparison system.} A \emph{refined step} consists
of balanced lengths with $b>0$ and $0<\ell\le1/5$, the direct majorant
\eqref{v:cmp:4} matched at $B(\ell)$, a moment parameter
$13/18\le\kappa\le1$ justified by a previously proved bound
$\beta_*(F)\le B_{\rm prev}$, that is, $\kappa\ge2B_{\rm prev}-1$ when
$\kappa<1$, the counts \eqref{v:cmp:10}, and, for every class with
$\delta\ge1/8$ and normalized amplitude $x$, the numerator average of
Proposition~\ref{v:prop:refined-numerator} with any shift
$0\le\theta\le g/(6\kappa)$, evaluated through
\eqref{eq:v-refined-full-exponent} at $d=h$. The step
yields $B(\ell)$ only if this endpoint exponent is
$\le0$ for all $(\delta,x)$ with
$1/8\le\delta\le2B(\ell)-1$, $0\le x\le1/2$ (bins up to
$a<\beta_*(F)$ occur in the contradiction argument). Conversely, a step
with $(b,\ell,\kappa)$ admissible (Remark~\ref{v:rem:refined-admissible})
yields $B(\ell)$; this scalar rule is realized by Lemma~\ref{v:lem:refined-closure}. As in Section~\ref{sec:v-detail-37}, a uniform
estimate must control the formal class $x=1/2$, whose nonemptiness is
not asserted. A parameter $\kappa>1$ is
dominated by $\kappa=1$ (Section~\ref{sec:v-detail-33}).

\begin{theorem}[Optimal refined comparison boundary]
\label{v:thm:refined-optimality}
In the refined comparison system, every finite chain of refined steps
whose initial bound is a bound $B_{\rm prev}\ge B_\flat$, for example
$11/12$, $7/8$ or $B_\infty$, yields only boundaries $B\ge B_\flat$, and
a chain whose initial bound exceeds $B_\flat$ never attains $B_\flat$.
Conversely, for every $\epsilon>0$ a finite chain of admissible steps
whose initial bound is $7/8$ yields a boundary below $B_\flat+\epsilon$.
Hence $B_\flat$ is the exact infimum of the system. These are scalar
statements about the system;
Lemma~\ref{v:lem:refined-closure} realizes each admissible step, and
$B_\flat$ is the boundary of Theorem~\ref{v:thm:refined-limit}.
\end{theorem}

\begin{proof}
Consider one refined step with boundary $B(\ell)$ and justifying bound
$B_{\rm prev}\ge B_\flat$. At every class the shift enters
\eqref{eq:v-refined-full-exponent} at $d=h$ only through
$\theta(h/2-l_y)=\theta(5\ell-1-b)/4$, whose coefficient is negative
because $\ell\le1/5<(1+b)/5$. Hence the least endpoint is attained at
$\theta=g/(6\kappa)$, where it equals $E_{b,\ell,\kappa}(\delta,x)$, and
the step requires $E_{b,\ell,\kappa}(\delta,x)\le0$ at every controlled
class. Since $\ell\le1/5$ gives $2B(\ell)-1\ge11/15>12/25$, the critical
class $(\delta_c(\kappa),1/2)$ of Lemma~\ref{v:lem:refined-critical} is
controlled, and part~(2) of that lemma gives $\ell\le\ell^*(\kappa)$,
that is,
\[
 B(\ell)\ge F_\flat(\kappa).
\]

If $\kappa<1$, then $\kappa\ge2B_{\rm prev}-1\ge2B_\flat-1=\kappa_\flat$;
if $\kappa=1$, then also $\kappa\ge\kappa_\flat$. Since
$\kappa_\flat>7497/10000$, part~(3) gives
$F_\flat(\kappa)\ge F_\flat(\kappa_\flat)=B_\flat$. Thus every step yields
$B(\ell)\ge B_\flat$, and the next justifying bound again satisfies
$B_{\rm prev}\ge B_\flat$. Induction over the chain proves the lower
bound. If $B_{\rm prev}>B_\flat$, then $\kappa>\kappa_\flat$ when $\kappa<1$,
and part~(3) gives the strict inequality
$B(\ell)\ge F_\flat(\kappa)>B_\flat$, as also for $\kappa=1$;
so such a chain never attains $B_\flat$.

Chains may also interleave whole steps of the pointwise-numerator
system. Such a step, justified by $B_{\rm prev}$, yields a boundary
$B\ge\min\{B_{\rm prev},B_\infty\}$: if $B<B_{\rm prev}$ then
$B<(1+\kappa)/2$. The step uses a moment parameter
$\kappa\ge\max\{13/18,2B_{\rm prev}-1\}$ in the range
$13/18\le\kappa\le1$ of Section~\ref{sec:v-detail-33}, and its two contour
duals give $B\ge F(\kappa)$ there; thus the compatibility
$F(\kappa)\le B\le(1+\kappa)/2$ used in the proof of
Theorem~\ref{thm:v-optimality} holds, and that proof gives $B\ge B_\infty$. Since $B_\infty>B_\flat$, the
induction is unchanged.

For the converse, let $0<\epsilon\le10^{-4}$, start from the initial bound
$7/8$, that is $\kappa_0=3/4$, and define inductively
$\kappa_{j+1}=2B_j-1$, where $B_j=B(\ell_j)$ with $\ell_j<\ell^*(\kappa_j)$
chosen by Lemma~\ref{v:lem:refined-family} so that
$(b^*(\kappa_j),\ell_j,\kappa_j)$ is admissible and
$B_j<F_\flat(\kappa_j)+\epsilon_j$, $\epsilon_j=4^{-j}\epsilon$.

Put
$\Delta_j=\kappa_j-\kappa_\flat$. The lower bound gives $B_j>B_\flat$, hence
$\Delta_{j+1}>0$. If $\kappa_j\in(\kappa_\flat,3/4]$, then
$0<F_\flat'\le0.0065/4$ on $[\kappa_\flat,3/4]$ (proof of
Lemma~\ref{v:lem:refined-critical}) gives
\[
 \Delta_{j+1}=2\bigl(B_j-F_\flat(\kappa_\flat)\bigr)<0.0033\,\Delta_j+2\epsilon_j .
\]
Since $\Delta_0=3/4-\kappa_\flat<2.27\times10^{-4}$, this gives
$\Delta_{j+1}<7.5\times10^{-7}+2\times10^{-4}<3/4-\kappa_\flat$, so by induction
every $\kappa_j$ lies in $(\kappa_\flat,3/4]$, the family range. The same
recursion gives $\Delta_j\to0$, and
$B_j<F_\flat(\kappa_j)+\epsilon_j\le B_\flat+0.0017\,\Delta_j+\epsilon_j\to B_\flat$.
Each step is justified by the bound proved at the previous step.
\end{proof}

\begin{remark}[Scope of the refinement]\label{v:rem:refined-scope}
Theorem~\ref{v:thm:refined-optimality} concerns the stated
system: balanced lengths with $\ell\le1/5$, the counts
\eqref{v:cmp:10}, and the averaged numerator at every class with
$\delta\ge1/8$. It does not cover a class-by-class choice between the
averaged and pointwise numerator within one step (at the
tangent parameters the pointwise endpoint at the critical class is
$\approx+4\times10^{-4}$, and such mixing does not appear to lower $B_\flat$),
unequal lengths, or stronger analytic inputs (a
smaller count at the critical class, a larger admissible prime length in M2, a
sharper Gram or reflected-energy estimate, a smaller amplification
exponent).

Theorem~\ref{thm:v-optimality} concerns
the pointwise-numerator system and is not contradicted, since
Proposition~\ref{v:prop:refined-numerator} lies
outside that system. $B_\flat$ is obtained by the
single step of Theorem~\ref{v:thm:refined-limit}; the
rational chain of Appendix~\ref{app:refined-chain} stops at
$B_\sharp$, with $B_\sharp-B_\flat<1.04\times10^{-12}$.
\end{remark}

\section{Optimality in the specified comparison system}
\label{v:sec:strong-parameter-barrier}
\label{sec:v-optimality}
In this section, we show that the pointwise-numerator comparison of
Section~\ref{v:sec:compensated-assembly} has the optimal boundary of
Theorem~\ref{thm:iv-system-optimality}: with free
total length, asymmetry, detector cutoffs and the local and global
numerator contours below, Theorem~\ref{thm:v-optimality} gives
$B\ge B_\infty$, attained only as a limit in
Theorem~\ref{v:thm:original-target}. The argument follows
Sections~\ref{sec:iv-count-obstruction}--\ref{sec:iv-further-parameters}
over $F$.
This system bounds the reflected numerator pointwise.
The refined boundaries $B_\flat<B_\sharp<B_\infty$ of
Section~\ref{v:sec:refined} and Appendix~\ref{app:refined-chain} use the
averaged numerator estimate of Proposition~\ref{v:prop:refined-numerator},
are not parameter choices in this system, and are not covered by
Theorem~\ref{thm:v-optimality}; Theorem~\ref{thm:v-optimality} is
transferred to arbitrary base fields in Section~\ref{vi:sec:optimality}.

\subsection{The estimate system}
\label{sec:v-detail-33}
In this subsection, we fix the estimate system to be optimized: the
normalization is that of the principal term in
Section~\ref{sec:v-detail-28}, with unequal lengths allowed and relaxed
length conditions, and the low theta average is bounded by the direct
majorant $B_{\rm low}$ of Proposition~\ref{prop:prev-free-low}. We also
fix the moment parameter $\kappa$, the root $\xi(\kappa)$ and
$F(\kappa)=(18\xi+7)/(18\xi+9)$; the critical
class of Section~\ref{sec:iv-count-obstruction} is retained for every
admissible choice, and its count is minimized in
Section~\ref{sec:v-detail-34}.

Put
\[
 M=l_x+l_y,\quad b=l_y-l_x,\quad h=1-l_x+\ell,\quad
 C(s)=s-5/6+l_x/3+\ell/6.
\]
The final principal term and low theta average use this normalization.
In the nonempty range of the defining theta sum we relax the conditions
$l_x,l_y>\ell>0$, $b>0$ of $\mathscr A$ in
Section~\ref{sec:iv-system-optimality} to $M,l_y,\ell\ge0$.
The length exponents are fixed and bounded; every remaining full-index
subset, coprimality indicator from the character zero extension and common profile is retained.
The supplied direct majorant is that of Proposition~\ref{prop:prev-free-low}
\cite[(22.4)]{LiuArXiv}:
\[
 B_{\rm low}=5/6+l_x/6-\ell/6+G+\tfrac12 k(M,\ell),
\]
\[
 G=\max(0,b/12,b-l_y/2),\qquad
 k=\max\left(0,\frac{1+3\ell-2M}{4},M+\ell-1,
                  \frac{3\ell-1}{2}\right).
\]

The empty subset controls the full subset estimate: a subset of length
$d$ changes the exponent by at most $-3d/8$. This majorant is the
completed mean square of the original theta sum with the actual additive
Gram, not canonical $Z^F$ energy.
Over $F$ and on the relaxed range it is part of the derivation
(Remark~\ref{v:rem:analytic-inputs}) and goes beyond
\Prev{Proposition}{22.2}, stated for $l_x,l_y>\ell$, $b>0$: the four-case calculation of
Lemma~\ref{lem:prev-free-reflection} is applied to the local length
variables of Section~\ref{sec:v-detail-27}, with the outer range $M'$ free
and the term $\max(H,v+\ell_b)$ retained, and the Gram contribution is that
of \eqref{v:cmp:3}. The continuation argument of
Section~\ref{v:sec:compensated-assembly} uses only the balanced case
\eqref{v:cmp:4}.
Every parameter-only improvement using this majorant has
$B\ge B_{\rm low}$.

A parameter $\kappa>1$ is dominated by the absolute-counting endpoint
$\kappa=1$, since
\[
 Z^{z+\epsilon}\le Z^{\kappa z+\epsilon},\qquad
 n_1+n_2+6z\le n_1+n_2+6\kappa z.
\]
Let $13/18\le\kappa\le1$ be the moment parameter,
$c=2/(3\kappa)$, $\alpha=5/6$, and let $\xi=\xi(\kappa)$ be the smaller
root of
\[
 (216\kappa-18)\xi^2-(468\kappa-57)\xi+150\kappa-20=0.
\]
Then $19/50<\xi<39/100$. Write
\[
 F(\kappa)=\frac{18\xi+7}{18\xi+9},\qquad a=(1+\xi)/2.
\]
The critical uniform class has $\delta=\xi$ and maximal normalized
prime amplitude $u=1/2$; it must be controlled, since no estimate here
proves that it is absent.

\subsection{The critical counting exponent}
\label{sec:v-detail-34}
In this subsection, we show that over all detector cutoffs and length
splits of Section~\ref{sec:iv-count-obstruction}, including amplified
long ranges, the minimal supplied count at the critical class is $2/3$
(the min--max lemma below, crossing $R_{\rm short}$ and $R_{\rm long}$
at $t_*$). This count enters the constraints of
Section~\ref{sec:v-detail-35}; the critical arithmetic class is not
assumed nonempty.

The detecting cutoff and detecting lengths obey
\[
 1\le t\le3/2,\quad t-1/2\le r\le t,\quad0\le m\le1/2,\quad r+m\ge t.
\]
At maximal amplitude the optimized separate supplied count exponents are
\[
 e_I(r)=1-\xi(1+r)/2\quad(r\le1),\qquad
 e_S(m)=1-\xi\{2m+(c/4)(1-2m)\},
\]
\[
 e_L(r)=1-\alpha+(\alpha-\xi)r\quad(r\ge1).
\]
These are exponents in \emph{upper bounds}, not lower bounds for the original sums.
The short inverse length bound is $(1-r)/2$; the second length inequality
for the inverse moment with prime-polynomial factors is weaker for
$r\ge1/2$. The ordinary estimate uses two copies of the same detected
element, giving $|S_m|^4$ and admissible length $(1-2m)/(6\kappa)$; one
copy, or interpolation, gives no better exponent at maximal amplitude.

The next lemma is the counterpart over $F$ of
Lemma~\ref{lem:iv-count-obstruction} and \Prev{Lemma}{34.2}, with the
short/long counts of \eqref{v:cmp:10}.

\begin{lemma}[Critical count min--max]
Let
\[
 R_{\mathrm{opt}}(t)=\sup_{(r,m)}\min(e_{\rm inv}(r),e_S(m)),
 \quad e_{\rm inv}=e_I\ (r\le1),\ e_L\ (r\ge1).
\]
With the above domain of detecting lengths,
\[
 \inf_{1\le t\le3/2}R_{\mathrm{opt}}(t)=2/3.
\]
\end{lemma}
\begin{proof}
For $D=(5-c)/2$, $P=1-c/4$, and $J=(\alpha-\xi)D+\xi P$,
the crossing on $m=t-r$ is
\[
 r_*(t)=\frac{(2-c/2)t+c/4-1/2}{D}.
\]
It is admissible and yields
\[
 R_{\rm short}(t)=1-\xi+\frac{\xi P}{D}(3/2-t).
\]
The second admissible split $(r,m)=(t,0)$ yields
$R_{\rm long}(t)=e_L(t)$ because
$e_S(0)-e_L(t)\ge7/120>0$. Limits cover the absent length-zero
annulus. The short/long upper-bound argument gives the reverse
inequality, so $R_{\mathrm{opt}}=\max(R_{\rm short},R_{\rm long})$.
The first decreases and the second increases; they cross at
\[
 t_*=1+\frac{\xi P}{2J},\qquad
 R_*=1-\xi+\frac{(\alpha-\xi)\xi P}{2J}=2/3,
\]
the last equality being the critical quadratic.
Finite prime supply, strict margins and mesh losses only weaken the
relaxed bounds, and at a given split a convex combination of marginal
exponents is at least their minimum.

Enlarging a moment's combined exponent of the row range and coprimality
modulus to $M'\ge1$ does not help: beyond short length bounds the costs have slopes
$2-\xi>0$ and $6\kappa-\xi>0$; for $r>1$ the unamplified inverse cost
$r(1-\xi)$ is at least $e_L(r)$.
Even granting sixth-power amplification while enlarging a selected
inverse row range $M'$, its count cost is
$1-\alpha+\alpha M'-\xi r-\xi z$. For $r>1$,
$M'\ge r+2z$ gives $e_L(r)+(2\alpha-\xi)z\ge e_L(r)$.
For the ordinary two-copy estimate the corresponding excess slope is
$6\alpha\kappa-\xi>0$. Hence prime-sum factors in the enlarged amplified
moment do not improve the adverse long split.
Different detector cutoffs cannot exclude both adverse splits without
an additional relation between the resulting detected elements.

The same detector outside $[1,3/2]$ does not improve the min--max either.
If $0\le t<1$, use the crossing split for $t=1$ when
$r_*(1)\le t$, and otherwise $(r,m)=(t,0)$; its supplied exponent is at
least $R_{\rm short}(1)>2/3$. If $t>3/2$, the split $(t,0)$ has both
$e_L(t)\ge e_L(3/2)>2/3$ and $e_S(0)>2/3$ (a relaxed scalar
observation, not a detector theorem).
\end{proof}

\subsection{Low constraints and the two numerator-contour duals}
\label{sec:v-detail-35}
In this subsection, we combine Sections~\ref{sec:v-detail-33} and
\ref{sec:v-detail-34} into the affine constraints
\eqref{v:param:1}--\eqref{v:param:3} and bound $B$ from below by
positive-weight duals as in Section~\ref{sec:iv-system-optimality}. Both
signs of the middle numerator slope $l_y-h/2$ are treated: the
positive-slope dual gives $B\ge F(\kappa)$ and the negative-slope dual
gives the larger bound $(11\xi+5)/(9\xi+7)$, so unequal lengths and $b<0$
do not evade the bound. The two duals are combined in the proof of
Theorem~\ref{thm:v-optimality}.

The direct estimate implies
\[
 B\ge L_0=\frac56+\frac{M}{12}-\frac\ell6,\quad
 B\ge L_K=\frac13+\frac{7M}{12}+\frac\ell3,\quad
 B\ge L_J=\frac{23}{24}-\frac M6+\frac{5\ell}{24}.
\]
These use $G\ge b/12$ and the first three branches of $k$.
Moreover $G\ge0$ gives the stronger asymmetry-sensitive pair
\begin{equation}
 B\ge L_0-b/12,\qquad B\ge L_K-b/12.                 \label{v:param:1}
\end{equation}
The positive identity
\[
 \frac{17}{30}L_0+\frac1{30}L_K+\frac25L_J=\frac{13}{15}
\]
forces $B\ge13/15$; with the compatibility $B\le(1+\kappa)/2$ this gives
$\kappa\ge11/15>13/18$, so smaller moment parameters do not help.
If $B<7/8$, then $L_K\le B$ and $\ell\ge0$ give $M<13/14$, hence
\[
 h-l_y=1-M+\ell>1/14,\qquad h>1/14,
\]
so the critical norm dyad $d=h$ is moderate even for $b\le0$ or
$l_x\le\ell$; the restriction $b>0$ is not needed.

At $d=h$ all $z$ terms cancel. For denominator line $s=a$ and numerator
conductor exponent $p(w)$, the optimistic critical error exponent is
\begin{equation}
 E(w)=a-B+h(1/2+p(w))+l_y(w-1)+(\xi-1)\ell/2.       \label{v:param:2}
\end{equation}
Positive contour margins are omitted; Euler losses outside a permitted
box are nonnegative, so omitting them only strengthens a parameter
choice. Buffered endpoint bounds, convexity and reflection give
\begin{equation}
 p(w)=
 \begin{cases}
  1/2-w,&w\le1-a,\\
  (a-w)/2,&1-a\le w\le a,\\
  0,&w\ge a .
 \end{cases}                                      \label{v:param:3}
\end{equation}
The left-exterior slope of \eqref{v:param:2} is $l_y-h<0$, and its right-exterior
slope is $l_y\ge0$, so no exterior numerator contour beats a buffered
endpoint; no Euler holomorphy outside the connecting boxes is needed.

At $w=1-a$ and $w=a$ the exponent is $H_\xi-B$, resp.\ $H_a-B$, where
\[
 H_\xi=1+\xi-(1+\xi)M/2+\xi\ell,\qquad
 H_a=H_\xi+\xi(l_y-h/2).
\]
The middle slope is $l_y-h/2$, so both endpoints must be considered;
the three-constraint dual of Proposition~\ref{prop:iv-dual-optimality}
alone does not exclude negative slopes.

\begin{lemma}[Positive-slope dual]
If $l_y\ge h/2$ and the critical comparison is nonpositive, then
$B\ge F(\kappa)$.
\end{lemma}
\begin{proof}
The minimizing contour is $w=1-a$, so $B\ge H_\xi$. With
\[
 (w_0,w_K,w_H)=\frac{(4+18\xi,2,3)}{9+18\xi},
 \qquad w_0,w_K,w_H>0,\qquad w_0+w_K+w_H=1,
\]
\[
 B\ge w_0L_0+w_KL_K+w_HH_\xi
   =\frac{18\xi+7}{18\xi+9}=F(\kappa).
\]
\end{proof}

\begin{lemma}[Negative-slope dual]
If $l_y\le h/2$ and the critical comparison is nonpositive, then
\[
 B\ge\frac{11\xi+5}{9\xi+7},\qquad
 \frac{11\xi+5}{9\xi+7}>7/8,\qquad
 \frac{11\xi+5}{9\xi+7}>F(\kappa)
\]
for $19/50<\xi<39/100$.
\end{lemma}
\begin{proof}
Put $S=2-3M+2\ell-b=-4(l_y-h/2)\ge0$. The minimizing contour
is $w=a$, so $B\ge H_a$. The positive weights
\[
 (v_0,v_K,v_a)=\frac{(11\xi+1,\,2(1-\xi),\,4)}{9\xi+7}
\]
sum to one, and
\[
 v_0(L_0-b/12)+v_K(L_K-b/12)+v_aH_a
 =\frac{11\xi+5}{9\xi+7}
   +\frac{1-\xi}{4(9\xi+7)}S.
\]
Both claims follow from
\[
 \frac{11\xi+5}{9\xi+7}-\frac78
   =\frac{25\xi-9}{8(9\xi+7)}>0,
\quad
 \frac{11\xi+5}{9\xi+7}-F
   =\frac{4(3\xi-1)(3\xi+1)}{9(2\xi+1)(9\xi+7)}>0.
\]

\end{proof}

\subsection{Global contours and self-consistency}
\label{sec:v-detail-36}
In this subsection, we compare \eqref{v:param:2}--\eqref{v:param:3}
with global reciprocal contours and positive margins, showing that they
do not bypass the local obstruction, and impose the prior zero-free bound
required for the moment parameter, i.e.\ $F(\kappa)\le B\le(1+\kappa)/2$,
as in Corollary~\ref{cor:iv-all-constructions}. Theorem~\ref{thm:v-optimality}
then gives $B\ge B_\infty$, with $\xi_\infty$ a root of
$1188\xi^3-2160\xi^2+171\xi+190$, for every admissible step and every
limit of such steps.

Let $\beta=\beta_*(F)$ be the family supremum and suppose $\beta\ge B$. Grant the strongest standard
zero-free-endpoint interpolation exponent
$p_\beta(w)=(\beta-w)/2$ on $[1-\beta,\beta]$ and retain the same
critical supplied count $2/3$ rather than the weaker trivial count.

On the global denominator line $s=\beta$ the two endpoint exponents are
\[
 E_{\rm glob}(1-\beta)
 =\beta-B+\beta(1-M)+[\beta+(\xi-1)/2]\ell>0,
\]
\[
 E_{\rm glob}(\beta)
 =\beta-B+(1-M)/2+(\beta-1/2)l_y+\xi\ell/2>0.
\]
Their affine interpolation stays positive; exterior slopes give no
improvement since $l_y<h$ on the left and $l_y\ge0$ on the right, and
positive shifts of the denominator contour only add loss. These
estimates are optimistic (actual correction losses only increase the
supplied bound), so global reciprocal contours do not bypass the local
obstruction; absolute contours are weaker still.

The moment parameter must come from a previously proved family bound:
$\beta_*(F)\le\widehat B$ gives $\kappa\ge2\widehat B-1$, as in
Section~\ref{sec:iv-system-optimality}, and in a supremum contradiction
$\kappa=2\beta-1$. The algebraic obstruction is the pair of necessary
conditions $F(\kappa)\le B\le(1+\kappa)/2$.

\begin{theorem}[Optimal comparison boundary]\label{thm:v-optimality}
Within the above uniform supplied-count, completed-low and
pointwise-numerator comparison, varying the lengths, cutoff, finite
mesh, asymmetry and these contours cannot give a boundary below
\[
 B_\infty=\frac{18\xi_\infty+7}{18\xi_\infty+9},
\quad
 1188\xi_\infty^3-2160\xi_\infty^2+171\xi_\infty+190=0,
\]
where $.38858<\xi_\infty<.38859$. 
\end{theorem}
\begin{proof}
The two contour duals force $B\ge F(\kappa)$.
The critical quadratic inverts to
\[
 \kappa(\xi)=\frac{18\xi^2-57\xi+20}{6(36\xi^2-78\xi+25)}.
\]
On the smaller-root interval its derivative is positive and
$36\xi^2-78\xi+25>0$. The fixed-point identity is
\[
 \kappa(\xi)-(2F-1)
 =-\frac{1188\xi^3-2160\xi^2+171\xi+190}
        {18(2\xi+1)(36\xi^2-78\xi+25)}.
\]
Compatibility $F\le B\le(1+\kappa)/2$ implies the numerator cubic is
nonpositive. Its derivative is negative on $[.38,.39]$, so
$\xi\ge\xi_\infty$ and $B\ge F(\xi_\infty)$.
The exact root bracket gives
\[
 .8749570194200989461286<B_\infty<
 .8749570194200989461287.
\]

The compensated assembly realizes the boundary as an infimum of finite
strict steps. The three dual branches $L_0,L_K,H_\xi$ of
Section~\ref{sec:iv-count-obstruction} meet there at
$\ell=(5-6\xi)/(9+18\xi)$, $M=1-\ell$; positive dual weights force these
equalities, so $\ell=\ell_\infty\notin\mathbb Q$ at the infimum. Hence,
as in \Prev{Section}{34.3}, rational parameters have the same infimum
but do not attain it.
\end{proof}

\subsection{Scope of the optimality statement}
\label{sec:v-detail-37}
In this subsection, we delimit the scope of
Theorem~\ref{thm:v-optimality}, as in Section~\ref{sec:iv-further-parameters}:
it optimizes the specified uniform counting
and comparison exponents; it gives no lower bound for the moments of the
original sums, and no impossibility theorem for zero-free regions. We
list changes that require stronger analytic inputs rather than other
parameters.

The uniform comparison allows formal critical classes without proving
them nonempty; an estimate excluding them, or exploiting dependence
between detected elements from different cutoffs, would change the min--max.
Additional mean cancellation in the coupled numerator/correction,
a sharper Gram-matrix estimate for the original sum or completed-energy
estimate, stronger length ranges for inverse moments with
prime-polynomial factors, or an amplification exponent $\alpha<5/6$
likewise changes an input; e.g.\ $\partial_\alpha R=\xi^2P^2/(2J^2)>0$
for the optimized count.

Changing the principal residues $w=1,z=1/6$, using oscillatory
prime profiles that kill their leading normalizer, combining theta averages
with cancellation, or extracting new local residues is a different
operator/estimate, as in item~(3) following
\eqref{eq:prev-critical-necessity}, taken from the proof of
\Prev{Lemma}{34.3}, not a variation of the current exponents.
A common power of $Z$ multiplying theta average and principal term
changes neither comparison; fixed smooth tests and norm twists change
constants and height orders, not these powers.
The cutoff $d=h$ cannot be omitted without an additional row estimate;
its $z$ dependence cancels.
A larger or point-dependent zero-class threshold cannot discard the
critical class without paying the trivial count or proving its absence:
by item~(1) following \eqref{eq:prev-critical-necessity}, a scalar
computation unchanged over $F$, a threshold $\delta_f\ge\xi$ requires
$B\ge H_\xi+h/3$.

The relaxed geometry above, over all norm scales and with the relaxed
length conditions of Section~\ref{sec:v-detail-33}, is a barrier for the
supplied majorant where it applies; it does not assert that every such
geometry has a nonempty theta average or a complete analytic
construction. Identically zero theta averages, or principal terms made
to vanish, do not supply the continuation criterion with a nonzero
principal term.

\section{Class-field descent and the arithmetic scope}
\label{sec:v-descent}
In this section, we descend the cyclotomic bound $\beta_*(F)\le B_\flat$ of
Theorem~\ref{v:thm:refined-limit} to the pairs $(K,\eta)$ of
\eqref{eq:v-class} by abelian class-field factorization, as in
\cite[Section 9.1]{Artin}: Proposition~\ref{prop:v-descent} treats a
fixed character and Corollary~\ref{cor:v-abelian} every abelian $K$.
Proposition~\ref{prop:v-metabelian}
characterizes this scope by the normal closure of $M_\eta$, completing
Theorem~\ref{thm:v-main}. Part~\ref{part:VI} removes the
restriction to $\mathscr C$ by descending through $K(\zeta_{12})$
(Section~\ref{vi:sec:descent}).

\subsection{Descent for a fixed character}
\label{sec:v-character-descent}
In this subsection, we prove Proposition~\ref{prop:v-descent}: Theorem~\ref{v:thm:refined-limit} applies to the Hecke factors
of $\zeta_{M_\eta F}$ over $F$, and
\eqref{eq:v-top-factorization}--\eqref{eq:v-cyclic-descent} descend
zero-freeness on $\Rea s>B_\flat$ to $L_K(s,\eta)$, with $s=1$ handled
by the zeta residues. Corollary~\ref{cor:v-abelian} applies it to every
abelian $K$ via Kronecker--Weber; the same three factorizations reappear
over $K(\zeta_{12})$ in Proposition~\ref{vi:prop:descent}.

\begin{proposition}\label{prop:v-descent}
Let $(K,\eta)\in\mathscr C$, and choose a cyclotomic field $F$
as in \eqref{eq:v-class}. Assume that every finite-order Hecke
function over this $F$ is nonzero on $\Rea s>B_\flat$, with
principal poles permitted. Then $L_K(s,\eta)$ has no zeros for
$\Rea s>B_\flat$, with principal poles permitted.
\end{proposition}
\begin{proof}
Choose $F$ as in \eqref{eq:v-class}, and put $M=M_\eta$, $E=MF$.
Since $E/F$ is abelian, class field theory gives
\begin{equation}\label{eq:v-top-factorization}
 \zeta_E(s)=\prod_{\psi\in\widehat{\operatorname{Gal}(E/F)}}L_F(s,\psi).
\end{equation}
By hypothesis this product is zero-free on $\Rea s>B_\flat$, $s\ne1$.
Since $F/\mathbb Q$ is Galois, $E/M$ is Galois and
\[
 \operatorname{Gal}(E/M)\simeq\operatorname{Gal}(F/F\cap M)
 \subseteq\operatorname{Gal}(F/\mathbb Q)
\]
is abelian, so
\begin{equation}\label{eq:v-first-descent}
 \zeta_E(s)=\zeta_M(s)
 \prod_{\substack{\omega\in\widehat{\operatorname{Gal}(E/M)}\\\omega\ne1}}
 L_M(s,\omega).
\end{equation}
The nonprincipal factors are entire, so $\zeta_M$ has no zeros there
either. For $d=\operatorname{ord}\eta$, the cyclic extension $M/K$ gives
\begin{equation}\label{eq:v-cyclic-descent}
 \zeta_M(s)=\prod_{j=0}^{d-1}L_K(s,\eta^j).
\end{equation}
All factors are holomorphic away from $1$, proving the assertion there.
If $d>1$, the simple zeta residues give
\[
 \prod_{j=1}^{d-1}L_K(1,\eta^j)
 =\frac{\operatorname*{Res}_{s=1}\zeta_M(s)}
 {\operatorname*{Res}_{s=1}\zeta_K(s)}>0,
\]
so each nonprincipal factor is nonzero at $1$; if $d=1$ the simple pole
is permitted. These are identities of full Hecke
$L$-functions, including ramified local factors; see \cite{MilneCFT}.
\end{proof}

\begin{corollary}[Abelian base fields]\label{cor:v-abelian}
For every finite abelian extension $K/\mathbb Q$ and every finite-order
Hecke character $\eta$ of $K$, one has $(K,\eta)\in\mathscr C$.
By Theorem~\ref{v:thm:refined-limit} and
Proposition~\ref{prop:v-descent}, $L_K(s,\eta)\ne0$ for
$\Rea s>B_\flat$, with principal poles permitted.
\end{corollary}
\begin{proof}
By Kronecker--Weber, $K\subset\mathbb Q(\zeta_m)$ for some $m$.
Take $F=\mathbb Q(\zeta_{\operatorname{lcm}(m,12)})$. Then
\[
 K\subseteq F,\qquad \mu_{12}\subset F,\qquad
 \operatorname{Gal}(M_\eta F/F)
 \simeq\operatorname{Gal}(M_\eta/M_\eta\cap F)
 \subseteq\operatorname{Gal}(M_\eta/K).
\]
The last group is cyclic, so $(K,\eta)\in\mathscr C$. By
Theorem~\ref{v:thm:refined-limit}, every finite-order Hecke $L$-function
over $F$ is nonzero on $\Rea s>B_\flat$, which is the hypothesis of
Proposition~\ref{prop:v-descent}; that proposition gives the assertion.
\end{proof}

\subsection{The normal-closure criterion}
\label{sec:v-metabelian}
In this subsection, we characterize the class \eqref{eq:v-class} without
choosing a cyclotomic overfield in advance: by
Proposition~\ref{prop:v-metabelian}, $(K,\eta)\in\mathscr C$ if and only
if the normal closure $N_\eta$ of $M_\eta/\mathbb Q$ has metabelian
Galois group. We then discuss examples, the cyclotomic restriction
removed in Part~\ref{part:VI}, and the scope of
Theorem~\ref{thm:v-optimality} after descent.

For a finite group $G$, write $G'=[G,G]$; the group is metabelian
if $G'$ is abelian.

\begin{proposition}\label{prop:v-metabelian}
Let $N_\eta$ be the normal closure of $M_\eta/\mathbb Q$. Then
\[
 (K,\eta)\in\mathscr C
 \quad\Longleftrightarrow\quad
 \operatorname{Gal}(N_\eta/\mathbb Q)\text{ is metabelian}.
\]
\end{proposition}
\begin{proof}
Suppose $M_\eta F/F$ is abelian for a cyclotomic $F$. Since
$F/\mathbb Q$ is Galois, every $\mathbb Q$-conjugate of $M_\eta F$ is
abelian over $F$, hence so is their compositum $N_\eta F$. In
\[
 1\longrightarrow\operatorname{Gal}(N_\eta F/F)
 \longrightarrow\operatorname{Gal}(N_\eta F/\mathbb Q)
 \longrightarrow\operatorname{Gal}(F/\mathbb Q)\longrightarrow1
\]
both end groups are abelian, so the middle group and its quotient
$\operatorname{Gal}(N_\eta/\mathbb Q)$ are metabelian.

Conversely, put $N=N_\eta$, $G=\operatorname{Gal}(N/\mathbb Q)$ and
suppose $G'$ is abelian. The fixed field $A=N^{G'}$ is abelian over
$\mathbb Q$. Choose a cyclotomic $F\supset A$ containing $\mu_{12}$.
Then
\[
 \operatorname{Gal}(NF/F)\simeq\operatorname{Gal}(N/N\cap F)
 \subseteq\operatorname{Gal}(N/A)=G'
\]
is abelian, hence so is $M_\eta F/F$, i.e.\ $(K,\eta)\in\mathscr C$.
\end{proof}

For $\eta=1$ this gives the Dedekind zeta bound whenever the normal
closure of $K/\mathbb Q$ has metabelian Galois group, e.g.\
$K=\mathbb Q(\sqrt[3]{2})$, with group $S_3$ and $S_3'=A_3\simeq C_3$.
For nonprincipal $\eta$ the criterion concerns $N_\eta$, not merely the
normal closure of $K$; no conclusion for all characters over an
arbitrary nonabelian field follows from the zeta assertion alone.

The metabelian restriction comes from requiring $F$ cyclotomic, for
which $M_\eta F/F$ need not be abelian. Part~\ref{part:VI}
uses $F=K(\zeta_{12})\supseteq K$, for which $M_\eta F/F$ embeds in the
cyclic group $\operatorname{Gal}(M_\eta/K)$ while $M_\eta(\zeta_{12})/M_\eta$
is abelian, and so reaches every pair $(K,\eta)$ under Theorem~\ref{vi:prop:transfer}.

Theorem~\ref{thm:v-optimality} concerns the analytic comparison over $F$,
and the factorizations
\eqref{eq:v-top-factorization}--\eqref{eq:v-cyclic-descent} preserve its
numerical boundary; this asserts no optimality of the true zero-free boundary of any
$L$-function and excludes no estimate outside the pointwise-numerator
system of Section~\ref{sec:v-optimality} (whose optimum is that of
Theorem~\ref{thm:iv-system-optimality}). The descended boundary
$B_\flat$ lies outside that system; its scalar counterpart is
Theorem~\ref{v:thm:refined-optimality}, and descent likewise preserves
its boundary.

\part{Arbitrary number fields}
\label{part:VI}

In this part, we remove the cyclotomic restriction of
Part~\ref{part:V}. Section~\ref{vi:sec:scope} states Theorem~\ref{vi:thm:main}, the zero-free half-plane $\Rea s>B_\flat$ for
every finite-order Hecke $L$-function, in particular every Dedekind zeta
function, over every number field. Section~\ref{vi:sec:generality}
checks that the derivation of Part~\ref{part:V} uses only
$\mu_{12}\subset F$ and runs it over $F=K(\zeta_{12})$
(Theorem~\ref{vi:prop:transfer});
Section~\ref{vi:sec:descent} gives the class-field descent
from $F$ to $K$ (Proposition~\ref{vi:prop:descent}); and
Section~\ref{vi:sec:optimality} records the scope of the optimality
statements. Remark~\ref{vi:rem:status} specifies the analytic inputs
and the relation to earlier work.

\section{The general-field family theorem and its boundary}
\label{vi:sec:scope}
In this section, we state Theorem~\ref{vi:thm:main}, which
extends the refined boundary $B_\flat$ of Part~II to every number
field and finite-order Hecke character in two steps: (i) the analytic
argument of Part~II runs over any number field $F\supset\mu_{12}$,
since its proofs use only $\mu_{12}\subset F$ and global facts valid for
every number field (Section~\ref{vi:sec:generality}, Theorem~\ref{vi:prop:transfer}); (ii) every pair
$(K,\chi)$ is reached from $F=K(\zeta_{12})$ by class-field
factorization (Section~\ref{vi:sec:descent},
Proposition~\ref{vi:prop:descent}). Remark~\ref{vi:rem:status} records
the analytic dependencies and the relation to \cite{CT26}.

Throughout Part~III,
\[
 \begin{gathered}
 B_\flat=0.874886647039964\ldots,\\
 133488B_\flat^4-376596B_\flat^3+384534B_\flat^2-167541B_\flat+26231=0,
 \end{gathered}
\]
is the boundary of \eqref{eq:v-refined-limit}, and for a number field
$F$ we write
\[
 \beta_*(F)=\sup\bigl(\{1/2\}\cup
 \{\Rea\rho:L_F(\rho,\nu)=0,\ \nu\text{ primitive of finite order}\}\bigr),
\]
with only nontrivial zeros included and principal poles excluded, as in
Section~\ref{sec:v-scope}.

\begin{theorem}[Zero-free half-plane for every number field]
\label{vi:thm:main}
For every number field $K$ and every finite-order Hecke character
$\chi$ of $K$,
\[
 L_K(s,\chi)\ne0\qquad(\Rea s>B_\flat),
\]
with the pole of the principal character at $s=1$ permitted. In
particular $\zeta_K(s)\ne0$ for $\Rea s>B_\flat$, $s\ne1$, for every
number field $K$. The boundary is independent of $K$ and $\chi$;
implied constants and thresholds in the underlying comparison may
depend on the fixed fields and characters. No assertion is made on
$\Rea s=B_\flat$.
\end{theorem}

\begin{proof}
Put $F=K(\zeta_{12})$. By Theorem~\ref{vi:prop:transfer},
the proof of Theorem~\ref{v:thm:refined-limit} gives
$\beta_*(F)\le B_\flat$; omitted Euler factors have zeros only on
$\Rea s=0$, so every finite-order Hecke $L$-function over $F$ is
nonzero on $\Rea s>B_\flat$.
Proposition~\ref{vi:prop:descent} gives the assertion for $(K,\chi)$.
\end{proof}

\begin{remark}[Analytic inputs and relation to earlier work]
\label{vi:rem:status}
Part~III uses the analytic estimates of
Sections~\ref{sec:v-packet}--\ref{v:sec:refined}, specified in
Remark~\ref{v:rem:analytic-inputs}, for $F=K(\zeta_{12})$, together with their extension to
fields containing $\mu_{12}$ in Theorem~\ref{vi:prop:transfer}.
Section~\ref{vi:sec:generality} gives the base-field argument and
Section~\ref{vi:sec:descent} the descent.
Cureton and Thorner
\cite[Theorem 1.1, Corollary 1.2 and Sections 5.1 and 14]{CT26} obtain
$\Rea s>31/32$ for every number field in version~2 of their report,
by an argument over
$K(\zeta_{12})$ of the same type; the boundary $B_\flat$ requires the
varying-parameter moments and the averaged numerator of Part~II.
\end{remark}

\section{Base-field generality of the Part~II argument}
\label{vi:sec:generality}
In this section, we check that every input of the analytic argument of
Part~II holds for any number field $F\supset\mu_{12}$:
Lemma~\ref{vi:lem:mu12} collects the consequences of $\mu_{12}\subset F$,
and Theorem~\ref{vi:prop:transfer} transfers the derivation of
Part~II, which gives $\beta_*(F)\le B_\flat$. This is the analytic input of
Theorem~\ref{vi:thm:main}.

Part~II calls its base field cyclotomic, following
\cite[Theorem 1.2 and Proposition 4.1]{Artin}. The proofs in Sections~3--8
and Appendix~A of that source are written for a fixed field and, by a line-by-line check, use ``cyclotomic'' only
through $\mu_{12}\subset F$ (the source's Section~5 opens by recording
only this consequence). The field enters only through
$\mu_{12}\subset F$ (total imaginarity, $\mathbb Q(\omega)\subset F$ in
its Lemma~3.4, the alternating pairing in its Lemma~3.2, and
$\chi_p(-1)=1$), through choices of $S$
valid for every field (principality of $\OO_S$, non-self-dual places,
small primes), and through theorems valid for every number field
(Kazhdan--Patterson, Goldmakher--Louvel, Tate, Weil); the source uses
cyclotomicity genuinely only for the auxiliary fields of
Section~9 of that source, not used here. In Part~II cyclotomicity enters
only the descent of Section~\ref{sec:v-descent} (through the abelian
group $\operatorname{Gal}(F/\mathbb Q)$) and the definitions of
$\beta_*(F)$ for $F=\mathbb Q(\zeta_m)$ and of $\mathscr C$; these are
replaced in Theorem~\ref{vi:prop:transfer} and
Section~\ref{vi:sec:descent}.

\begin{lemma}[Consequences of $\mu_{12}\subset F$]\label{vi:lem:mu12}
Let $F$ be a number field containing $\mu_{12}$, and $r_F=[F:\mathbb Q]/2$.
\begin{enumerate}
\item $F$ is totally complex, contains $\mathbb Q(\omega)$, $\omega^2+\omega+1=0$,
      and $-1=\zeta_{12}^6\in F^{\times6}$.
\item Every finite-order Hecke character of $F$ has trivial infinity type,
      and for primitive $\psi$ the completed $L$-function is
      $\mathcal Q^{s/2}\Gamma(s)^{r_F}L(s,\psi)$ with
      $\mathcal Q=|\Delta_F|q_{\mathfrak f}/(2\pi)^{2r_F}$.
\item After enlarging a finite set $S$ of places by finitely many finite
      places generating the class group, $R=\OO_{F,S}$ is a principal
      ideal domain; $U_S=R^\times$ is finitely generated of rank
      $|S|-1$, and $F_S^1/U_S^6$ is compact.
\item For $\ell\in\{2,3\}$ and $S$ containing the archimedean places and
      the places above $6$,
      $\dim_{\mathbb F_\ell}F_S^\times/F_S^{\times\ell}=2(|S_f|+r_F)$ and
      $\dim_{\mathbb F_\ell}U_S/U_S^\ell=|S_f|+r_F$; the $S$-Hilbert pairing
      of exponent $\ell$ is nondegenerate and alternating.
\item An element of $F^\times$ that is a cube in every completion is a
      cube in $F$.
\end{enumerate}
\end{lemma}

\begin{proof}
(1) A real place would embed $\zeta_3$ into $\mathbb R$; and
$-1=\zeta_{12}^6$. (2) The infinity type of a finite-order character is
trivial on the connected group $(F\otimes\mathbb R)^\times=(\mathbb C^\times)^{r_F}$
(by (1));
the functional equation is Hecke's. (3) The class group is finite, and
Dirichlet's $S$-unit theorem applies; the quotient of the norm-one
subgroup $F_S^1$ by a finite-index subgroup of $U_S$ is compact.

(4) The local dimensions
are the standard counts of $F_v^\times/F_v^{\times\ell}$ when
$\mu_\ell\subset F_v$, including the places above $\ell$, which lie in
$S$; the unit dimension is the $S$-unit theorem together with
$\mu_\ell\subset U_S$. Local nondegeneracy is local Kummer duality, and
the symbol is alternating because $(x,-x)_\ell=1$ and
$-1\in F_v^{\times\ell}$ by~(1). (5) Otherwise $F(x^{1/3})/F$ is a cyclic
cubic extension, and Chebotarev gives a prime outside any given finite
set, unramified and inert in it, with $v(x)=0$, at which $x$ is not a
local cube (Grunwald--Wang for odd exponent).
\end{proof}

\begin{theorem}[Transfer of the Part~II derivation]\label{vi:prop:transfer}
Let $F$ be any number field containing $\mu_{12}$, and choose $S$ as in
Section~\ref{sec:v-scope}; in particular $S$ contains the archimedean
places, the places above $6$ and the support of the different of $F$,
so that $R$ is self-dual under the trace character. The derivations of
Sections~\ref{sec:v-scope} (apart from the class $\mathscr C$ and
Theorem~\ref{thm:v-main}),
\ref{sec:v-packet}--\ref{v:sec:compensated-assembly} and
\ref{v:sec:refined} apply verbatim to $F$, the words ``cyclotomic'' and
``$F=\mathbb Q(\zeta_m)$'' being replaced by ``containing $\mu_{12}$'',
including in the definition of $\beta_*(F)$ in Section~\ref{sec:v-scope}.
In particular the proofs of Theorems~\ref{v:thm:original-target}
and \ref{v:thm:refined-limit} give, for $F$,
\[
 \beta_*(F)\le B_\flat .
\]
The scalar exponent functions, admissibility inequalities and
certificate verifications are the same for every such $F$.
For each fixed $F$, the contradiction geometry and its numerical
choices follow Remark~\ref{v:rem:parameter-ledger} and are common
to all targets over $F$; they may depend on $\beta_*(F)$.
Implied constants and thresholds depend on $F$, $S$ and the fixed
target, and no uniformity in $F$ is used. The only
degree-dependent length is the range $m\le1/2+r_F/100+\eps$ in the proof
of Proposition~\ref{v:prop:refined-numerator}, which is harmless because
M2 without prime sums allows arbitrary bounded lengths.
\end{theorem}

\begin{proof}
For each part of the argument we list the facts used and their source.

\emph{Allowed generators, reciprocity and unit characters}
(Lemmas~\ref{v:lem:allowed-reciprocity} and \ref{v:lem:unit-products}).
The proofs use Lemma~\ref{vi:lem:mu12}(3),(4), Hilbert reciprocity, and
Kummer theory for $\mu_6\subset F$: injectivity of $U_S/U_S^\ell$ uses
only that $R$ has trivial class group and that the places above $\ell$
lie in $S$, so that a Kummer extension unramified outside $S$ and split
at $S$ lies in the trivial $S$-Hilbert class field. Symplectic linear
algebra gives the isotropic complement $B_S$. No step uses
$\operatorname{Gal}(F/\mathbb Q)$.

\emph{Theta datum and coefficient characters}
(Proposition~\ref{v:prop:checked-inputs}, Lemma~\ref{v:lem:theta-characters}).
The Kazhdan--Patterson metaplectic theory
\cite[Introduction and Sections I--II]{KP}, with its correction
\cite{KPerr}, is developed over an arbitrary number field containing $\mu_n$; here $n=3$, and
$\mu_6\subset\mu_{12}\subset F$ also gives total imaginarity. Tame local uniqueness comes from
\cite[Corollary I.3.6]{KP}, whose criterion $2(c+1)\equiv0\pmod n$ depends
only on $n$ and $c$, and uniqueness at the remaining places, including
those above $3$, from the global-to-local argument of
\cite[Theorem II.2.5]{KP}; only soft local choices are made at the fixed
places of $S$. The global compatibility uses
Lemma~\ref{vi:lem:mu12}(5).

The Eisenstein character $\Lambda$ is
defined through the norm map to $\mathbb Q(\omega)\subset F$, which
exists for every finite extension; its construction and the Jacobi-sum
and Davenport--Hasse identities use Lemma~\ref{vi:lem:mu12}(1), the class
number one of $\mathbb Z[\omega]$ (to define the primary generator),
that exterior primes lie over primes of $\mathbb Q(\omega)$ prime to
$3$, and that residue extensions are cyclic. No Galois action of
$F/\mathbb Q(\omega)$ or $F/\mathbb Q$ is used: there is no Galois
equivariance of Gauss sums, no Stickelberger-type factorization and no
uniform splitting of primes; Gauss sums are evaluated prime by prime
from the residue field, and the Weil index product formula holds over
every global field.

\emph{The unitary character in the reflected theta identity.} The
marked terminal lemma uses a unitary character $L$ of $F_S^\times$, fixed
after the classes of $U_S/U_S^6$, taken from \cite[Section 6]{Artin}, and
the finite-place support bounds of \cite[Section 5]{Artin}. Only the
multiplicativity of $L$ on allowed generators and $|L|=1$ are used; such
a character exists over any $F$, and the support bounds are local
statements at the fixed places of $S$.

\emph{Radial functional equation} (Proposition~\ref{v:prop:radial}).
The archimedean computation is done at complex places, which are all
places at infinity by Lemma~\ref{vi:lem:mu12}(1); the gamma multiplier
has exponent $r_F$, as already written in
Section~\ref{sec:v-detail-26}.

\emph{Sieve inputs} (Section~\ref{v:sec:primary-kummer-transport}). The quadratic large sieve of Goldmakher and Louvel
\cite[Definition 1 and Theorem 1.1]{GL} is stated for a quadratic Hecke
family over an arbitrary number field $k$ (Definition~1 requires
$\mu_n\subset k$, here $n=2$), and the cubic large sieve
\cite[Theorem 1.3]{BGL} for a number field $K$ containing $\mu_n$,
$n\ge3$ (hence totally imaginary); here $n=3$.
Lemma~\ref{v:lem:primary-kummer-transport} verifies the required
primitive character families, conductor conditions and finite twists,
using Lemma~\ref{vi:lem:mu12}(2) for trivial infinity type and
Lemma~\ref{v:lem:allowed-reciprocity} for reciprocity.

\emph{Poisson summation over unit orbits, inverse and product moments}
(the sections from Section~\ref{v:sec:primary-kummer-transport} to
Section~\ref{v:sec:compensated-assembly}). These use the
compact quotient $F_S^1/U_S^6$ of Lemma~\ref{vi:lem:mu12}(3), the trace
character and self-dual lattice, Gauss sums and their correlation
identities, and Hecke functional equations with polynomial height
bounds for the fixed field. Constants depend on $F$ and $S$, with no
uniformity in the degree; the numerical mesh
of Lemma~\ref{v:lem:arithmetic-independent-mesh} is chosen before the
field.

\emph{Comparison, detector and refined numerator}
(Sections~\ref{v:sec:compensated-assembly} and \ref{v:sec:refined}).
These use the preceding inputs, the prime ideal theorem in fixed
ray classes, logarithmic control of Hecke $L$-functions over the fixed
field, the local Euler factors, which are computed at exterior primes
from the theta coefficients, and the approximate functional equation of
Lemma~\ref{v:lem:refined-afe}, which requires only
Lemma~\ref{vi:lem:mu12}(2). The exponent functions, the tangent family of
Lemma~\ref{v:lem:refined-family} and the optimality of
Section~\ref{v:sec:refined-optimality} are scalar and field
independent, as are the verifications of
Proposition~\ref{v:prop:refined-parameters} in
Appendix~\ref{app:refined-chain}.

\emph{The initial bound.} The bound $\beta_*(F)\le11/12$ of Theorem~\ref{v:thm:unmarked-seed} and then $7/8$ in
Section~\ref{sec:v-detail-31} concern the finite-order family over the
same field $F$ and use only the preceding estimates; no other field
enters before the descent.

These items and the proofs of Part~II give the assertion. The source-level facts used are the field-generality
of \cite[Sections 3--8 and Appendix A]{Artin} recorded above and the
Kazhdan--Patterson results over a global field containing $\mu_n$
\cite[Introduction]{KP}, $n=3$.
\end{proof}

\section{Descent to every number field}
\label{vi:sec:descent}
In this section, we prove the descent
Proposition~\ref{vi:prop:descent}: zero-freeness on $\Rea s>B$ of all
finite-order Hecke $L$-functions over $F=K(\zeta_{12})$ passes to those
over $K$. Unlike Proposition~\ref{prop:v-descent}, the descent from
$K(\zeta_{12})$ does not require the auxiliary field to be abelian over
$\mathbb Q$. Corollary~\ref{vi:cor:consequences} records the
consequences of Theorem~\ref{vi:thm:main} for Dedekind zeta
functions, monomial Artin $L$-functions and quotients $\zeta_L/\zeta_K$.

\begin{proposition}[Descent through $K(\zeta_{12})$]\label{vi:prop:descent}
Let $B\in[1/2,1)$, let $K$ be a number field, and put $F=K(\zeta_{12})$.
Assume that every finite-order Hecke $L$-function over $F$ is nonzero on
$\Rea s>B$, with principal poles permitted. Then every finite-order
Hecke $L$-function over $K$ is nonzero on $\Rea s>B$, with principal
poles permitted. In particular $\zeta_K(s)\ne0$ for $\Rea s>B$, $s\ne1$.
\end{proposition}

\begin{proof}
Let $\chi$ be a finite-order Hecke character of $K$, of order $d$, and
let $M/K$ be the cyclic extension cut out by $\ker\chi$. Put
$E=MF=M(\zeta_{12})$.

Restriction embeds $\operatorname{Gal}(E/F)$ into the cyclic group
$\operatorname{Gal}(M/K)$, so $E/F$ is abelian and
\[
 \zeta_E(s)=\prod_{\psi\in\widehat{\operatorname{Gal}(E/F)}}L_F(s,\psi)
\]
is a product of finite-order Hecke $L$-functions over $F$; hence
$\zeta_E$ has no zeros in $\Rea s>B$.
Since $\operatorname{Gal}(E/M)\hookrightarrow(\mathbb Z/12\mathbb Z)^\times$,
$E/M$ is abelian and
\[
 \zeta_E(s)=\zeta_M(s)\prod_{\substack{\omega\in\widehat{\operatorname{Gal}(E/M)}\\
 \omega\ne1}}L_M(s,\omega)
\]
with entire nonprincipal factors, so $\zeta_M$ has no zeros in
$\Rea s>B$, $s\ne1$.

Finally, $M/K$ is cyclic and
\[
 \zeta_M(s)=\prod_{j=0}^{d-1}L_K(s,\chi^j)
\]
with factors holomorphic away from $s=1$, so each is zero-free on
$\Rea s>B$, $s\ne1$. If $d>1$, the simple zeta residues give
\[
 \prod_{j=1}^{d-1}L_K(1,\chi^j)
 =\frac{\operatorname*{Res}_{s=1}\zeta_M(s)}
       {\operatorname*{Res}_{s=1}\zeta_K(s)}>0,
\]
so every nonprincipal factor is nonzero at $s=1$. Taking $j=1$ gives
the primitive $L$-function of $\chi$, and $\chi=1$ the zeta assertion.
An imprimitive $L$-function differs by finitely many factors
$1-\chi(\mathfrak p)q_{\mathfrak p}^{-s}$, with zeros on $\Rea s=0$ and
nonzero at $s=1$. These are identities of full Hecke $L$-functions, including the ramified local
factors, valid also when $\chi$ is ramified at real places of $K$; see
\cite{MilneCFT}.
\end{proof}

\begin{corollary}[Consequences]\label{vi:cor:consequences}
Under Theorem~\ref{vi:thm:main}:
\begin{enumerate}
\item every Dedekind zeta function $\zeta_K$ is zero-free on
      $\Rea s>B_\flat$, $s\ne1$;
\item for every finite Galois extension $L/K$, every Artin
      $L$-function $L(s,\rho,L/K)$ of a direct sum $\rho$ of monomial
      representations, in particular of a monomial $\rho$, is zero-free on
      $\Rea s>B_\flat$, with a pole at $s=1$ permitted
      only when $\rho$ contains the trivial representation;
\item for every finite extension $L/K$, the quotient $\zeta_L/\zeta_K$ is
      holomorphic and zero-free on $\Rea s>B_\flat$, including $s=1$.
\end{enumerate}
\end{corollary}

\begin{proof}
(1) is the case $\chi=1$. For (2), a monomial representation is induced
from a linear character $\lambda$ of a subgroup fixing an intermediate
field $K'$, and $L(s,\rho,L/K)=L_{K'}(s,\lambda)$ by inductivity; apply the
theorem to $(K',\lambda)$; the same holds for direct sums of monomial
representations. For (3), both $\zeta_L$ and $\zeta_K$ are holomorphic and
zero-free there away from $s=1$ by (1), and at $s=1$ the quotient is the
ratio of the two positive residues.
\end{proof}

\section{Optimality and scope over arbitrary fields}
\label{vi:sec:optimality}
In this section, we record that the optimality statements of Part~II,
Theorems~\ref{v:thm:refined-optimality} and \ref{thm:v-optimality}, hold
verbatim over every $F\supset\mu_{12}$ and are not affected by the
descent of Proposition~\ref{vi:prop:descent}, and we delimit their scope.

The refined comparison system of
Section~\ref{v:sec:refined-optimality} has the same exponent functions
over every $F\supset\mu_{12}$, since its endpoint, counts and
compatibility conditions are scalar. Hence Theorem~\ref{v:thm:refined-optimality}
holds verbatim over every such $F$: no chain of refined steps gives a
boundary below $B_\flat$, which is the exact infimum of the refined
system. The descent of Proposition~\ref{vi:prop:descent} transfers the
same strict half-plane and loses nothing.

As in Part~II, this optimality
concerns the specified comparison system; it does not exclude improvements from a
class-by-class combination of the averaged and pointwise numerators,
from unequal lengths, or from stronger analytic inputs. Likewise, the
pointwise comparison of Section~\ref{v:sec:compensated-assembly} gives
$\beta_*(F)\le B_\infty$ over every $F\supset\mu_{12}$, and
Theorem~\ref{thm:v-optimality} describes its optimum there. None of
these statements asserts optimality of the true zero-free boundary of
any individual $L$-function.

\appendix
\section{Index of cited results}
\label{app:inputs-index}
In this appendix, the table lists the cited results used in Parts~\ref{part:V}
and~\ref{part:VI} and their roles: the results of \cite{LiuArXiv}
summarized in Section~\ref{sec:prev-iv}, the cyclotomic construction,
descent and scope (Proposition~\ref{prop:v-metabelian}) of
Part~\ref{part:V}, and the base-field generality and descent of
Part~\ref{part:VI}.

\small
\begin{longtable}{@{}>{\raggedright\arraybackslash}p{0.31\textwidth}>{\raggedright\arraybackslash}p{0.64\textwidth}@{}}
\toprule Cited result & Use in Parts II--III\\\midrule
\endfirsthead
\toprule Cited result & Use in Parts II--III\\\midrule
\endhead
\bottomrule\endfoot
\cite[Proposition~3.1, Theorem~27.1, Proposition~18.2, (18.4)--(18.5)]{LiuArXiv}; Section~\ref{sec:prev-limit} &
Continuation criterion, the limiting bound $B_{\mathrm{new}}=B_\infty$, and the comparison exponent $E_B(d)$ used throughout Section~\ref{v:sec:compensated-assembly} and Section~\ref{v:sec:strong-parameter-barrier}.\\
\cite[(4.5), (4.12)--(4.14), (5.1), (7.1)--(7.2), (8.1)--(8.5), Proposition~17.2]{LiuArXiv}; Section~\ref{sec:prev-average} &
Modified average, allowed prime windows, principal correction, detector thresholds and principal savings, as models for Section~\ref{v:sec:compensated-assembly}.\\
\cite[Proposition~10.1, Lemma~11.1, Section~11.3, (12.1), Lemma~19.2]{LiuArXiv}; Section~\ref{sec:prev-initial} &
Initial $7/8$ comparison, its count and endpoint saving, and the row bounds for the remaining norm ranges, in Section~\ref{sec:v-detail-31}.\\
\cite[Lemma~16.1, Proposition~16.2]{LiuArXiv}; Section~\ref{sec:prev-free-direct} &
Retained dual lengths of the reflected numerator and the balanced direct bound, as models for \eqref{v:cmp:2}--\eqref{v:cmp:4} in Section~\ref{sec:v-detail-27}.\\
\cite[Lemma~22.1, Proposition~22.2]{LiuArXiv}; Section~\ref{sec:prev-free-direct} &
The free-length direct majorant $B_{\rm low}$ in Section~\ref{sec:v-optimality}.\\
\cite[Propositions~28.1, 29.1, 32.1; Lemmas~30.1, 31.1, 31.2, 33.1; Section~33]{LiuArXiv}; Section~\ref{sec:prev-iv} &
Models for the moments and count; balanced geometry and fixed-point parameters, endpoint and margin inequalities, one-step bound and iteration, used as stated.\\
\cite[Lemmas~34.2, 34.3; Propositions~34.1, 34.6; Theorem~34.4; Corollary~34.5]{LiuArXiv}; Sections~\ref{sec:iv-count-obstruction}--\ref{sec:iv-further-parameters} &
Critical counting exponent, necessary constraints, dual lower bound and optimality, reused for the pointwise-numerator system of Section~\ref{sec:v-optimality}.\\
\cite[Proposition 0.1.1; Theorems I.2.9 and I.4.3; Corollary I.3.6; Theorems II.2.1 and II.2.5]{KP}; \cite{KPerr} &
Exceptional representation, rational splitting, nonzero Fourier functional and local uniqueness: Proposition~\ref{v:prop:checked-inputs}.\\
\cite[Sections 3--6 and Appendix A]{Artin}; Proposition~\ref{v:prop:radial} &
Allowed generators; the coefficient characters of Lemma~\ref{v:lem:theta-characters}; local Whittaker factors and the radial functional equation with polynomial height bounds.\\
\cite[Definition 1 and Theorem 1.1]{GL}; \cite[Theorem 1.3 and Appendix II]{BGL} &
Primitive quadratic and cubic sieve families, transported in Lemma~\ref{v:lem:primary-kummer-transport}; shared-prime coprimality indicators are separated before applying the bounds.\\
Theorems~\ref{v:thm:unmarked-inverse}, \ref{v:thm:marked-inverse-new} and \ref{v:thm:M2-assembled} &
Inverse and unweighted product moments, retaining coefficient independence, natural zeros, finite test seminorms and the uniformity in the number of prime sums.\\
Section~\ref{v:sec:compensated-assembly}; Sections~\ref{sec:iv-iteration} and \ref{sec:iv-system-optimality} &
Completed theta average, nonzero principal correction, simultaneous count of detected elements, common exponent comparison and the limiting boundary $B_\infty$ of the pointwise-numerator comparison.\\
\cite{XS26}; \cite[Theorem~5.3]{IK}; Section~\ref{v:sec:refined} &
Amplified first moment of the reflected numerator, amplitude-dependent numerator contours, the tangent-family verification, the limiting boundary $B_\flat$ and its optimality in the refined comparison system.\\
Appendix~\ref{app:refined-chain} &
The independent rational chain of four stages verified by exact identities, reaching $B_\sharp$; its first two stages recover the boundary of \cite{XS26}.\\
Theorem~\ref{thm:v-optimality} &
Optimality within the stated pointwise-numerator moment/count and contour comparison system.\\
\cite{MilneCFT}; \cite[Section 9.1]{Artin}; \cite[Sections 5.1 and 14]{CT26} &
Abelian class-field factorization and descent in Proposition~\ref{prop:v-descent}; Kronecker--Weber gives Corollary~\ref{cor:v-abelian}.\\
Proposition~\ref{prop:v-metabelian} &
Equivalence between the cyclotomic compositum condition and metabelian normal closure of the cyclic class field.\\
\cite[Sections I--II]{KP}; \cite{KPerr}; \cite[Definition 1 and Theorem 1.1]{GL}; \cite[Theorem 1.3]{BGL}; \cite[Sections 3--8 and Appendix A]{Artin}; Lemma~\ref{vi:lem:mu12} &
Base-field generality: the Part~II argument over an arbitrary $F\supset\mu_{12}$, Theorem~\ref{vi:prop:transfer}.\\
\cite{MilneCFT}; \cite[Sections 5.1 and 14]{CT26} &
Descent through $K(\zeta_{12})$ to every number field, Proposition~\ref{vi:prop:descent}, and Theorem~\ref{vi:thm:main}.\\
\end{longtable}
\normalsize

\section{An explicit rational chain}
\label{app:refined-chain}

\paragraph{\it Reproducibility of the scalar certificates.}
The scalar replay bundle accompanying this manuscript is in the
directory \path{certificates/}. The entry point
\path{certificates/run_certificates.py} runs the three checkers:
\path{check_refined_certificates.py} checks the four exact rational stages
below, \path{optimality_proof_check.py} checks
Lemma~\ref{v:lem:refined-critical}, and \path{family_certificate.py}
checks Lemma~\ref{v:lem:refined-family}. The interval checkers share
\path{interval_support.py}, which encloses exact rational endpoints,
checks overlap of consecutive boxes, and accepts candidate root brackets
only after interval sign verification on the whole parameter box.
Acceptance comparisons use interval bounds and rational thresholds;
floating-point grids and printed decimal extrema are diagnostics only.
The checks were run with Python~3.11.4, SymPy~1.14.0, mpmath~1.3.0 and
NumPy~2.4.6, pinned in \path{certificates/requirements.txt}; the
accompanying \path{README.md} describes the commands and the checked
scope, and the runner records the versions, exit statuses and SHA-256
hashes of the scripts and logs.
These certificates check the scalar inequalities only; they do not verify the analytic estimates.

\label{v:sec:refined-certificates}
In this appendix, we verify admissibility of four explicit rational
choices (Proposition~\ref{v:prop:refined-parameters}) and run the
corresponding finite chain \eqref{eq:v-refined-order}, which gives
Theorem~\ref{v:thm:refined-target} with
$B_\sharp=0.874886647041$. The proof of Theorem~\ref{v:thm:refined-limit} does not use this appendix; the chain
is an independent exact-arithmetic check of the refined comparison of
Section~\ref{v:sec:refined}, whose first two stages recover over $F$ the boundary
$B_1=8748867/10^7$ obtained by Xu and Sang \cite{XS26} over
$\mathbb Q(\sqrt{-3})$ under the analytic inputs of \cite{QRH}.
The stages are
\begin{equation}\label{eq:v-refined-order}
 \begin{gathered}
 \beta_*(F)\le\frac78
 \Rightarrow\kappa_0=\frac34
 \Rightarrow\beta_*(F)\le B_0
 \Rightarrow\kappa_1=2B_0-1
 \Rightarrow\beta_*(F)\le B_1\\
 \Rightarrow\kappa_2=2B_1-1
 \Rightarrow\beta_*(F)\le B_2
 \Rightarrow\kappa_3=2B_2-1
 \Rightarrow\beta_*(F)\le B_\sharp,
 \end{gathered}
\end{equation}
with
\[
 B_0=\frac{69991}{80000},\quad
 B_1=\frac{8748867}{10^7},\quad
 B_2=\frac{87488664725}{10^{11}},\quad
 B_3=B_\sharp=\frac{874886647041}{10^{12}} .
\]
The stages $j=0,1$ reproduce that construction; the stages $j=2,3$
continue the fixed-point iteration $\kappa_{j+1}=2B_j-1$, whose fixed
point is $\kappa_\flat=2B_\flat-1$; each later $\kappa$ is used only
after its justifying bound is proved.
The exponent functions in Remark~\ref{v:rem:refined-admissible} are
field independent and coincide with those in \cite[Appendix~D]{XS26}.
The verifications for $j=0,1$ are those of that appendix,
rechecked by exact computation; those for $j=2,3$ are new.

\begin{proposition}[Verification of the refined parameters]
\label{v:prop:refined-parameters}
The triples
\[
\begin{aligned}
 (\kappa_0,b_0,\ell_0)&=\Bigl(\frac34,\frac{173}{2000},\frac{10027}{60000}\Bigr),\\
 (\kappa_1,b_1,\ell_1)&=\Bigl(\frac{29991}{40000},\frac{17307}{200000},
                          \frac{1253399}{7500000}\Bigr),\\
 (\kappa_2,b_2,\ell_2)&=\Bigl(\frac{3748867}{5000000},\frac{1081687}{12500000},
                          \frac{501360233}{3000000000}\Bigr),\\
 (\kappa_3,b_3,\ell_3)&=\Bigl(\frac{1499546589}{2000000000},\frac{1081687}{12500000},
                          \frac{125340058877}{750000000000}\Bigr)
\end{aligned}
\]
are admissible. Their boundaries and compatibility identities are
\begin{equation}\label{eq:v-refined-license}
 \begin{gathered}
 B(\ell_0)=B_0,\quad B(\ell_1)=B_1,\quad B(\ell_2)=B_2,\quad
 B(\ell_3)=B_\sharp,\\
 \frac{1+\kappa_0}2=\frac78,\quad\frac{1+\kappa_1}2=B_0,\quad
 \frac{1+\kappa_2}2=B_1,\quad\frac{1+\kappa_3}2=B_2 .
 \end{gathered}
\end{equation}
\end{proposition}

\begin{proof}
The lengths $(l_x,l_y,h)$ are
\[
\begin{gathered}
 \Bigl(\frac{44783}{120000},\frac{55163}{120000},\frac{31757}{40000}\Bigr),\qquad
 \Bigl(\frac{11195177}{30000000},\frac{13791227}{30000000},
                   \frac{7939473}{10000000}\Bigr),\\
 \Bigl(\frac{2239034887}{6000000000},\frac{2758244647}{6000000000},
       \frac{1587895193}{2000000000}\Bigr),\\
 \Bigl(\frac{559758721123}{1500000000000},\frac{689561161123}{1500000000000},
       \frac{396973798877}{500000000000}\Bigr),
\end{gathered}
\]
and substitution gives \eqref{eq:v-refined-license}.

\emph{The endpoint inequality.} Put $y=1/2-x\ge0$. The positive
denominator is
\[
 \mathcal J_{\kappa_j}=(5/6-\delta)\{5/2-c_j+(1+2c_j)y\}
          +\delta\{1-c_j/2+2y+2c_jy^2\},\qquad c_j=c_{\kappa_j}.
\]
Let $E_j$ be \eqref{eq:v-refined-endpoint} at the $j$th choice, and put
$r_{0,0}=10^{-6}$, $r_{0,1}=10^{-7}$, $r_{0,2}=r_{0,3}=0$. The
following are exact polynomial identities in $(\delta,y)$, to which
Lemma~\ref{v:lem:refined-square} applies:
\begin{equation}\label{eq:v-refined-certificate}
 2\mathfrak L_j\mathcal J_{\kappa_j}(-E_j-r_{0,j})
 =\mathcal A_j(y)\delta^2+\mathscr B_j(y)\delta+\mathcal C_j(y).
\end{equation}
For $j=0$, $\mathfrak L_0=2916000000$ and
\[
\begin{aligned}
 \mathcal A_0&=1209956550+2343372300y+1740394800y^2+794068800y^3,\\
 \mathscr B_0&=-1135891798-1173235148y+258347816y^2,\qquad
 \mathcal C_0=266598135+244982070y,
\end{aligned}
\]
and $4\mathcal A_0\mathcal C_0-\mathscr B_0^2=\sum_{k=0}^4d_{0,k}y^k$ with
\[
\begin{gathered}
 d_{0,0}=38461880464196,\quad d_{0,1}=1019289016601843792,\quad
 d_{0,2}=3362710436662080432,\\
 d_{0,3}=3158456603732169536,\quad d_{0,4}=711386879353694144 .
\end{gathered}
\]
For $j=1$, $\mathfrak L_1=1457125331220000000$ and
\[
\begin{aligned}
 \mathcal A_1={}&604514027207496297+1171060420363058502y\\
   &+869869372328057808y^2+396910005991680000y^3,\\
 \mathscr B_1={}&-567451842451072679-586296072701709724y
              +128972962522474840y^2,\\
 \mathcal C_1={}&133165863887059575+122393853878983830y,
\end{aligned}
\]
and $4\mathcal A_1\mathcal C_1-\mathscr B_1^2=\sum_{k=0}^4d_{1,k}y^k$ with
\[
\begin{aligned}
 d_{1,0}&=937158609712709801654383338059,\\
 d_{1,1}&=254346722868854219618398387158475448,\\
 d_{1,2}&=839298823187405114798749944189373584,\\
 d_{1,3}&=788518797558558572150007715039666880,\\
 d_{1,4}&=177683356043985413856187219704574400 .
\end{aligned}
\]
For $j=2$, $\mathfrak L_2=9106994451830472000000000$ and
{\small
\[
\begin{aligned}
 \mathcal A_2={}&3778193352391630914419166+7319095681689997745021676y\\
   &+5436669952149040932279984y^2+2480686164452655840000000y^3,\\
 \mathscr B_2={}&-3546551271962712879464750-3664330633401557418650380y\\
   &+806077231528623210233320y^2,\\
 \mathcal C_2={}&832277812324211460172875+764954585997834384069150y,
\end{aligned}
\]}
and $4\mathcal A_2\mathcal C_2-\mathscr B_2^2=\sum_{k=0}^4d_{2,k}y^k$ with
{\small
\[
\begin{aligned}
 d_{2,0}&=66805203898160358310056061958443800726500,\\
 d_{2,1}&=9935196157823724391857486958898069806530078059600,\\
 d_{2,2}&=32784651795347338089666584660352720939349257313200,\\
 d_{2,3}&=30801169650622232041005296697113854276386680897600,\\
 d_{2,4}&=6940688528488898576468164144697994010575161777600 .
\end{aligned}
\]}
For $j=3$, $\mathfrak L_3=364279675558208657202000000000000$ and
{\small
\[
\begin{aligned}
 \mathcal A_3={}&151127684185839530059810623477846
       +292763748278090742265981699601436y\\
   &+217466752914884442117530106885744y^2
       +99227432036821163838528000000000y^3,\\
 \mathscr B_3={}&-141862002928748992478778925491110
       -146573191082894512715025382770940y\\
   &+32243069682538330286141755738120y^2,\\
 \mathcal C_3={}&33291100806961283141506470318375
       +30598175656855954778442588127350y,
\end{aligned}
\]}
and $4\mathcal A_3\mathcal C_3-\mathscr B_3^2=\sum_{k=0}^4d_{3,k}y^k$ with
{\scriptsize
\[
\begin{aligned}
 d_{3,0}&=857181592314604050391326353715228847575410698827848900,\\
 d_{3,1}&=15896302610172371576342617451594213067657269893353110483242811600,\\
 d_{3,2}&=52455409304147673504470788182123404503762081489232854346746225200,\\
 d_{3,3}&=49281844619233302438837896759897078971214470337076103370658459200,\\
 d_{3,4}&=11105098039212537679396267721796581276981913925108233253978865600 .
\end{aligned}
\]}
All displayed coefficients of $\mathcal A_j$ and of
$4\mathcal A_j\mathcal C_j-\mathscr B_j^2$ are positive. For $y\ge0$ and
real $\delta$,
\[
 4\mathcal A_j(\mathcal A_j\delta^2+\mathscr B_j\delta+\mathcal C_j)
 =(2\mathcal A_j\delta+\mathscr B_j)^2+4\mathcal A_j\mathcal C_j-\mathscr B_j^2>0,
\]
and the multiplier $2\mathfrak L_j\mathcal J_{\kappa_j}$ is positive on the
rectangle. Hence $E_j<-r_{0,j}\le0$ for $0\le\delta\le5/6$,
$0\le x\le1/2$, which is condition~(3).

\emph{The scalar conditions.} The remaining conditions are rational
comparisons. At the four choices (truncated decimals of exact
rationals):
{\scriptsize\[
\begin{array}{c|cccc}
 &j=0&j=1&j=2&j=3\\\hline
 l_x-\ell&0.2060&0.2060&0.2060&0.2060\\
 l_y-\ell-11b/6&0.1339&0.1339&0.1339&0.1339\\
 \ell/h&0.21049&0.21049&0.21049&0.21049\\
 C_0+A_{\rm sm}/8&-0.00418&-0.00417&-0.00417&-0.00417\\
 C_0+(1/2+3\ell/2)/50&-0.01167&-0.01166&-0.01166&-0.01166\\
 13h/75-l_y/2+1/50&-0.07223&-0.07223&-0.07223&-0.07223
\end{array}
\]}
For $j=0,1$ the exact values of the last three rows are
\[
\begin{gathered}
 -\frac{191530333}{45730080000},\quad-\frac{5837}{500000},\quad
 -\frac{433393}{6000000},\\
 -\frac{47785783297087}{11432841120000000},\quad
 -\frac{8748053}{750000000},\quad-\frac{108354377}{1500000000}.
\end{gathered}
\]
Moreover $C_0<0$, $A_{\rm sm}>0$,
$6\kappa_j\ell_j>3/4$, $\ell_j/h_j<1/(6\kappa_j)$,
$h_j/2<l_{y,j}<1/2$, $B(\ell_j)>437/500$, $1/6<\ell_j<1/5$ and
$7497/10000\le\kappa_j\le3/4$, all checked in exact rational arithmetic
by \path{certificates/check_refined_certificates.py}.
\end{proof}

The four stages give a bound weaker than Theorem~\ref{v:thm:refined-limit}, without the critical class or the
tangent family of Section~\ref{v:sec:refined-closure}.

\begin{theorem}[The rational chain]
\label{v:thm:refined-target}
For every fixed cyclotomic field $F\supset\mu_{12}$ and every
finite-order Hecke character $\eta$ of $F$,
\[
 L_F(s,\eta)\ne0\qquad(\Rea s>B_\sharp),\qquad
 B_\sharp=\frac{874886647041}{10^{12}}=0.874886647041.
\]
Principal poles are permitted; no assertion is made on the boundary.
Constants and onsets may depend on the fixed target, whereas this
boundary is common.
\end{theorem}

\begin{proof}
The proof of Theorem~\ref{v:thm:original-target} gives
$\beta_*(F)\le7/8$ (indeed $\beta_*(F)\le B_\infty$). By
\eqref{eq:v-refined-license} this justifies $\kappa_0=3/4$, and
Proposition~\ref{v:prop:refined-parameters} with
Lemma~\ref{v:lem:refined-closure} gives $\beta_*(F)\le B_0$. This
justifies $\kappa_1$, and the next application gives
$\beta_*(F)\le B_1$; two further applications give
$\beta_*(F)\le B_2$ and $\beta_*(F)\le B_3=B_\sharp$. Each stage uses
its own theta average and normalizer, and each $\kappa$ only after its
justifying bound. Finite omitted Euler factors have zeros only on $\Rea s=0$, and
principal poles give zeros of the reciprocal.
\end{proof}

\section{Weighted moment estimates and uniformity}
\label{app:weighted-edges}
In this appendix, we make explicit the weighted estimates used in
Lemmas~\ref{v:lem:inverse-returned-family} and
\ref{v:lem:M2-returned-family}. Section~\ref{app:separate-weights}
gives the common integral estimate, and
Sections~\ref{app:inverse-return-measures} and
\ref{app:product-return-measures} apply it to the transformed inverse and product moments
of Sections~\ref{sec:v-detail-10} and~\ref{sec:v-detail-22}.
Section~\ref{app:M2-recurrence-derivation} combines the M2
transformations and weighted estimates to derive the recurrence.
Section~\ref{app:numerical-edge-budget} bounds the resulting exponent
loss uniformly in the arithmetic data, as required in
Lemma~\ref{v:lem:arithmetic-independent-mesh}.

\subsection{Weighted moment integrals}
\label{app:separate-weights}
In this subsection, we bound the two weighted integrals occurring
after Cauchy--Schwarz. These estimates justify the application of
Lemma~\ref{v:lem:finite-depth-loss} in
Sections~\ref{app:inverse-return-measures} and
\ref{app:product-return-measures}; the test-function bounds follow
from Section~\ref{sec:v-detail-8}.

For nonnegative measurable functions $\mathcal E_i,\mathfrak S_i$,
put
\[
 \mathcal I_i=\int\mathcal E_i\,d\nu_i,\qquad
 \mathcal T_i=\int\mathfrak S_i\,d\nu_i<\infty\qquad(i=1,2).
\]
If $\mathcal E_i\le C_0Z^\lambda\mathfrak S_i$ for $i=1,2$, then
\begin{equation}\label{eq:app-separate-integrals}
 (\mathcal I_1\mathcal I_2)^{1/2}
 \le C_0Z^\lambda(\mathcal T_1\mathcal T_2)^{1/2}.
\end{equation}
In particular, if $(\mathcal T_1\mathcal T_2)^{1/2}
\ll\mathfrak S Z^\delta$, the right side is
$\ll C_0\mathfrak S Z^{\lambda+\delta}$.
For $b_\gamma=|a_1(\gamma)a_2(\gamma)|w_\gamma\ge0$, use
\[
 \sum_\gamma b_\gamma(\mathcal E_1(\gamma)
                          \mathcal E_2(\gamma))^{1/2}
 \le\prod_{i=1}^2
       \left(\sum_\gamma b_\gamma\mathcal E_i(\gamma)\right)^{1/2}.
\]
Thus the paired extraction coefficient $b_\gamma$ is retained in
\emph{both} separate integrals.

Normalize the original test seminorms to be at most one. For the
Euler seminorms of Section~\ref{sec:v-detail-8},
\[
 p_J(W(y)y^{it})\ll_J(1+|t|)^Jp_J(W).
\]
Consequently, for prescribed derivative orders $J_f$ and height
degree $H$,
\begin{equation}\label{eq:app-child-profile-weight}
 \begin{aligned}
 W_{f,\zeta}(y)&=W_f(y)y^{i(t_f+L_f(\zeta))},\\
 \mathfrak S_i
 &=\prod_f p_{J_f}(W_{f,\zeta})^2(1+|\mathbf t_i|)^H
 \ll(1+|\mathbf t|+|\zeta|)^B,
 \qquad B=H+2\sum_fJ_f.
 \end{aligned}
\end{equation}
Increase $B$ for reflected tests when necessary. For the product moment, the new norm
character is common to both ordinary factors, all remaining prime
sums and both rectangles, by Lemma~\ref{v:m2profile:lem:whole}.
Bounded subunit dilations preserve these estimates.

For box dimension $n$, choose $2r>B+n$ in
\eqref{v:eq:H2-common-density}. Uniformly in the divisor parameters
and actual scale ratios,
\begin{equation}\label{eq:app-weighted-density}
 \int(1+|\mathbf t|+|\zeta|)^B
       |\widehat P_{\gamma,Z}(\zeta)|\,d\zeta
 \ll(1+\log Z)^{d_F}(1+|\mathbf t|)^{H_0}.
\end{equation}
For the two successive separations, first integrate in $\zeta$ and
then in $\xi$:
\[
 \iint(1+|\mathbf t|+|\xi|+|\zeta|)^B
 |\widehat P_{1,Z}(\xi)|
 |\widehat P_{2,\gamma,\xi,Z}(\zeta)|\,d\zeta\,d\xi
 \ll(1+\log Z)^{2d_F}(1+|\mathbf t|)^{H_1}.
\]
This proves absolute integrability. Restore the original homogeneous
seminorm factors after integration. The orders depend on the fixed
data, whereas the logarithmic exponent is independent of those orders.

\subsection{The transformed inverse moment}
\label{app:inverse-return-measures}
In this subsection, we apply
\eqref{eq:app-separate-integrals}--\eqref{eq:app-weighted-density}
to the inverse moment of Section~\ref{sec:v-detail-10}. We identify
the transformed coefficient and bound its two integrals, completing
the weighted estimates in Lemma~\ref{v:lem:inverse-returned-family}
used for Theorems~\ref{v:thm:unmarked-inverse}
and~\ref{v:thm:marked-inverse-new}.
Section~\ref{app:product-return-measures} treats the product moment.

With the notation of Sections~\ref{sec:v-detail-10}
and~\ref{sec:v-detail-11},
\begin{align*}
 k'&=d_kd_2k'',& f'&=JCd_2v',&
 \gamma&=(q_0,t',g'/d_2),\\
 \rho_\gamma(n)&=\rho(n)1_{(n,\rad(q_0t'g'/d_2))=1},&&&
 a(n)&=a_0(n)\nu_2(n)\rho_\gamma(n)\chi_n(k')\chi_n(f')^4.
\end{align*}
The actual column and label scales are
\[
 Y_N=Y_{s_i}/(Y_gY_v),\qquad Y_V=Y_CY_\theta Y_vY_J.
\]
Retain $|A(v')|^2$, including its zeros, before using $|A(v')|\le1$;
restore the marked prime sums by Section~\ref{sec:v-detail-13}.
At fixed $(\gamma,f',k')$, the number of reconstructions is
divisor-bounded. The old $f,h$ are not counted again.

Define the weighted sums and integrals by
\[
 \int G\,d\nu_a
 :=Z^{-C_c}\sum_\gamma A_\gamma
       \int G(\gamma,\zeta)
            |\widehat P_{2,a,\gamma,Z}(\zeta)|\,d\zeta,
 \quad
 \sum_\gamma A_\gamma\ll Z^{C_c+11\eta/2+\epsilon_1}.
\]
Any remaining first Fourier integration is performed as in
Section~\ref{app:separate-weights}. The divisor weight
$w_\gamma(f')$ remains in the label sum. Separately for $a=1,2$,
\[
 \mathcal T_a\ll Z^{11\eta/2+\epsilon_1}
 (1+\log Z)^{2d_F}(1+|\mathbf t|)^{H_a}.
\]
Thus the geometric mean contains one projection count. The
normalization and exponent estimates are
\[
 \kappa_i+\lambda_c+C_c+2F_c=F,\qquad
 2|F'-F_c|\le26\eta,\qquad
 81\eta+\frac{11\eta}{2}<256\eta.
\]
These give the allowance of Section~\ref{sec:v-detail-11}, after
subdividing the subsidiary losses to absorb the logarithmic factors.

\subsection{The transformed product moment}
\label{app:product-return-measures}
In this subsection, we apply the weighted integral estimate to
\eqref{v:m2:eq:M2-second-bridge}, retaining the complete divisor
expansions of Section~\ref{sec:v-return-3}. We specify the
transformed moments and their measures. Section~\ref{app:M2-recurrence-derivation}
then derives their recurrence without invoking the closure statement
of Lemma~\ref{v:lem:M2-returned-family};
Section~\ref{app:numerical-edge-budget} bounds the exponent loss.

Fix the first Cauchy norm with allowance $\Lambda_c$ in
\eqref{v:m2:eq:M2-allowance}; its exchange is treated separately.
The substitutions are as follows.
\begin{longtable}{@{}p{.18\textwidth}p{.77\textwidth}@{}}
\toprule
Parameters&Transformed moment\\\midrule
Common divisors&$D_2,E_2$, with their prime-power decompositions and
unit/nonunit partition; $\rad(D_2E_2)=sB$, $(s,B)=1$,
$q_B\ll Z^{p_2}/q_s$. The first divisor $s$ is fixed.\\
Summation and character&$j=G_cV_{\rm id}h'$,
$q_{h'}\le Z^{m'+\delta_2}$, $0\le\delta_2<\xi$;
$\psi_{h'}(n)=\tau_2(n)\chi_n(G_cV_{\rm id}h')$,
with the unit twist included in $\tau_2$. Conjugate the second
column as a whole, retaining its zeros.\\
Moduli&$q'=\widetilde q+w_o+t_2+V$, before character cancellation.
The separate sixth-power coprimality modulus and inherited $R_*$
are completely expanded when omitted from this exponent.\\
Test functions&$X_i'=X_i/(q_{g_i}q_{d_i})$,
$Y_i'=Y_i/(q_{g_i}q_{d_i})$; $X_1'X_2'=Y_1'Y_2'$.
Artificial divisor quotients have no additional coprimality restriction.
Remaining prime sums retain their full annuli; incompatible
assignments contribute zero.\\
Remaining terms&Retain the normalized correlation, extracted
character and test-function values, and M\"obius signs before
absolute estimates. The origins, second diagonal and characters in
$\Theta$ are estimated directly; the other terms are moments with
$M'_{\rm act}\le M-\sigma+\xi$.\\\bottomrule
\end{longtable}
The finite correlation coefficient satisfies
\[
 \left|Z^{-g_2+t_2}
 \mathcal F(D_2,E_2;G_cV_{\rm id}h')\right|\le1
\]
on its specified partition, and is zero elsewhere. Retaining the
central extraction factors, put
\begin{equation}\label{eq:app-product-normalization}
 \begin{aligned}
 \lambda_0&=K+g-\ell-a_0+p_2-s_0+g_2-t_2-b_2,\\
 \Lambda_c-\lambda_0&=M'+\Delta',\qquad \Delta'\ge w+\ell,\\
 Z^{-(p_2-s_0)}\sum_{\gamma_2}a_{\gamma_2}&\ll Z^{\epsilon_a}.
 \end{aligned}
\end{equation}
The factor $Z^{p_2-s_0}$ has already been included in $Z^{\lambda_0}$.

For $T\le q_{\mathfrak t}<2T$, let $t_-=\log_ZT$ and use
$\widetilde r_j,e_j,\omega_j,\pi_j,\pi_{0j}$ from
Section~\ref{sec:v-return-4}. Define
\[
 b_{\mathfrak t,\alpha}=|c_{\mathfrak t,\alpha}|
 Z^{(-r_1^++e_1+\pi_{01}-r_2^++e_2+\pi_{02})/2},\qquad
 r_j^+=\widetilde r_j-(\pi_j-\pi_{0j}).
\]
The coefficient $c_{\mathfrak t,\alpha}$ contains the extracted test
values and the divisor-decomposition multiplicity. By
$\widetilde r_j+\omega_j\ge t_-$,
\begin{align*}
 t_-+\frac{-r_1^++e_1+\pi_{01}-r_2^++e_2+\pi_{02}}2
 &\le\frac{e_1+\omega_1+e_2+\omega_2+\pi_1+\pi_2}2
 \le2\theta_N,\\
 \sum_{\mathfrak t,\alpha}b_{\mathfrak t,\alpha}
 &\ll Z^{2\theta_N+\epsilon_a}.
\end{align*}
Here the sum is over the fixed divisor dyad. Both Cauchy integrals
use this same paired weight:
\begin{align*}
 \int G\,d\nu_i
 &:=Z^{-(p_2-s_0)}\sum_{\gamma_2}a_{\gamma_2}
       \sum_{\mathfrak t,\alpha}b_{\mathfrak t,\alpha}
       \int G\,
       |\widehat P_{2,\gamma_2,\mathfrak t,\alpha,Z}(\zeta)|\,d\zeta,\\
 \mathcal T_i
 &\ll(1+|\mathbf t|)^{H_i}(1+\log Z)^{2d_F}
                  Z^{2\theta_N+\epsilon_b}\qquad(i=1,2).
\end{align*}
The successive sums and integrals are taken in the displayed order;
\eqref{eq:app-weighted-density} gives uniform bounds for the inner
integral. Removal of an initial segment of prime sums gives
\[
 \kappa d_z\le\Delta'+\delta_1/6+\theta_N/3+\eta,
 \qquad
 \eta+\delta_1/6+\theta_N/3+\delta_2+2\theta_N+\epsilon_b
 \le\eta+\frac72\xi+\epsilon_b.
\]
Use \eqref{eq:app-product-normalization} to absorb $\Delta'$.
The first Cauchy norm uses the normalized divisor count
$Z^{-(p+s_0)}\sum_{\gamma_1}a_{\gamma_1}$ and the exponent inequality
following \eqref{v:m2:eq:M2-Gauss-target}; its Fourier integration
is the outer integration of Section~\ref{app:separate-weights}.

For the amplified main term, $Z^{-\ell}$ already occurs in
\eqref{v:m2:eq:M2-second-bridge}. The error terms instead have
$\ell=0$ and the normalized prime sum
$(\#\mathcal A)^{-1}\sum_{p\in\mathcal A}$.
For reflection, with absolute Euler and lower-annulus weights
$a_\alpha\ge0$, Hilbert-space triangle and Cauchy give
\[
 \left\|\sum_\alpha a_\alpha G_\alpha\right\|^2
 \le L\sum_\alpha a_\alpha\|G_\alpha\|^2,
 \qquad L=\sum_\alpha a_\alpha.
\]
Thus the integrated weight has total mass $L^2$. The two-scale
Sobolev estimate of Lemma~\ref{v:m2:lem:M2-reflection} retains
\[
 |G|^2+|\partial_{\lambda_1}G|^2+
 |\partial_{\lambda_2}G|^2+
 |\partial_{\lambda_1}\partial_{\lambda_2}G|^2
\]
on unit logarithmic boxes of total volume $O((1+\log Z)^2)$.
The subsidiary losses absorb these masses; differentiating the tests
increases seminorm orders only.

\subsection{Derivation of the two-transform recurrence}
\label{app:M2-recurrence-derivation}
In this subsection, we derive the recurrence used in the finite-width
induction of Theorem~\ref{v:thm:M2-assembled}. We combine the two
Poisson transformations of Section~\ref{sec:v-detail-22}, the
exceptional-row estimates of Section~\ref{sec:v-detail-23}, and the
weighted measures of Sections~\ref{app:separate-weights}
and~\ref{app:product-return-measures}.

In particular, we use
the signed identities
\eqref{v:m2:eq:M2-first-bridge} and
\eqref{v:m2:eq:M2-second-bridge}, the local Gauss and correlation
formulas and the fixed-character mean estimate. This derives the
conclusion summarized in Lemma~\ref{v:lem:M2-returned-family}
without assuming that lemma. Moment estimates are applied only
at earlier stages of the induction.

\indent{\it The first transform and its allowance.}
Put $M=m+q$, $A=n_1+n_2+z$ and use the complete common-support
extraction $U=Ca$, $V=Db$ of Section~\ref{sec:v-return-1}.
The zero-frequency term satisfies
\[
 \mathcal O_1\ll
 Z^{m-A}\#\{UV\text{ powerful}:q_{UV}\ll Z^{2A}\}
 Z^{\eps_0/2}\ll Z^{m+\eps_0},
\]
after decreasing the input divisor loss. Here a nonzero character
mean forces $UV$ to be powerful, and its allocations contribute
only a fixed divisor power.

For nonzero frequencies the Poisson factor and finite Gauss
normalization multiply to
\[
 \frac{Z^{m-A}}{q_eq_rq_aq_b}\sqrt{q_rq_aq_b}
 =\frac{Z^{m-A}}{q_e\sqrt{q_r}\sqrt{q_aq_b}}.
\]
The inverse roots in \eqref{v:m2:eq:M2-first-bridge} are therefore
retained. Cauchy--Schwarz gives separate positive norms with
characters $\tau_C$ and $\tau_D$; these characters need not agree.

Put
\[
 \widetilde q=q+R+E,\qquad K_0=2A-c-d+R+E-m,
 \qquad K\le K_0+\delta_1,\quad 0\le\delta_1<\xi.
\]
The omitted frequency dyads are estimated on the coprime sum
before the M\"obius extension. Fourier separation in the norms
of the whole products preserves the centered difference and
introduces one common norm power on each side.

For $s\mid a,b$, $s_0=\log_Zq_s$, the required $C$-side exponent
is $A-c+\widetilde q+B_c-s_0$, with
$B_c=(3c-5d-R)_+/6$. There are at most
$Z^{p+s_0+\eps_0}$ choices of the common factors and $s$.
Combining the two Cauchy bounds gives
\[
 M+\tfrac12\{B_c+B_d-(c+d-2p-R)\}\le M,
\]
by the primewise inequality following
\eqref{v:m2:eq:M2-Gauss-target}. Thus the count occurs once.
We now estimate each positive Gauss norm.

\indent{\it The amplification identity and its errors.}
For one such norm, before amplification, put $\alpha_0=A-c$ and
\[
 \mathscr P(h)=Z^{-\alpha_0/2}\sum_u
 B(u)\tau_1(u)q_u^{it}\mathcal G(u,h),
\]
where $B$ includes the condition $s\mid u$ and the centered
difference when present.

For an amplifying prime $p\nmid h$ outside the fixed supports,
put $P=q_p$ and
$\Delta_{p,i}(h)=\mathcal G(p^i,h)-\mathcal G(p^i,hp^6)$.
The exact identity is
\begin{align*}
 \mathscr P(h)-\mathscr P(hp^6)
 =Z^{-\alpha_0/2}\sum_{i\in\{1,6,7\}}\ \sum_{p\nmid v}
 &B(p^iv)\tau_1(p^iv)q_{p^iv}^{it}\chi_p(v)^{2i}\\
 &\times\Delta_{p,i}(h)\mathcal G(v,h).
\end{align*}
Indeed CRT gives the factor $\chi_p(v)^{2i}$ and
$\mathcal G(v,hp^6)=\mathcal G(v,h)$; all other valuations cancel
or vanish by Lemma~\ref{v:m2:lem:full-Gauss-new}.

After normalizing the residual sum at length
$\alpha_0-i\log_ZP$, the squared local coefficients are
\[
 P^{-i}|\Delta_{p,i}(h)|^2=
 \begin{cases}
 P^{-1},&i=1,7,\\
 (1-P^{-1})^2,&i=6.
 \end{cases}
\]
The amplifying primes lie beyond all remaining prime annuli, so
the powers $p^i$ are allocated only to the two ordinary sums.
Both rectangle
scales are divided by the same powers. For $i=1,7$ the residual
character acquires the additional factor supported at $p$; for
$i=6$ the condition $p\nmid v$ is an extra common coprimality condition.

Take $\ell_*=\sigma/3$ and let $\mathcal A$ be the amplifying
primes in Section~\ref{sec:v-return-2}. For $q_h\ll Z^K$, put
\[
 \mathcal A_h=\{p\in\mathcal A:p\nmid hs
       \text{ and }p\text{ avoids the fixed supports}\}.
\]
Then
\[
 \#\mathcal A=Z^{\ell_*+o(1)},\qquad
 \#(\mathcal A\setminus\mathcal A_h)=O(\log Z),\qquad
 \#\mathcal A_h\ge\tfrac12\#\mathcal A.
\]
Averaging the signed identity over $\mathcal A_h$ and applying
Jensen costs a fixed constant.

For each transformed index $H=hp^6$, the number of possible $p$
is bounded: $q_H\le Z^{O(1)}$ and $q_p\ge Z^{\ell_*/2}$.
After this count, enlarge the positive sum over $H$.
The averaged main term has factor $Z^{-\ell_*+o(1)}$.
For the errors, retain $(\#\mathcal A)^{-1}\sum_{p\in\mathcal A}$;
there is no additional factor $\#\mathcal A$.

In each case put
\[
 a_0=\alpha_0-w,\qquad J=d-c+(K_0-K)-2w+w_o.
\]
The choices are
\[
\begin{array}{c|cccc}
 &w&w_o&\ell&g\\\hline
 \text{amplified main term}&0&0&\ell_*&J_++2\sigma\\
 \text{error }i=1,7&i\log_ZP&\log_ZP&0&J_++\sigma\\
 \text{error }i=6&6\log_ZP&0&0&J_++\sigma\\
 \text{without prime sums}&0&0&0&J_++\sigma
\end{array}
\]
The errors are estimated without further amplification.
After removing the squared coefficient $Z^{-w_o}$, the required
bound has exponent
\[
 \Lambda_c=a_0+\widetilde q+w_o+w+B_c-s_0,
\]
as in \eqref{v:m2:eq:M2-allowance}.

\indent{\it The second transform and the transformed moments.}
The resulting positive norm has central factor $Z^{-a_0}$,
row scale $Z^{K+g}$ and density $Z^{-\ell}$.
Finite Fourier inversion supplies the factor
$Z^{K+g-\ell-a_0}/\sqrt{q_uq_v}$ in the second bridge.
At $j=0$, $\mathcal F(u,v;0)$ vanishes unless $u=v$.
The condition $s\mid u$ then gives exponent
$K+g-\ell-s_0+\theta_N+\eps_0$. Its excess over $\Lambda_c$ is
at most $A-M+5\sigma/3+O(\xi+\eps_0)$, which is estimated directly.

For $j\ne0$, remove the genuine high-frequency tail, extract
$D_2,E_2$, and put $j=G_cV_{\rm id}h'$ as in
Section~\ref{sec:v-return-3}. On its partition,
$|\mathcal F(D_2,E_2;j)|\le Z^{g_2-t_2}$.
The transformed character is
$\tau_2(n)\chi_n(G_cV_{\rm id}h')$, including its natural zeros,
and
\[
 \begin{aligned}
 m'&=2a_0-K-g-g_2-V,&
 q'&=\widetilde q+w_o+t_2+V,\\
 M'&=M+J-g-g_2+t_2\le M-\sigma.
 \end{aligned}
\]
Here $g_2\ge t_2$ and $g-J_+\ge\sigma$.
The enclosing ball adds $\delta_2<\xi$ when the sum is nonempty.
Since $s$ divides the genuine common radical, the number of choices
of the common factors is $\ll Z^{p_2-s_0+\eps_0}$.
Multiplication by the Poisson prefactor and the inverse-root and
column-normalization factors gives $Z^{\lambda_0}$ in
\eqref{eq:app-product-normalization}. In particular
\[
 \Lambda_c-\lambda_0=M'+\Delta',\qquad
 \Delta'=b_2-p_2+w+B_c+\ell\ge w+\ell.
\]

Insert the complete last-divisor expansion before positivity.
Its union identity fixes a prime position or divides a plain
variable, with unrestricted quotient; both rectangle scales are
divided by the same factor. Thus unfixed prime factors retain their
full original annuli and coefficients, and the equal-product
identity survives. Extracted character values, including zeros,
are bounded row scalars. Pure sixth-power masks omitted from $q'$
are completely expanded before applying a moment estimate on
nonexceptional rows. The radical of the previously displayed character factors
remains counted in $q'$, even when these factors cancel.

The paired weights $b_{\mathfrak t,\alpha}$ and the normalized
sum over the common factors are those constructed in
Section~\ref{app:product-return-measures}. Their count is bounded
there by $Z^{2\theta_N+\epsilon_a}$.
The excess in the transformed length condition is at most
$6(w+\ell)+\delta_1+2\theta_N$. Removing an initial segment of
prime factors costs
\[
 \kappa d_z\le\Delta'+\delta_1/6+\theta_N/3+\eta.
\]
The pointwise bound of Lemma~\ref{v:m2:lem:prime-capacity} applies
after the residual character has been formed and its exceptional
rows separated. This cost is absorbed by $\Delta'$, apart from
the displayed errors and one final prime of length at most $\eta$.
The transformed moments use their exact rectangle scales.
Nonempty scales below one are converted separately by bounded dilations.

\indent{\it Exceptional rows and the centered saving.}
At an equal multiplicity-one common prime in the nonunit partition,
$G_cV_{\rm id}$ contributes $p^2$. Both columns retain the natural
zeros of the earlier character factors, so this prime lies outside
their support. The earlier finite character is unramified there.
Consequently an exceptional row satisfies
\[
 2+v_p(h')\equiv0\pmod6,\qquad v_p(h')\ge4.
\]
If these primes have radical length $v_1$ and $f=2v_1$, counting
the prescribed sixth-power-free part times a sixth power gives
\[
 \#\{h'\text{ exceptional}\}
 \ll Z^{(m'_{\rm act}-2f)/6+\eps_0}
 \le Z^{(m'-f)/6+\delta_2/6+\eps_0}.
\]
This count applies on the enlarged ball by its character condition.

With raw paired volume the excess is
\eqref{v:m2:eq:M2-exceptional-excess}. Its local estimates give
$F_1\ge2(c+w)/3-3\sigma-O(\xi)$ and $F_2\ge2b_2/3$.
They suffice when $A\le5M/6$.

In the centered range all four initial ordinary scales are at least
$Z^L$. Put
\[
 L=M/4,\qquad v=c+w+\min(c_2,d_2),\qquad r=(L-v)_+.
\]
One side retains this lower scale before the last-divisor removals.
After the Gauss coefficient has disappeared, both its rectangles
have the same fixed character, common coprimality condition and norm
power. Lemma~\ref{v:m2profile:lem:mean} therefore cancels their
leading products exactly and saves $(r-\widetilde r_1)_+$.

Including fixed prime factors in each total removal is essential:
$\widetilde r_j+\omega_j\ge t_-$ implies
\[
 t_--\widetilde r_1-\widetilde r_2
       -(r-\widetilde r_1)_+
 \le\omega_2-r.
\]
Thus the paired volume, including the last-divisor count, is
$Z^{a_0-b_2-r+2\theta_N+\eps_0}$ times the retained seminorm and
height factors. Since $b_2\ge\min(c_2,d_2)$, its excess is at most
\[
 A-5M/6-2v/3-(M/4-v)_++3\sigma+O(\xi+\eps_0)
 \le(A-M)_++3\sigma+O(\xi+\eps_0).
\]
Indeed,
\[
 \sup_{v\ge0}\{-2v/3-(M/4-v)_+\}=-M/6.
\]

Exceptional signed rectangles keep their original scales even below
one; in that case the raw-volume bound of
Proposition~\ref{v:m2profile:prop:center} applies.
Thus the exceptional contribution is bounded directly.

\indent{\it Separate weighted integrals and the recurrence.}
For either application of Cauchy--Schwarz put
\[
 \mathcal I_i=\int E_i\,d\nu_i,\qquad
 \mathcal T_i=\int\mathfrak S_i\,d\nu_i\quad(i=1,2).
\]
The measures $\nu_i$ retain the same paired extraction weight in
both integrals and the normalized sum over the common factors.
For prescribed derivative orders,
\eqref{eq:app-child-profile-weight} bounds the transformed seminorms by
a polynomial in the original and Fourier heights.
Equation~\eqref{eq:app-weighted-density}, integrated successively
in the two Fourier variables, gives
\[
 \mathcal T_i\ll\mathfrak S_{\rm in}
 Z^{2\theta_N+\epsilon_b}(1+\log Z)^{2d_F}
 (1+|\mathbf t|)^{H_i}.
\]
Here $\mathfrak S_{\rm in}$ is a homogeneous seminorm of the initial
test functions. Orders are chosen backward through the finite induction.
The fixed-ratio boxes make the logarithmic exponent independent of
these orders; no factor $Z^{J\eta}$ is introduced.

Ordinary reflection contributes the absolute Euler and annular
masses and the common two-scale Sobolev measure displayed above,
including its mixed derivative. These have subsidiary exponent
loss and finite seminorm and polynomial height bounds.

If the induction hypothesis gives $E_i\le A_*Z^\lambda\mathfrak S_i$,
the required implication is precisely
\begin{equation}\label{eq:app-M2-separate-propagation}
 (\mathcal I_1\mathcal I_2)^{1/2}
 \le A_*Z^\lambda(\mathcal T_1\mathcal T_2)^{1/2}.
\end{equation}
This bounds the two integrals separately and applies also when
$\nu_1\ne\nu_2$. For comparison, when both measures equal $\nu$,
Cauchy--Schwarz gives
\[
 \int(E_1E_2)^{1/2}\,d\nu
 \le\left(\int E_1\,d\nu\int E_2\,d\nu\right)^{1/2}.
\]
Thus a pointwise geometric mean cannot replace the right-hand
side. With different measures, no such pointwise comparison is
made. Apply \eqref{eq:app-M2-separate-propagation} separately to
the two first-transform norms, including when their characters
differ. The resulting bound contains $A_*Z^\lambda$ to the first power.

Denote the moment under consideration, divided by $Z^M$, by $E_v$.
After restoring homogeneous seminorms and absorbing the reserved
logarithmic and divisor losses, the preceding estimates imply
\begin{equation}\label{eq:app-M2-derived-recurrence}
 E_v\le C\mathfrak S_vZ^\tau
       +CA_*\mathfrak S_vZ^{\lambda+2\delta_{\rm edge}},
 \qquad \delta_{\rm edge}\le C_*(\eta+\xi+\eps_0),
\end{equation}
whenever all transformed moments of smaller rank obey the stated
$A_*Z^\lambda$ bound with their required seminorms.
The first zero frequency, second diagonal and exceptional rows
are included in the first term, with
\[
 \tau\le\rho+\delta+14\sigma/3+C_*(\xi+\eps_0).
\]

After two transformations,
\[
 M'_{\rm act}\le M-\sigma+C_*\xi\le M-\sigma/2.
\]
For reflection or centering at the same width, the other indices
in $\mathcal R$ decrease, as shown in the theorem's proof.
Every estimate on the right of \eqref{eq:app-M2-derived-recurrence}
therefore belongs to an earlier stage of that induction.
Ordinary reflection restores the balanced length condition;
the error $\delta$ is not accumulated at each step.

For $M\le\rho$, removal of the remaining prime factors costs at most
$Z^{\kappa z}\le Z^{\rho/6}$. As in the theorem's proof, set
\[
 \widetilde E_v=Z^{-t\rho/6}E_v,\qquad
 t=\begin{cases}0,&\text{without prime sums},\\
                 1,&\text{with prime sums}.
   \end{cases}
\]
This incurs $Z^{\rho/6}$ only once when passing from the ordinary
estimate to the estimate with prime sums. Induction over at most
$D_\sharp$ steps gives total exponent
\[
 \tau+\rho/6+2D_\sharp\delta_{\rm edge}<\eps
\]
with the choices in Theorem~\ref{v:thm:M2-assembled}.
This proves the recurrence and its propagation;
Lemma~\ref{v:lem:M2-returned-family} summarizes the result.

\subsection{An explicit bound for the exponent loss}
\label{app:numerical-edge-budget}
In this subsection, we give a numerical constant for
Lemma~\ref{v:lem:arithmetic-independent-mesh}. We collect the
exponents from Sections~\ref{sec:v-return-1}--\ref{sec:v-return-6}
and Sections~\ref{app:separate-weights},
\ref{app:product-return-measures} and~\ref{app:M2-recurrence-derivation}.
Their uniform bound is used in
Theorem~\ref{v:thm:M2-assembled}, and hence in the compensated
comparison of Section~\ref{v:sec:compensated-assembly} and the
refined remainder in Proposition~\ref{v:prop:refined-high}.

Take $C_*=2^{20}$. Put
\[
 M=m+q,\qquad D_0=d+K_0-K,\qquad b_2=(c_2+d_2)/2,
 \qquad c_{2,1}=c_2,\quad c_{2,2}=d_2,
\]
and let $F_1,F_2$ be as in Section~\ref{sec:v-return-5}.
The following family $\mathscr P$ consists of 32 expressions.
\begin{longtable}{@{}p{.12\textwidth}p{.83\textwidth}@{}}
\toprule
Indices&Exponent expressions\\\midrule
1--5&$\widetilde q$, $K_0$, $K-K_0$, $B_c$, $B_d$.\\
6--8&$m-A-R/2-E+p+s_0$; the two exponents
$A-c+\widetilde q+B_c-s_0$ and its exchange.\\
9--14&$a_0$, $\Lambda_c$, $J$, $m'$, $q'$, $M'$.\\
15--19&$\lambda_0$, $\Delta'$;
$a_0-c_{2,j}+e_j-\widetilde r_j+\pi_j$, $j=1,2$;
$t_-+(-\widetilde r_1+\pi_1+e_1-\widetilde r_2+\pi_2+e_2)/2$.\\
20--24&$w-D_0-w_o+g+2\theta_N$; $F_1$, $F_2$;
the two reflected lengths $M-n_i$.\\
25--28&$3M/2-A+2z$, $3M/2-A+(6\kappa+1)z$;
$r=(M/4-v)_+$, $v=c+w+\min(c_2,d_2)$;
$A-5M/6-2v/3-r$.\\
29--32&$A-M+5\sigma/3+\delta_1$;
$A-5M/6-F_1-F_2+2\theta_N+\delta_2/6+2\eps_0$;
$\delta_1/6+\theta_N/3+\eta$; $M'+\delta_2$.\\\bottomrule
\end{longtable}
For an affine expression, define its coefficient sum by
\[
 c\left(a+\sum_j a_jx_j\right)=\sum_j|a_j|.
\]
For expressions containing positive parts or minima, take the maximum
over their affine pieces. Substitute all displayed definitions,
including $g=J_++\sigma$ or $J_++2\sigma$; the fixed $\sigma$
terms contribute to the terminal allowance. Uniformly for
$13/18\le\kappa\le1$,
\[
 c(J)\le11,\quad c(B_c)\le\frac32,\quad
 c(F_1)<32,\quad c(F_2)=\frac{31}{6},\quad
 \max_{P\in\mathscr P}c(P)\le64.
\]
The two first Cauchy norms, two second transformed factors, two
rectangles and ordinary/comparison estimates require at most
$2^4$ occurrences. Hence
\[
 2^4\sum_{P\in\mathscr P}c(P)
 \le16\cdot32\cdot64=32768<C_*.
\]
The remaining margin includes the two reflection endpoint allowances
and the single excess prime factor in the initial-segment removal.

All prime removals are summed into the displayed total lengths before
this calculation. Fixed support constants and logarithmic factors
satisfy, after the arithmetic data have been fixed,
\[
 \theta_N=\frac{H_N}{\log Z}<\xi,\qquad
 (1+\log Z)^A\ll_{A,\epsilon_0}Z^{\epsilon_0}.
\]
Divide the reserved $\epsilon_0$ among the finitely many divisor,
test-function and logarithmic factors. By
\eqref{eq:app-weighted-density}, derivative orders affect constants
and height polynomials only. This proves
\[
 \delta_{\rm edge}\le C_*(\eta+\xi+\epsilon_0),\qquad C_*=2^{20},
\]
with $C_*$ and the prime-length mesh fixed before
$F,S,\Theta$, the tests and the number of prime sums, as required
in Lemma~\ref{v:lem:arithmetic-independent-mesh}.
\label{app:weighted-edges-end}

\end{document}